%% file: dissertation.tex
\documentclass[final,11pt]{thesis}

\usepackage{url}
\usepackage{amssymb, amsmath, amsthm} 
\usepackage{epsfig,graphicx} 
\usepackage{tikz}
\usepackage{multirow}
\usetikzlibrary{shapes,backgrounds}
\usepackage{subcaption}
\usepackage{xspace}
\usepackage{cite}
\usepackage{arydshln}
\usepackage{color}
\usepackage{verbatim}
\usepackage{tikz}
\usetikzlibrary{shapes,backgrounds}
\usepackage{morefloats}
\usepackage{rotating}
\usepackage{pifont} 
\usepackage{color}
\usepackage{listings}
\usepackage{mathtools} 
\input{macros}

\usepackage{etoolbox}
\usepackage{lipsum}

\BeforeBeginEnvironment{equation}{\begin{singlespace}}
\AfterEndEnvironment{equation}{\end{singlespace}\noindent\ignorespaces}
\BeforeBeginEnvironment{align}{\begin{singlespace}}
\AfterEndEnvironment{align}{\end{singlespace}\noindent\ignorespaces}
\BeforeBeginEnvironment{align*}{\begin{singlespace}}
\AfterEndEnvironment{align*}{\end{singlespace}\noindent\ignorespaces}
\BeforeBeginEnvironment{eqnarray*}{\vspace{-0.75cm}\begin{singlespace}}
\AfterEndEnvironment{eqnarray*}{\end{singlespace}\noindent\ignorespaces}
\BeforeBeginEnvironment{eqnarray}{\vspace{-0.75cm}\begin{singlespace}}
\AfterEndEnvironment{eqnarray}{\end{singlespace}\noindent\ignorespaces}
\newcommand\numberthis{\addtocounter{equation}{1}\tag{\theequation}}
\usepackage{pdfpages}

\title{Prediction in Projection: A new paradigm in delay-coordinate reconstruction}

\author{J.}{Garland}

\otherdegrees{B.S., Colorado Mesa University, 2009 \\
	      M.S., University of Colorado, 2011}

\degree{Doctor of Philosophy}		
	{Ph.D., Computer Science}		

\dept{Department of}			
	{Computer Science}		

\advisor{Prof.}				
	{Elizabeth Bradley}			

\reader{Prof.~James D. Meiss}		
\readerThree{Prof.~James  Crutchfield}		
\readerFour{Prof.~Aaron Clauset}
\readerFive{Prof.~Sriram Sankaranarayanan}
\readerSix{Prof.~ Robert Easton}

\abstract{  \OnePageChapter	

\begin{singlespace}
Delay-coordinate embedding is a powerful, time-tested mathematical
framework for reconstructing the dynamics of a system from a series of
scalar observations.  Most of the associated theory and heuristics are
overly stringent for real-world data, however, and real-time use is
out of the question due to the expert human intuition needed to use
these heuristics correctly.  The approach outlined in this thesis
represents a paradigm shift away from that traditional approach.  I
argue that perfect reconstructions are not only unnecessary for the
purposes of delay-coordinate based forecasting, but that they can
often be less effective than reduced-order versions of those same
models.  I demonstrate this using a range of low- and high-dimensional
dynamical systems, showing that forecast models that employ {\it
imperfect} reconstructions of the dynamics---{\it i.e.,} models that
are not necessarily true embeddings---can produce surprisingly
accurate predictions of the future state of these systems.  I develop
a theoretical framework for understanding why this is so.  This
framework, which combines information theory and computational
topology, also allows one to quantify the amount of predictive
structure in a given time series, and even to choose which forecast
method will be the most effective for those data.
\end{singlespace}	
	}

\dedication[Dedication]{	
To Mom.
\newline\newline
This dissertation---and my entire academic journey---would have been abandoned multiple times, and likely never would have begun, if it were not for your endless love and support.}

\acknowledgements{	\OnePageChapter	
	\begin{singlespace}
		
	\noindent\textbf{Liz}, I have been told that as you look back on your life there will be one or maybe two people that radically alter 
	your entire life course. I have no doubt you are that person in my life. You personally introduced me to complexity science and reinvigorating me as a scientist at a time in my life  
	I was ready to leave academics for good. You have not only shaped me as an academic but as a better human being.  You opened more doors than I can count and did everything in your power to give me every possible 
	opportunity to succeed.  You are the epitome of advisor, mentor, colleague, champion and friend. \newline\newline
	\noindent \textbf{My Family}, without your love and support I am nothing.\newline\newline
	\noindent \textbf{Jim Meiss,} for all the trident visits and willingness to let me bounce ideas off of you over the years, from the beginning to the end of my dissertation. All our conversations about computational topology were invaluable in providing my thesis with a sound theoretical foundation. \newline\newline
	\noindent \textbf{Ryan James,} for introducing me to, and patiently teaching me information theory. 	\newline\newline
	\noindent \textbf{Jim Crutchfield,} for numerous conversations, at the Santa Fe Institute and at UC Davis, throughout the research process. Your guidance in information theory was invaluable in finalizing the story of my thesis. \newline\newline
	\noindent \textbf{Holger Kantz,} for hosting me at Max-Planck-Institut F\"ur Physik Komplexer Systemer, where I had multiple interactions with you and your students that drastically matured the overall picture of my thesis.   \newline\newline
	\noindent\textbf{Bob Easton,} for all the Trident seminars.\newline\newline
	\noindent \textbf{Aaron Clauset}, \textbf{Sriram Sankaranarayanan}, \textbf{Mark Muldoon}, \textbf{Simon DeDeo}, and \textbf{Zach Alexander,} for interesting and useful conversations throughout this process.\newline\newline
	 \noindent \textbf{Santa Fe Institute,} my academic home and sanctuary every summer. Your interdisciplinary halls recharged me every year, gave me a place to forge interesting collaborations with people (and in fields) I never would have imagined and took my career to heights I could have never dreamed. 
	 \end{singlespace}
	}

 \emptyLoF	

 \emptyLoT	

\begin{document}

\begin{singlespace}
\input{overview} 
\input{background} 
\input{caseStudies} 
\input{pnp} 
\input{explanation} 
\input{pecomplexity} 
\input{concl} 
\end{singlespace}
\bibliographystyle{plain}	
\bibliography{master-refs}		


\end{document}

%% file: macros.tex
\newcommand{\alert}[1]{{\color{red}#1}}
\newcommand{\espace}[2]{\ensuremath{\mathbb{E}[#1,#2]}}

\newcommand{\mytau}{\ensuremath{\mathcal{A}_\tau}\xspace}
\newcommand{\smytau}{\ensuremath{\mathcal{A}_\mathcal{S}}\xspace}

\newtheorem*{mydef}{Definition}
\newtheorem{thm}{Theorem}
\newtheorem*{remark}{Remark}

\newtheorem{example}[thm]{Example}

\newcommand{\gcc}{{\tt 403.gcc}\xspace}
\newcommand{\sphinx}{{\tt 482.sphinx}\xspace}
\newcommand{\eig}{{\tt dgeev}\xspace}

\newcommand{\svd}{{\tt dgesdd}\xspace}
\newcommand{\svdone}{{\tt dgesdd}$_1$\xspace}
\newcommand{\svdtwo}{{\tt dgesdd$_2$}\xspace}
\newcommand{\svdthree}{{\tt dgesdd$_3$}\xspace}
\newcommand{\svdfour}{{\tt dgesdd$_4$}\xspace}
\newcommand{\svdfive}{{\tt dgesdd$_5$}\xspace}
\newcommand{\svdsix}{{\tt dgesdd$_6$}\xspace}
\newcommand{\naive}{na\"ive\xspace}
\newcommand{\row}{{\tt row\_major}\xspace}
\newcommand{\col}{{\tt col\_major}\xspace}
\newcommand{\arima}{ARIMA\xspace}

\newcommand{\citeNTau}{\cite{josh-tdAIS,Gibson92,fraser-swinney,Olbrich97,kantz97,Buzug92Comp,liebert-wavering,Buzugfilldeform,Liebert89,rosenstein94}\xspace}

\newcommand{\citeNM}{\cite{liebert-wavering,Cao97Embed,Kugi96,Buzugfilldeform,KBA92,Hegger:1999yq,kantz97}\xspace}

\newcommand{\citeDCEFORECASTING}{\cite{weigend-book,casdagli-eubank92,Smith199250,pikovsky86-sov,sugihara90,lorenz-analogues}}

\newcommand{\fnnLMA}{{\tt fnn-LMA}\xspace}
\newcommand{\roLMA}{{\tt ro-LMA}\xspace}
\newcommand{\ipc}{IPC\xspace}

\usetikzlibrary{calc}
\tikzset{filled/.style={fill=blue, opacity=0.2}}

\def \setupdiagrams {
  \def \delt {2cm}
  \def \radius {2cm}

  \coordinate (A) at (0, 0);
  \coordinate (B) at (240:\delt);
  \coordinate (C) at (-60:\delt);

  \coordinate (center) at (barycentric cs:A=1/3,B=1/3,C=1/3);

  \def \Acirc { (A) circle (\radius) }
  \def \Bcirc { (B) circle (\radius) }
  \def \Ccirc { (C) circle (\radius) }

  \def \drawdiagram {
    \draw[very thick] \Acirc;
    \draw[very thick] \Bcirc;
    \draw[very thick] \Ccirc;

    \node at ($ (center) ! 3   ! (A) $) {$H[X_{j+h}]$};
    \node at ($ (center) ! 3.5 ! (B) $) {$H[X_{j-\tau}]$};
    \node at ($ (center) ! 3.5 ! (C) $) {$H[X_{j}]$};
  }
}

\tikzset{filled/.style={fill=blue, opacity=0.2}}

\def \setupdiagramonecircle {
  \def \delt {2cm}
  \def \radius {2cm}

  \coordinate (A) at (0, 0);

  \def \Acirc { (A) circle (\radius) }

  \def \drawdiagram {
    \draw[very thick] \Acirc;

    \node at ($ (center) ! 3   ! (A) $) {$H[Q]$};
  }
}

\tikzset{filled/.style={fill=blue, opacity=0.2}}

\def \setupdiagramstwocircles {
  \def \delt {2cm}
  \def \radius {2cm}

  \coordinate (B) at (240:\delt);
  \coordinate (C) at (-60:\delt);

  \def \Bcirc { (B) circle (\radius) }
  \def \Ccirc { (C) circle (\radius) }

  \def \drawdiagram {
    \draw[very thick] \Bcirc;
    \draw[very thick] \Ccirc;

    \node at ($ (center) ! 3.5 ! (B) $) {$H[Q]$};
    \node at ($ (center) ! 3.5 ! (C) $) {$H[R]$};
  }
}

\usetikzlibrary{calc}
\tikzset{filled/.style={fill=blue, opacity=0.2}}

\def \setupdiagramsgeneral {
  \def \delt {2cm}
  \def \radius {2cm}

  \coordinate (A) at (0, 0);
  \coordinate (B) at (240:\delt);
  \coordinate (C) at (-60:\delt);

  \coordinate (center) at (barycentric cs:A=1/3,B=1/3,C=1/3);

  \def \Acirc { (A) circle (\radius) }
  \def \Bcirc { (B) circle (\radius) }
  \def \Ccirc { (C) circle (\radius) }

  \def \drawdiagram {
    \draw[very thick] \Acirc;
    \draw[very thick] \Bcirc;
    \draw[very thick] \Ccirc;

    \node at ($ (center) ! 3   ! (A) $) {$H[Q]$};
    \node at ($ (center) ! 3.5 ! (B) $) {$H[R]$};
    \node at ($ (center) ! 3.5 ! (C) $) {$H[S]$};
  }
}

%% file: overview.tex

\chapter{Overview and Motivation}\label{ch:overview}

Complicated nonlinear dynamics are ubiquitous in natural and
engineered systems.  Methods that capture and use the state-space
structure of a dynamical system are a proven strategy for forecasting
the behavior of systems like this, but use of these methods is not always
straightforward.  The governing equations and the state variables are
rarely known; rather, one has a single (or perhaps a few) series of
scalar measurements that can be observed from the system.  It can be a
challenge to model the full dynamics from data like this, especially
in the case of {\it forecast} models, which are only really useful if
they can be constructed and applied on faster time scales than those
of the target system.  While the traditional state-space
reconstruction machinery is a good way to
accomplish the task of modeling the dynamics, it is problematic in
real-time forecasting because it generally requires input from and
interpretation by a human expert.  This thesis argues that that roadblock can be sidestepped
by using a reduced-order variant of delay-coordinate embedding to
build forecast models: I show that the resulting
forecasts can be as good as---or better than---those obtained using complete
embeddings, and with far less computational and human effort. I then explore the underlying reasons for this using a novel combination of techniques from computational topology and information theory.

 Modern approaches to modeling a time series for forecasting arguably began with Yule's work on predicting the annual number of sunspots~\cite{Yule27} through what is now known as {\it autoregression}. Before this, time-series forecasting was done mostly through simple global extrapolation~\cite{weigend-book}. Global linear
methods, of course, are rarely adequate when one is working with
nonlinear dynamical systems; rather, one needs to model the details of
the state-space dynamics in order to make accurate predictions.  The
usual first step in this process is to reconstruct that dynamical
structure from the observed data.  The state-space reconstruction
techniques proposed by Packard {\it et al.} \cite{packard80} in 1980
were a critical breakthrough in this regard.  
In 1981, Takens showed that this method, {\it delay-coordinate embedding},  provides a topologically correct representation of a nonlinear
dynamical system if a specific set of theoretical assumptions are
satisfied. I discuss this in detail in Section~\ref{sec:dce} alongside the appropriate citations.
  
A large number of creative strategies have been developed for using
the state-space structure of a dynamical system to generate
predictions, as discussed in depth in Section~\ref{sec:lma}.
Perhaps the most simple of these is the ``Lorenz Method of Analogues''
(LMA), which is essentially nearest-neighbor
prediction~\cite{lorenz-analogues}. Even this simple strategy, which builds predictions by
looking for the nearest neighbor of a given point and taking that
neighbor's observed path as the forecast---works quite well for
forecasting nonlinear dynamical systems.  LMA and similar methods have
been used successfully to forecast measles and chickenpox
outbreaks~\cite{sugihara90}, marine phytoplankton
populations \cite{sugihara90}, performance dynamics of a running
computer(\begin{it}e.g.,\end{it}~\cite{josh-IDA11,josh-IDA13}), the fluctuations in a
far-infrared laser~\cite{sauer-delay,weigend-book}, and many more.

The reconstruction step that is necessary before any of these methods
can be applied to scalar time-series data, however, can be
problematic.  Delay-coordinate embedding is a powerful piece of
machinery, but estimating good values for its two free parameters, the
time delay $\tau$ and the dimension $m$, is not trivial.  A large
number of heuristics have been proposed for this task, but
these methods, which I cover in depth in Sections~\ref{sec:numericaltau} and \ref{sec:numericalM}, are computationally intensive and they require input
from---and interpretation by---a human expert.  This can be a real
problem in a prediction context: a millisecond-scale forecast is not
useful if it takes seconds or minutes to produce.  If effective
forecast models are to be constructed and applied in a manner that
outpaces the dynamics of the target system, then, one may not be able
to use the full, traditional delay-coordinate embedding machinery to
reconstruct the dynamics.
And the hurdles of delay-coordinate reconstruction are even more
of a problem in nonstationary systems, since the reconstruction
machinery is only guaranteed to work for an infinitely long noise-free
observation of a single dynamical system.  
This means that no matter how much effort and human intuition is put into estimating $m$, or how precise a heuristic is developed for that process, {\it the theoretical constraints of delay-coordinate embedding can never be satisfied in practice}. This means that an experimentalist can never guarantee, on any theoretical basis, the correctness of their embedding, no matter their choice of $m$. In Section~\ref{sec:NLD}, I provide an in-depth discussion of these issues.

The conjecture that forms the basis for this thesis is that a 
formal embedding, although mandatory for detailed dynamical
analysis, {\it is not necessary for the purposes of
prediction}---in particular, that reduced-order variants of delay-coordinate
reconstructions are adequate for the purposes of forecasting, even
though they are not true embeddings~\cite{joshua-pnp}.  As a first step towards
validating that conjecture, I construct two-dimensional
time-delay reconstructions from a number of different time-series
data sets, both simulated and experimental, and then build forecast
models in those spaces.  
I find that forecasts produced using the Lorenz method of analogues
on these reduced-order models of the dynamics are roughly as accurate
as---and often even {\it more} accurate than---forecasts produced by
the same method working in the complete embedding space of the corresponding 
system. This exploration is detailed in Chapter~\ref{ch:pnp}.

\begin{figure*}[ht!]
        \centering
        \begin{subfigure}[b]{0.5\textwidth}
                \includegraphics[width=\textwidth]{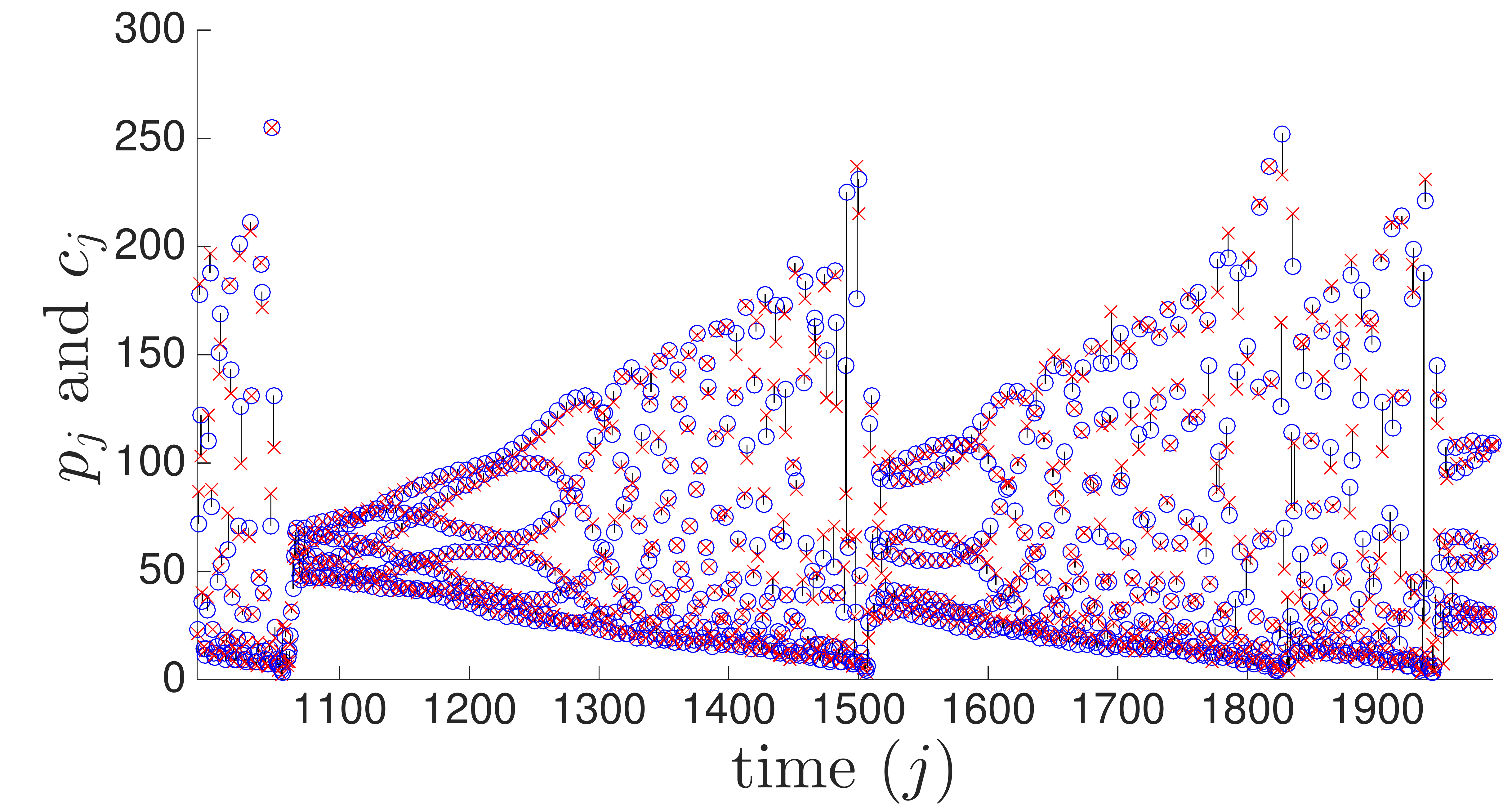}
                \caption{LMA on a complete embedding}
                \label{fig:sfiA12D500tspred}
        \end{subfigure}%
        \begin{subfigure}[b]{0.5\textwidth}
                \includegraphics[width=\textwidth]{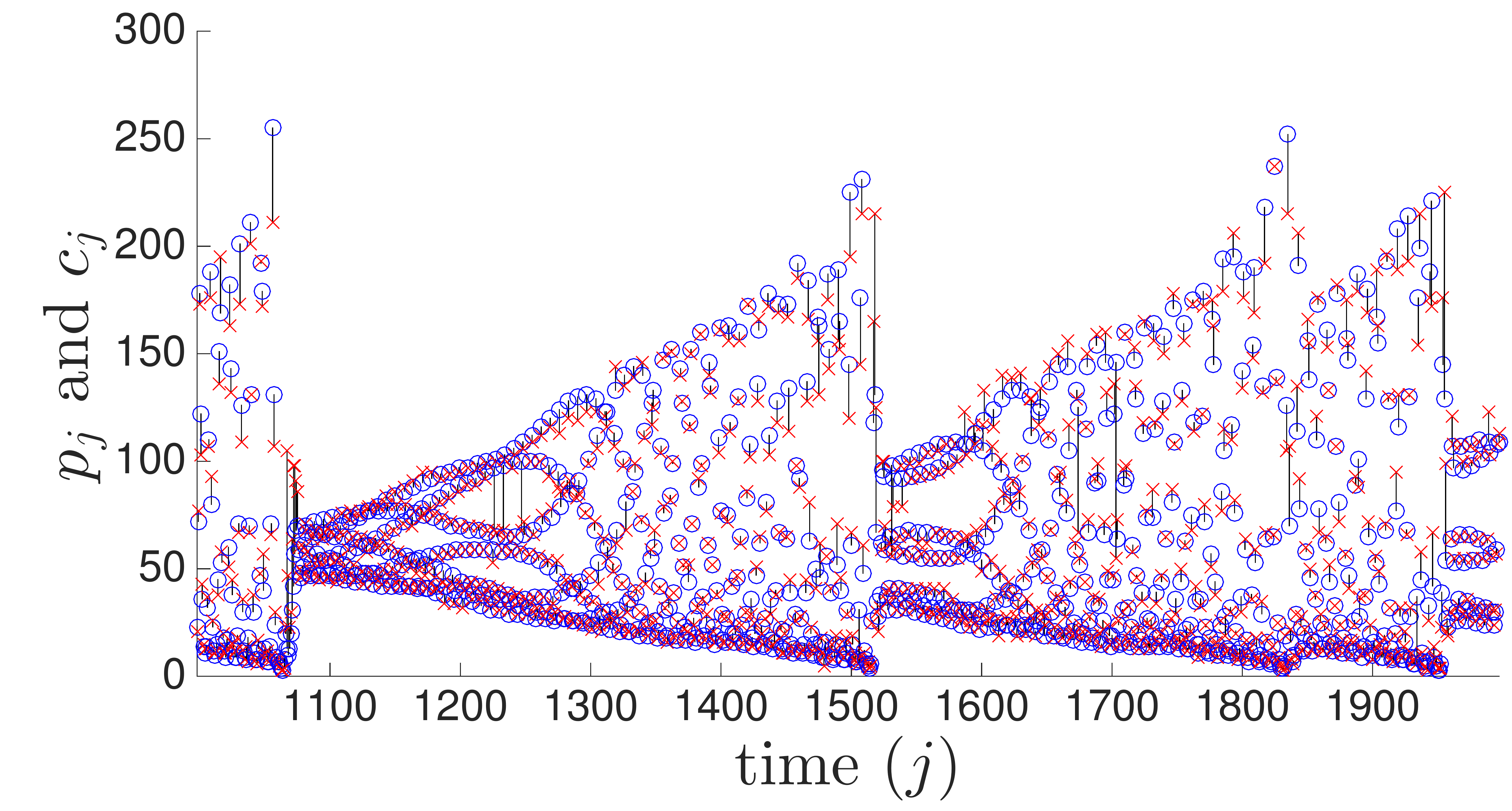}
                \caption{LMA on a two-dimensional delay reconstruction}
                \label{fig:sfiA2D500tspred}
        \end{subfigure}

\caption{Forecasts of SFI Dataset A using
Lorenz's method of analogues in {(a) a delay-coordinate embedding of the
state-space dynamics and (b) a $2D$ delay reconstruction of those
dynamics}.  Blue {\color{blue}o}s are the true continuation $c_j$ of
the time series and red {\color{red}x}s ($p_j$) are the forecasts;
black error bars are provided if there is a discrepancy between the
two.  Reconstruction parameter values for (a) were estimated using
standard techniques: the first minimum of the average mutual
information \cite{fraser-swinney} for the delay in both images and the
false-near neighbor (FNN) method of Kennel {\sl et al.}~\cite{KBA92},
with a threshold of 10\%, for the dimension in the left-hand image.
Even though the $2D$ reconstruction used in (b) is not faithful to the
underlying topology, it enables successful forecasting of the time
series.}
\label{fig:projExample}
\end{figure*}

Figure~\ref{fig:projExample} shows a quick proof-of-concept example: a pair of
forecasts of the so-called ``Dataset A," a time series from a
far-infrared laser from the Santa Fe Institute prediction
competition~\cite{weigend-book}.
Even though the low-dimensional reconstruction used to generate the
forecast in the right panel of the figure is not completely faithful to
the underlying dynamics of this system, it appears to be good enough
to support accurate short-term forecast models of nonlinear dynamics.
While this example is encouraging, 
Dataset A is only one time series and it was drawn from a
comparatively simple system---one that is well-described by a
first-return map (or, equivalently, a one-dimensional surface of
section). The examples presented in Chapter~\ref{ch:pnp} offer a
broader validation of this thesis's central claim by constructing 
forecasts using two-dimensional delay reconstructions
of ensembles of data sets from a number of different systems whose
dynamics are far more complex than the time series in
Figure~\ref{fig:projExample}.  Uniformly, the results indicate that
the full complexity (and effort) of the delay-coordinate `unfolding'
process may not be strictly necessary to the success of forecast
models of real-world nonlinear dynamical systems.  Finally, I want to
emphasize that this reduced-order strategy is intended as a {\it short}-term forecasting
scheme.  Dimensional reduction is a double-edged sword; it enables
on-the-fly forecasting by eliminating a difficult estimation step, but
it effects a guaranteed memory loss in the model.    
I explore this limitation experimentally in Section~\ref{sec:time-scales} and theoretically in Section~\ref{sec:dataandhorizon}.

While the results in Chapter~\ref{ch:pnp} are interesting from a practical perspective, in that they allow delay-coordinate reconstruction to be used in real time, they are perhaps even more interesting from a theoretical perspective. The central premise of this thesis is a heresy, according to the dogma of delay-coordinate embedding, but regardless, it works. This naturally leads to the need for a deeper exploration into why such a deviation from theory provides so much practical traction. 

That exploration is precisely the focus of Chapter~\ref{ch:explain}, where I provide two disjoint explanations of why prediction in projection---my reduced-order strategy---works. The first is from an information theoretic perspective; the second utilizes computational topology. These two disjoint branches of mathematics offer two very different, but quite complementary, tools for exploring this discontinuity between theory and practice. The prior, the subject of Section~\ref{sec:tdAIS}, provides a framework for understanding how information is stored and transmitted from past to future in delay-coordinate reconstructions. 
Building upon ideas from this field, I develop a novel method called {\it time-delayed active information storage} (\mytau) that can be used to select forecast-optimal parameters for delay-coordinate reconstructions~\cite{josh-tdAIS}.
Using \mytau, I show that for noisy finite length time series, a two-dimensional projection ({\it i.e.,} $m=2$) often provides as much---or more---information about the future than a traditional embedding. This further corroborates the central premise of this thesis. This counter-intuitive result, its source, and its implications are discussed in depth in Section~\ref{sec:dataandhorizon}. Section~\ref{sec:compTopo} offers an alternative view of the reconstruction process---one based on topology. As I discuss in Section~\ref{sec:dce}, the theoretical restrictions of delay-coordinate embedding are intended to ensure a {\it diffeomorphic} reconstruction, something that is required for analysis of dynamical invariants but that is excessive for reconstruction of topology. I conjecture that one of the reasons prediction in projection works is that topology (which is preserved by a homeomorphism) becomes correct at much lower embedding dimensions than what one would expect from the associated theorems. The results in Section~\ref{sec:compTopo} confirm this, providing insight into the mechanics of prediction in projection, explaining why this approach exhibits so much accuracy while being theoretically wrong, and offering a deeper understanding into delay-coordinate reconstruction theory in general.

Of course, no forecast model will be ideal for {\it every} task. In fact, as a corollary of the undecidability of the halting
problem~\cite{halting-problem}, no single forecasting schema
is ideal for all noise-free deterministic
signals~\cite{weigend-book}---let alone all real-world time-series data sets.  
I do not want to give the impression that this reduced-order model will be effective for {\it every} time series, but I have shown that it is effective for a broad spectrum of signals. Following this line of reasoning, it is important to be able to determine when prediction is possible at all, and, if so, what forecast model is the best one to use. To this end, I have developed a Shannon information-theory based heuristic for quantifying precisely when a given time-series is predictable given the correct model~\cite{josh-pre}. This heuristic---the focus of Chapter~\ref{ch:wpe}---allows for {\it a priori} evaluation of when prediction in projection will be effective.

The rest of this thesis is organized as follows. Chapter~\ref{ch:background} reviews all the necessary background and related work, including the theory and practice of delay-coordinate embedding, information theory, and the forecast methods, as well as the figure of merit that I use for assessing 
forecast accuracy. In Chapter~\ref{ch:systems}, I introduce the case studies used in this thesis. In Chapter~\ref{ch:pnp}, I demonstrate the effectiveness of this reduced order forecast strategy on a range of different examples, comparing it to traditional linear and nonlinear forecasting strategies, and exploring some of its limitations.  In Chapter~\ref{ch:explain}, I provide a mathematical foundation for {\it why} prediction in projection works.
In Chapter~\ref{ch:wpe}, I describe a measure for quantifying time-series predictability to understand {\it when} my reduced-order method---or  {\sl any} forecasting strategy---will be effective. 
At the end of Chapters~\ref{ch:pnp}-\ref{ch:wpe} I discuss specialized avenues of future research directly associated with the specific contribution of that chapter. In Chapter~\ref{ch:concandfuture}, I conclude and outline the next frontier of this work: developing strategies for grappling with nonstationary time series in the context of delay coordinate based forecasting---which, I believe, will require a combination of all aspects of this thesis to solve.

%% file: background.tex
\chapter{Background and Related Work}\label{ch:background}

\section{Reconstructing Nonlinear Deterministic Dynamics}\label{sec:NLD}
The term \begin{it}nonlinear deterministic dynamical system\end{it} describes a set $\mathbb{X}$ combined with a deterministic nonlinear evolution or update rule $\Phi$, also called the {\it generating equations}. The set $\mathbb{X}$ could be as simple as $\mathbb{R}^n$ or a similar geometric manifold, or as abstract as a set of symbols \cite{meissDynamcis}. Elements of the set $\mathbb{X}$ are referred to as {\it states} of the dynamical system; the set $\mathbb{X}$ is generally referred to as the {\it state space}.  The update or evolution rule is a fixed mapping that gives a unique image to any particular element of the set. In the problems treated in this thesis, this update rule is deterministic and fixed: given a particular state, the next state of the system is completely determined. 
 The theory of dynamical systems is both vast and rich.
This section of this dissertation is intended to review the subset of this field that is needed to understand the core ideas of my thesis. It is not intended as a general review of this field. For more complete reviews, see \cite{meissDynamcis,Holger-and-Liz, kantz97}. 

Dynamical systems can be viewed as falling into one of two categories: those that are discrete in time and those that are continuous in time. The former are referred to as maps and denoted by
\begin{equation}
\vec{y}_{n+1} = \Phi(\vec{y}_n), n \in \mathbb{N}
\end{equation} The latter are referred to as flows and are represented by a system of first-order ordinary differential equations 
\begin{equation}
\frac{d}{dt}\vec{y}(t) = \Phi(\vec{y}(t)), t \in \mathbb{R}^{+}
\end{equation}
When the generating equations $\Phi$ of a dynamical system are known, the future state of any particular initial condition can be completely determined. Unfortunately, knowledge of the generating equations (or even the state space) is a luxury that is very rarely afforded to an experimentalist.  In practice, the dynamical system under study is a black box that is observed at regular time intervals. In such a situation, one can reconstruct the underlying dynamics using so-called {\it delay-coordinate embedding}, the topic of the following section.

\subsection{Delay-Coordinate Embedding}\label{sec:dce}

The process of collecting a {\it time series} $\{x_{j}\}_{j=1}^{N}$ or trace is formally the evaluation of an {\it observation function}~\cite{sauer91} $h: \mathbb{X} \rightarrow \mathbb{R}$ at the {\it true} system state $\vec{y}(t_j)$ at time $t_j$ for $j=1,..,N$, {\it i.e.}, $x_j = h(\vec{y}(t_j))$ for $j=1,\dots,N$.
 Specifically, $h$ smoothly observes the path of the dynamical system through state space at regular time intervals, {\it e.g.}, measuring the angular position of a pendulum every 0.01 seconds or measuring the average number of instructions executed by a computer per cycle~\cite{mytkowicz09,zach-IDA10}.  Provided that the underlying dynamics $\Phi$ and the observation function $h$---are both smooth and generic,
Takens \cite{takens} formally proves that the delay coordinate map 
\begin{equation}\label{eqn:takens}
F(h,\Phi,\tau,m)(\vec{y}(t_j)) = (  [x_j~x_{j-\tau} ~ \dots ~ x_{j-(m-1)\tau}])=\vec{x}_j
\end{equation}
 from an $d$-dimensional smooth compact manifold $M$ to $\mathbb{R}^{m}$ is almost always a diffeomorphism on $M$ whenever $\tau>0$ and
$m$ is large enough, {\it i.e.,} 
$m>2d$. 
\begin{mydef}[Diffeomorphism, Diffeomorphic]
A function $f: M \rightarrow N$ is said to be a diffeomorphism if it is a $C^1$ bijective correspondence whose inverse is also $C^1$. Two manifolds $M$ and $N$ are said to be diffeomorphic if there exists a diffeomorphism $F$ that maps $M$ onto $N$. 
\end{mydef}

What all of this means is that, given an observable deterministic dynamical system---a computer for example, a highly
complex nonlinear dynamical system \cite{mytkowicz09} with no obvious $(\mathbb{X},\Phi)$---I can
measure a single quantity ({\it e.g.}, instructions executed per cycle or  L2 cache misses) and use that time series to
faithfully reconstruct the underlying dynamics {\it up to
  diffeomorphism}. In other words, the true unknown dynamics $(\mathbb{X},\Phi)$ and the
dynamics reconstructed from this scalar time series {\it have the same
  topology}.  Though this is less information than one might like, it
is still very useful, since many important dynamical properties ({\it e.g.},
the Lyapunov exponent that parametrizes chaos) are invariant under
diffeomorphism.  It is also useful for the purposes of prediction---the goal of  this thesis.

The delay-coordinate embedding process involves two parameters: the
time delay $\tau$ and the embedding dimension $m$. For notational convenience, I denote the embedding space with dimension $m$ and time delay $\tau$ as \espace{m}{\tau}.
To assure topological
conjugacy, the Takens proof~\cite{takens} requires that the embedding dimension $m$ must be
at least twice the dimension $d$ of the ambient space; a tighter bound
of $m>2d_{cap}$, the capacity dimension of the original dynamics, was
later established by Sauer {\it et al.}~\cite{sauer91}.
\begin{mydef}[Capacity Dimension \cite{meissDynamcis}]
Let $N(\epsilon)$ denote the minimum number of open sets ($\epsilon$-balls) of diameter less than or equal to $\epsilon$ that form a finite cover of a compact metric space $X$. Then the \underline{capacity dimension} of $X$ is a real number $d_{cap}$ such that:  $N(\epsilon)\approx\epsilon^{-d_{cap}}$ as $\epsilon \rightarrow 0$, explicitly 
\begin{equation}d_{cap} \equiv -\lim_{\epsilon \rightarrow 0^+}\frac{\ln N(\epsilon)}{\ln \epsilon}
\end{equation} 
 if this limit exists. 
\end{mydef}
Operationalizing either of these theoretical constraints can be a real
challenge.  $d$ is not known and accurate $d_{cap}$ calculations are
not easy with experimental data.  And besides, one must first embed
the data before performing those calculations.  

Apropos of the central claim of this thesis, it is worth considering the intention behind these bounds on $m$. The worst-case bound of
$m>2d_{cap}$ is intended to eliminate
{\it all} projection-induced trajectory crossings in the reconstructed 
dynamics.  For most systems, and most projections, the dimensions of
the subspaces occupied by these false crossings are far smaller than
those of the original systems~\cite{sauer91}; often, they are sets of
\label{page:crossings-are-a-small-set}
measure zero.  For the delay-coordinate map to be a diffeomorphism,
all of these crossings must be unfolded by the embedding process.
This is necessary if one is interested in calculating dynamical
invariants like Lyapunov exponents.  However, the near-neighbor
relationships that most state-space forecast methods use in making
their predictions are {\it not} invariant under diffeomorphism, so it
does not make sense to place that strict condition on a model that
one is using for those purposes.  False crossings will, of course,
cause incorrect predictions, but that is not a serious problem in
practice if the measure of that set is near zero, particularly when
one is working with noisy, real-world data.

My reduced-order strategy explicitly fixes $m=2$. This choice takes care of one of the two free parameters in the delay-coordinate reconstruction process, but selection of a value for
the delay, $\tau$, is still an issue.  The theoretical constraints in
this regard are less stringent: $\tau$ must be greater than zero and
not a multiple of the period of any orbit~\cite{sauer91,takens}.  In
practice, however, the noisy and finite-precision nature of digital
data and floating-point arithmetic combine to make the choice of
$\tau$ much more delicate~\cite{kantz97}. It is to this issue that I will turn next.  

\subsection{Traditional Methods for Estimating the Embedding Delay $\tau$}\label{sec:numericaltau}

The $\tau$ parameter defines the amount of time separating each coordinate of the delay vectors: $\vec{x}_j=[x_j,~x_{j-\tau}, ~ \dots, ~ x_{j-(m-1)\tau}]^T$. 
The theoretical constraints on the time delay are far from stringent and this parameter does not---in theory \cite{sauer91,takens}---play  a role in the correctness of the embedding.  However, that assumes an infinitely long noise-free time series \cite{sauer91,takens}, a luxury that is rare in practice. As a result of this practical limitation, the time delay $\tau$ plays a crucial role in the {\it usefulness}\footnote{Here by {\it usefulness} I mean that not only are the dynamical invariants ({\it e.g.}, Lyapunov exponents and fractal dimension) and topological properties, ({\it e.g.}, neighborhood relations)  preserves, but also that those quantities are attainable from the reconstructed dynamics. 
} of the embedding \cite{fraser-swinney,kantz97,Buzug92Comp,liebert-wavering,Buzugfilldeform,Liebert89,rosenstein94}. 

The fact that the time delay does not play into the underlying mathematical framework is a double-edged sword. Because there are no theoretical constraints, there is no practical way to derive an ``optimal" lag or even know what criterion an ``optimal" lag would satisfy \cite{kantz97}. Casdagli {\it et al.} \cite{Casdagli:1991a} provide a theoretical discussion of this, together with some treatment of the impacts of $\tau$ on reconstructing an attractor using a noisy observation function.  Unfortunately no practical methods for estimating $\tau$ came from that discussion, but it does nicely outline a range of $\tau$ between {\it redundancy} and {\it irrelevance}.
For very small $\tau$, especially with noisy observations, $x_j$ and $x_{j-\tau}$ are effectively indistinguishable. 
 In this situation, the reconstruction coordinates are highly {\it redundant}~\cite{Casdagli:1991a,Gibson92}, {\it i.e.}, they contain nearly the same information about the system.\footnote{This is made more rigorous in Section~\ref{sec:ReviewInfoTheory}, where I discuss information theory.} This is not a good choice for  $\tau$ because additional coordinates add almost nothing new to the model. Choosing an arbitrarily {\it large} $\tau$ is undesirable as well. On this end of the spectrum, the coordinates of the reconstruction become causally unrelated, {\it i.e.}, the measurement of $x_{j-\tau}$ is {\it irrelevant} in understanding $x_j$ \cite{Casdagli:1991a}. Useful $\tau$ values lie somewhere between these two extrema. In practice, selecting useful $\tau$ values can be quite challenging, as demonstrated in the following example.

\begin{example}\label{ex:tau}
To explore the effects of $\tau$ on an embedding, I first construct an artificial time series by integrating the R\"{o}ssler system~\cite{rossler76}  
\begin{eqnarray}
 \dot{x} &=& -y-z  \\
  \dot{y} &=& x +ay \\
  \dot{z} &=& b+z(x-c)
\end{eqnarray}
with $a= 0.15$, $b = 0.20$, and $c = 10.0$, using a standard fourth-order Runge-Kutta integrator starting from $[x(0),y(0),z(0)]^T =[10,0,0]^T$ for 100,000 time steps with a time step of $\pi/100$.
This results in a trajectory of the form $\vec{y}(t_j)=[x(t_j), y(t_j), z(t_j)]^T$, where $t_j = j(\pi/100)$ for $j = 1,\dots,100,000$. This trajectory is plotted in Figure~\ref{fig:rossler}(a). To discard transient behavior, I remove the first 1,000 points of this trajectory. I define the observation function as $h(\vec{y}(t_j))=x(t_j)={x_j}$, resulting in the time series: $\{{x}_j\}_{j=1001}^{100,000}$. The first 5,000 points of this time series can be seen in Figure~\ref{fig:rossler}(b). 
\begin{figure}[ht!]
        \centering
                \vspace*{0.5cm}
        \begin{subfigure}[b]{0.5\textwidth}
                \includegraphics[width=\textwidth]{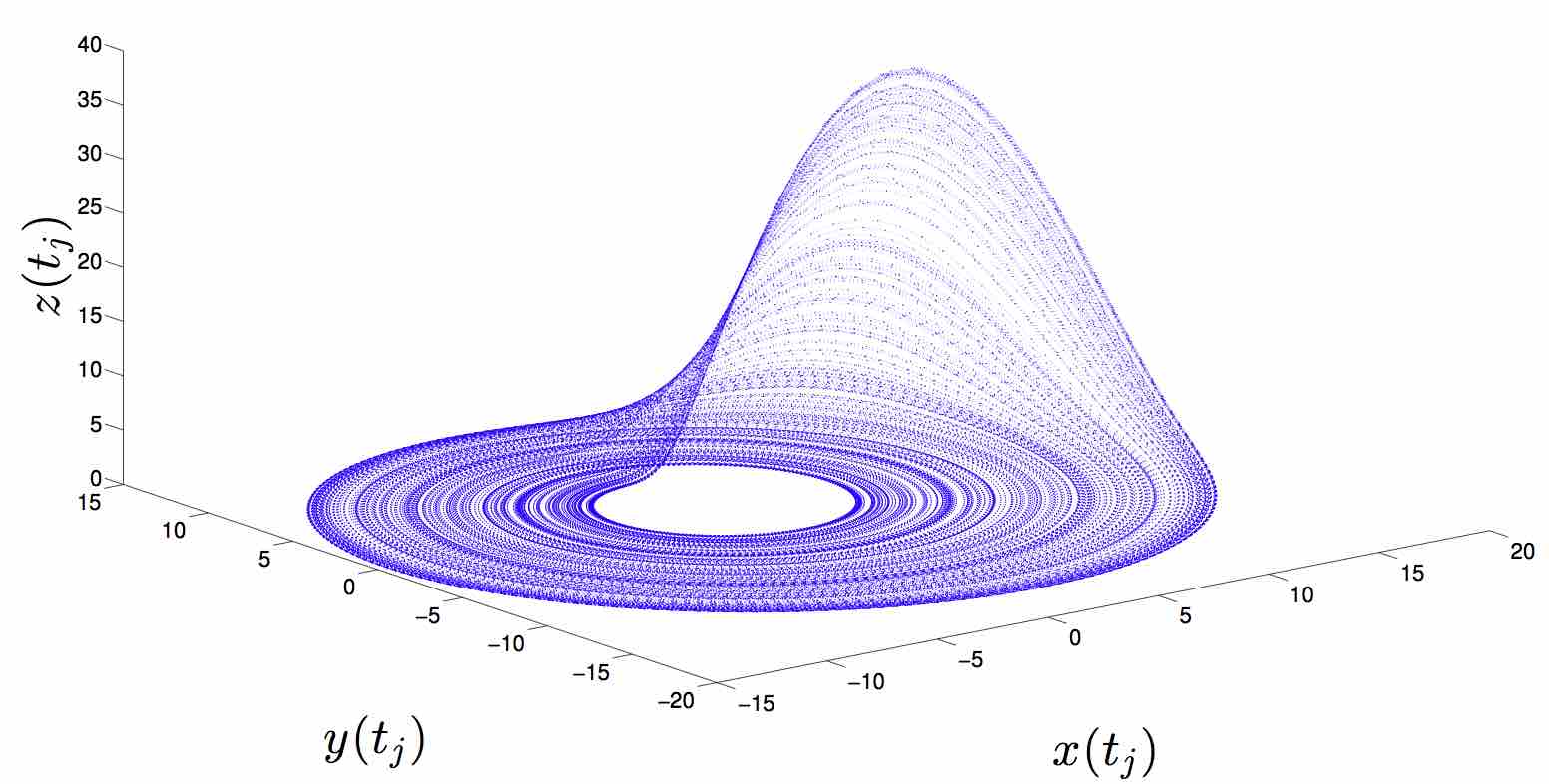}
                \caption{R\"ossler trajectory}
                \label{fig:rosslerAttractor}
        \end{subfigure}%
        ~ 
        \begin{subfigure}[b]{0.5\textwidth}
                \includegraphics[width=\textwidth]{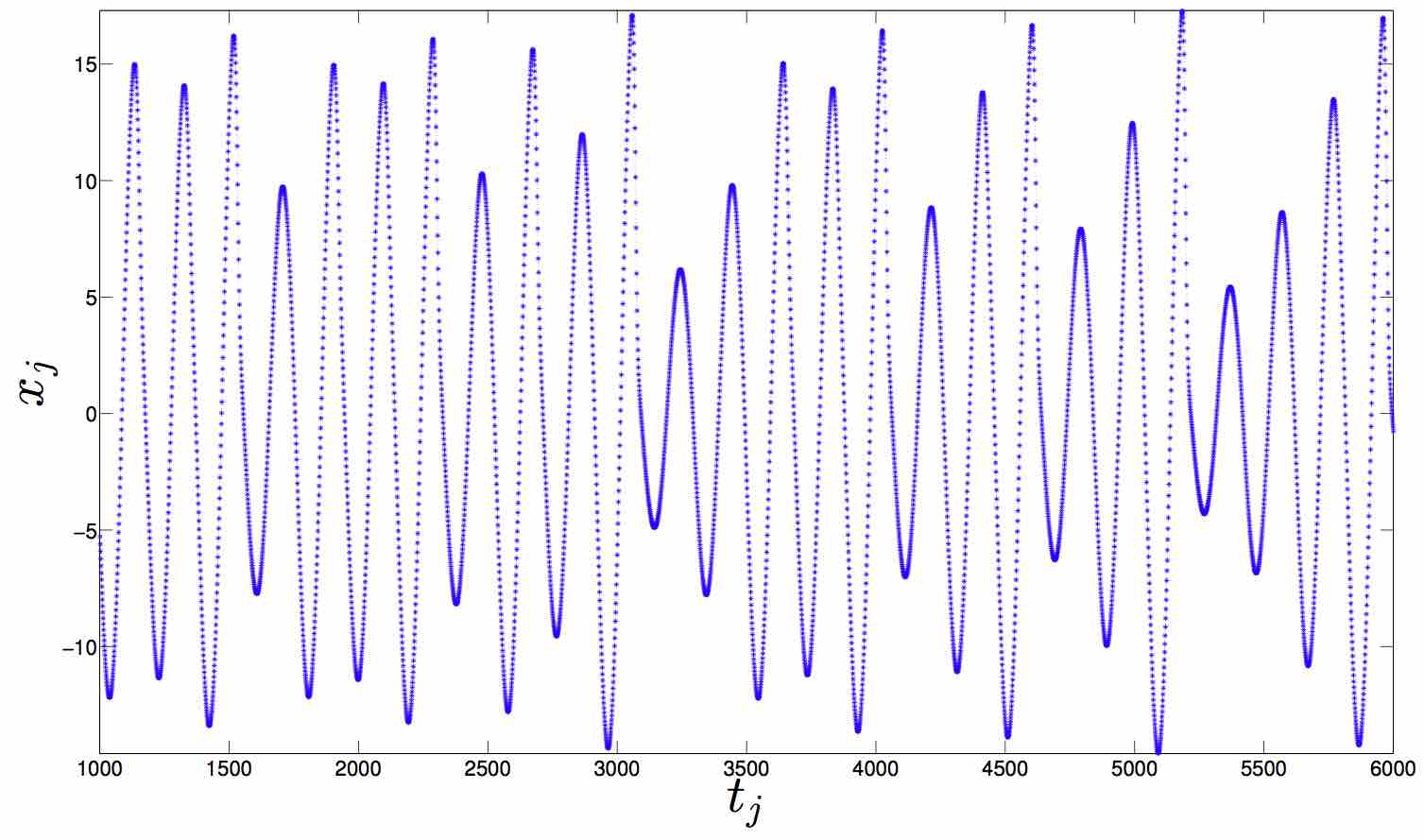}
                \caption{R\"ossler time series.}
                \label{fig:rosslerTS}
        \end{subfigure}
        \caption{The R\"ossler attractor and a segment of the time series of its $x$ coordinate.}\label{fig:rossler}
\end{figure}

\begin{figure}[ht!]
\begin{center}
                \vspace*{0.5cm}
 \includegraphics[width=\textwidth]{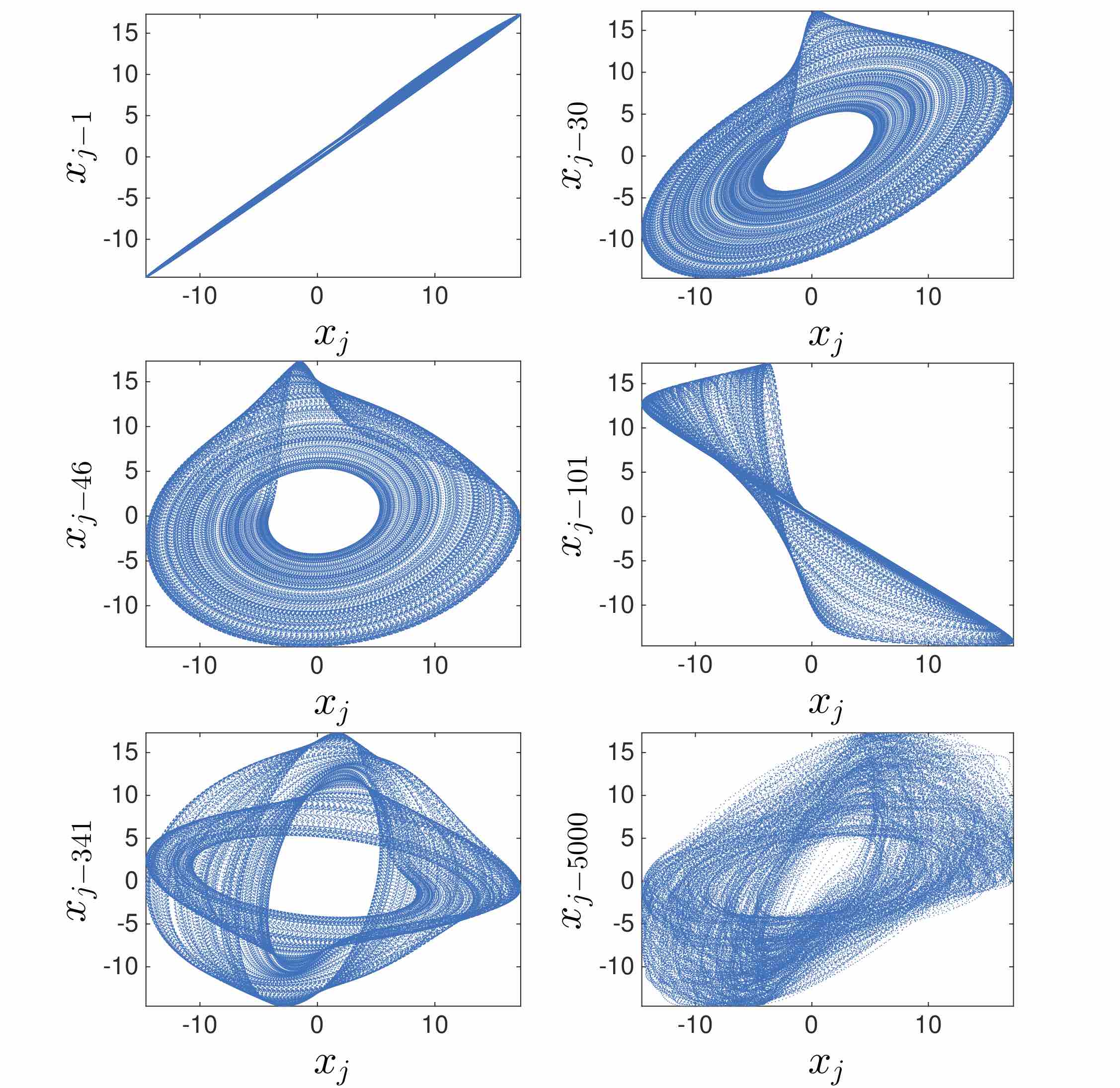} 
 
\caption{Delay-coordinate reconstruction of the R\"ossler time series in Figure~\ref{fig:rossler}(b) with $m=2$ and varying $\tau$.}\label{fig:timeDelayExample}
\end{center}
\end{figure}

To illustrate the role of $\tau$ in the delay-coordinate embedding process, I embed $\{{x}_j\}_{j=1001}^{100,000}$ using $m=2$ and several different choices of $\tau$. These embeddings are shown in Figure \ref{fig:timeDelayExample}. 
In theory, each of the choices of $\tau$ in Figure \ref{fig:timeDelayExample} should yield correct, topologically equivalent embeddings---given the right choice of $m$. In practice, however, that is not the case.

First consider the top-left panel of Figure~\ref{fig:timeDelayExample} where $\tau=1$. Here, the axes are spread apart so little that the embedding appears to be a noisy line. This is because $x_j$ and $x_{j-\tau}$ are effectively indistinguishable at this small $\tau$.
In the embedding in the bottom-right panel of Figure~\ref{fig:timeDelayExample}, 
 the reconstruction appears to be a ball of noise with only traces of underlying structure. At this large $\tau$, the coordinates of the reconstruction are causally unrelated.
This is known as ``overfolding." To visualize this concept, consider the progression in Figure~\ref{fig:timeDelayExample} from $\tau = 30$ to $\tau=341$. 
As $\tau$ increases, the reconstruction expands away from the diagonal and begins to resemble the original attractor. 
However, as $\tau$ increases past this point,
the top corner of the
reconstruction is slightly folded over.

This ``melting" effect is called folding in the literature.  ``Overfolding" occurs when the reconstructed attractor folds back on itself, completely collapsing back to the diagonal, (as can be seen for $\tau = 101$) and then re-expanding away from the diagonal, (as can be seen for $\tau=341$). Overfolding produces an unnecessarily (and in the case of noise, often incorrect) reconstruction \cite{Casdagli:1991a,rosenstein94,kantz97}. Compare, for example, the bottom-left panel of Figure~\ref{fig:timeDelayExample} with the actual attractor in Figure~\ref{fig:rossler}(a). From a theoretical standpoint, given the right choice of $m$, these two objects are topologically equivalent; from a practical standpoint, however, the embedding is overly complex. 

If the time series were noisy, this overfolding would likely introduce additional error. 
With knowledge of the true attractor, it is easy to say that the $\tau=30$ and $\tau=46$ embeddings most closely match its shape; without that knowledge, however, the choice is not obvious.
The situation is even more delicate than this. If one knew, somehow, that $\tau=30$ and $\tau=46$ were both good reconstructions, how would one know which of these two choices was optimal? With $\tau=30$, no folding has occurred, which is beneficial because with noisy or projected dynamics (choosing $m$ too small), foldings may cause {\it false crossings.}\footnote{A {\it false crossing} is when two trajectories intersect due to projection or measurement error, a phenomenon that cannot happen in a theoretical deterministic dynamical system.} But the trajectory in \espace{m}{30} is not as ``space filling" as $\tau=46$ or as spread apart from the diagonal, so the coordinates are most likely more redundant. 
Weighing the importance of these kinds of criteria is non-trivial and, I believe, application specific. In the rest of this section, I review heuristics aimed at optimizing the estimation of $\tau$ by weighing these different attributes against one another.
\end{example}

There are dozens of methods for estimating $\tau$---{\it e.g.}, \citeNTau. This is a central issue in my thesis, so the following section surveys this literature in some depth.
Choice of $\tau$ is application- and system-specific \cite{Buzug92Comp,kantz97,rosenstein94}; a $\tau$ that works well for Lyapunov exponent calculation may not work well for forecasting. For this reason, Kantz \& Schreiber \cite{kantz97} suggest that it may be necessary to perform additional system- and application-specific tuning of $\tau$ after using any generic selection heuristic. In my first set of examples, 
I use the method of {\it mutual information}~\cite{fraser-swinney, Liebert89}---described below in detail. While this is the standard $\tau$-selection method, I will show in Section~\ref{sec:varyingproj} that this choice is almost always {\it suboptimal} for forecasting. In Section~\ref{sec:tdAIS}, I provide a solution to this: an alternative $\tau$ selection method that leverages ``active information storage" to select a $\tau$ that is optimal for forecasting specific reconstructions~\cite{josh-tdAIS}. 

\subsubsection{Linear Independence and Autocorrelation}
A \naive strategy for selecting the time delay would be to choose a $\tau$ that forces the coordinates of the delay vectors to be {\it linearly} independent. This is equivalent to choosing the first zero of the autocorrelation function $R(\tau)$ %
\begin{equation}R(\tau) = \frac{1}{N-\tau}\frac{\sum_j(x_j-\mu_x)(x_{j-\tau}-\mu_x)}{\sigma_x^2}\end{equation}
where $N$, $\mu_x$ and $\sigma_x$ are respectively the length, average and standard deviation of the time series~\cite{kantz97,fraser-swinney}
. Several other methods have been proposed that suggest instead choosing $\tau$ where the autocorrelation function first drops to a particular fraction of its initial value, or at the first inflection point of that function \cite{kantz97,rosenstein94}. 

An advantage to this class of methods is that its members are extremely computationally efficient; the autocorrelation function, for instance, can be calculated with the fast Fourier transform~\cite{rosenstein94}. However, autocorrelation is a linear statistic that ignores nonlinear correlations. This often yields $\tau$ values that work well for some systems and not well for others~\cite{kantz97,rosenstein94,fraser-swinney}.

\subsubsection{General Independence and Mutual Information}\label{sec:mutual}

Instead of seeking linear independence between delay coordinates, it may be more appropriate to seek {\it general} independence---{\it i.e.}, coordinates that share the least amount of {\it information} (also called ``redundancy") with one another. The following discussion requires some methods from information theory; for a review of these concepts, please refer to Section~\ref{sec:ReviewInfoTheory}.
Fraser \& Swinney argue that selecting the first minimum of the {\it time-delayed mutual information} will minimize the redundancy of the embedding coordinates, maximizing the information content of the overall delay vector \cite{fraser-swinney}. In that approach, one obtains generally independent delay coordinates by symbolizing the two time series $X_j=\{x_j\}_{j=1}^N$ and $X_{j-\tau}=\{x_{j-\tau}\}_{j=1+\tau}^N$ by binning, discussed in Section~\ref{sec:binning}, and then computes the mutual information between $X_j$ and $X_{j-\tau}$ for a range of $\tau$, call this $I[X_j,X_{j-\tau}; \tau]$.  Then for each $\tau$, $I[X_j,X_{j-\tau}; \tau]$ is the amount of information shared between the coordinates $x_j$ and $x_{j-\tau}$, {\it i.e,} $I[X_j,X_{j-\tau}]$ quantifies how redundant the second axis is~\cite{fraser-swinney}.
According to \cite{fraser-swinney}, choosing a $\tau$ that minimizes $I[X_j,X_{j-\tau}]$, results in generally independent delay coordinates, {\it i.e,} delay coordinates that are minimally redundant.

The argument for choosing $\tau$ in this way applies {\it strictly} to two-dimensional embeddings~\cite{fraser-swinney,kantz97}, but was extended to work in $m$ dimensions in \cite{Liebert89}. To accomplish this, Liebert \& Schuster rewrote mutual information in terms of second-order Reny\'i entropies. This transformation allowed them to show that the minima of $I[X_j,X_{j-\tau}; \tau]$ agreed with the minima of the correlation sum \cite{GrassbergerPhysicaD}, $C(\epsilon;m,\tau)$, defined as
\begin{equation}C(\epsilon;m,\tau) = \frac{1}{N^2}\sum_{{\stackrel{i,j=1}{i\neq j}}}^N\Theta [\epsilon- ||\vec{x_i} - \vec{x_j}||]\end{equation}
where $N$ is the length of the time series, $\Theta$ is the Heavyside function, and $\vec{x_i},\vec{x_j}$ are the $i^{th}$ and $j^{th}$ delay vectors in \espace{m}{\tau}. 
In addition to extending the argument of \cite{fraser-swinney} to $m$ dimensions, the modification of \cite{Liebert89} allowed for much faster approximations of $\tau$ by simply finding the minimum of $C(\epsilon;m,\tau)$, which can be done quickly with the Grassberger-Procaccia algorithm \cite{GrassbergerPhysicaD,Liebert89}.

The choice of the {\it first} minimum of $I[X_j,X_{j-\tau}; \tau]$ is intended to avoid the kind of overfolding of the reconstructed attractor and irrelevance between coordinates that was demonstrated in Figure~\ref{fig:timeDelayExample}. 
This choice is discussed and empirically verified in \cite{Liebert89} by showing that the first minimum of $C(\epsilon;m,\tau)$ (so in turn $I[X_j,X_{j-\tau}; \tau]$) corresponded to the most reliable calculations of the {\it correlation dimension}~\cite{GrassbergerPhysicaD}.
\begin{mydef}[Correlation Dimension] 
If the correlation sum, $C(\epsilon)$,   
%
%
 decreases like a power law, $C(\epsilon)\sim\epsilon^D$, then $D$ is called the \underline{correlation dimension}. Formally 
 \begin{equation}D=\lim_{\epsilon\to 0} \frac{\log C(\epsilon)}{\log \epsilon}\end{equation}
if this limit exists. The Grassberger-Procaccia algorithm\footnote{The term Grassberger-Procaccia algorithm is used generically for any algorithm that estimates the correlation dimension (and more generally the correlation integral) from the small-$\epsilon$ behavior of the correlation sum $C(\epsilon;m,\tau)$.}\cite{GrassbergerPhysicaD} allows the correlation dimension to be approximated for \espace{m}{\tau} as
 \begin{equation}D=\lim_{\epsilon\to 0} {\frac{\log C(\epsilon;m,\tau)}{\log \epsilon}}\end{equation} 
\end{mydef}

It was {\it not} shown in \cite{Liebert89}, however, that this choice corresponds to the {\it best choice of $\tau$ for the purposes of forecasting}. In Section~\ref{sec:varyingproj} I show that this is in fact {\it not} the case for all time series. 
Even so, it is a reasonable starting point, as this method is the gold standard in the associated literature.  In my first round of experiments, and as a point of comparison, I select $\tau$ at the first minimum of the mutual information~\cite{fraser-swinney,Liebert89} as calculated by {\tt mutual} in the {\tt TISEAN} package~\cite{Hegger:1999yq}. 
There are a few possible drawbacks to this method. For example, there is no guarantee that $I[X_j,X_{j-\tau}; \tau]$ will ever achieve a minimum; a first-order autoregressive process, for example, does not~\cite{kantz97}. Rosenstein {\it et al.} \cite{rosenstein94} argue that calculating mutual information is too computationally expensive. Several papers~\cite{hasson2008influence,martinerie92} have argued that mutual information can give inconsistent results, especially with noisy data.

\subsubsection{Geometric and Topological Methods}

There are  several geometric and topological methods for approximating $\tau$ that address some of the shortcomings of mutual information, including: {\it wavering product} \cite{liebert-wavering}, {\it fill factor}, {\it integral local deformation}~\cite{Buzugfilldeform}, and {\it displacement from diagonal}~\cite{rosenstein94}, among others. Most of these methods have the distinct advantage of attempting to solve for both $m$ and $\tau$ simultaneously, albeit at the cost of being more complicated and less computationally efficient. (This additional computational overhead is not a factor in my reduced-order framework as I explicitly fix $m=2$.)

\subsubsection{Wavering Product}

The wavering product of Liebert {\it et al.} \cite{liebert-wavering} is a topological method for simultaneously determining embedding dimension and time delay. This approach focused on detecting when the attractor is properly unfolded, {\it i.e.}, the situation in which projection-induced overlap disappears.

Liebert {\it et al.} focused on preserving neighborhood relations of points in \espace{m}{\tau}. When transitioning from \espace{m}{\tau} to \espace{m+1}{\tau}, an embedding preserves neighborhood relations of every point in \espace{m}{\tau}, {\it i.e.,}  inner points remain inner points, and analogously with the boundary points. If these neighborhood relations are preserved, then $m$ is a sufficient embedding dimension. The so-called ``direction of projection"~\cite{liebert-wavering} that mitigates false crossings is associated with the best choice of $\tau$, {\it i.e.}, the $\tau$ that yields (for a fixed dimension) the smallest amount of overlap. To this end, they defined two quantities
\begin{equation} Q_1(i,k,m,\tau) = \frac{dist^{\tau}_{m+1}(i,j(k,m))}{dist^{\tau}_{m+1}(i,j(k,m+1))},\,\,\, Q_2(i,k,m,\tau) = \frac{dist^{\tau}_{m}(i,j(k,m))}{dist^{\tau}_{m}(i,j(k,m+1))}\end{equation}
 where, $dist^{\tau}_{m+1}(i,j(k,m))$ is the standard Euclidean distance measured in \espace{m+1}{\tau} between an $i^{th}$ reference point $\vec{x_i}$ in \espace{m}{\tau} and its $k^{th}$ nearest neighbor $\vec{x}_{j(k,m)}$ in \espace{m}{\tau} or similarly for $dist^{\tau}_{m+1}(i,j(k,m+1))$, the $k^{th}$ nearest neighbor of $\vec{x}_i$ in \espace{m+1}{\tau}. 
To determine if the neighborhood relations are preserved in the embedding, they defined the {\it wavering product}
\begin{equation}W_i(m,\tau) = \left(( \prod_{k=1}^{N_{nb}} Q_1(i,k,m,\tau) Q_2(i,k,m,\tau)\right)^{1/(2N_{nb})}\end{equation}
where $N_{nb}$ is the number of neighbors used in each neighborhood. If $W_i(m,\tau) \approx 1$, then the topological properties are preserved {\it locally} by the embedding~\cite{liebert-wavering}. In order to compute this {\it globally}, Liebert {\it et al.} defined the average wavering product as
\begin{equation}\overline{W}(m,\tau) = \ln\left[ \frac{1}{N_{ref}}\sum_{i=1}^{N_{ref}}W_i(m,\tau)\right]\end{equation}
where $N_{ref}$ is the number of 
reference points, typically chosen to be about 10\% of the signal. A minimum of $\overline{W}(m,\tau)/\tau$ as a function of $
\tau$ yields an optimal $\tau$ for that choice of $m$. They also showed that a 
sufficient embedding dimension can be found when $\overline{W}(m,
\tau)/\tau$ converged to zero. Choosing the 
embedding parameters in this way guarantees that the embedding faithfully preserves neighborhood 
relations. This is particularly important when forecasting based on neighbor relations.

This technique works very well on many systems, including the R\"ossler system and the Mackey-Glass system. In particular, Liebert {\it et al.}  showed that choosing $m$ and $\tau$ in this way allowed for accurate estimation of the information dimension~\cite{GrassbergerPhysicaD}. Noise, however, is a serious challenge for this heuristic~\cite{KBA92}, so it may not be useful for real-world datasets.

\subsubsection{Integral Local Deformation}
Integral local deformation, introduced by Buzug \& Pfister in \cite{Buzugfilldeform}, attempts to maintain continuity of the dynamics on the reconstructed attractor: {\it viz.,} neighboring trajectories remain close for small evolution times.  
The underlying rationale is that false crossings will cause what look like neighboring trajectories to separate exponentially in very short evolution time. 
Integral local deformation quantifies this. Buzug \& Pfister show that choosing $m$ and $\tau$ to minimize this quantity gives an approximation of $\tau$ that minimizes false crossings created by projection.

In my work, I rely strongly on the continuity of the reconstructed dynamics, since I use the image of neighboring trajectories for forecasting. Integral local deformation seems useful at first glance for choosing a $\tau$ that helps to preserve the continuity of the underlying dynamics in the face of projection. However, the computational complexity of this measure makes it ineffective for on-the-fly adaptation or selection of $\tau$.

\subsubsection{Fill Factor}

In \cite{Buzugfilldeform}, Buzug \& Pfister introduced a purely geometric heuristic for estimating $\tau$. This method attempts to maximally fill the embedding space by {\it spatially} spreading out the points as far as possible. To accomplish this, Buzug \& Pfister calculate the average volume of a large number of $m$-dimensional parallelepipeds, spanned by a set of $m+1$ arbitrarily chosen $m$-dimensional delay vectors.
They then show that the first maximum of the average of these volumes as a function of $\tau$ (for a fixed $m$) maximizes the distance between trajectories. This method is computationally efficient, as no near-neighbor searching is required. However, for any attractor with multiple unstable foci, there is no significant maximum of the fill factor as a function of $\tau$ \cite{rosenstein94,Buzugfilldeform}. In addition, this method cannot take into account overfolding, as an overfolded embedding may be more space-filling than the ``properly" unfolded counterpart~\cite{rosenstein94}. This consideration is corrected (at the cost of additional computational complexity) in the method described next. 

\subsubsection{Average Displacement / Displacement from Diagonal}

The average displacement method introduced by Rosenstein {\it et al.} \cite{rosenstein94}, which is also known as the displacement from diagonal method \cite{kantz97}, also seeks a $\tau$ that causes the embedded attractor to fill the space as much as possible, while mitigating error caused by overfolding and also addressing some other concerns  \cite{Buzugfilldeform}. Rosenstein {\it et al.} define the average displacement (from diagonal) for  \espace{m}{\tau} as
 \begin{equation}\left<S_m(\tau)\right> = \frac{1}{N}\sum_{i=1}^N\sqrt{\sum_{j=1}^{m-1}(x_{i+j\tau}-x_i)^2}\end{equation}
For a fixed $m$, $\left<S_m(\tau)\right>$ increases with increasing $\tau$ (at least initially; the attractor may collapse for large $\tau$ due to overfolding).  Rosenstein {\it et al.} suggest choosing $\tau$ and $m$ where the slope between successive $\left<S_m(\tau_i)\right>$ drops to around 40\% of the slope between $\left<S_m(\tau_1)\right>$ and $\left<S_m(\tau_2)\right>$, where $\tau_1$ and $\tau_2$ are the first and second choices of $\tau$. In noisy data sets, this leads to consistent and accurate computation of the correlation dimension. However, this---like most heuristics---was developed to correctly approximate dynamical invariants ({\it e.g.,} correlation dimension), and comes with no guarantees about forecast accuracy.

\begin{remark}
Several papers ({\it e.g.}, \cite{Kugi96,McNames98anearest,rosenstein94,Small2004283,PhysRevE.87.022905}) have claimed that the emphasis should be placed on the window size $\tau_w = \tau m$ rather than $\tau$ or $m$ independently. The basic premise behind this idea is that it is more important to choose $\tau_w$ to span an important time segment ({\it e.g.}, mean orbital period) than the actual choice of either $\tau$ or $m$ independently. This is something I have not found to be the case when choosing parameters for delay reconstruction-based forecasting. 
\end{remark}

\subsection{Traditional Methods for Estimating the Embedding Dimension $m$}\label{sec:numericalM}

As the embedding dimension $m$ is not a parameter in my reduced-order algorithm, I only review a few important methods for estimating it.  This discussion is important mainly because these conventions are the point of departure (and comparison) for my work.

\begin{figure}[tb!]
        \centering
        \begin{subfigure}[b]{0.5\textwidth}
                \includegraphics[width=\textwidth]{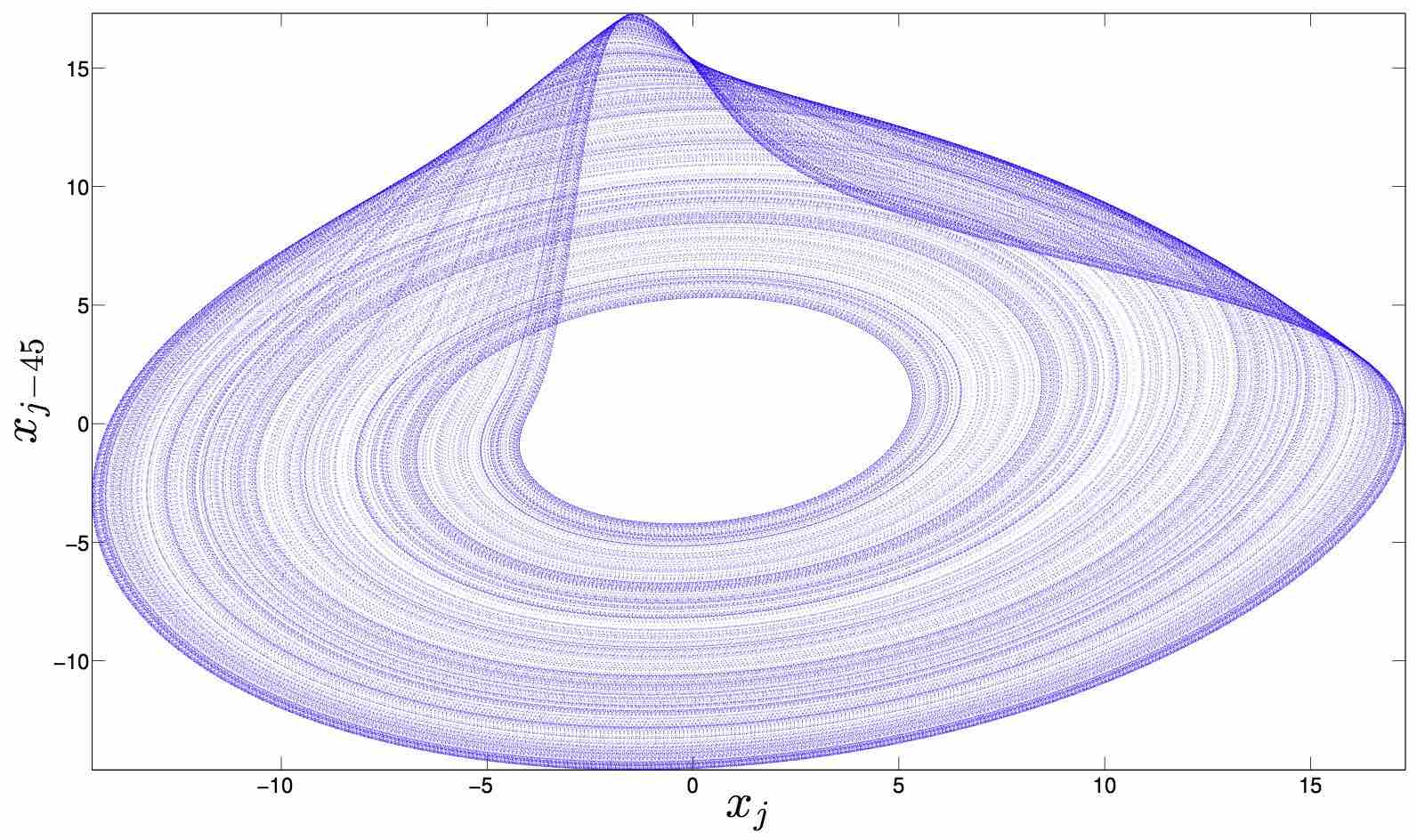}
                \caption{R\"ossler time series in \espace{2}{45}}
                \label{fig:rossler2D}
        \end{subfigure}%
        ~ 
        \begin{subfigure}[b]{0.5\textwidth}
                \includegraphics[width=\textwidth]{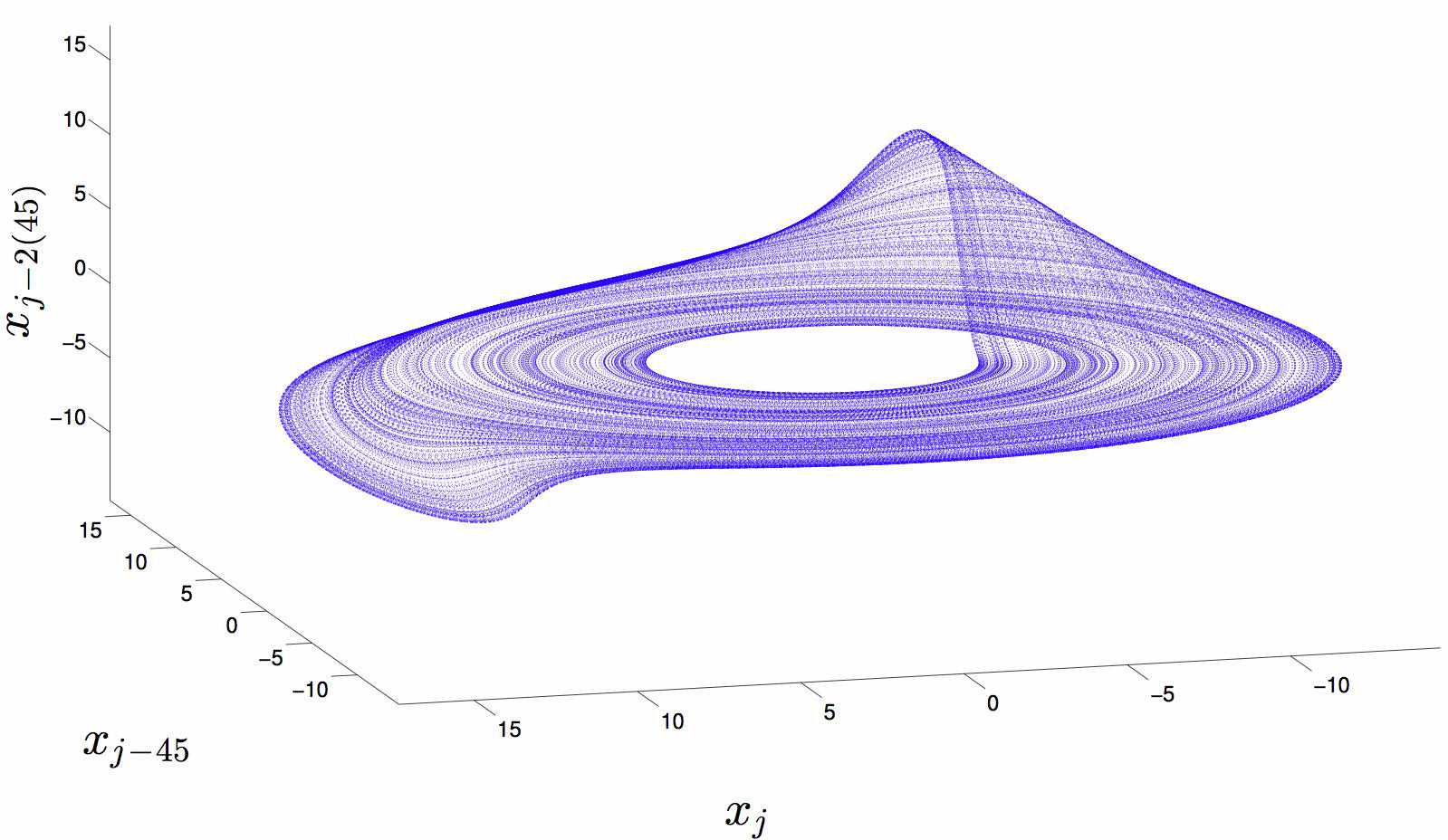}
                \caption{R\"ossler time series in \espace{3}{45}}
                \label{fig:rossler3D}
        \end{subfigure}
        \caption{An illustration of the utility of higher embedding dimensions to eliminate false crossings in the dynamics.}\label{fig:2to3D}
\end{figure}

A scalar time series $\{x_j\}_{j=1}^N$ measured from a dynamical system is a projection of the original state space onto a one-dimensional sub-manifold. 
A fundamental concern in the theoretical embedding dimension requirement $m>2d_{cap}$  is to ensure that the embedding has enough dimensions to ``spread out" and thus avoid false crossings. Such crossings violate several properties of deterministic dynamical systems, {\it e.g.,} determinism, uniqueness and continuity. In Figure~\ref{fig:2to3D}(a), for example, the \espace{2}{45} embedding of the $x$-coordinate of the R\"ossler system contains trajectory crossings.  In Figure~\ref{fig:2to3D}(b), however, the top-right region of Figure~\ref{fig:2to3D}(a) has folded under the attractor and the intersections on the top left of Figure~\ref{fig:2to3D} have become a ``tunnel."
The issue here is that the dynamics do not have enough space to spread apart in two dimensions. However when this dimension is increased, the attractor can spread out and the intersections disappear.  According to \cite{sauer91}, choosing $m>2d_{cap}$ ensures that the attractor has enough space to spread out and false crossings will not occur. More precisely, the probability of a crossing occurring in a ball of size $\epsilon$ is $p_c\propto \epsilon^{m-2d_{cap}}$. Recall from Section~\ref{sec:dce}, however, that this is for an infinitely long noise free time series and may not hold in practice, as noise can easily cause false crossings and violate the assumptions that went into this estimation.
It should also be noted that, even if I knew the capacity dimension $d_{cap}$ of the system---which is generally not the case---I do not necessarily want to choose $m$ to be $2d_{cap}+1$. This is a (generally loose) sufficient bound that should ensure the correctness of the embedding. But it is often the case that the embedding unfolds completely before $m=2d_{cap}+1$. The Lorenz system~\cite{lorenz} has $d_{cap} = 2.06\pm 0.01$, for example. \cite{sauer91} would suggest using $m=5$, but in fact this system can be embedded properly using $m=3$ \cite{KBA92}.  

Na\"ively, it may seem that simply choosing an ``extremely large" $m$ would be a simpler and completely reasonable choice. This is not true, in practice. First, the complexity of many of the algorithms that deduce information about dynamical systems scale exponentially with $m$~\cite{Liebert89}. Worse yet, each noisy data point creates $m$ noisy embedding points in the reconstruction ~\cite{Casdagli:1991a}. This amplification of noise quickly destroys the usefulness of an embedding. In light of both of these concerns, good values for the {\it minimal} $m$ are highly sought after. For a noisy real-valued time series, this is still an open problem, but there exist several heuristic approximations ({\it e.g.}, \citeNM).  Recall, too, that several of the methods presented in the previous section for estimating $\tau$---{\it e.g.},{\it wavering product}~\cite{liebert-wavering} and {\it integral local deformation} \cite{Buzugfilldeform}---simultaneously estimate both the delay and the dimension, $m$---the other free parameter in the embedding process.

There are two standard classes of methods for estimating the minimal $m$, the {\it method of dynamical invariants} and the {\it method of false neighbors}. In the following sections, I review the basics of these two families.

\subsubsection{Method of Dynamical Invariants}

Dynamical invariants, such as correlation dimension, are topological measures of a system that persist under diffeomorphism. In theory, this means that once a particular choice of embedding dimension, say $\hat{m}$, yields a topologically valid reconstruction, increasing $m$ should have no impact on these dynamical invariants. This is the case because in theory every \espace{m>\hat{m}}{\tau} will be topologically conjugate, to one another and to the original dynamics. This implies that dynamical invariants will become {\it persistent} for increasing $m$, once $\hat{m}$ has been reached. Hence, choosing the first $m$ for which dynamical invariants stop changing is a good way to estimate the minimal dimension needed to obtain a topologically valid reconstruction.
The class of methods that is the topic of this section follows directly from this logic:
to choose $m$, one approximates some dynamical invariant ({\it e.g.}, dominant Lyapunov exponent or correlation dimension) for a range of embedding dimensions, choosing the first embedding dimension for which it becomes persistent, and then corroborates with other dynamical invariants.

For example, one can approximate the correlation dimension for a range of embedding dimensions using the Grassberger \& Procaccia algorithm \cite{GrassbergerPhysicaD}, choosing the first $m$ for which that approximation stops changing. Then one corroborates this choice by approximating the dominant Lyapunov exponent for a range of $m$ (using for example the algorithm in \cite{wolf}), then choosing the first $m$ where this result stops changing. Finally, one ensures these two estimates of $m$ are consistent with each other.

Recall, though, that noise in the data can impact any dynamical invariant calculation, and that that impact increases with $m$~\cite{Casdagli:1991a}. It is more often the case that there is a range of embedding dimensions for which the dynamical invariant being approximated stays ``fairly consistent." Ascertaining this is computationally expensive and requires time-intensive post processing and human interpretation. For these reasons, it is common to use an alternative heuristic, such as those covered in the next section, to narrow down the search to a smaller range of embedding dimensions and then select from this range using the method of dynamical invariants.


\subsubsection{The Method of False Neighbors}

The method of false neighbors was proposed by Kennel {\it et al.} in  \cite{KBA92}. This heuristic searches for points that appear close only because the embedding dimension is too small. Consider a 
point on the top of the tunnel in Figure~\ref{fig:2to3D}(b) and a point directly below this 
point on the planar part of the R\"ossler attractor. These two points are near neighbors in \espace{2}{45} 
because the tunnel collapses down on the planar region; however, they are not near neighbors in \espace{3}{45} 
because the embedded attractor inflates, separating points on the 
top of the tunnel from the points on the planar region. 
Since these two points are neighbors in  \espace{2}{45} and {\it not} neighbors in \espace{3}{45}, they are {\it false near(est) neighbors} at $m=2$.
Consider, in contrast, two neighboring points on the top of the tunnel in \espace{3}{45}. If the space is projected down to \espace{2}{45}, these points would still be neighbors.  

More formally, Kennel {\it et al.} define the $k^{th}$ nearest neighbor $\vec{x}_{j(k,m)}\in$ \espace{m}{\tau} of $\vec{x}_i\in$ \espace{m}{\tau} as a {\it false near(est) neighbor} if
\begin{equation}\left(\frac{ dist^{\tau}_{m+1}(i,j(k,m))^2-dist^{\tau}_{m}(i,j(k,m))^2}{dist^{\tau}_{m}(i,j(k,m))^2}\right)^{1/2}  > R_{tol}\end{equation}
where $R_{tol}$ is some tolerance. Recall that the $dist^{\tau}_{m}(i,j(k,m))$ is the distance between the $i^{th}$  point $\vec{x}_i\in$ 
\espace{m}{\tau} and its $k^{th}$ nearest neighbor $\vec{x}_{j(k,m)}\in$ \espace{m}{\tau}. But notice for delay vectors that $dist^{\tau}_{m+1}(i,j(k,m))^2 = |x_{i-m\tau} -  x_{j(k,m)-m\tau}|^2+ dist^{\tau}_{m}(i,j(k,m))^2$, so this condition simplifies to
\begin{equation}\label{eq:fnn}\frac{|x_{i-m\tau} -  x_{j(k,m)-m\tau}|}{dist^{\tau}_{m}(i,j(k,m))}>R_{tol}\end{equation}
In particular, a neighbor is a false neighbor if the distance between the two points in \espace{m+1}{\tau} is significantly more ({\it viz.,} $R_{tol}$) than the distance between the two neighbors in \espace{m}{\tau}. Kennel {\it et al.} claim that choosing a single nearest neighbor is sufficient ({\it i.e.}, $k=1$) \cite{KBA92}.
  In addition, they claim that empirically $R_{tol}\geq 10$ seems to give robust results. This tolerance can be interpreted as defining false neighbors as points that are 10 times farther apart in \espace{m+1}{\tau} than in \espace{m}{\tau}.


This heuristic alone is not enough to distinguish chaos from uniform noise and can incorrectly classify time series constructed from a uniform distribution as having low-dimensional dynamics. Kennel {\it et al.} found that for a uniformly distributed random time series, on average, the nearest neighbor of a point is not near at all. Rather, $dist^{\tau}_{m}(i,j(k,m))\approx R_A$, where $R_A = \sqrt{\frac{1}{N}\sum_{j=1}^N(x_j - \mu_x)^2}$. That is, the average distance to the nearest neighbor is the size of the attractor.
To handle this, they defined a secondary heuristic $\frac{dist^{\tau}_{m+1}(i,j(k,m))}{R_A}>A_{tol}$, where $A_{tol}$ is another free parameter chosen as 2.0 without justification. I want to note that this heuristic is added to distinguish pure-uniform noise from chaotic dynamics, {\it not} to aid in estimating embedding dimension for noisy observations of a chaotic system. 

For a time series with noise, near-neighbor relations---which are the basis for this class of heuristics---can cause serious problems in practice. For well-sampled, noise-free data, it makes sense to choose $m$ as the first embedding dimension for which the ratio of true to false neighbors goes to zero \cite{KBA92}. For noisy data however, this is unrealistic; in practice, the standard approach is to choose the first $m$ for which the percentage of false near(est) neighbors drops below 10\%. If topological correctness is vitally important for the application, a range of embedding dimensions for which the percentage of false near(est) neighbors drops to around $\approx10$\% is typically chosen and then this range is refined using the method of dynamical invariants described above. 
This 10\% is an arbitrary threshold, however; depending on the magnitude of noise present in the data, it may need to be adjusted, as may $R_{tol}$ and $A_{tol}$. For example, in the computer performance data presented in Section \ref{sec:compPerfExperiments}, the percentage of false near(est) neighbors rarely dropped below even 20\%.

Recently an extension of the false near(est) neighbor  method was proposed by Cao \cite{Cao97Embed}, which attempts to get around the three different tolerances ($R_{tol}$, $A_{tol}$ and the percentage of false neighbors) in \cite{KBA92}. Cao points out that the tolerances---in particular, $R_{tol}$---need to be specified on a per-time-series and even per-dimension basis. Assigning these tolerances universally is inadvisable and in many cases will lead to inconsistent estimates. In \cite{Cao97Embed}, he illustrates that different choices of these three tolerances result in very different estimates for $m$.  To get around this, he defines an alternative heuristic that is tolerance free 
\begin{equation}E(m)=\frac{1}{N-m\tau}\sum_{i=1}^{N-m\tau}\frac{dist^{\tau}_{m+1}(i,j(k,m))}{dist^{\tau}_{m}(i,j(k,m))}\end{equation}
 While the inside of the sum is very similar to Equation \ref{eq:fnn} of \cite{KBA92}, the numerator is slightly different: $dist^{\tau}_{m+1}(i,j(k,m))$ instead of $|x_{i-m\tau} -  x_{j(k,m)-m\tau}|$. That is, the former measures the distance between an element and its $k^{th}$-nearest neighbor in \espace{m}{\tau} measured in \espace{m+1}{\tau}, whereas the latter measures the {\it change} in distance between $\vec{x}_i, \vec{x}_{j(k,m)}\in$ \espace{m}{\tau}, and the same vectors extended in \espace{m+1}{\tau}. Cao then defines $E1(m) =\frac{E(m+1)}{E(m)}$ and shows that when $E1(m)$ stops changing, a sufficient embedding dimension has been found. He also claims that if $E1(m)$ does not stop changing, then one is observing noise and not deterministic dynamics~\cite{Cao97Embed}. 
Cao does admit that it is sometimes hard to determine if the $E1(m)$ curve is just slowly growing but will plateau eventually (in the case of high dimensional dynamics) or just constantly growing (in the case of noise). To deal with this, he defines a secondary heuristic to help distinguish these two cases. As this method has been shown to give more consistent $m$, I hoped that this method could provide a more accurate comparison point. However, I was never able to successfully  replicate the results in \cite{Cao97Embed} on any experimental data, so I chose to use the traditional version of this algorithm proposed by Kennel {\it et al.} in  \cite{KBA92}.

The astute reader may have noticed a similarity between the method of false neighbors \cite{KBA92}, the method of wavering products ~\cite{liebert-wavering}, and the methods of Cao~\cite{Cao97Embed}. It is true that these methods are quite similar. In fact, almost all the methods for determining minimum embedding dimension \citeNM are based in some way on minimizing the number of false crossings. 
As this parameter is not important in my work, I do not go into all of these nuances but simply use the standard false near(est) neighbor approach to which the rest of these methods are fundamentally related. In particular, I use the {\tt TISEAN}~\cite{Hegger:1999yq} implementation of this algorithm (\verb|false_nearest|) to choose $m$ with a $\approx20\%$ threshold on the percentage of neighbors and the $R_{tol}$ and $A_{tol}$ selected by the {\tt TISEAN} implementation. In my later discussion, I refer to the reconstruction produced in this manner as an embedding of the data. This is by no means perfect, but since it is the most widely used method for estimating $m$, it is the most useful for the purposes of comparison.

\subsection{Delay-Coordinate Embedding Reality Check}\label{sec:predinproj}

As discussed in Section~\ref{sec:dce}, the theory of delay-coordinate embedding~\cite{sauer91,takens} outlines beautiful machinery to reconstruct---up to diffeomorphism---the dynamics of a system from a scalar observation. 
Unfortunately this theoretical machinery requires both infinitely-long and noise-free observations of the dynamical system: luxuries that are never afforded to a time-series analyst in practice. 
While there has been a tremendous amount of informative literature on estimating the free parameters of delay-coordinate embedding, 
at the end of the day these heuristics are just that: empirical estimates with no theoretical guarantees. 
This means that, even if the most careful rigorous in-depth analysis is used to estimate $\tau$ and $m$, there is no way to guarantee, in the experimental context, that the reconstructed dynamics are in fact diffeomorphic to the observed dynamics.

Even worse, overestimating $m$ has drastic impacts on the usefulness of the embedding, as it exponentially amplifies the noise present in the reconstruction.  
If little usable structure is present in a time series in the first place, perverting this precious structure by amplifying noise is something a practitioner can ill afford to do. Moreover, the methods that are most commonly used for estimating $m$ are based on neighbor relations, which are easily corrupted by noisy data. As a result, these heuristics tend to overestimate $m$.

In addition to noise amplification concerns and the lack of theoretical guarantees, the methods for estimating minimal embedding dimension are highly subjective, dependent on the estimate of $\tau$, and require a great deal of human intuition to interpret correctly. This time-consuming, error-prone human-intensive process makes it effectively impossible to use delay-coordinate embedding for automated or `on-the-fly' forecasting. As stated in Chapter~\ref{ch:overview}, this is unfortunate because delay-coordinate embedding is such a powerful modeling framework. My reduced-order framework---the foundation of this thesis---will, I hope, at least partially rectify this shortcoming.

\section{Information Theory Primer}\label{sec:ReviewInfoTheory}

In this section, I provide a basic overview of notation and concepts from Shannon information theory~\cite{shannon64}, as well as a review of some more-advanced topics that are utilized throughout the thesis. I will first cover the basics; an expert in this field can easily skip this part. I will then move on to non-traditional topics {\it viz.,} multivariate information theory (Section~\ref{sec:multiinfoandidiagrams}), methods for computing information measures on real-valued time series (Section~\ref{sec:symbolize}), and measures to quantify the predictability of a real-valued time series (Section~\ref{sec:meaComplex-intro}). 

\subsection{Entropy}\label{sec:entropyprimer}
Perhaps the most fundamental concept or building block in information theory is the concept of Shannon Entropy.
\begin{mydef}[(Shannon) Entropy \cite{shannon64}]
Let $Q$ be a discrete random variable with support $\{q_1,...,q_n\}$ and a probability mass function $p$ that maps a possible symbol to the probability of that symbol occurring, {\it e.g.}, $p(q_i) = p_i$, where $p_i$ is the probability that an observation $q$ is measured to be $q_i$. 
The average amount of information gained by taking a measurement of $Q$ and thereby specifying an observation $q$ is the \underline{Shannon Entropy} (or simply entropy) $H$ of $Q$, defined by
\begin{equation}H[Q] = - \sum_{i=1}^np(q_i)\log(p(q_i))\label{eq:entropy}\end{equation}
\end{mydef}
\noindent Throughout this thesis, $\log$ is calculated with base two, so that the information is in bits. The entropy $H[Q]$ can be interpreted as the amount of ``surprise" in observing a measurement of a discrete random variable $Q$, or equivalently the average uncertainty in the outcome of a process, or the amount of ``information" in each observation of a process.  
\begin{example}[Entropy of fair and biased coins]
First consider a fair coin: $Q = \{h,t\}$ and $p(h) = p(t) = 1/2$. 
\begin{eqnarray}
H[Q] &=& -[p(h)\log(p(h)) + p(t)\log(p(t))]   \\ 
  &=& - [\frac{1}{2}\log(\frac{1}{2}) + \frac{1}{2}\log(\frac{1}{2})] \\
  &=& - (\frac{-1}{2} +\frac{-1}{2}) = 1
\end{eqnarray} 
At every flip of the coin there is one bit of ``new" information, or one bit of surprise. In contrast, consider an extremely biased coin with heads on both sides: $Q = \{h,t\}$ and $p(h) = 1$, $p(t) = 0$. Then $H[Q] = 0$, {\it i.e.}, there are zero bits of ``new" information at each toss, as the coin always gives heads. 
 \end{example}
To gain an intuitive understanding of what phrases like `one bit of ``new" information', or `one bit of surprise' mean, it is sometimes easier to interpret Equation~(\ref{eq:entropy}) as the average number of (optimal) yes-no questions one needs to ask in order to determine what the outcome of observing a system will be. Returning to the coin-flip example above, since the fair coin had $H[Q]=1$, on average, one (optimal) question needs to be asked to determine the outcome of the coin flip: ``Was the coin a head?" With the biased coin, however, the entropy was zero, which means on average no questions were needed in order to infer the observation was a head (it always is!). 
The following example clarifies this.
 \begin{example}[Entropy of Animal-Vegetable-Mineral~\cite{simon-it}]
You may have, at some point during your childhood, played the game
``Animal-Vegetable-Mineral." If not, the rules are simple: player one thinks of an object, and by a series of yes-no questions, and
the other players attempt to guess the object. ``Is it bigger than a breadbox?" No. ``Does it have fur?" Yes.
``Is it a mammal?" No. This continues until the players can guess the object.

As anyone who has played this game will attest, some questions are better than others---for example, you usually try to focus on general
categories first (hence, the name of the game itself---is it an animal?) then get more specific within that category.
Asking on the first round ``is it a dissertation?" is likely to waste time---unless, perhaps, you are
playing with a graduate student who is about to defend.

If a game lasts too long, you may begin to wonder if there exists an optimal set of questions to
ask. ``Could I have gotten the answer sooner, if I had skipped that useless question about the fur?"
A moment's reflection shows that, in fact, the optimal set of questions depends upon the player: if
someone is biased towards automobiles, it would be sensible to focus on questions that specify
make, model, year, etc. You could then imagine writing down a script for this player: ``first ask if it is a car;
then if yes, ask if it is domestic, if no, ask if it is a Honda...;" or for my nieces: ``first ask if it's Elsa
from Frozen." (It almost always is.)

For each player and their preferences ({\it i.e.,} for every probability distribution over the things that player
might choose), there is an optimal script. And for
each optimal script for a given person, the game will last five rounds, or ten rounds, or
seven, or twenty, depending on what they choose that time. Profoundly, the number of questions you have to ask on average for a particular person and
optimal script pair is given by Equation~(\ref{eq:entropy}). In particular, we are measuring information (and uncertainty): the
average number of yes-no questions we'll need to ask to find out an answer.
 \end{example}
 
  Understanding entropy is important to the rest of the discussion in this chapter as it is the fundamental building block of all other information-theoretic quantities. 
 \subsection{Mutual Information}\label{sec:intromi}
It is often interesting to consider how knowledge about something informs us about something
else. People carrying umbrellas, for example, tells us something about the weather; it is not perfect
but (informally) if you tell me something about the weather, you also reduce my uncertainty about
umbrella-carrying~\cite{simon-it}. To constructively introduce the so-called {\it mutual information}, I will adapt the next example from \cite{simon-it}.
 
How can we quantify the information between the weather $W$ and umbrella-carrying $U$? For simplicity, I will assume you only get to see one
person---who is either carrying ($u_1$) or not carrying ($u_2$) an umbrella, {\it i.e.,} $U=\{u_1,u_2\}$ with probability $p(u_i)$. 
 
Now assume that there are some finite number of weather types (say ``rain", ``cloudy", ``windy" etc., labeled with $w_j$), each with a probability of occurring, $p(w_j)$. Then from Section~\ref{sec:entropyprimer}, the uncertainty in the weather is simply 
\begin{equation}H[W] = -\sum_{i=1}^Np(w_j)\log(p(w_j))\end{equation}
 We are interested in the probability of seeing a particular weather type given that we see the person carrying an umbrella. For this, consider the conditional probability of weather type $i$ given that you see someone carrying
an umbrella---$p(w_i|u_1)$. Generally, $p(w_i|u_1)$ will be higher than $p(w_j|u_1)$ when $i$ is labeling weather with precipitation and $j$ is not, 
so the uncertainty of the weather given that someone is carrying an umbrella is then  
\begin{equation}H[W|u_1] = -p(u_1)\sum_{j}p(w_j|u_1)\log(p(w_j|u_1))\end{equation}
 or in words, ``the uncertainty about the weather, given that the person who walked in was carrying an umbrella." Similarly, for the reverse
case, we could compute the associated uncertainty $H[W|u_2]$,  to determine ``the uncertainty about the weather, given that the person who walked in was not carrying an umbrella." Combining (summing) these two we get the  {\it conditional entropy} between two variables (weather type and state of umbrella-carrying in this example).
\begin{mydef}[Conditional Entropy \cite{shannon64}]
Define $Q$ and $R$ to be discrete random variables with support $\{q_1,...,q_n\}$ and $\{r_1,...,r_m\}$ respectively. Then the \underline{conditional entropy} is defined as \begin{equation}H[Q|R] = -\sum_ip(r_i)\sum_{j}p(q_j|r_i)\log(p(q_j|r_i))\end{equation} 
 where $p(q_j|r_i)$ is the conditional probability of $q_j$ given $r_i$.
\end{mydef}

We can then quantify the ``reduction in uncertainty" in the weather given that someone is carrying an umbrella by $H[W]- H[W|u_1]$, and the reverse with $H[W]- H[W|u_2]$. Note that the reduction can be positive or negative---in some climates, seeing your colleague not
carrying an umbrella will make you more uncertain about the weather. Consider, for example, an
extremely rainy climate; it is either sunny, cloudy, or rainy, but most often rainy. You are generally
quite certain about the weather before you see your colleague (it is raining). So when they walk
through the door without their umbrella, you think it is less likely to be raining, and so you are more
uncertain (the options sunny, cloudy, or rainy are now more evenly balanced).

Now consider the ``average reduction in uncertainty" of the weather given the state of umbrella carrying
\begin{eqnarray}
I[W,U] &=& p(u_1)(H[W] - H[W|u_1]) +p(u_2)(H[W] - H[W|u_2])  \\
&=& H[W] - (p(u_1)H[W|u_1]+p(u_2)H[W|u_2]) \\
&=& H[W]-H[W|U]
\end{eqnarray}
This is called the {\it mutual information}; it tells us how much
less uncertain we are, on average, about $W$ given that we know $U$.

 \begin{mydef}[Mutual Information]
Define $Q$ and $R$ to be discrete random variables with support $\{q_1,...,q_n\}$ and $\{r_1,...,r_m\}$ respectively, and let $H(Q)$ be the entropy of $Q$ and $H(Q|R)$ be the conditional entropy. Then the \underline{mutual information} $I$ between $Q$ and $R$ is defined as 
\begin{eqnarray}
I[Q,R] &=&  -\sum_{i,j}p(q_j,r_i)\log\frac{p(q_j,r_i)}{p(q_j)p(r_i)}  \\
&=& H[Q]-H[Q|R] 
\end{eqnarray}
Note: $I[Q,R] = I[R,Q]$ \cite{fraser-swinney}.
\end{mydef}

In the next section, I extend this discussion to information shared between {\it more than} two variables. In the language of {\it this} section, that is equivalent to the situation where I have two or more colleagues with umbrellas $U_{C1}$ and $U_{C2}$ and I want to know the average reduction in uncertainty of the weather given the state of $U_{C1}$ and $U_{C2}$, {\it i.e.,} $I[W,U_{C1},U_{C2}]$. This is unfortunately not a straightforward generalization and there is little agreement in the literature about interpreting or even defining multivariate mutual information.

\subsection{$I$-Diagrams and Multivariate Mutual Information}\label{sec:multiinfoandidiagrams}
The mathematical definitions of multivariate mutual information can get quite confusing to interpret, especially when comparing and contrasting the difference in these definitions. To clarify this discussion, I will use $I$-Diagrams of
Yeung~\cite{yeung2012first}---a highly useful visualization technique for interpreting information theoretic quantities. 
\subsubsection{$I$-Diagrams}

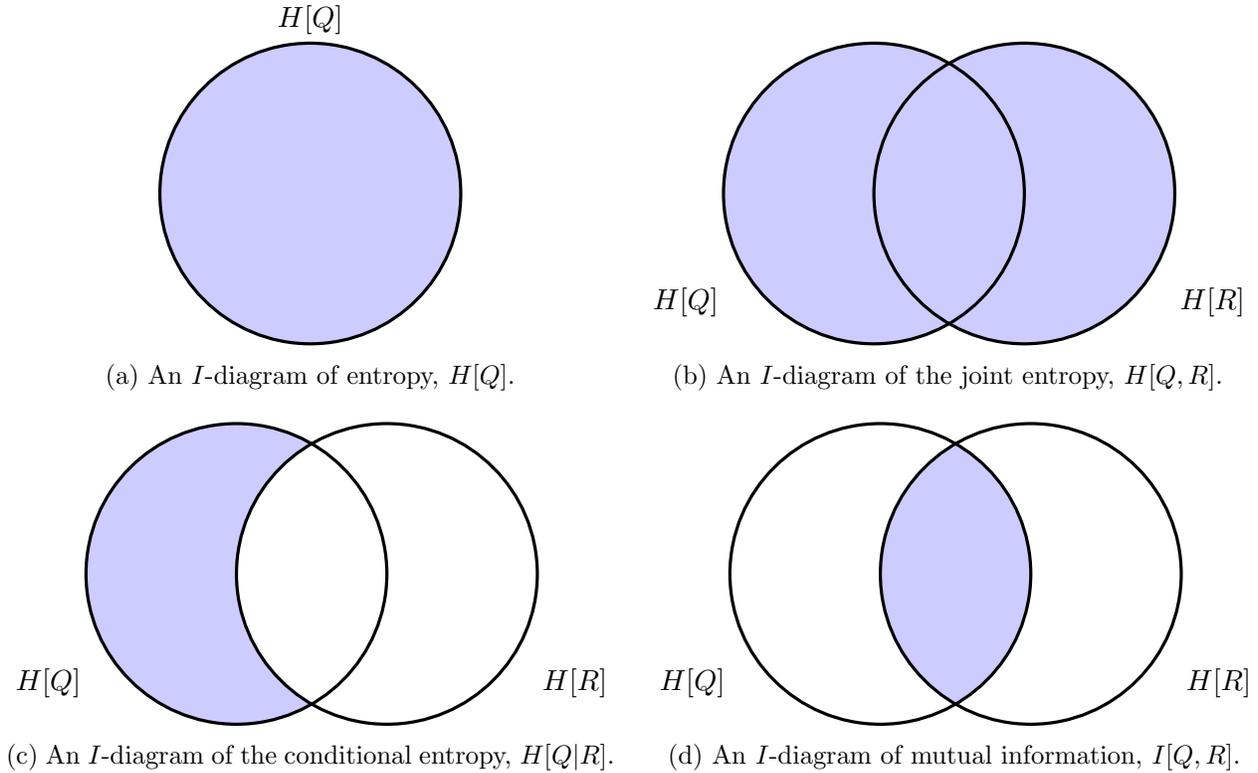
\begin{figure}[tb!]
  \begin{subfigure}[b]{0.49\textwidth}
  	\centering
    \begin{tikzpicture}[baseline=0,label=$m$]
    \setupdiagrams
      \setupdiagramonecircle
      \begin{scope}
        \fill[filled] \Acirc;
      \end{scope}

      \drawdiagram
    \end{tikzpicture}
    \caption{An $I$-diagram of entropy, $H[Q]$.}
    \label{fig:entropy}
  \end{subfigure}%
  \quad
    \begin{subfigure}[b]{0.49\textwidth}
    \begin{tikzpicture}[baseline=0]
      \setupdiagramstwocircles

      \begin{scope}[even odd rule]
        \clip  \Bcirc;
      \end{scope}
      \begin{scope}
        \fill[filled] \Bcirc;
         \fill[white] \Ccirc;
        \fill[filled] \Ccirc;
         \end{scope}

      \drawdiagram
    \end{tikzpicture}
    \caption{An $I$-diagram of the joint entropy,
      ${H}[Q,R]$.}
    \label{fig:joint-entropy}
  \end{subfigure}
  \begin{subfigure}[b]{0.49\textwidth}
    \begin{tikzpicture}[baseline=0]
      \setupdiagramstwocircles

      \begin{scope}[even odd rule]
        \clip  \Bcirc;
      \end{scope}
      \begin{scope}
        \fill[filled] \Bcirc;
        \fill[white] \Ccirc;
      \end{scope}

      \drawdiagram
    \end{tikzpicture}
    \caption{An $I$-diagram of the conditional entropy,
      ${H}[Q|R]$.}
    \label{fig:cond-entropy}
  \end{subfigure}
  \quad
  \begin{subfigure}[b]{0.49\textwidth}
    \begin{tikzpicture}[baseline=0]
      \setupdiagramstwocircles

      \begin{scope}[even odd rule]
        %
        \clip  \Bcirc;
      \end{scope}

      \begin{scope}[even odd rule]
	      \clip \Ccirc;
         \clip \Ccirc ;
      \end{scope}


      \begin{scope}
        \clip \Bcirc;
        \fill[filled] \Ccirc;
      \end{scope}

      \drawdiagram
    \end{tikzpicture}
    \caption{An $I$-diagram of mutual information,
      ${I}[Q,R]$.}
  \end{subfigure}
  \caption{$I$-Diagrams of $H[Q]$, H[Q,R], $H[Q|R]$ and $I[Q,R]$ }\label{fig:i-diagramintro}
\end{figure}

Figure~\ref{fig:i-diagramintro} shows
$I$-diagrams of some of the important quantities introduced in Sections~\ref{sec:entropyprimer} and ~\ref{sec:intromi}: (a) entropy $H[Q]$, (b) joint entropy $H[Q,R]$ (c) conditional entropy $H[Q|R]$ and (d) mutual information $I[Q,R]$. In $I$-Diagrams, each circle represents the uncertainty in a
particular variable and the shaded region is the information quantity of interest, {\it e.g.,} in (a) we are interested in $H[Q]$---the uncertainty in $Q$---so the entire circle is shaded. Figure~\ref{fig:i-diagramintro}(b) introduces a new measure: joint entropy $H[Q,R] = \sum_{q,r}p(q,r)\log(p(q,r))$. $H[Q,R]$ is uncertainty about processes $Q$ {\it and} $R$; this is easily depicted in an $I$-Diagram by simply shading both circles.

The real magic of $I$-Diagrams comes from their ability to depict more-complex information theoretic measures by simply manipulating shaded regions. For example, recall from Section~\ref{sec:intromi} that conditional entropy $H[Q|R]$---Figure~\ref{fig:i-diagramintro}(b)---is the uncertainty about process $Q$ given knowledge of $R$. One way of writing this is $H[Q|R] = H[Q,R] - H[Q]$:  {\it i.e.,} subtracting the shaded regions in (a) and (b) produces the shaded region in (c). The same can be done with mutual information. Recall from Section~\ref{sec:intromi} that $I[Q,R]$ is the shared uncertainty between $Q$ and $R$ or  $I[Q,R] = H[Q] -H[Q|R]$: {\it i.e.,} subtracting the shaded region in (a) from the shaded region in (c) produces the shaded region in (d). While obviously not a proof, this kind of approach allows us to easily build intuition about more complicated identities, {\it e.g.,} symmetry of mutual information:
$I[Q,R] = H[Q] -H[Q|R] =H[R] -H[R|Q] = I[R,Q].$

In the next section, I will use $I$-diagrams to explore three common interpretations of multivariate mutual information,  interaction information~\cite{McGill-1954} (also commonly called the co-information~\cite{Bell03theco-information}),  the binding information~\cite{binding} (also called the dual total correlation~\cite{Han1978133}),  and total correlation~\cite{total-cor} (also commonly called multi-information~\cite{Studen1998}).

\subsection{Multivariate Mutual Information}
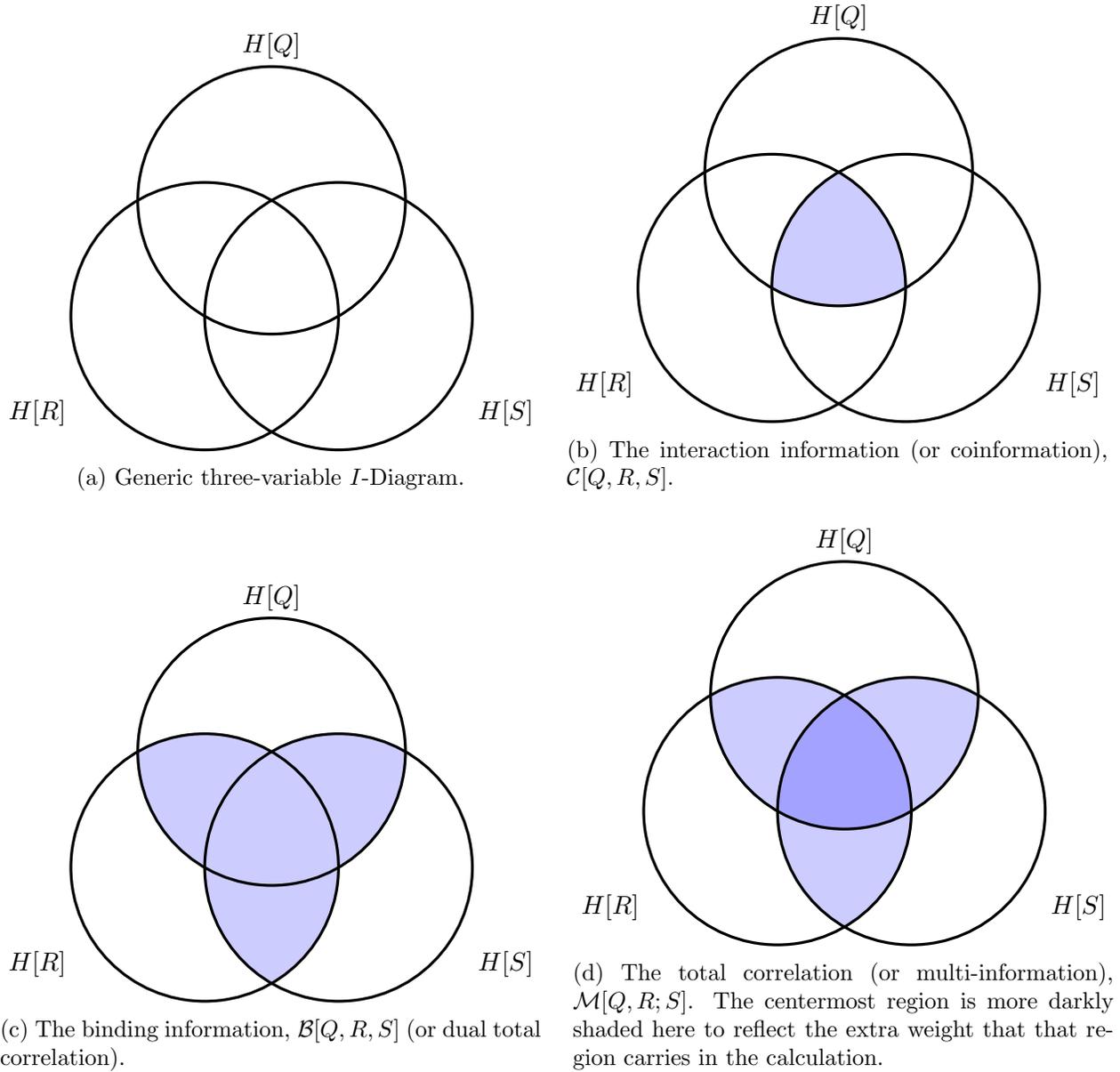
\begin{figure}[tb!]
  \begin{subfigure}[b]{0.49\textwidth}
    \begin{tikzpicture}[baseline=0,label=$m$]
      \setupdiagramsgeneral

      \begin{scope}
        \clip \Acirc;
        \clip \Bcirc;
      \end{scope}

      \drawdiagram
    \end{tikzpicture}
    \caption{Generic three-variable $I$-Diagram.}
    \label{fig:empty}
  \end{subfigure}%
  \quad
    \begin{subfigure}[b]{0.49\textwidth}
    \begin{tikzpicture}[baseline=0,label=$m$]
      \setupdiagramsgeneral

      \begin{scope}
        \clip \Acirc;
        \clip \Bcirc;
        \fill[filled] \Ccirc;
      \end{scope}

      \drawdiagram
    \end{tikzpicture}
    \caption{The interaction information (or coinformation), $\mathcal{C}[Q,R,S]$.}
    \label{fig:co-information}
  \end{subfigure}%
  
  \begin{subfigure}[b]{0.49\textwidth}
    \begin{tikzpicture}[baseline=0]
      \setupdiagramsgeneral

      \begin{scope}[even odd rule]
	      \clip \Acirc;
        \clip \Acirc \Bcirc;
        \fill[filled] \Ccirc;
      \end{scope}

      \begin{scope}[even odd rule]
	      \clip \Ccirc;
        \clip \Ccirc \Acirc;
        \fill[filled] \Bcirc;
      \end{scope}

      \begin{scope}[even odd rule]
	      \clip \Bcirc;
        \clip \Bcirc \Ccirc;
        \fill[filled] \Acirc;
      \end{scope}

      \begin{scope}
        \clip \Acirc;
        \clip \Bcirc;
        \fill[filled] \Ccirc;
        \fill[filled] \Ccirc;
      \end{scope}

\begin{scope}
        \clip \Acirc;
        \clip \Bcirc;
        \fill[white] \Ccirc;
      \end{scope}
\begin{scope}
        \clip \Acirc;
        \clip \Bcirc;
        \fill[filled] \Ccirc;
      \end{scope}

      \drawdiagram
    \end{tikzpicture}
    \caption{The binding information, $\mathcal{B}[Q,R,S]$ (or dual total correlation).}
  \end{subfigure}
  \quad
  \begin{subfigure}[b]{0.49\textwidth}
    \begin{tikzpicture}[baseline=0]
      \setupdiagramsgeneral

      \begin{scope}[even odd rule]
	      \clip \Acirc;
        \clip \Acirc \Bcirc;
        \fill[filled] \Ccirc;
      \end{scope}

      \begin{scope}[even odd rule]
	      \clip \Ccirc;
        \clip \Ccirc \Acirc;
        \fill[filled] \Bcirc;
      \end{scope}

      \begin{scope}[even odd rule]
	      \clip \Bcirc;
        \clip \Bcirc \Ccirc;
        \fill[filled] \Acirc;
      \end{scope}

      \begin{scope}
        \clip \Acirc;
        \clip \Bcirc;
        \fill[filled] \Ccirc;
        \fill[filled] \Ccirc;
      \end{scope}

      \drawdiagram
    \end{tikzpicture}
    \caption{The total correlation (or multi-information),
      $\mathcal{M}[Q,R;S]$. The centermost
      region is more darkly shaded here to reflect the extra weight
      that that region carries in the calculation.}
    \label{fig:multi-information}
  \end{subfigure}
  \caption{Generalizations of the mutual information to the multivariate case.} \label{fig:multi-information-panel}

\end{figure}

When interpreting $I[Q,R]$ using $I$-Diagrams, the situation is quite simple, as there is exactly one region of ``shared uncertainty;" when generalizing even to three variables $I[Q,R,S]$, the situation becomes much more confusing---and this is reflected in the mathematical uncertainties as well. Consider the generic three-variable $I$-Diagram in Figure~\ref{fig:multi-information-panel}. Instead of having one region of overlap as in Figure~\ref{fig:i-diagramintro}(d), there are now four. There are three standard ways of shading each these regions to quantify $I[Q,R,S]$.

One interpretation is the so-called {\it interaction information}~\cite{McGill-1954,Bell03theco-information}
\begin{align*}
  \mathcal{C}[Q,R,S]\equiv I[Q,R,S] \equiv & - ( H[Q] + H[R] + H[S] ) \\
               & + ( H[Q,R] + H[Q,S] + H[R,S] ) \\
               & - H[Q,R,S] \numberthis
\end{align*}
As depicted in Figure~\ref{fig:multi-information-panel}(b), this is the intersection of $H[Q]$, $H[R]$ and $H[S]$. It describes the reduction in uncertainty that any {\it two} processes ({\it e.g.,} $Q$ and $R$), together, provide regarding the third process ({\it e.g.,} $Q$ and $R$). While this may seem like the natural extension of mutual information, it does not take into account the information that is shared between the two process but {\it not with the third}. One common criticism of this interpretation is that $\mathcal{C}[Q,R,S]$ is quite often negative. For example, when the shared information between $\{Q,R\}$ is due entirely to information in $S$, the interaction information can be negative as well as positive. Many interpretations of negative information have been provided---{\it e.g.,} that the variable $S$ inhibits ({\it i.e.,} accounts for or explains some of) the correlation between \{Q,R\}---but in general negative information is frowned upon~\cite{Abdallah}. 

The next obvious step is to take into account the information that is shared between any two process but {\it not shared with the third}, as well as the information shared between all three processes. This is called the binding information~\cite{binding,Han1978133}  
\begin{align}
\mathcal{B}[Q,R,S]\equiv I[Q,R,S] \equiv H[Q,R,S] + \left[ \sum_{i=1}^3 H[X_{(i-1)\%3},X_{(i+1)\%3]}) - H[Q,R,S]\right]
\end{align}%
\noindent where $X_0= Q$, $X_1= R$, and $X_2= S$, and \% is the modulus operator. This quantity is depicted in Figure~\ref{fig:multi-information-panel}(c). $\mathcal{B}[Q,R,S]$ has the nice feature that it is always positive, but it equally weights information contained in two variables as information contained in all three. 
The total correlation~\cite{total-cor,Studen1998}
\begin{align}
  \mathcal{M}[Q,R,S] \equiv I[Q,R,S] \equiv I[X_0,X_1,X_2] \equiv \sum_{\mathclap{i=0 }}^2 \left( H[X_i] \right) - H[X_0,X_1,X_2]
  \label{eq:multi-information}
\end{align}
depicted in Figure~\ref{fig:multi-information-panel}(d) addresses this shortcoming, but is equally criticized for over emphasizing information that is shared by all three variables. 

The total correlation and binding information are both always positive but their relative merits are a subject of contention. For a nice comparison of these and discussion of the associated issues,  please see \cite{Abdallah}. The takeaway of this section should be that extending mutual information as defined in Section~\ref{sec:intromi} to even the three-variable case, let alone beyond that, is non-trivial and not well understood at all. This will become very important in Section~\ref{sec:tdAIS}, where I propose a new information-theoretic method for selecting delay reconstruction parameters.

\subsection{Estimating Information from Real-Valued Time-Series Data}\label{sec:symbolize}
Note that all the information theory discussed thus far has been on {\it discrete} random variables. The topic of this thesis, however, involves {\it real-valued} time series. To compute any information measure on {\it real-valued} time series, one must ``symbolize" that data, {\it i.e.,} map the real values to a set of discrete symbols. Ideally, this symbolization should preserve the information and/or dynamics of the original time series, but this can be hard to accomplish in practice.   The processes by which this is accomplished, and the issues that make it difficult, are the focus of this section.
\subsubsection{Simple Binning}\label{sec:binning}

A common (and by far the simplest) symbolization method is {\it binning}. To symbolize a real-valued time series $\{x_j\}_{j=1}^N$ by binning, one breaks the time series support into $n$ bins, which need not be equally spaced. Then one defines a discrete random variable $Q$ to have a symbol for each bin $b_i$, {\it i.e.,} $Q$ has support $\{b_i\}_{i=1}^n$. The associated probability mass function is then computed using \begin{equation}p(b_i) = \frac{|\{j|x_j \in b_i\}|}{N}\end{equation}
For example, consider a time series with support on $[0,1]$ and bins $b_1 = [0,0.5)$ and $b_2 = [0.5,1]$. Then one simply estimates the probability mass function associated with $b_1$ and $b_2$ by counting the number of time-series elements that appear in each subinterval $[0,0.5)$ and $[0.5,1]$. 

This method is extremely simple, but simplicity is often a double-edged sword. Binning is a very fast and efficient symbolization, but it is known to introduce severe biasing and 
spurious dynamics if the bin boundaries do not happen to create a so called {\it generating partition} of the dynamics~\cite{KSG, bollt2001}. 
\begin{mydef}[Generating Partition]
Given a dynamical  system $f: \mathcal{M}\rightarrow\mathcal{M}$ on a measure space $(\mathcal{M},F,\mu)$, a finite partition $P = \{b_i\}_{i=1}^n$ is said to be generating if the union of all images and preimages of $P$ gives the set of all $\mu$-measurable sets $F$. In other words, the ``natural" tree of partitions always generates some sub-$\sigma$-algebra, but if it gives the full $\sigma$-algebra of all measurable sets $F$, then $P$ is called \underline{generating}~\cite{rudolph-measurable-dynamics}.
\end{mydef}

Unfortunately, even for most canonical dynamical systems, let alone all real-valued time series, the generating partition is not known or computable. Even when the partition is known it can be fractal, as is the case with the H\'enon map, for example, and thus useless for creating a {\it finite} partition. For a good review of the difficulties in finding a generating partition see \cite{PhysRevE.61.1353}. I review this in more detail in Section \ref{sec:meaComplex-intro}.

\subsubsection{Kernel estimation methods}\label{sec:KSG}
\label{sec:KSG}

A useful alternative to simple binning is a class of methods known as {\it kernel estimation} \cite{PhysRevLett.85.461,PhysRevLett.99.204101}, in which the relevant probability density functions are estimated via a function $\Theta$ with a resolution or bandwidth $\rho$ that measures the similarity between two points in $Q \times R$ space.\footnote{ In the case of delay-coordinate reconstruction, $Q \times R = X_j \times X_{j-\tau}$} Given points $\{q_i,r_i\}$ and $\{q'_i,r'_i\}$ in $Q \times R$, one can define
%
\begin{equation}
  \hat{p}_\rho(q_i,r_i) = \frac{1}{N} \sum_{i'=1}^N\Theta
  \left( \begin{array}{|c|}
    q_i-q'_i \\%
    r_i-r'_i
  \end{array} -\rho \right)
  \end{equation}
where $\Theta(x>0) = 0$ and $\Theta(x\le 0 ) = 1$.  That is,
$\hat{p}_\rho(q_i,r_i)$ is the proportion of the $N$ pairs of points in
$Q\times R$ space that fall within the kernel bandwidth $\rho$ of
$\{q_i,r_i\}$, i.e., the proportion of points similar to
$\{q_i,r_i\}$. When $|\cdot|$ is the max norm, this is the so-called box
kernel.  This too, however, can introduce bias~\cite{jidt} and is obviously 
dependent on the choice of bandwidth $\rho$.  After these estimates, and/or the
analogous estimates for $\hat{p}(q)$, are produced, they are
then used directly to
compute local estimates of entropy or mutual information for each point in space,
which are then averaged over all samples to produce the entropy or mutual
information of the time series.  For more details on this procedure,
see~\cite{jidt}.

A less biased method to perform kernel estimation when one is interested in computing mutual information is the
Kraskov-St\"ugbauer-Grassberger (KSG) estimator~\cite{KSG}.  This
approach dynamically alters the kernel bandwidth to match the density
of the data, thereby smoothing out errors in the probability density
function estimation process.  In this approach, one first finds the
$k^{th}$ nearest neighbor for each sample $\{q,r\}$ (using max norms
to compute distances in $q$ and $r$), then sets kernel widths $\rho_q$
and $\rho_r$ accordingly and performs the pdf estimation.  There are two
algorithms for computing $I[Q,R]$ with the KSG estimator~\cite{jidt}.
The first is more accurate for small sample sizes but more biased; the
second is more accurate for larger sample sizes.  I use the second of
the two in the results reported in this dissertation, as I have fairly long time series.  This
algorithm sets $\rho_q$ and $\rho_r$ to the $q$ and $r$ distances to the
$k^{th}$ nearest neighbor.  One then counts the number of neighbors
within and on the boundaries of these kernels in each marginal space,
calling these sums $n_q$ and $n_r$, and finally calculates
\begin{equation}
  I[Q,R] = \psi(k) - \frac{1}{k}-\langle \psi(n_q) +\psi(n_r) \rangle + \psi(n)
\end{equation}
where $\psi$ is the digamma function\footnote{The formula for the
  other KSG estimation algorithm is subtly different; it sets $\rho_q$
  and $\rho_r$ to the maxima of the $q$ and $r$ distances to the $k$
  nearest neighbors.}.  This estimator has been demonstrated to be
robust to variations in $k$ as long as $k\ge4$~\cite{jidt}.

In this thesis, I employ the Java Information Dynamics Toolkit (JIDT)
implementation of the KSG estimator~\cite{jidt}.  The computational
complexity of this implementation is $\mathcal{O}(kN\log N)$, where
$N$ is the length of the time series and $k$ is the number of
neighbors being used in the estimate.  While this is more expensive
than traditional binning $(\mathcal{O}(N)$), it is bias corrected,
allows for adaptive kernel bandwidth to adjust for under- and
over-sampled regions of space, and is both model and parameter free
(aside from $k$, to which it is very robust).

\subsection{Estimating Structural Complexity and Predictability}\label{sec:meaComplex-intro}
An understanding of the predictive capacity of a real-valued time series---{\it i.e.,} whether or not it is even predictable---is essential to any forecasting strategy. In joint work with Ryan James, I propose to quantify the complexity of a signal by approximating the entropy production of the system that generated it. In general,   
estimating the entropy (production) of an arbitrary, real-valued time series is a
challenging problem, as discussed above, but recent advances in Shannon information theory---in particular, permutation entropy~\cite{bandt2002per,fadlallah2013}---have reduced this challenge. I review this class of methods in this section.  

For the purposes of this thesis, I view the Shannon entropy---in
particular, its growth rate with respect to word length (the
{\it Shannon entropy rate})---as a measure of the complexity and hence the 
predictability for a time series.  Time-series data consisting of i.i.d. random variables, such as white noise, have high entropy rates,
whereas highly structured time-series---for example, those that are periodic---have very
low (or zero) entropy rates. A time series with a high entropy rate is
almost completely unpredictable, and conversely. This can be made more rigorous: Pesin's
relation \cite{pesin1977characteristic} states that in chaotic dynamical systems, the
Kolmogorov-Sinai (KS) entropy is equal to the sum of the positive
Lyapunov exponents $\lambda_i$.  These exponents directly
quantify the rate at which nearby states of the system diverge with
time: $\left| \Delta x(t) \right| \approx e^{\lambda t} \left| \Delta
x(0) \right|$.  The faster the divergence, the larger the entropy.
The KS entropy is defined as the supremum of the Shannon entropy rates
of all partitions---{\it i.e.,} all possible choices for binning~\cite{petersen1989}. As an aside, an alternative definition of the generating partition defined above is a partition that achieves this supremum.

From a different point of view, I can consider the information (as
measured by the Shannon entropy) contained in a single observable of the system at a given point in time. This information can be partitioned
into two components: the information shared with past
observations---{\it i.e.}, the mutual information between the past and
present---and the information in the present that is not contained in
the past ({\it viz.,} ``the conditional entropy of the present given the
past'').  The first part is known as the {\it redundancy};  the second 
is the aforementioned {\it Shannon entropy rate}.  Again working with R.~G.~James, I establish that 
the more redundancy in a signal, the more predictable it is~\cite{josh-pre,ISIT13}. This is discussed in more detail in Chapter~\ref{ch:wpe}.

Previous approaches to measuring temporal complexity via the Shannon
entropy rate \cite{Shannon1951, mantegna1994linguistic} required
categorical data: $x_i \in \mathcal{S}$ for some finite or countably
infinite {\it alphabet} $\mathcal{S}$.  Data taken from real-world
systems are, however, effectively\footnote{Measurements from
  finite-precision sensors are discrete, but data from modern
  high-resolution sensors are, for the purposes of entropy
  calculations, effectively continuous.}  real-valued. So for this reason I need to symbolize the time series, as discussed above. The methods discussed above however, are generally biased or fragile in the face of noise.

Bandt and Pompe introduced the {\it permutation entropy} (PE) as a
``natural complexity measure for time series''~\cite{bandt2002per}. Permutation entropy involves a method for
symbolizing real-valued time series that follows the intrinsic
behavior of the system under examination.  This method  has many advantages,
including robustness to observational noise, and its application
does not require any knowledge of the underlying mechanisms of the system.  Rather
than looking at the statistics of sequences of values, as is done when
computing the Shannon entropy, permutation entropy looks at the
statistics of the {\it orderings} of sequences of values using
ordinal analysis. Ordinal analysis of a time series is the process of
mapping successive time-ordered elements of a time series to their
value-ordered permutation of the same size.  By way of example, if
$(x_1, x_2, x_3) = (9, 1, 7)$ then its {\it ordinal pattern},
$\phi(x_1, x_2, x_3)$, is $231$ since $x_2 \leq x_3 \leq x_1$.  The
ordinal pattern of the permutation $(x_1, x_2, x_3) = (9, 7, 1)$ is
$321$.

\begin{mydef}[Permutation Entropy]

  Given a time series $\{x_i\}_{i = 1,\dots,N}$, define $\mathcal{S}_\ell$ as all $\ell!$ permutations $\pi$ of order $\ell$. For each $\pi \in \mathcal{S}_\ell$, define the relative frequency of that permutation occurring in $\{x_i\}_{i = 1,\dots,N}$
  \begin{equation}
    p(\pi) = \frac{\left|\{i|i \leq N-\ell,\phi(x_{i+1},\dots,x_{i+\ell}) = \pi\}\right|}{N-\ell+1}
  \end{equation}
  where $p(\pi)$ quantifies the probability of an ordinal and
  $|\cdot|$ is set cardinality. The \underline{permutation entropy} of
  order $\ell \ge 2$ is defined as
  \begin{equation}
    \textrm{PE}(\ell) = - \sum_{\pi \in \mathcal{S}_\ell} p(\pi) \log_2 p(\pi)
  \end{equation}

\end{mydef}

Notice that $0\le PE(\ell) \le \log_2(\ell!)$ \cite{bandt2002per}.
With this in mind, it is common in the literature to normalize
permutation entropy as follows: $\frac{PE(\ell)}{\log_2(\ell!)}$.  With
this convention, ``low'' PE is close to 0 and ``high'' PE is close to
1. Finally, it should be noted that the permutation entropy has been
shown to be identical to the Kolmogorov-Sinai entropy for many large
classes of systems \cite{amigo2012permutation}, as long as
observational noise is sufficiently small. As mentioned before, PE
 is equal to the Shannon entropy rate of a generating partition of the
system. This transitive chain of equalities, from permutation entropy
to Shannon entropy rate via the KS entropy, allows one to approximate
the redundancy of a signal---being the dual of the Shannon entropy rate---by $1 - \frac{PE(\ell)}{\log_2(\ell!)}$.

In this thesis, I utilize a variation of the basic permutation entropy
technique, the {\it weighted permutation entropy} (WPE), which was
introduced in~\cite{fadlallah2013}.  The intent behind the weighting
is to correct for observational noise that is larger than the trends
in the data, but smaller than the larger-scale features.  Consider,
for example, a signal that switches between two fixed points and
contains some additive noise. The PE is dominated by the noise about
the fixed points, driving it to $\approx 1$, which in some sense hides
the fact that the signal is actually quite structured.  
 To
correct for this, the {\it weight} of a permutation is taken into
account
\begin{equation}
  w(x_{i+1}^\ell) = \frac{1}{\ell} \sum_{j = i}^{i+\ell}
                      \left( x_j - \bar{x}_{i+1}^\ell \right)^2
\end{equation}
where $x_{i+1}^\ell$ is a sequence of values $x_{i+1}, \ldots,
x_{i+\ell}$, and $\bar{x}_{i+1}^\ell$ is the arithmetic mean of
those values.

The weighted probability of a permutation is defined as
\begin{equation}
  p_w(\pi) = \frac{\displaystyle \sum_{i \le N - \ell} w(x_{i+1}^\ell) \cdot \delta(\phi(x_{i+1}^\ell), \pi) }{\displaystyle \sum_{i \le N - \ell} w(x_{i+1}^\ell)}
\end{equation}
where $\delta(x, y)$ is 1 if $x = y$ and 0 otherwise. Effectively,
this weighted probability emphasizes permutations that are involved in
``large'' features and de-emphasizes permutations that are small in
amplitude, relative to the features of the time series.  The
standard form of weighted permutation entropy is
\begin{equation}
  \textrm{WPE}(\ell) = - \sum_{\pi \in \mathcal{S}_\ell} p_w(\pi) \log_2 p_w(\pi),
\end{equation}
which can also be normalized by dividing by $\log(\ell!)$, to make $0
\le \textrm{WPE}(\ell) \le 1$. 

In practice, calculating permutation entropy and weighted permutation
entropy involves choosing a good value for the word length $\ell$. The
primary consideration in that choice is that the value be large enough
that forbidden ordinals are discovered, yet small enough that
reasonable statistics over the ordinals are gathered.  If an average
of 100 counts per ordinal is considered to be sufficient, for
instance, then $\ell = \operatornamewithlimits{argmax}_{\hat{\ell}} \{
N \gtrapprox 100 \hat{\ell}! \}$.  In the literature, $3 \le \ell \le
6$ is a standard choice---generally without any formal justification.
In theory, the permutation entropy should reach an asymptote with
increasing $\ell$, but that can require an arbitrarily long time
series. In practice, what one should do is calculate the
{\it persistent} permutation entropy by increasing $\ell$ until the
result converges, but data length issues can intrude before that
convergence is reached.  I use this approach to choose $\ell = 6$
for the experiments presented in this thesis.  This value represents a
good balance between accurate ordinal statistics and finite-data
effects.

\section{Forecast Methods}\label{sec:forecastmodels}

Any discussion of new prediction technology is incomplete, of course,
without a solid comparison to traditional techniques.   
In this section, I describe the four different forecasting methods
used in my thesis as points of comparison. These forecast methods include:
\begin{itemize}
\item The {\it random-walk} method, which uses the previous value in
  the observed signal as the forecast,

\item The {\it \naive} method, which uses the mean of the
  observed signal as the forecast,

\item The {\it ARIMA} (auto-regressive integrated moving average)
  method, a common linear forecast strategy for scalar time-series data, instantiated via the
  {\it \arima} procedure \cite{autoARIMA}, and


\item The {\it LMA} (Lorenz method of analogues) method, which uses a
  near-neighbor forecast strategy on a delay-coordinate reconstruction of the
  signal.
\end{itemize}
ARIMA, as its name suggests, is based on standard autoregressive techniques.  LMA, introduced in Chapter~\ref{ch:overview}, is designed
to capture and exploit the deterministic structure of a signal from a
nonlinear dynamical system.  The \naive ~and random-walk methods,
somewhat surprisingly, often outperform these more-sophisticated
prediction strategies in the case of highly complex signals, as
discussed briefly below and in depth in Chapter~\ref{ch:wpe}.

\subsection{Simple Prediction Strategies}
\label{sec:simple}

A random-walk predictor simply uses the last observed measurement as
the forecast: that is, the predicted value $p_i$ at time $i$ is
calculated using the relation 
\begin{equation}p_i \equiv x_{i-1}\end{equation}
where $\{x_j\}_{j=1}^N$ is the time-series data. The prediction
strategy that I refer to using the term ``na\"ive" averages all prior
observations to generate the forecast 
\begin{equation}p_i \equiv \mu_{x,i-1} =
\sum_{j=1}^{i-1}\frac{x_j}{i-1}\end{equation}
 While both of these methods are
simplistic, they are not without merit.  For a time series that possess very
little predictive structure (WPE $\approx 1$), these two methods can actually be the
best choice.  In forecasting currency exchange rates, for instance,
sophisticated econometrics-based prediction models fail to
consistently outperform the random-walk method~\cite{rwMeese,rwCCE}.
These signals are constantly changing, subject to jump processes, noisy, and possess very little predictive
structure, but their variations are not---aside from jump processes---very large, so
the random-walk method's strategy of simply guessing the last known
value is not a bad choice.  If a signal has a unimodal distribution
with low variance, the \naive ~prediction strategy will perform quite
well, even if the signal is highly complex, simply because the mean
is a good approximation of the future behavior.  Moreover, the \naive
~prediction strategy's temporal average effects a low-pass filtering
operation, which mitigates the complexity in signals with
very little predictive structure.

Both of these methods have significant weaknesses, however.  Because
they do not model the temporal patterns in the data, or even the
distribution of its values, they cannot track changes in that
structure.  This causes them to fail in a number of important
situations.  Random-walk strategies are a particularly bad choice for
time series that change significantly at every time step.  In the
worst case---a large-amplitude square wave whose period is equivalent
to twice the sample time---a random-walk prediction would be exactly
180 degrees out of phase with the true continuation.  The \naive
~method would be a better choice in this situation, since it would
always predict the mean.  It would, however, perform poorly when a
signal had a number of long-lived regimes that have significantly
different means.  In this situation, the inertia of the \naive
~method's accumulating mean is a liability and the agility of the
random-walk method is an advantage, since it can respond quickly to
regime shifts.

Of course, methods that can capture and exploit the geometry of the
data, or its temporal patterns, can be far more effective in the
situations described in the previous two paragraphs.  The ARIMA and LMA
methods covered in the following sections are
designed to do exactly that. 
However, if a signal contains little predictive structure, forecast
strategies like ARIMA and LMA have nothing to work with and thus will
often be outperformed by the two simple strategies described in this
section.  This contrast is explored further in Sections~\ref{sec:accuracy}
and~\ref{sec:compPerfProj}.

\subsection{(ARIMA) A Regression-Based Prediction Strategy}
\label{sec:arima}

A simple and yet powerful way to capture and exploit the structure of
data is to fit a hyperplane to the dataset and then use it to make
predictions.  The roots of this approach date back to the original
autoregressive schema of Yule~\cite{Yule27}, which forecasts the value at the next time
step through a weighted average of the past $q$ observations 
\begin{equation}p_i \equiv \sum_{j=i-q}^{i-1} a_j x_j\end{equation}
The weighting coefficients $a_j$ are
generally computed using either an ordinary least-squares approach, or
with the method of moments using the Yule-Walker equations.  To
account for noise in the data, one can add a so-called ``moving
average'' term to the model; to remove nonstationarities, one can
detrend the data using a differencing operation.  A strategy that
incorporates all three of these features is called a {\it nonseasonal
  ARIMA model}.  If evidence of periodic structure is present in the
data, a {\it seasonal ARIMA model}, which adds a sampling operation
that filters out periodicities, can be a good choice.

There is a vast amount of theory and literature regarding the
construction and use of models of this type; please see \cite{davislinearts} for an in-depth exploration.  
For the purposes of
this thesis, I choose 
seasonal ARIMA models to serve as a good exemplar for a broad class of
linear predictors and a useful point of comparison for my work. 
Fitting such a model to a time series involves
choosing values for the various free parameters in the autoregressive,
detrending, moving average, and filtering terms.  I employ the
automated fitting techniques described in~\cite{autoARIMA} to
accomplish this.  This procedure uses sophisticated methods---KPSS
unit-root tests~\cite{KPSSunit}, a customization of the Canova-Hansen
test~\cite{Canova1995}, and the Akaike information
criterion~\cite{akaike1974}, conditioned on the maximum likelihood
of the model fitted to the detrended data---to select good values for
these free parameters.

ARIMA forecasting is a common and time-tested procedure.  Its
adjustments for seasonality, nonstationarity, and noise make it an
appropriate choice for short-term predictions of time-series data
generated by a wide range of processes.  If information is being
generated and/or transmitted in a nonlinear way, however, a global
linear fit is inappropriate and ARIMA forecasts can be inaccurate.
Another weakness of this method is prediction horizon: an ARIMA
forecast is guaranteed to converge to a constant or linear trend after some number of
predictions, depending on model order.  To sidestep this issue, and make the comparison as fair as possible, I
build ARIMA forecasts in a stepwise fashion: {\it i.e.}, fit the model to the
existing data, use that model to perform a one-step prediction, rebuild it
using the latest observations, and iterate until the desired
prediction horizon is reached.  For consistency, I take the same
approach with the other models in this proposal as well, even
though doing so amounts to artificially hobbling LMA, the method that is the topic of the next section.

\subsection{Lorenz Method of Analogues}\label{sec:lma}
 The dynamical systems community has
developed a number of methods that leverage delay-coordinate reconstruction for the purposes of forecasting dynamical systems ({\it e.g.},~\citeDCEFORECASTING).   Since the goal of this
thesis is to show
that {incomplete reconstructions}---those that are not true embeddings---can give these kinds of
methods enough traction to generate useful predictions, I choose one
of the oldest and simplest members of that family to use in my analysis: Lorenz's method of analogues (LMA), which is essentially
nearest-neighbor forecasting in reconstruction space. 

To apply LMA to a scalar time-series data set $\{x_j\}_{j=1}^n$, one
begins by performing a delay-coordinate {reconstruction} to
produce a trajectory of the form
 \begin{equation}\{\vec{x}_j=[x_j~x_{j-\tau} ~ \dots ~ x_{j-(m-1)\tau}]^T \}_{j=1-(m-1)\tau}^{n}\end{equation}
  using one or more of the heuristics presented in Sections~\ref{sec:numericaltau} and \ref{sec:numericalM} to choose $m$ and $\tau$.
Forecasting the next point in the time series, $x_{n+1}$, amounts to
reconstructing the next delay vector $\vec{x}_{n+1}$ in the
trajectory.  Note that, by the form of delay-coordinate vectors, all
but the first coordinate of $\vec{x}_{n+1}$ are known.  To choose that first coordinate, LMA finds the nearest neighbor of $\vec{x}_{n}$ in
the reconstructed space\footnote{$\vec{x}_{n}$ should not be chosen as its own neighbor as it has no forward image. In some cases, a longer Theiler exclusion may be useful~\cite{theiler-window}.}---namely $\vec{x}_{j(1,m)}$---and maps that
vector forward using the delay map, obtaining 
\begin{equation} \vec{x}_{j(1,m)+1}=[x_{j(1,m)+1}~x_{j(1,m)+1-\tau} ~ 
\dots ~ x_{j(1,m)+1-(m-1)\tau}]^T\end{equation}
Using the image of the neighbor, one defines \begin{equation}\vec{p}_{n+1} \equiv
[x_{j(1,m)+1}~x_{n+1-\tau} ~ \dots ~ x_{n+1-(m-1)\tau}]^T\end{equation} The LMA
forecast of $x_{n+1}$ is then $p_{n+1}\equiv x_{j(1,m)+1}$.  If performing
multi-step forecasts, one appends the new delay
vector \begin{equation}\vec{p}_{n+1}=[x_{j(1,m)+1}~x_{n+1-\tau} ~ \dots ~
x_{n+1-(m-1)\tau}]^T\end{equation} to the end of the trajectory and repeats this
process as needed.

In my work, I use the LMA algorithm in two ways: first---as a baseline for
comparison purposes---on an embedding of each time series, with $m$
chosen using the false near(est) neighbor method~\cite{KBA92}; second,
with $m$ fixed at 2.  In the rest of this thesis, I will refer to
these as \fnnLMA and \roLMA, respectively.  The experiments reported
in Chapter~\ref{ch:pnp}, unless stated otherwise, use the same $\tau$ value for
both \fnnLMA and \roLMA, choosing it at the first minimum of the
time-delayed mutual information of the time
series~\cite{fraser-swinney}.  In Section~\ref{sec:time-scales}, I
explore the effects of varying $\tau$ on the accuracy of both methods. In Section~\ref{sec:tdAIS}, I show that a time-delayed version of the so-called {\it active information storage} is a highly effective method for selecting $\tau$, and $m$ as well, when forecasting is the end goal.

Dozens---indeed, hundreds---of more-complicated variants of 
the LMA algorithm have appeared in the literature ({\it
e.g.},~\cite{weigend-book,casdagli-eubank92,Smith199250,sugihara90}),
most of which involve building some flavor of local-linear model
around each delay vector and then using it to make the prediction of
the next point.  
For the purposes of this thesis, I chose to use the basic original LMA
because it is dynamically the most straightforward and thus provides a
good baseline assessment.  While I believe that the claims stated
here extend to other state space-based forecast methods, the
pre-processing steps involved in some of those methods make a careful
analysis of the results somewhat problematic.  One can use GHKSS-based
techniques, for instance, to project the full dynamics onto linear
submanifolds \cite{ghkss} and then use those manifolds to build
predictions.  While it might be useful to apply a method like that to
an incomplete reconstruction, the results would be some nonlinear
conflation of the effects of the two different projections and it
would be difficult to untangle and understand the individual effects.
(Note that the careful study of the effects of projection in
forecasting that are undertaken in this thesis may suggest why
GHKSS-based techniques work so well; this point is discussed further in Chapter~\ref{ch:pnp}.)  

Since LMA does not rest on an assumption of linearity (as ARIMA models
do), it can handle both linear and nonlinear processes.  If the
underlying generating process is nondeterministic, however, it can
perform poorly. For an arbitrary real-valued time series, without any knowledge of the generating process and with all of the attendant
problems (noise, sampling issues, and so on), answers to the question as to which forecast model is best should, ideally, be derived from the data, but that is a difficult task.  By quantifying the
balance between redundancy, predictive structure, and entropy for
these real-valued time series---as I describe in Chapter~\ref{ch:wpe}---I can begin to answer these questions in
an effective and practical manner.

\section{Assessing Forecast Accuracy}\label{sec:accuracy}

To assess and compare the prediction methods studied here, I calculate a figure of merit in the following way. I split each $N$-point time series into two pieces: the first 90\%, referred to as the ``initial training"
signal and denoted $\{x_j\}_{j=1}^{n}$, and the last 10\%, known as
the ``test" signal $\{c_j\}_{j=n+1}^{(k+n+1)=N}$. Following the
procedures described in Section~\ref{sec:forecastmodels}, I build a model from the initial training signal, use that model to
generate a prediction $p_{n+1}$ of the value of $x_{n+1}$, and compare
$p_{n+1}$ to the true continuation, $c_{n+1}$. I then rebuild the model using $\{c_{n+1}\}\cup\{x_j\}_{j=1}^{n}$ and repeat the process $k$ times, out to the end of the observed time series.  This
``one step prediction'' process is not technically necessary in the
\fnnLMA or \roLMA methods, which can generate arbitrary-length\footnote{Although the accuracy of these predictions will
  degrade with prediction horizon, in the presence of positive Lyapunov
  exponents.}predictions, but the performance of the other three methods used here
will degrade severely if the associated models are not periodically
rebuilt.  In order to make the comparison fair, I use the iterative
one-step prediction schema {\it for all five methods}.  This has the
slightly confusing effect of causing the ``test'' signal to be used
both to assess the accuracy of each model and for periodic refitting.

As a numerical measure of prediction accuracy, for each $h$-step forecast, I calculate the $h$-step mean
absolute scaled error ($h$-MASE) between the true and predicted signals, defined as
\begin{equation}
h\mathrm{-MASE} = \sum_{j=n+1}^{k+n+1}\frac{|p_j-c_j|
}{\frac{k}{n-h}\sum^n_{i=1}\sqrt{\frac{\sum^h_{\iota=1}(x_{i}-x_{i+\iota})^2}{h}}}\label{eq:hMASE}
\end{equation}
%
$h$-MASE is a normalized measure: the
scaling term in the denominator is the average $h$-step in-sample forecast error for a random-walk prediction
over the initial training signal $\{x_i\}^n_{i=1}$.  That is, $h$-MASE$<1$ means that the
prediction error in question was, on the average, smaller than the
error of an $h$-step random-walk forecast on the same data.  Analogously,
$h$-MASE$>1$ means that the corresponding prediction method did
{\it worse}, on average, than the random-walk method.  I choose this
error metric because it allows for fair comparison across varying
methods, prediction horizons, and signal scales, and is a standard error measure in the forecasting literature~\cite{MASE}. 

To provide insight into interpreting $h$-MASE values, I will refer back to the proof-of-concept example presented in Chapter~\ref{ch:overview}.  
  The one-step forecasts in Figure~\ref{fig:projExample}, for instance, had 1-MASE values of 0.117 and
0.148---{\it i.e.,} the \fnnLMA and \roLMA forecasts of the SFI dataset A
were, respectively $\frac{1}{0.117}=8.5$ and $\frac{1}{0.148}=6.5$
times better than a one-step random-walk forecast of the initial training
portion of the same signal.

For any non-constant signal, $h$-step forecasting with random walk will degrade as $h$ increases. In general, then, it is to be expected that $h$-MASE will decrease drastically with increasing prediction horizon. Thus, $h$-MASE scores should not be compared for different $h$. For example, 10-MASE can be compared to 10-MASE for two different methods or signals but should not be compared to 1-MASE or 100-MASE, even on the same signal. While its comparative nature may seem odd, this error metric allows
for fair comparison across varying methods, prediction horizons, and
signal scales, making it a standard error measure in the forecasting
literature---and a good choice for this thesis, which involves a number of very different signals.



%% file: caseStudies.tex
\chapter{Case Studies}\label{ch:systems}

I use several different dynamical systems as case studies throughout this thesis, both real and synthetic. 
Two of them---the Lorenz 96 model and sensor data from a laboratory experiment on computer performance dynamics---persist across all chapters of this document; I use a number of others as well to drive home different points in different chapters.  Each is described in more depth in the following sections.

\section{Synthetic Case Studies}
When developing any new mathematical theory or method it is important to first explore it in the context of well-understood synthetic examples. This gives me a controlled environment where I can test the boundaries of my theory, {\it e.g.,} increasing data length or adding a (controlled) signal-to-noise ratio.
\subsection{The Lorenz-96 Model}\label{sec:lorenz96}

\begin{figure}[b!]
        \centering
        
\includegraphics[width=0.48\textwidth]{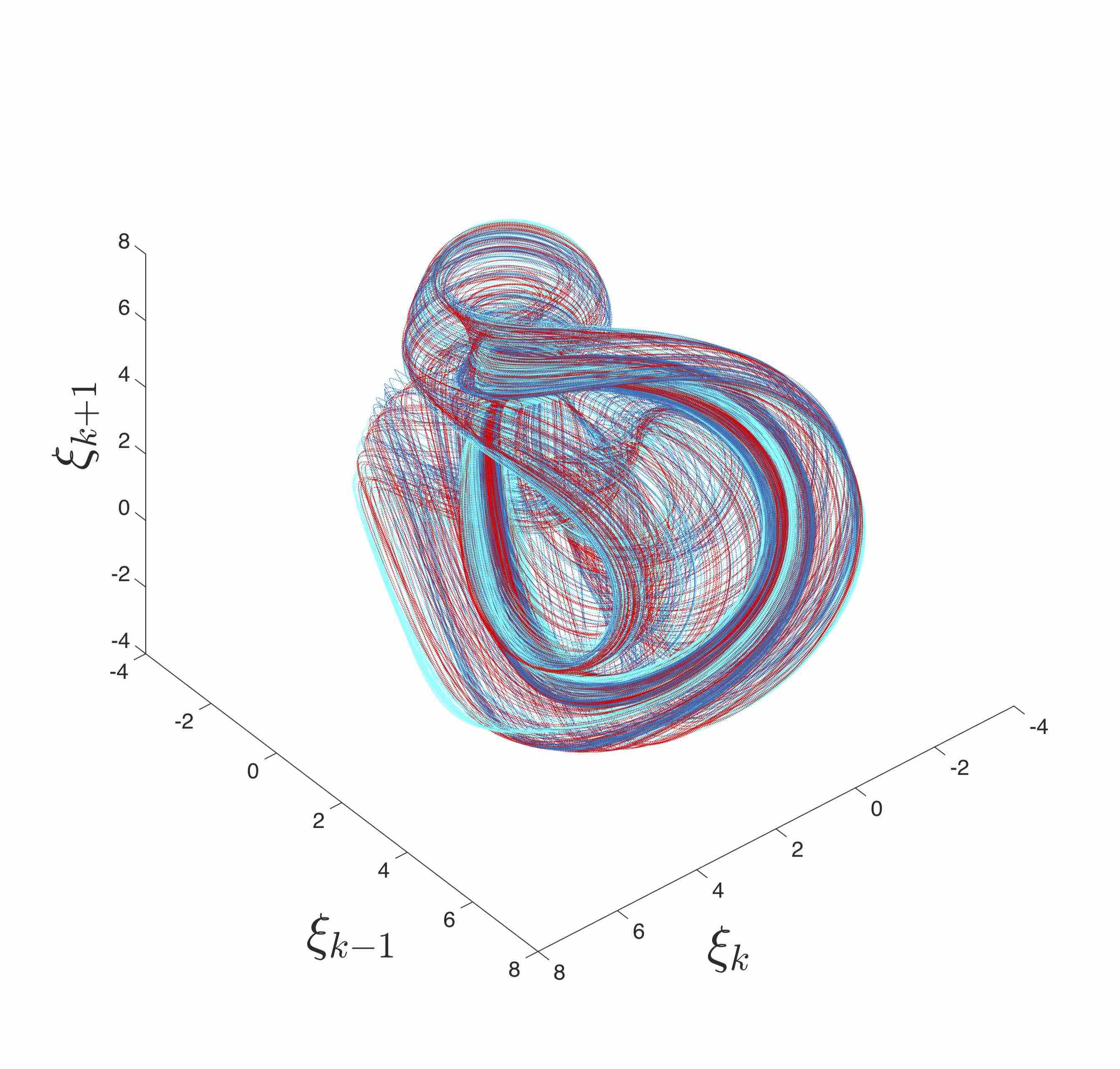}%
\quad
\includegraphics[width=0.48\textwidth]{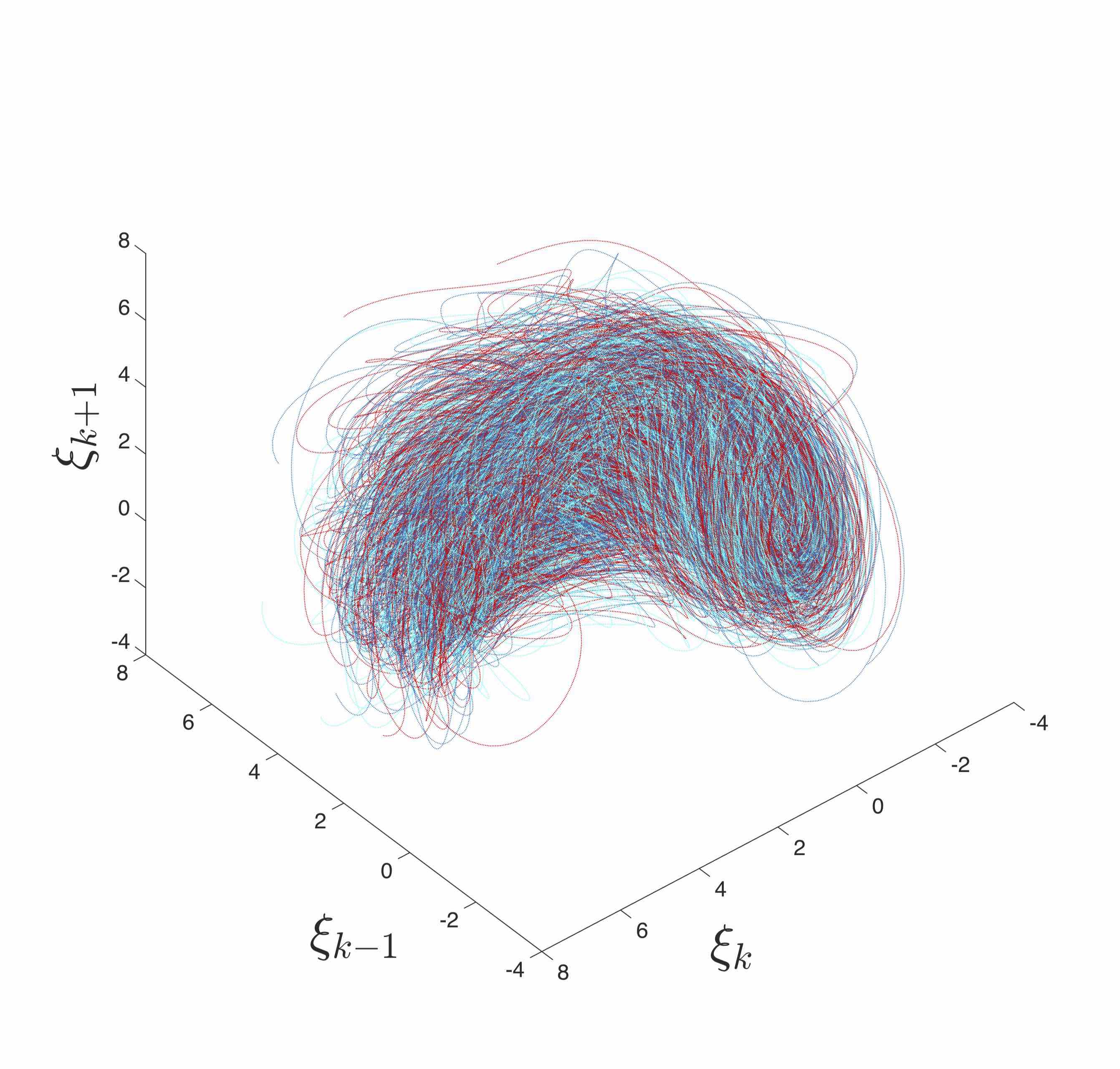}%
        
\caption{The Lorenz 96 attractor ($F=5$) with (left) $K=22$ and (right) $K=47$. Since 22 and 47 dimensional plots are not possible, I plot 3D projections of these systems. In particular, each differently colored trajectory represents different projections or equivalently choices of $k$: $k$=2 is aqua, $k=6$ is blue and $k=18$ is red.}\label{fig:L96Traj}
\end{figure}
The Lorenz-96 model was introduced by Edward Lorenz in \cite{lorenz96Model} to study atmospheric predictability. Lorenz-96 is defined by a set of $K$ first-order differential equations relating the $K$ state variables $\xi_1\dots\xi_K$
\begin{equation}\label{eq:lorenz96}
\dot{\xi}_k= (\xi_{k+1}-\xi_{k-2})(\xi_{k-1})-\xi_k + F
\end{equation}
for $k=1,\dots,K$, where $F\in \mathbb{R}$ is a constant forcing term that is  independent of $k$. In this model, each $\xi_k$ is some atmospheric quantity (such as temperature or vorticity) at a discrete location on a periodic lattice representing a latitude circle of the earth. Following standard
practice~\cite{KarimiL96}, I enforce periodic boundary conditions
and solve Equation~(\ref{eq:lorenz96}) from several randomly chosen
initial conditions using a standard fourth-order Runge-Kutta
solver for 60,000 steps with a step size of $\tfrac{1}{64}$
normalized time units.  I then discard the first 10,000 points of
each trajectory in order to eliminate transient behavior.  Finally, I
create scalar time-series traces by individually ``observing'' each
of the $K$ state variables of the trajectory: {\it i.e.,}
$h_i(\xi_i(t_j)) = x_{j,i}$ for $j\in\{10,000,\dots,60,000\}$ and for
$i\in \{1,\dots,K\}$.  I repeat all of this from a number of
different initial conditions---seven for the $K=47$ Lorenz-96 system
and 15 for the $K=22$ case---producing a total of 659 traces for my
forecasting study. 

In \cite{KarimiL96}, Karimi \& Paul studied this model extensively, performing and analyzing numerous parameter sweeps showing that it exhibits a vast array of possible dynamics: everything from fixed points and periodic attractors to low- and high-dimensional chaos. One particularly interesting feature of this dynamical system is the relationship between state-space dimension and how much of that space is occupied by dynamics. For different choices of the parameter values, the Lorenz-96 system yields dynamics with low fractal dimensions in large state spaces as well as large fractal dimensions in large state spaces. 
All of these features make this model an ideal candidate for testing and evaluating \roLMA.
For my initial investigation, I fix $F=5$  and choose $K=22$ and $K=47$---choices that yield chaotic trajectories with low and
high~\cite{KarimiL96} Kaplan-Yorke (Lyapunov) dimension~\cite{kydimension} respectively: $d_{KY}\lessapprox3$ for the $K=22$
dynamics and $d_{KY}\approx19$ for $K=47$.  Projections of trajectories on these attractors can be seen in Figure~\ref{fig:L96Traj}. 

\begin{table}[tb]
\centering
\begin{tabular}{ccccc}
\hline\hline 
Grid Points & $m$-fnn & $\tau$ & $m$-Embedology & $m$-Takens \\ 
\hline 
$K=22$ & 8 & 26 & $\approx 7$ & 45 \\ 
$K=47$ & 10 & 31 & $\approx 41$ & 95 \\ 
\hline\hline 
\end{tabular}
\caption{Estimated and theoretical embedding parameter values for the Lorenz-96 model. $m$-fnn is the embedding dimension produced by the false-nearest neighbor method. $\tau$ is chosen as the first minimum of the mutual information curve. $m$-Embedology is derived following \cite{sauer91} and $m$-Takens is derived from \cite{takens}.}\label{tab:embedL96}
\end{table}

  For each of these time series, I use the procedures outlined
in Section~\ref{sec:dce} to estimate values for the free parameters of
the embedding process, obtaining $m=8$ and $\tau=26$ for all traces in
the $K=22$ case, and $m=10$ and $\tau=31$ for the $K=47$
traces. Table~\ref{tab:embedL96} tabulates the estimated and theoretical embedding parameter values for these two test cases, derived using the methodologies described in Sections~\ref{sec:dce}-\ref{sec:numericalM}.

It has been shown in~\cite{sprottBook} that
$d_{KY} \approx d_{cap}$ for typical chaotic systems.  This suggests
that embeddings of the $K=22$ and $K=47$ time series would require
$m\gtrapprox6$ and $m\gtrapprox 38$, respectively.  The values
suggested by the false-near neighbor method for the $K=22$ traces are
in line with this, but the $K=47$ false-near neighbor values are far smaller than $2
d_{KY}$. For $K=47$ there are two potential causes for this disparity. First, the false-near neighbor method does not guarantee $m>2
d_{KY}$, it is simply a heuristic to mitigate false crossings in the dynamics. Second $m>2
d_{KY}$ is a {\it sufficient} bound---it could very well be the case that false crossings are eliminated before this bound is reached.

\subsection{Lorenz 63}\label{sec:lorenz63}
The now canonical Lorenz-63 model was introduced by Lorenz in 1963 as a first example of ``Deterministic Nonperiodic Flow," what is now known as {\it chaos.}\footnote{Although Lorenz did not coin this term.} Lorenz 63 is defined by a set of three first-order differential equations system~\cite{lorenz}
\begin{eqnarray}
  \dot{x} &=& \sigma(y-x)\\
  \dot{y} &=& x(\rho -z)-y \\
  \dot{z} &=& xy -\beta z
\end{eqnarray}
with the typical chaotic parameter selections: $\sigma=10$, $\rho=28$,
and $\beta=8/3$. 
\begin{figure}[tb!]
\begin{centering}
        \begin{subfigure}[b]{0.32\textwidth}
                \includegraphics[width=\textwidth]{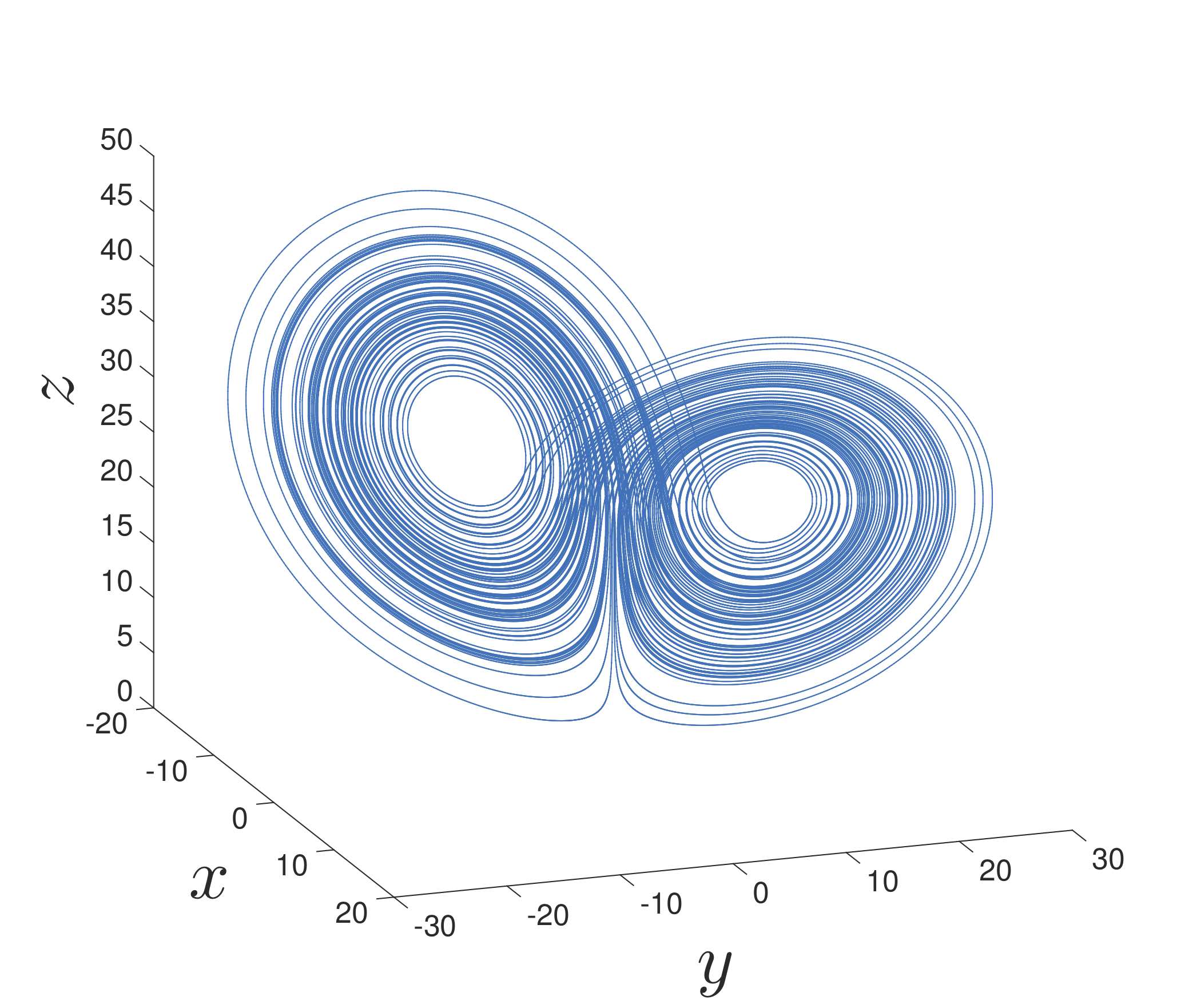}
                \caption{}
        \end{subfigure}
        \begin{subfigure}[b]{0.32\textwidth}
                \includegraphics[width=\textwidth]{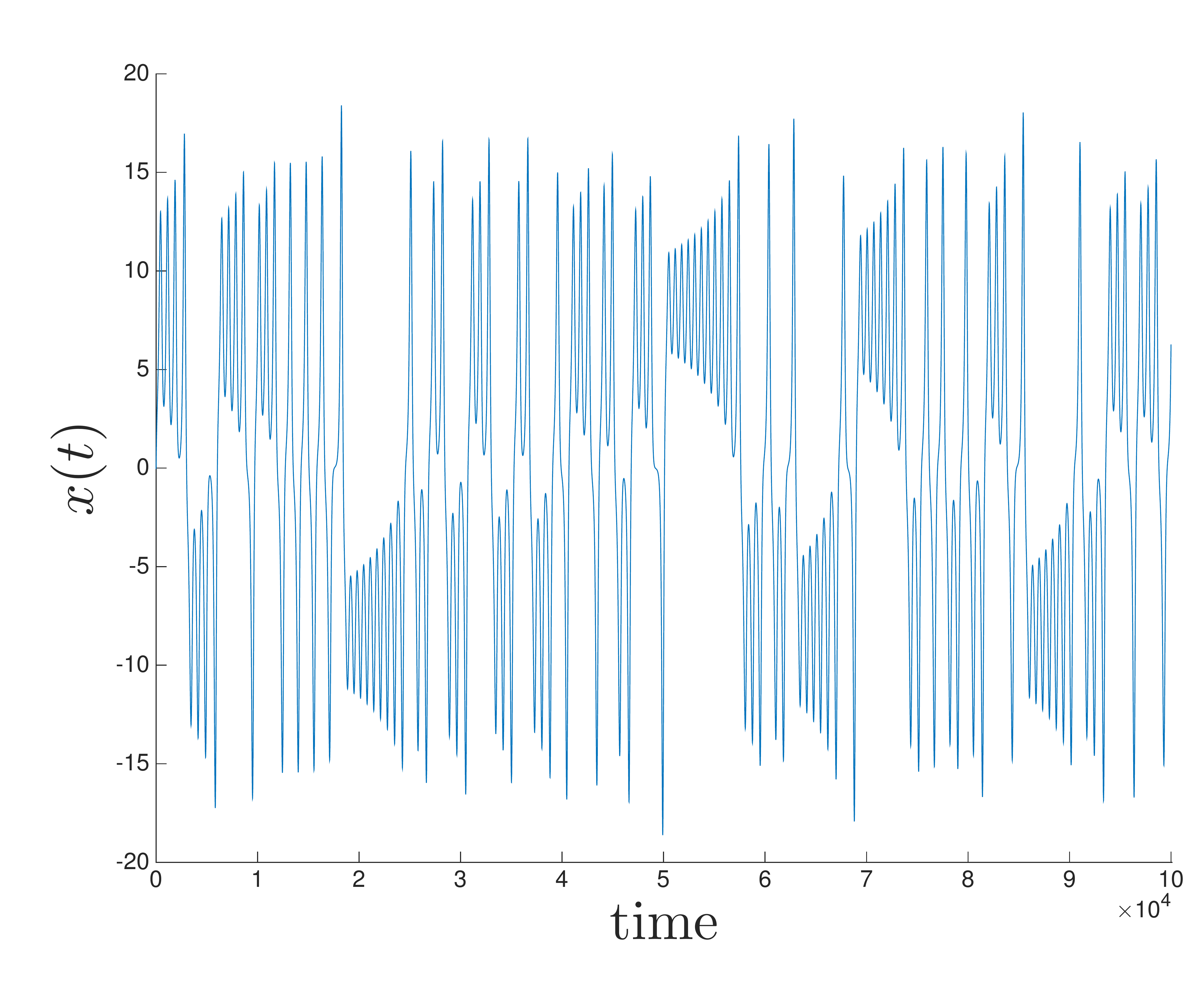}
                \caption{}
        \end{subfigure}
        \begin{subfigure}[b]{0.32\textwidth}
                \includegraphics[width=\textwidth]{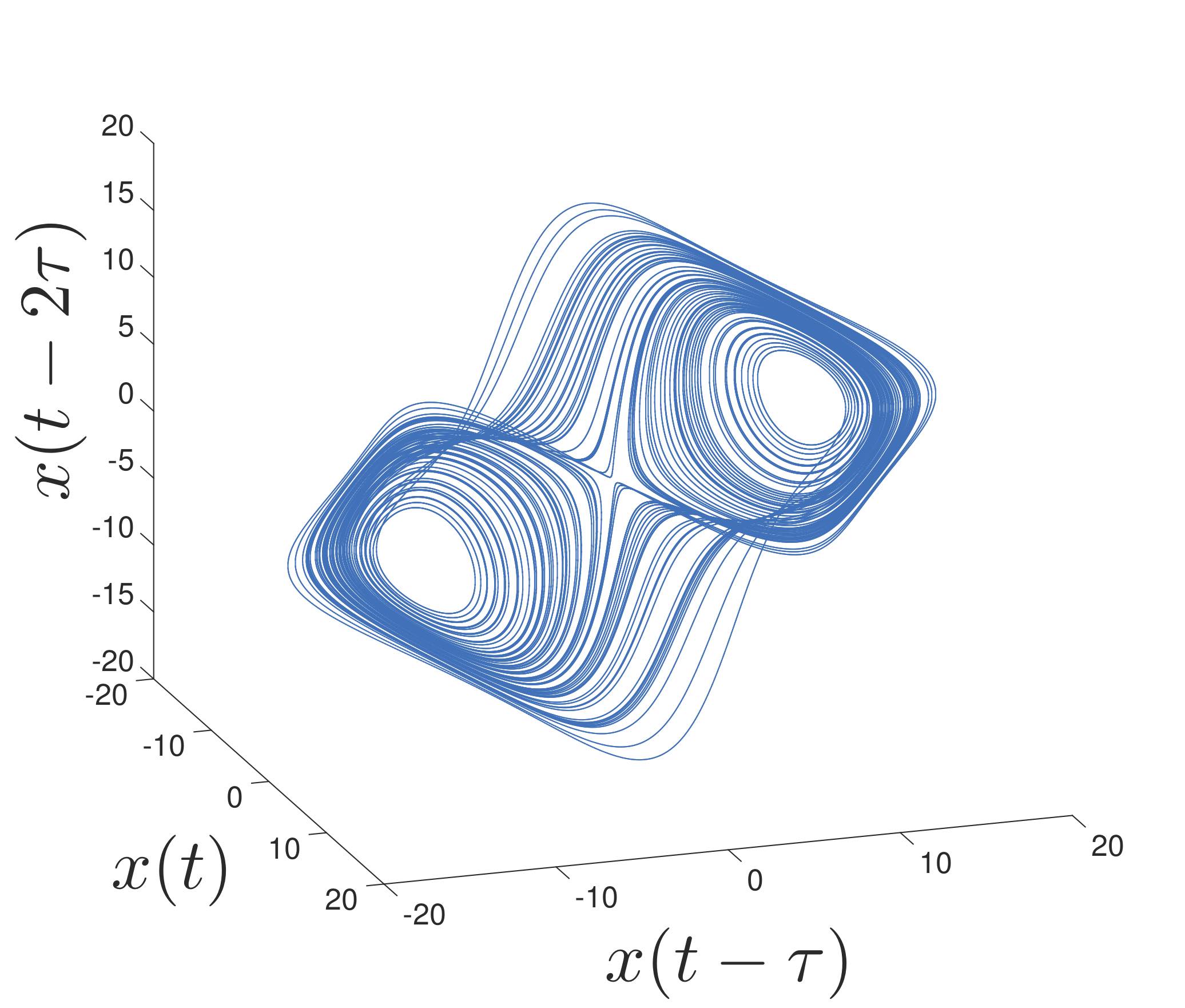}
                \caption{}
        \end{subfigure}
        \end{centering}
\caption{Classic Lorenz attractor ($\sigma=10$, $\rho=28$, $\beta=8/3$): (a)
  A 50,000-point trajectory in $\mathbb{R}^3$ generated using fourth-order
  Runga-Kutta with a time step of $\frac{1}{64}$.  (b) A time-series trace
  of the $x$ coordinate of that trajectory.  (c) A $3D$ projection of
  a delay-coordinate embedding of the trajectory in (b) with dimension $m=5$ and delay $\tau =
  12$.}
\label{fig:Lorenz63Example}
\end{figure}
  Figure~\ref{fig:Lorenz63Example}(a) shows a 50,000-point trajectory in $\mathbb{R}^3$ generated using fourth-order
  Runga-Kutta on those equations with a time step of $T=\frac{1}{64}$, as well as a time-series trace
  of the $x$ coordinate of that trajectory---Figure~\ref{fig:Lorenz63Example}(b)---and a $3D$ projection of
  a delay-coordinate embedding of that trace with dimension $m=5$ and delay $\tau =
  12$, Figure~\ref{fig:Lorenz63Example}(c). This canonical example is used in Chapter~\ref{ch:explain} to establish an explanation of why \roLMA is able to get traction even though it works with a reconstruction that does not meet the theoretical conditions of an embedding.  
  Table~\ref{tab:embedL63} tabulates the estimated and theoretical embedding parameter values for this system.
  \begin{center}
\begin{table}[tb]
\centering
\begin{tabular}{ccccc}
\hline\hline  & $m$-fnn & $\tau$ & $m$-Embedology & $m$-Takens \\ 
\hline  
& 5 & 174 & $\approx 5$ & 7 \\ 
\hline\hline 
\end{tabular}
\caption{Estimated and theoretical embedding parameter values for the Lorenz 63 model, chosen as described in the caption of Table~\ref{tab:embedL96}.}\label{tab:embedL63}
\end{table}
\end{center}

\section{Experimental Case Studies}

Validation with synthetic data is an important first step in
evaluating any new theory, as it provides a controlled, well-defined and well-understood environment. However, these are luxuries, rarely if ever afforded 
to an experimentalist. Furthermore, experimental data often misbehaves more---and in different ways---than synthetic data. For these reasons, it is vital to test a method with experimental time-series
data if one is interested in real-world applications.

\subsection{Computer Performance}\label{sec:compPerfCaseStudy}

As has been established in prior work by our group, it is highly effective to treat computers
as nonlinear dynamical systems \cite{todd-phd,zach-IDA10,mytkowicz09,josh-IDA11,josh-IDA13,joshua-pnp,josh-tdAIS,ISIT13, josh-pre}. In this view, register and memory contents and
physical variables like the temperature of different regions of the processor chip define the state of
the system. The logic hardwired into the computer, combined with the software executing on that
hardware, defines the dynamics of the system. Under the influence of these dynamics, the
state of the processor moves on a trajectory through its high-dimensional state space as the clock cycles progress
and the program executes.
Like Lorenz-96, this
system has been shown to exhibit a range of interesting deterministic
dynamical behavior, from periodic orbits to low- and high-dimensional
chaos~\cite{zach-IDA10,mytkowicz09}, making it a good test case for
this thesis.  It also has important practical implications; these
dynamics, which arise from the deterministic, nonlinear interactions
between the hardware and the software, have profound effects on
execution time and memory use. 

\subsubsection{Theoretical Description}\label{sec:theoryofComputerPerformance}

For the purposes of this thesis, I will consider a ``stored-program
computer,'' {\it i.e.,} a standard von Neumann architecture, as a
deterministic nonlinear dynamical system.  In a stored-program
computer, the current state---both instructions and data---are stored
in some form of addressable memory.
%
%
The contents of this memory are, as established in~\cite{mytkowicz09},
the state space $\mathbb{X}$ of the computer.  Other components of the
computer, such as external memory and video cards, also play roles in
its state.  Those roles depend on the decisions made by the computer
designers---how things are implemented and connected---almost all of
which are proprietary.  In order to distinguish known and unknown
effects, I follow~\cite{mytkowicz09} and define the state space $\mathbb{X}$ as a composition of the
addressable memory elements $\vec{m}$ and the unknown implementation
variables $\vec{u}$:
 \begin{equation}\label{eqn:statespace} \mathbb{X} = \{\vec{\xi} ~|~\vec{\xi} = 
[\vec{m},~ \vec{u}]\} \end{equation} 
The distinction between $\vec{m}$ and $\vec{u}$ is important
because the dynamics of a running computer have two distinct sources:
a map $\vec{F}_{code}$ that acts on the addressable memory $\vec{m}$
directly, as dictated by the program instructions, and a map
$\vec{F}_{impl}$ that captures how the implementation affects the
evolution of the computer state.  The overall dynamics of the
computer---that is, the mapping from its state at the $j^{th}$ clock
cycle to its state at the ${j+1}^{st}$ clock cycle---is a composition
of these two maps:

\begin{equation}\label{eqn:fperf}\vec{\xi}(t_{j+1})
   =\Phi(\vec{\xi}(t_j))=\vec{F}_{P}(\vec{\xi}(t_j)) = \vec{F}_{impl} \circ
   \vec{F}_{code}(\vec{\xi}(t_j))\end{equation} 

\noindent where $\vec{F}_{P}$ is the performance dynamics of the
computer.  An improved design for the processor, for instance---that
is, a ``better'' $\vec{F}_{P}$---is a change in $\vec{F}_{impl}$.  The
form of the map $\vec{F}_{code}$ is dictated by the combination of the
computer's formal specification
({\tt x86\_64}, for the Intel i7 used in the experiments here)
and the software that it is running.  Both $\vec{F}_{impl}$ and
$\vec{F}_{code}$ are nonlinear and deterministic, and their composed
dynamics must be modeled together in order to predict future computer performance.

The framework outlined in the previous paragraph lets me use the
methods of nonlinear dynamics---in particular, delay-coordinate embedding---to model $\vec{F}_{P}$, as long as I observe those dynamics in a way that satisfies the associated theorems (see Section~\ref{sec:dce}). The hardware performance monitor registers (HPMs) that are
built into modern processors can be programmed to count events on
the chip: the
total number of instructions executed per cycle (\ipc), for instance,
or the total number of references to the data cache.  These are some
of the most widely used and salient metrics in the computer
performance analysis
literature~\cite{alameldeen,lebeck02,nussbaum01,tippdirk,sherwood02}.  \ipc is a
good proxy for processor efficiency because most modern microprocessors can execute
more than one instruction per clock cycle.  
While this
metric may not be an element of the state vector $\vec{\xi}$, the
fundamental theorems of delay-coordinate embedding only require that
one measures a quantity that is a smooth, generic function of at least
one state variable.  It was shown in~\cite{zach-IDA10} that the
transformation performed by the HPMs in sampling the state\footnote{This process entails
subtracting successive HPM readings, checking for overflow and
adjusting accordingly.}$\vec{\xi}$
is indeed smooth and generic unless those registers overflow---an
unlikely event, given that they are 64 bits long and that I read them
every $100,000$ instructions. 

The choice of that sample interval is
important for another reason as well.  The HPMs are part of the system under
study, so accessing them can disturb the very dynamics that they are sampling.  This potential  {\it observer problem} was addressed in~\cite{todd-phd} by
varying the sample interval and testing to make sure that the sampling
was not affecting the dynamics.  To further reduce
perturbation, the measurement infrastructure used to gather the data for the experiments reported here only monitors events when the
target program is running, and not when the operating system (or the
monitoring tool itself) have control of the microprocessor. I have completed a careful examination of the impact of interrupt rate on prediction results that corroborates the discussion above; these results are reported in \cite{josh-IDA13}. Finally,
I follow best practices from the computer performance analysis
community~\cite{georges07} when measuring the system: I only use
local disks and limit the number of other processes that are running
on the machine ({\it i.e.}, Linux \texttt{init} level 1). 

The next section describes the experimental observation of this system and the different choices of $\vec{F}_{code}$. For an
in-depth description of this custom-measurement infrastructure,
including a deeper discussion of the implications of the sampling interval,
please see \cite{zach-IDA10,mytkowicz09,todd-phd,josh-IDA13}.

\subsubsection{Experimental methods}\label{sec:compPerfExperiments}

\begin{figure}[tb!]
        \centering
       \vspace{0.5cm} 
\resizebox{.99\textwidth}{!}{%
\includegraphics[height=3cm]{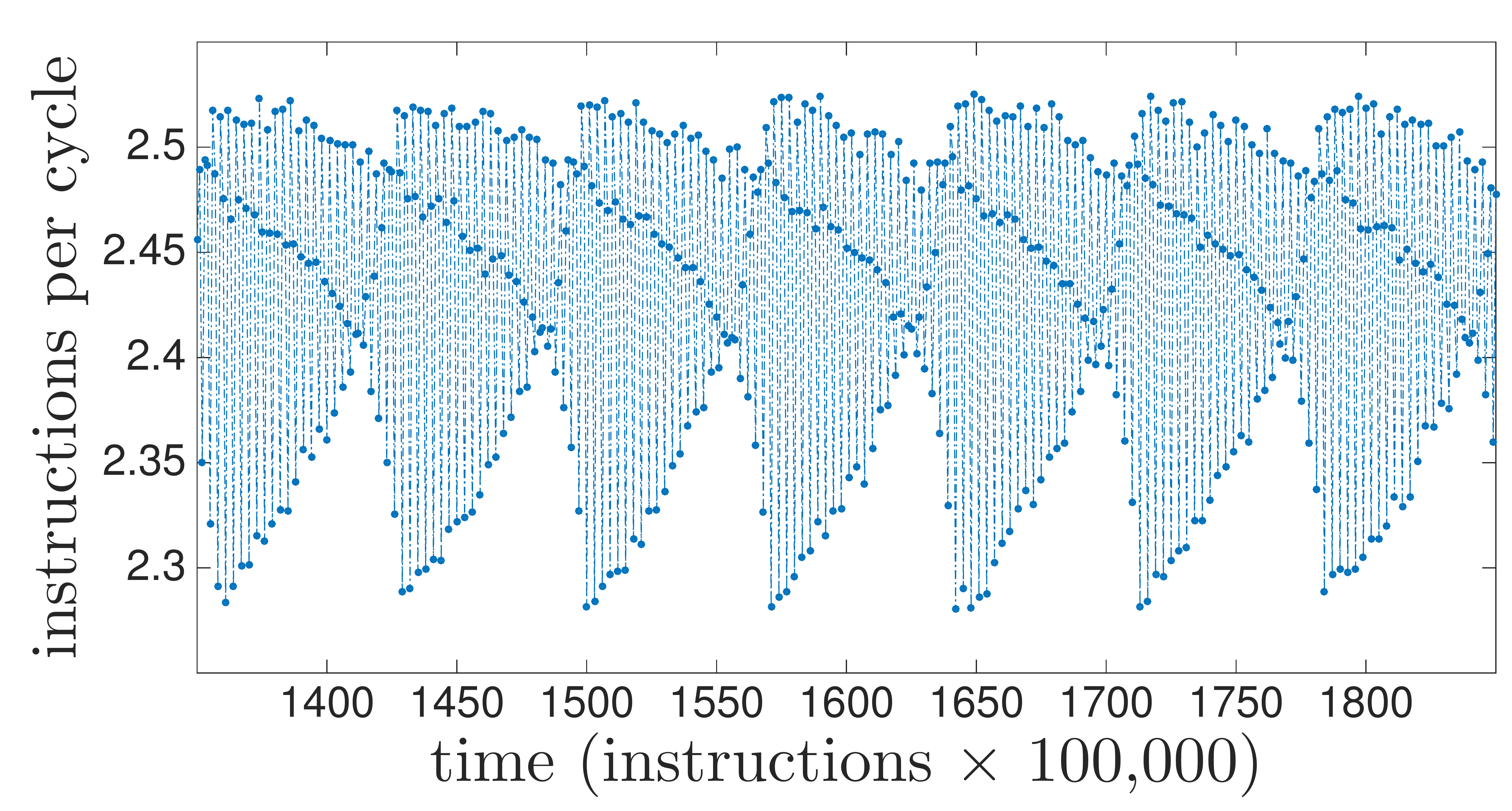}%
\quad
\includegraphics[height=3cm]{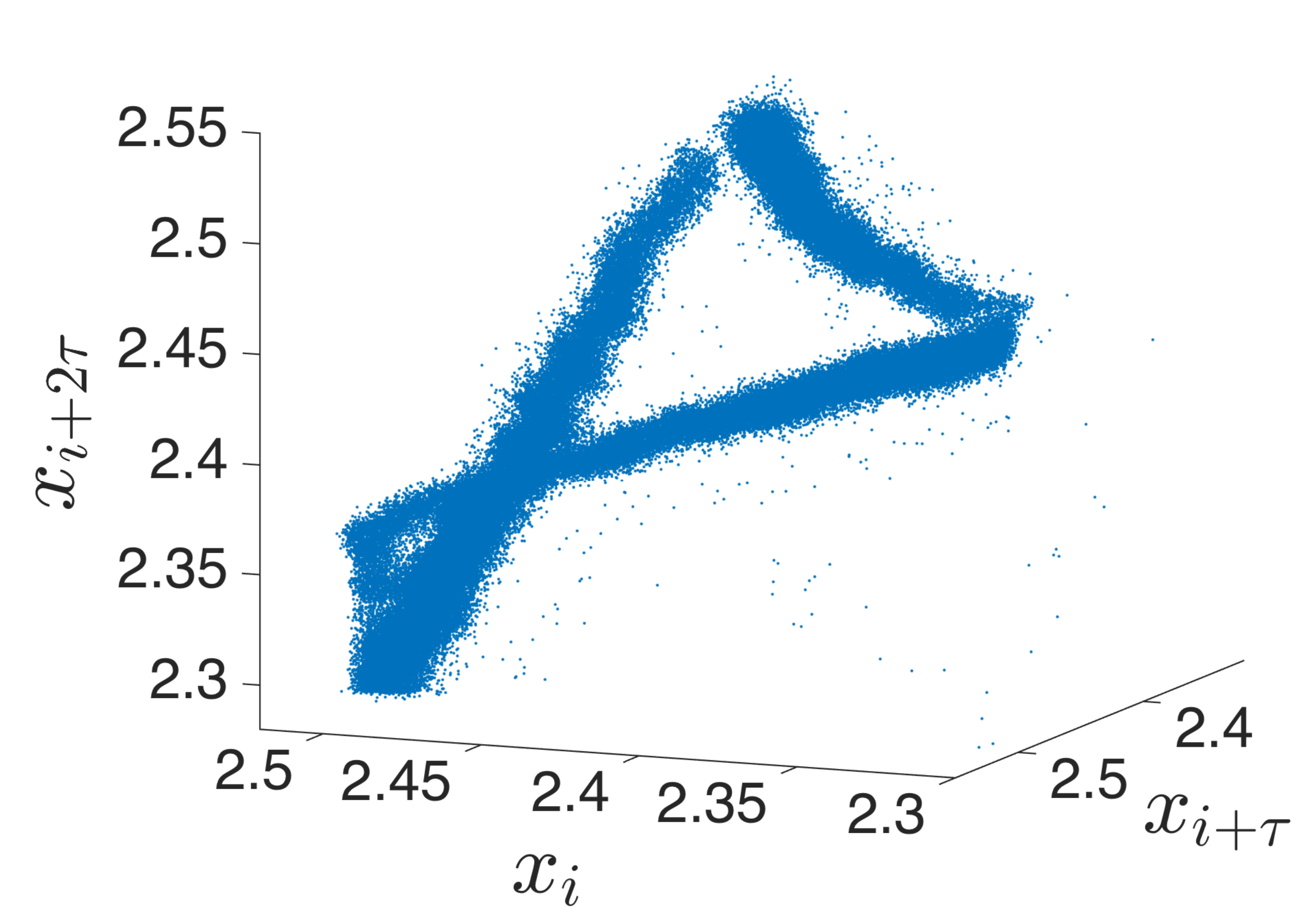}%
}        
        
\caption{ (Left) Time-series data from a computer performance experiment:
processor load traces, in the form of instructions executed per cycle
(\ipc) of a simple program (\col) that repeatedly initializes a
256 $\times$ 256 matrix. Each point is the average \ipc over a 100,000
instruction period.  (Right) A 3D projection of a 12D embedding of this time series.}\label{fig:col12DEmbedding}
\end{figure}

The computer performance time-series data sets for the experiments presented in this thesis were collected on an
Intel Core\textsuperscript{\textregistered} i7-2600-based machine
running the 2.6.38-8 Linux kernel. I also carried out experiments on an Intel Core2 Duo. Those Core2 results, reported in \cite{josh-IDA11} but omitted here, are consistent with the results reported in this dissertation.  This i7 microprocessor
chip has eight processing units, a clock rate of 3.40 GHz, and a cache
size of 8192 KB.  The experiments in this thesis involve performance traces gathered during the execution
of several different programs, beginning with the simple \col loop whose performance
is depicted in the left panel of Figure~\ref{fig:col12DEmbedding}, 
as well as a more-complex
program from the SPEC 2006CPU benchmark suite (\gcc) \cite{spec2006}.  In addition, I also carried out careful analysis of standard computer performance benchmarks programs such as \sphinx \cite{spec2006}, linear algebra software from LAPACK (Linear Algebra PACKage) such as \svd and \eig~\cite{lapack} and  \row (the row-major analogue of \col). Many of these experiments are omitted here for brevity; please see \cite{josh-IDA11,josh-IDA13,ISIT13, josh-pre} for these companion results. I select \col and \gcc from this larger constellation of experiments for the discussion in my thesis because they are informative in their own right and representative of the other results I encountered; \col is a simple highly-structured chaotic time series, while \gcc is a chaotic time series where almost all structure has been consumed by noise. 

In all of these experiments,
the scalar observation $x_j$ is a measurement of the processor
performance at time $j$ during the execution of each program.  To
record these measurements, I use the {\tt libpfm4} library, via PAPI (Performance Application Programming Interface)  
5.2~\cite{papi}, to stop program execution at 100,000-instruction
intervals---the unit of time in these experiments---and read off the contents of
the CPU's onboard hardware performance monitors, which I programmed to count how
many instructions are executed in each clock cycle(\ipc). I also recorded and analyzed other metrics including total L2 cache misses, missed branch predictions, and L2 instruction cache hits. Description of these metrics, as well as corresponding analysis, are published in \cite{josh-IDA11,josh-IDA13}. In this thesis, I only report results on \ipc , as it is representative of all of these results.    For statistical
validation, I collect 15 performance traces from each of the 
programs.  These traces, and the processes that generated them, are
described in more depth in the rest of this section.

\col is a simple C program that repeatedly initializes the upper
triangle of a 2048 $\times$ 2048 matrix in column-major order by
looping over the following three lines of code:
\begin{lstlisting}
  for(i=0;i<2048;i++)
  	for(j=i;j<2048;j++)
  		data[j][i] = 0;
\end{lstlisting}
%
%
%
\noindent As mentioned in Chapter~\ref{ch:overview} and shown in Figure~\ref{fig:col12DEmbedding}, this simple program
exhibits surprisingly complicated behavior. I also collected data from the row-major analogue to \col. These time series were very different than \col, but the forecasting results were largely the same, so they are omitted here but can be found in \cite{ISIT13,josh-IDA11}.


\begin{figure}[tb!]
  \centering
  \vspace{0.5cm}
  \resizebox{.99\textwidth}{!}{%
\includegraphics[height=3cm]{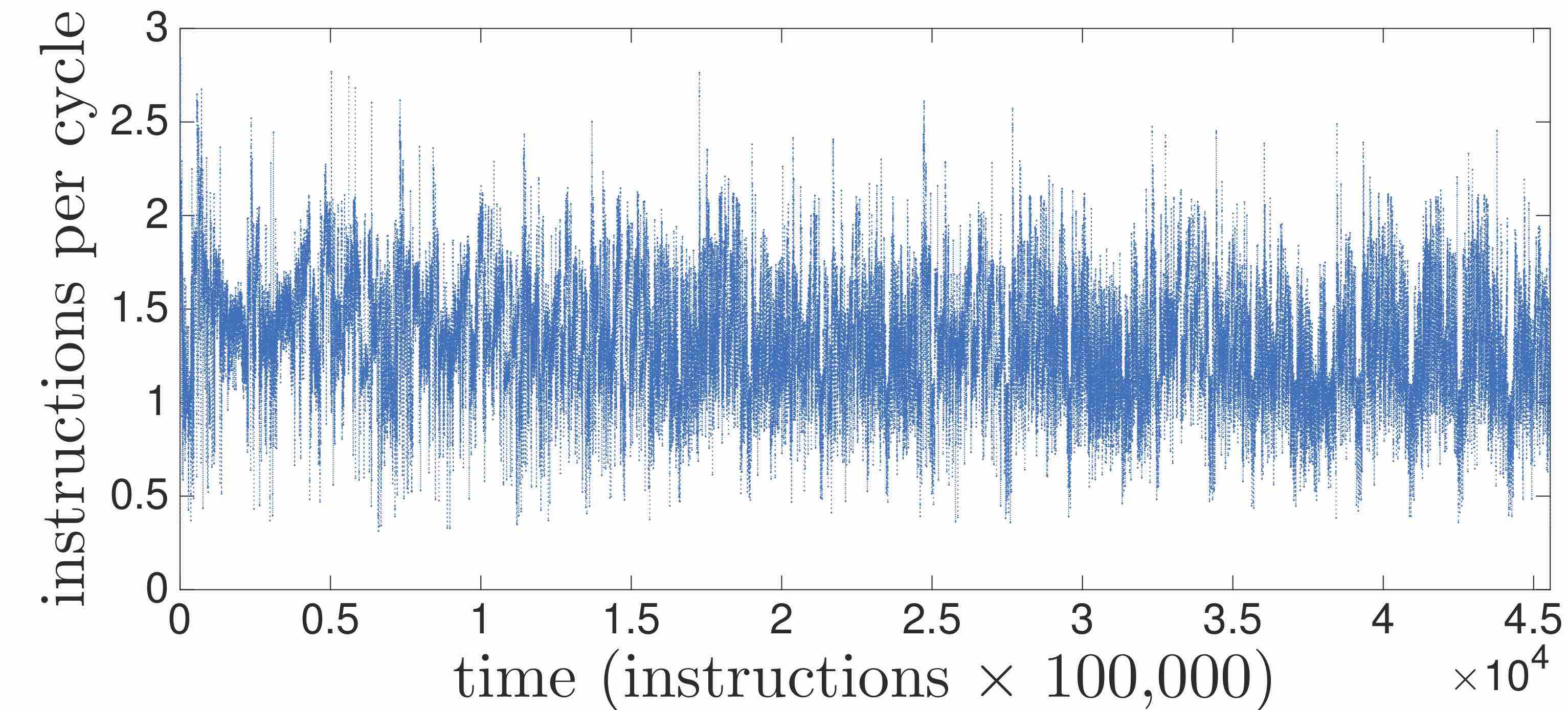}%
\quad
\includegraphics[height=3cm]{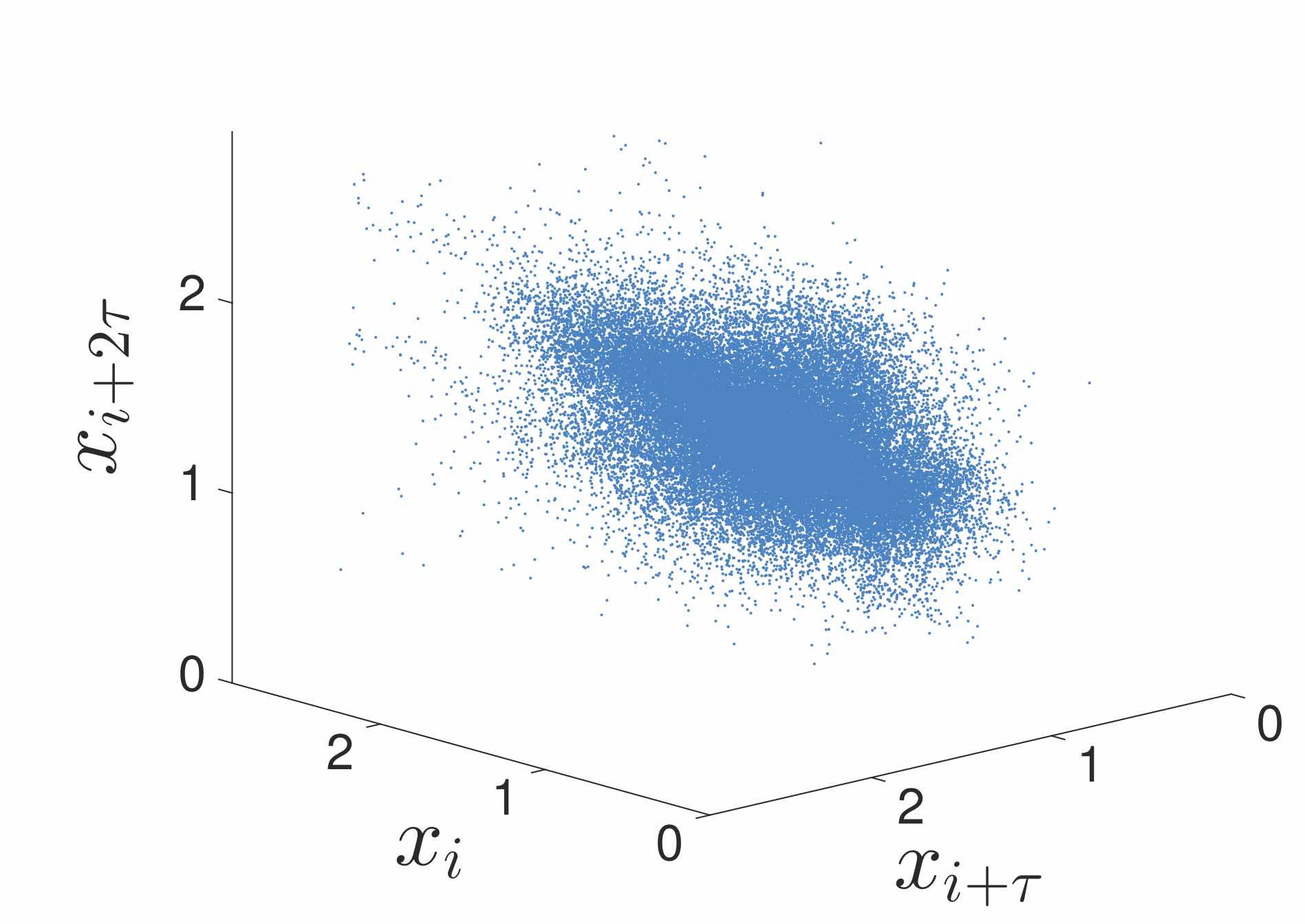}%
}        
  \caption{Processor load traces (\ipc) of the SPEC benchmark \gcc. Each point is the average IPC in a 100,000
      instruction period. }\label{fig:gccipc}
\end{figure}
The SPEC CPU2006 benchmark suite \cite{spec2006} is a collection of
complicated programs that are used in the computer-science community
to assess and compare the performance of different computers.  \gcc is
a member of that suite.
It is a {\it compiler}: a program that translates code written in a
high-level language (C, in the case of \gcc) into a lower-level format
that can be executed by the processor chip.  Its behavior is far more
complicated than that of \col, as is clear from
Figure~\ref{fig:gccipc}. 
Unlike \col, where the processor utilization is quite structured, the performance of \gcc appears almost random. In addition to \gcc, I also studied \sphinx from this benchmark suite: a speech-recognition tool~\cite{spec2006}. \sphinx and the associated results are covered in more depth in~\cite{ISIT13,josh-IDA13}.

Table~\ref{tab:embedcompperf} tabulates the estimated embedding parameters for \col and \gcc. Notice that because these dynamical systems are not understood from a theoretical perspective, {\it i.e.}, the governing equations or knowledge of the state space dimension are unknown, I must rely on the heuristics presented in Section~\ref{sec:numericalM}. I should note that it has been estimated that the state space of these systems is at least $2^{32}$ dimensions~\cite{mytkowicz09}, which would mean that $m$-Takens$>2^{33}$. However, the same paper suggests the actual fractal dimension of these dynamics are much smaller, due in part to standard programming and design principles, which have the effect of reducing the dimension of the dynamics, and that $m$-Embedology is probably less that ten.
\begin{table}[h!]
\centering
\begin{tabular}{ccccc}
\hline\hline  & $m$-fnn & $\tau$ & $m$-Embedology & $m$-Takens \\ 
\hline  
\col & 12 & 2 & $\alert{**}$ & \alert{**} \\ 
\gcc & 13 & 10 & $\alert{**}$ & \alert{**} \\
\hline\hline 
\end{tabular}
\caption{Estimated embedding parameter values for the computer performance experiments.  $\tau$ and $m$-fnn chosen as in Table~\ref{tab:embedL96}. $m$-Embedology and $m$-Takens are not provided as these dimensions are unknown for a experimental system like this one.}\label{tab:embedcompperf}
\end{table}

%% file: pnp.tex

\chapter{Prediction in Projection}\label{ch:pnp}

In this chapter, I demonstrate that the accuracies of forecasts
produced by \roLMA---Lorenz's method of analogues, operating on a
two-dimensional time-delay reconstruction of a trajectory from a
dynamical system---are similar to, and often better than, forecasts
produced by \fnnLMA, which operates on an embedding of the same
dynamics. While the brief example in Chapter~\ref{ch:overview} is a
useful first validation of that statement, it does not support the
kind of exploration that is necessary to properly evaluate a new
forecast method, especially one that violates the basic tenets of
delay-coordinate embedding.  The SFI dataset A is a single trace
from a single system---and a low dimensional system at that.  My goal in this chapter is to show that \roLMA is comparable to or better than \fnnLMA for a {\it range} of systems and parameter values---and to repeat each
experiment for a number of different trajectories from each system. This exploration serves as an experimental validation of the central premise of this thesis. And of course, any discussion of new forecasting strategies is incomplete without a solid comparison with traditional methods. To this end, I present results for two dynamical systems, one simulated and one
real: the Lorenz-96 model and sensor data from a laboratory experiment
on computer performance dynamics. I produce \roLMA forecasts of these systems and compare them to forecasts using the four traditional strategies presented in Section~\ref{sec:forecastmodels}.

\section{A Synthetic Example: Lorenz-96}\label{sec:roLMALorenz96}

In this example, I perform two sets of forecasting experiments with ensembles of
traces from the Lorenz-96 model~\cite{lorenz96Model}, introduced in Section~\ref{sec:lorenz96}: one with $K=22$ and the other with
$K=47$.

\subsection{Comparing \roLMA and \fnnLMA}
\begin{figure}[tb!]
        \centering
        \begin{subfigure}[b]{\columnwidth}
                \includegraphics[width=0.98\textwidth]{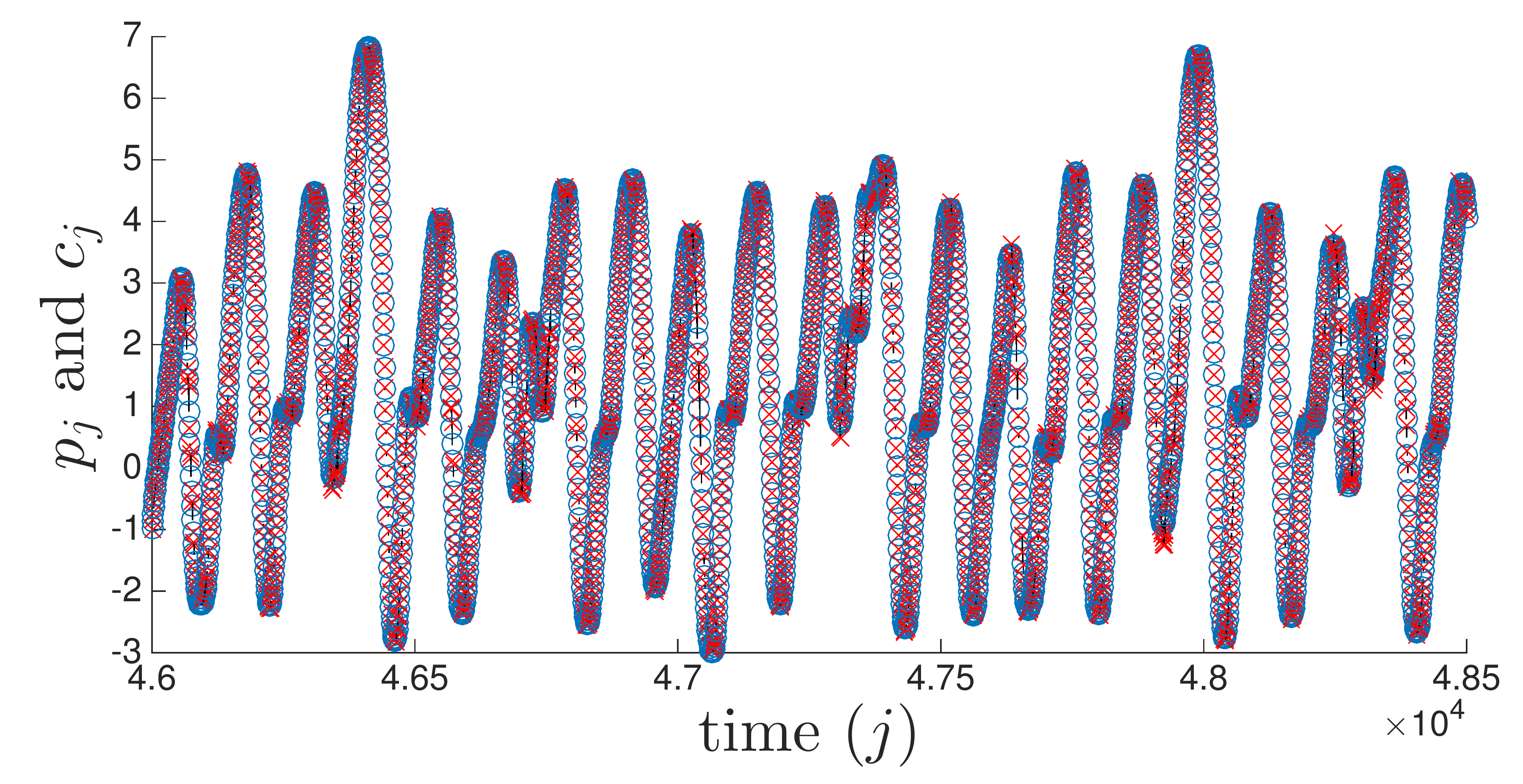}

                \caption{2,500-point forecast using the reduced-order
                forecast method \roLMA. Blue circles and red
                {\color{red}$\times$}s are the true and predicted
                values, respectively; vertical bars show where these
                values differ.}

                \label{fig:K22m2k1per10tspred}
        \end{subfigure}%
        ~ 
          
        \begin{subfigure}[b]{0.49\columnwidth}
                \includegraphics[width=\columnwidth]{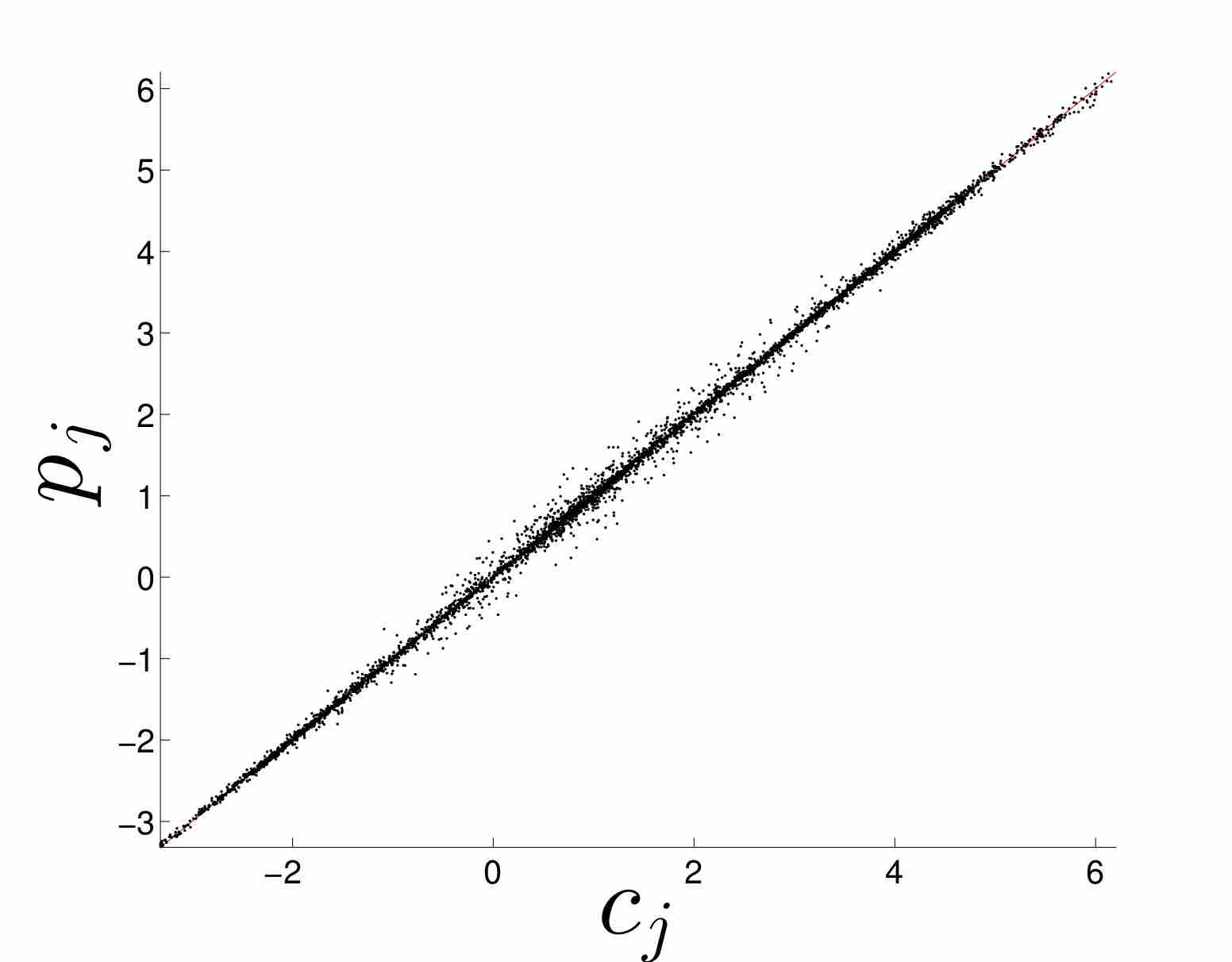}
                \caption{\roLMA forecast}
                \label{fig:K22m2k1per10pjcjpred}
        \end{subfigure}
                \begin{subfigure}[b]{0.49\columnwidth}
                \includegraphics[width=\columnwidth]{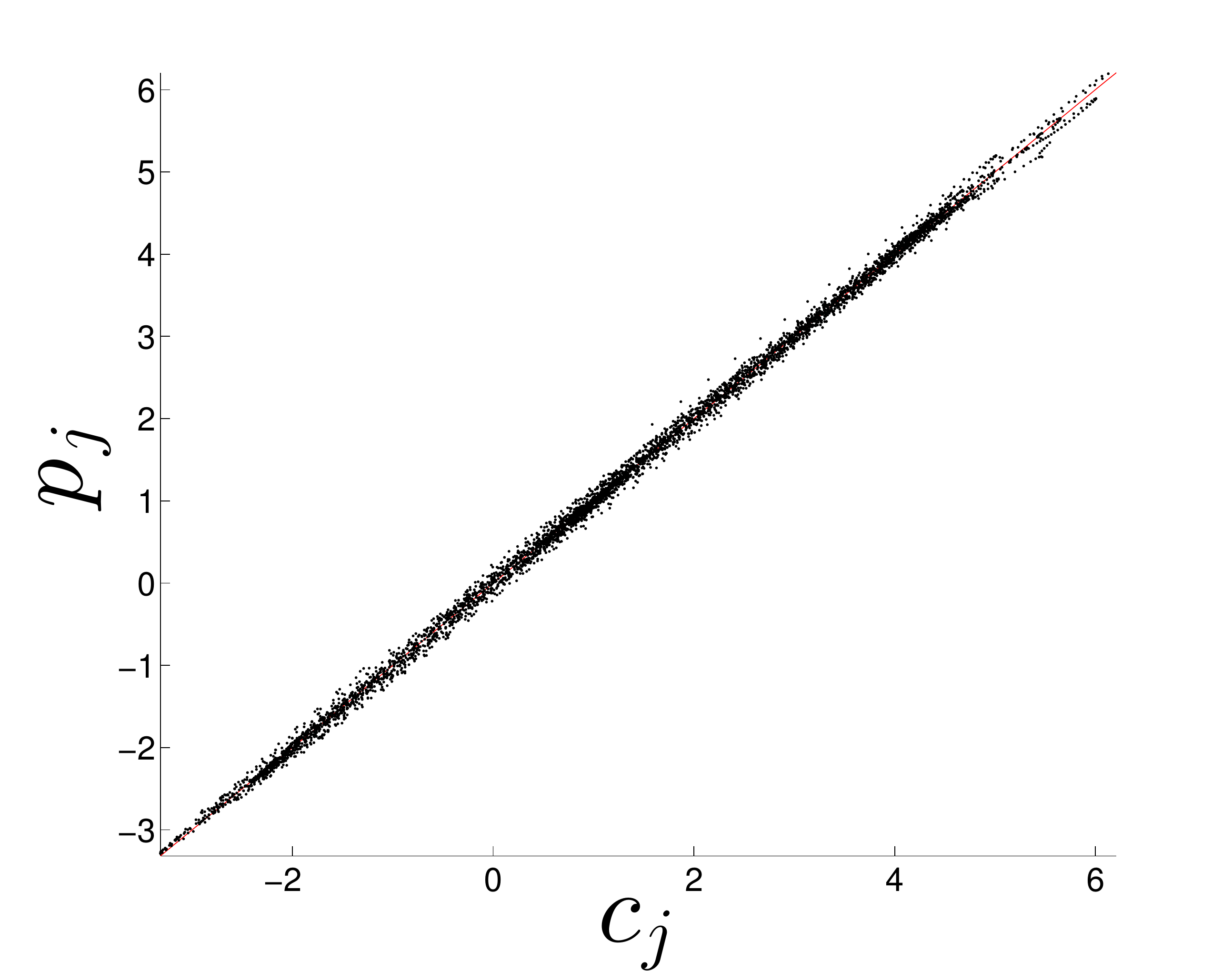}
                \caption{\fnnLMA forecast}
                \label{fig:K22m8k1per10pjcjpred}
        \end{subfigure}

\caption{\roLMA and \fnnLMA forecasts of a representative trace from the
Lorenz-96 system with $K=22$ and $F=5$.  Top: a time-domain plot of
the first 2,500 points of the \roLMA forecast.  Bottom: the predicted
($p_j$) vs true ($c_j$) values for forecasts of that trace generated
by (b) \roLMA and (c) \fnnLMA.  On such a plot, a perfect prediction
would lie on the diagonal.  The 1-MASE scores of the forecasts in (b)
and (c) were 0.392 and 0.461, respectively.
}\label{fig:K22predictions}
\end{figure}

As I will illustrate in the following discussion, both \roLMA and \fnnLMA worked quite well for the $K=22$ dynamics.
See Figure~\ref{fig:K22predictions}(a) for a time-domain plot of
an \roLMA forecast of a representative trace from this system and
Figures~\ref{fig:K22predictions}(b) and (c) for graphical
representations of the forecast accuracy on that trace for both
methods. In Figures~\ref{fig:K22predictions}(b) and (c), the vertical axis is the prediction $p_j$ and the horizontal axis is the true continuation $c_j$. On this type of plot, a perfect prediction would lie on the diagonal. 
The diagonal structure on the $p_j$ vs. $c_j$ plots in the Figure
indicates that both of these LMA-based methods perform very well on
this trace.  More importantly---from the standpoint of evaluation of
my primary claim---the LMA forecasting strategy worked {\it better} on a
two-dimensional reconstruction of these dynamics than on a full
embedding, and by a statistically significant margin: the 1-MASE scores\footnote{one-step ahead Mean Absolute Scaled Error}of \roLMA and \fnnLMA forecasts,
computed following the procedures described in
Section~\ref{sec:accuracy}, are $0.391 \pm 0.016$ and $0.441 \pm
0.033$, respectively, across the 330 traces at this parameter value.
 This is somewhat startling, given that the {two-dimensional delay
reconstruction used by \roLMA} falls far short of the 
requirement for topological conjugacy~\cite{sauer91} in this system.
Clearly, though, it captures {\it enough} structure to allow LMA to
generate good predictions.

The $K=47$ case is a slightly different story:  here, \roLMA still
outperforms \fnnLMA, but not by a statistically significant margin.
The 1-MASE scores across all 329 traces were $0.985 \pm 0.047$ and
$1.007 \pm 0.043$ for \roLMA and \fnnLMA, respectively.  In view of
the higher complexity of the state-space structure of the $K=47$
version of the Lorenz-96 system, the overall increase in 1-MASE scores over
the $K=22$ case makes sense.
Recall that $d_{KY}$ is far higher for the $K=47$ case: this attractor
fills more of the state space and has many more manifolds that are
associated with positive Lyapunov exponents.  

This has obvious
implications for predictability. Since I use the same traces for both methods, one might be 
tempted to think that the better performance or \roLMA is a predictable
consequence of data length---simply because filling out a
higher-dimensional object like the reconstruction used by the \fnnLMA model requires more data.
When I re-run the experiments with longer traces, the 1-MASE scores
for \roLMA and \fnnLMA did converge, but not until the traces are over
$10^6$ points long, and at the (significant) cost of near-neighbor
searches in a space with many more dimensions.  Note, too, that the
longer delay vectors used by \fnnLMA span far more of the training
set, which at first glance would seem to be a serious advantage from an information-theoretic
standpoint (although, as shown later in Section~\ref{sec:tdAIS}, this is not always an advantage).  In view of this, the comparable performance of \roLMA is quite impressive. All of these issues are explored at more length in
Section~\ref{sec:time-scales}.

\subsection{Comparing \roLMA with Traditional Linear Methods} 

For the $K=22$ time series, the LMA-based methods do significantly better than the \naive and ARIMA methods and about twice as well as the random walk method and this makes sense. Each point in the time series of Figure~\ref{fig:K22predictions}(a) is very close to its predecessor and successor, which plays to the strengths of random walk. In contrast, the oscillations of the signal and the inertia of the mean make the \naive method ineffective.  The fact that the LMA-based methods outperform the random walk at all, let alone twice as well, is quite impressive. Successive points of this signal are so close together that it is really an ideal candidate for random walk forecasting, leaving little room for another method to be more successful. However, both of the LMA-based techniques successfully meet this challenge: about 2.5 times and 1.8 times better than random walk, respectively, for \roLMA and \fnnLMA. See Figures~\ref{fig:K22predictions}(b) and ~\ref{fig:K22predictions}(c) for a visual comparison. 

$K=47$ is a very similar story; all of the LMA-based methods outperform \naive and ARIMA by several orders of magnitude for the same reasons discussed in the previous paragraph. The comparison between the LMA methods and random walk is more interesting. Both \roLMA and \fnnLMA MASE scores are almost identical to those of random walk forecasts. For the reasons discussed above, this is not surprising; random walk is very well suited for this signal leaving a very small margin to be outperformed. Table~\ref{tab:lorenz96error} compares the forecast accuracy of all Lorenz-96 time series with each of the methods discussed above.

\begin{table*}[tb]
\caption{The average 1-MASE scores of all four forecast methods for the two ensembles of Lorenz-96 time series.}
  \begin{center}
  \begin{tabular*}{\textwidth}{@{\extracolsep{\fill} } cccccc}
  \hline\hline  Parameters & \roLMA  & \fnnLMA  &ARIMA& na\"{i}ve \\ 
  \hline
$\{K=22, F=5\}$    &$0.391\pm0.016$& $0.441\pm0.033$  & $17.031\pm0.310$ &$17.006\pm  0.233$ \\
 $\{K=47, F=5\}$   &$0.985\pm0.047$ & $1.007\pm0.043$  & $18.330\pm0.583$  & $17.768\pm 0.765$ \\
   \hline\hline
  \end{tabular*}
  \end{center}
 \label{tab:lorenz96error}
  \end{table*}%

\section{Experimental Data: Computer Performance Dynamics}
\label{sec:compPerfProj}

Validation with synthetic data is an important first step in
evaluating any new forecast strategy, but experimental time-series
data are the acid test if one is interested in real-world
applications.  My second set of tests of \roLMA, and comparisons of
its accuracy to that of traditional forecast strategies,  involves data from the laboratory
experiment on computer performance dynamics that was introduced in Section~\ref{sec:compPerfCaseStudy}.  

I have tested \roLMA on traces of many different processor and memory
performance metrics gathered during the execution of a variety of
programs on several different computers (see {\it e.g.,} \cite{josh-IDA11,josh-IDA13, my-masters,ISIT13, josh-pre}).  Here, for conciseness, I
focus on {\it processor} performance traces from two different
programs, one simple (\col) and one complex (\gcc), running on the same Intel
i7-based computer. As discussed in Section~\ref{sec:compPerfCaseStudy}, computer performance dynamics result from a composition of
hardware and software. These two programs represent two different
dynamical systems, even though they are running on the same computer.
The dynamical differences are visually apparent from the traces in
Figures~\ref{fig:col12DEmbedding} and~\ref{fig:gccipc}; they are also mathematically apparent from
nonlinear time-series analysis of embeddings of those
data~\cite{mytkowicz09}, as well as in calculations of the information
content of the two signals.  Among other things, \gcc has much less
predictive structure than \col and is thus much harder to
forecast~\cite{josh-pre}.  These attributes make this a useful pair of
experiments for an exploration of the utility of reduced-order
forecasting.

\begin{figure}[tb!]
  \centering
  \vspace*{0.2cm}
    \includegraphics[width=\columnwidth]{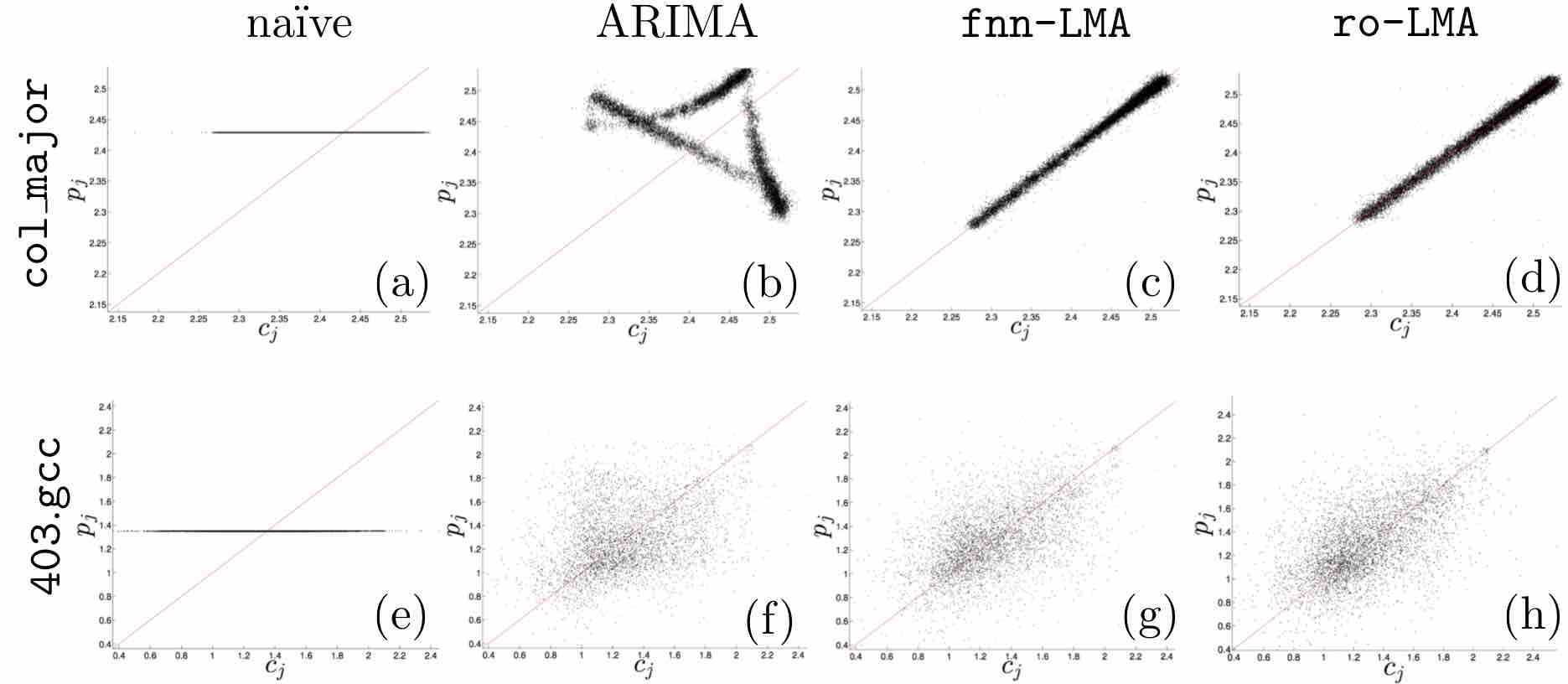}
    \caption{Predicted ($p_j$) versus true values ($c_j$) for \col and \gcc generated with four of the forecast methods considered in this thesis. 
}\label{fig:forecast-example}
\end{figure}

For statistical validation, I collect 15 performance traces from
the computer as it ran each program, calculated embedding parameters as described in Section~\ref{sec:dce}, and generated forecasts
of each trace using \roLMA and the traditional methods outlined in Section~\ref{sec:forecastmodels}.
Figure~\ref{fig:forecast-example} shows some representative examples.
Recall that on such a plot, a perfect prediction would lie on the diagonal.  Horizontal lines result when a constant predictor ({\it e.g.},
\naive) is used on a non-constant signal. In the case of Figure~\ref{fig:forecast-example}, \fnnLMA and \roLMA
both generate very accurate predictions of the \col trace, while ARIMA
does not. Note that the
shape of Figure~\ref{fig:forecast-example}(b)
(ARIMA on \col) is reminiscent of the
projected embedding in the right panel of Figure~\ref{fig:col12DEmbedding}. This structure is also present in a $p_j$ vs. $c_j$ plot of a random-walk forecast (not shown) on this same signal.    Indeed, for a
random-walk predictor, a $p_j$ vs. $c_j$ plot is technically
equivalent to a two-dimensional embedding with $\tau=1$.  For ARIMA,
the correspondence is not quite as simple, since the $p_j$ values are
linear combinations of a number of past values of the $c_j$, but the
effect is largely the same.

 The 1-MASE scores
for \roLMA and \fnnLMA across all 15 trials in this set of experiments were $0.050 \pm 0.002$
and $0.063 \pm 0.003$, respectively; ARIMA scored much worse ($0.599\pm0.211$). This difference in performance is not surprising; the \col time series 
contains plenty of nonlinear structure that the LMA-based methods can capture and utilize, whereas ARIMA can not. 
These 1-MASE scores mean that both \fnnLMA and \roLMA perform roughly 20 times better on \col than a
random-walk predictor, while ARIMA only outperform random walk by a
factor of 1.7. This is in accordance with the visual appearance of
the corresponding images in Figure~\ref{fig:forecast-example}.
 For \gcc, however, \roLMA is somewhat more accurate: 1-MASE scores of $1.488 \pm 0.016$
versus \fnnLMA's $1.530 \pm 0.021$.  Note that the \gcc 1-MASE scores
are higher for both forecast methods than for \col, simply because
the \gcc signal contains less predictive structure~\cite{josh-pre}.
This actually makes the comparison somewhat problematic, as discussed
at more length in Section~\ref{sec:time-scales}.

Comparing \roLMA to the the \naive method is illustrative. 
\roLMA does significantly better on all signals but \gcc. This is reassuring as \gcc has very high complexity, almost no redundancy, and very little predictive structure\cite{josh-pre}. With signals like this, simple forecast methods that do not rely on predictive structure tend to do very well; this is discussed in more depth in Chapter~\ref{ch:wpe}.

Table~\ref{tab:comp-perf-error} summarizes all of the computer performance experiments presented in this discussion.  Overall, these results are consistent with the Lorenz-96 example in
the previous section: prediction accuracies of \roLMA and \fnnLMA are
quite similar on all traces, despite the former's use of a
theoretically incomplete reconstruction.  This amounts to a
validation of the conjecture on which this thesis is based.  And in
both numerical and experimental examples, \roLMA
actually {\it outperform} \fnnLMA on the more-complex traces (\gcc,
$K=47$).  I believe that this is due to the noise mitigation that is
naturally effected by a lower-dimensional reconstruction.

\begin{table}[tb!]
\caption{
The average 1-MASE scores of all four forecast methods for the 15 trials of \col and \gcc.}
\centering
  \begin{tabular*}{\textwidth}{@{\extracolsep{\fill} } cccccc}
  \hline\hline Signal & \fnnLMA MASE & \roLMA MASE  & ARIMA MASE & na\"{i}ve MASE  \\ 
  \hline
 \col           & $ 0.050 \pm0.002  $ &$0.0625\pm0.0032$& $0.599  \pm 0.211 $ & $0.571\pm0.002$ \\

\gcc           & $ 1.530\pm 0.021$ &$1.4877\pm0.016$& $1.837 \pm0.016 $ & $0.951 \pm 0.001$ \\
  \hline\hline
  \end{tabular*}
 \label{tab:comp-perf-error}
  \end{table}%

%
%

\section{Time Scales, Data Length and Prediction Horizons}
\label{sec:time-scales}
In this Section, I explore the effects of the values of the $\tau$ parameter (Section~\ref{sec:varyingproj}), prediction horizon (Section~\ref{sec:hforecasts}) and data length (Section~\ref{sec:datalength}) on \roLMA. For the remainder of this chapter, I discontinue comparing \roLMA to traditional linear methods, as that comparison would not add anything to the discussion, and instead focus on the direct comparison between \roLMA and \fnnLMA. 
%


\subsection{The $\tau$ Parameter}\label{sec:varyingproj}

The embedding theorems require only that $\tau$ be greater than zero
and not a multiple of any period of the dynamics.  In practice,
however, $\tau$ can play a critical role in the success of
delay-coordinate reconstruction---and any nonlinear time-series
analysis that follows~\cite{fraser-swinney,kantz97,rosenstein94,joshua-pnp}. It
follows naturally, then, that $\tau$ might affect the accuracy
of an LMA-based method that uses the structure of a {time-delay
reconstruction} to make forecasts.

\begin{figure}[bt!]
        \centering

        ~ 
          
        \begin{subfigure}[b]{0.49\textwidth}
                \includegraphics[width=\textwidth]{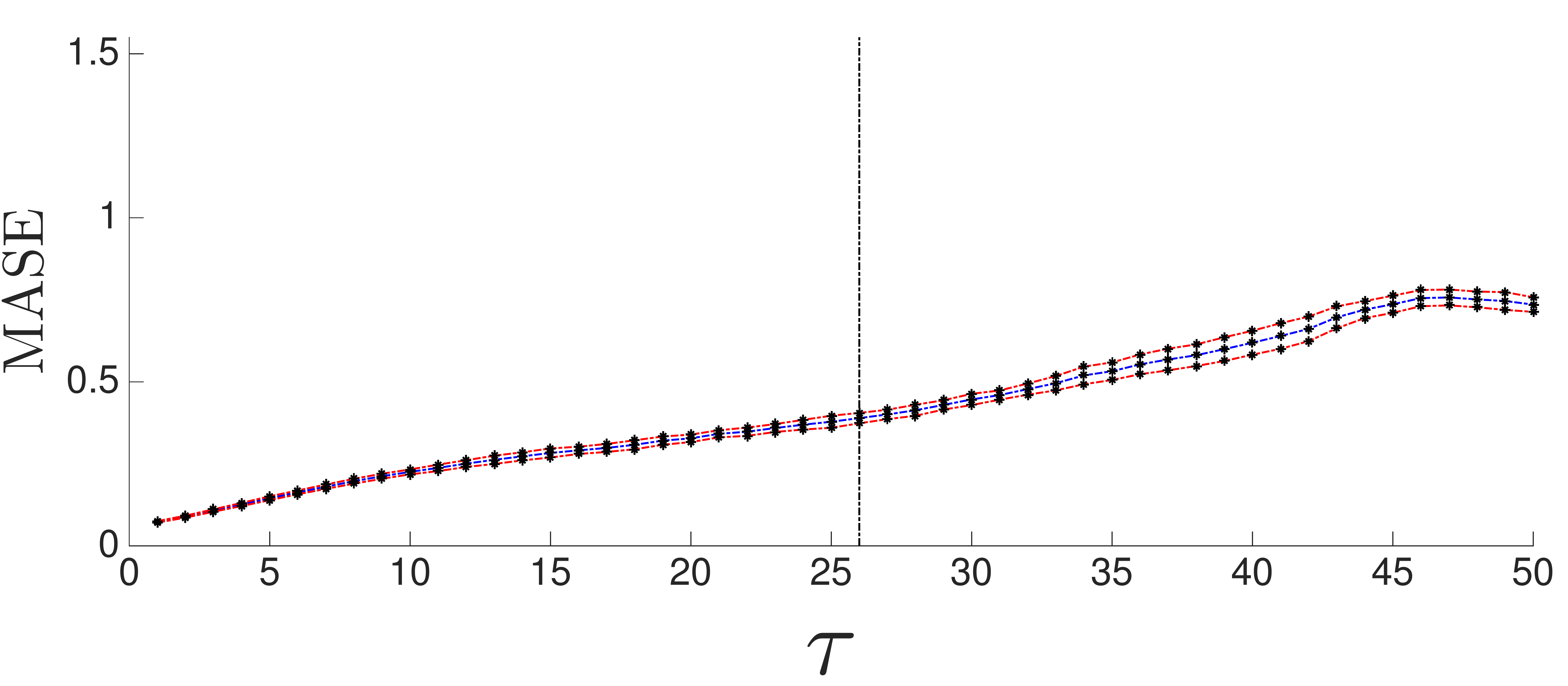}
                \caption{\roLMA on Lorenz-96 with $K=22$. }
        \end{subfigure}
\begin{subfigure}[b]{0.49\textwidth}
                \includegraphics[width=\textwidth]{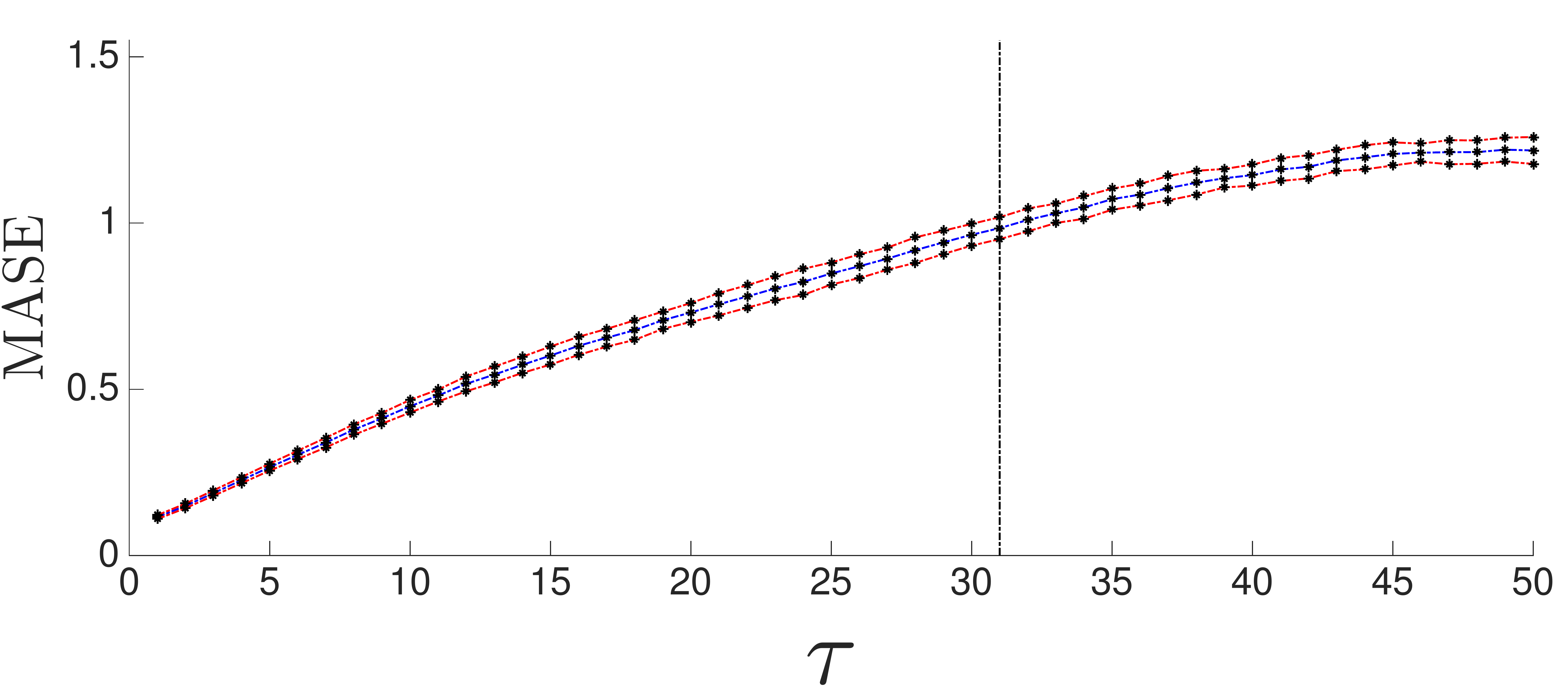}
                \caption{\roLMA on Lorenz-96 with $K=47$. }
        \end{subfigure}%
        
\begin{subfigure}[b]{0.5\textwidth}
                \includegraphics[width=\textwidth]{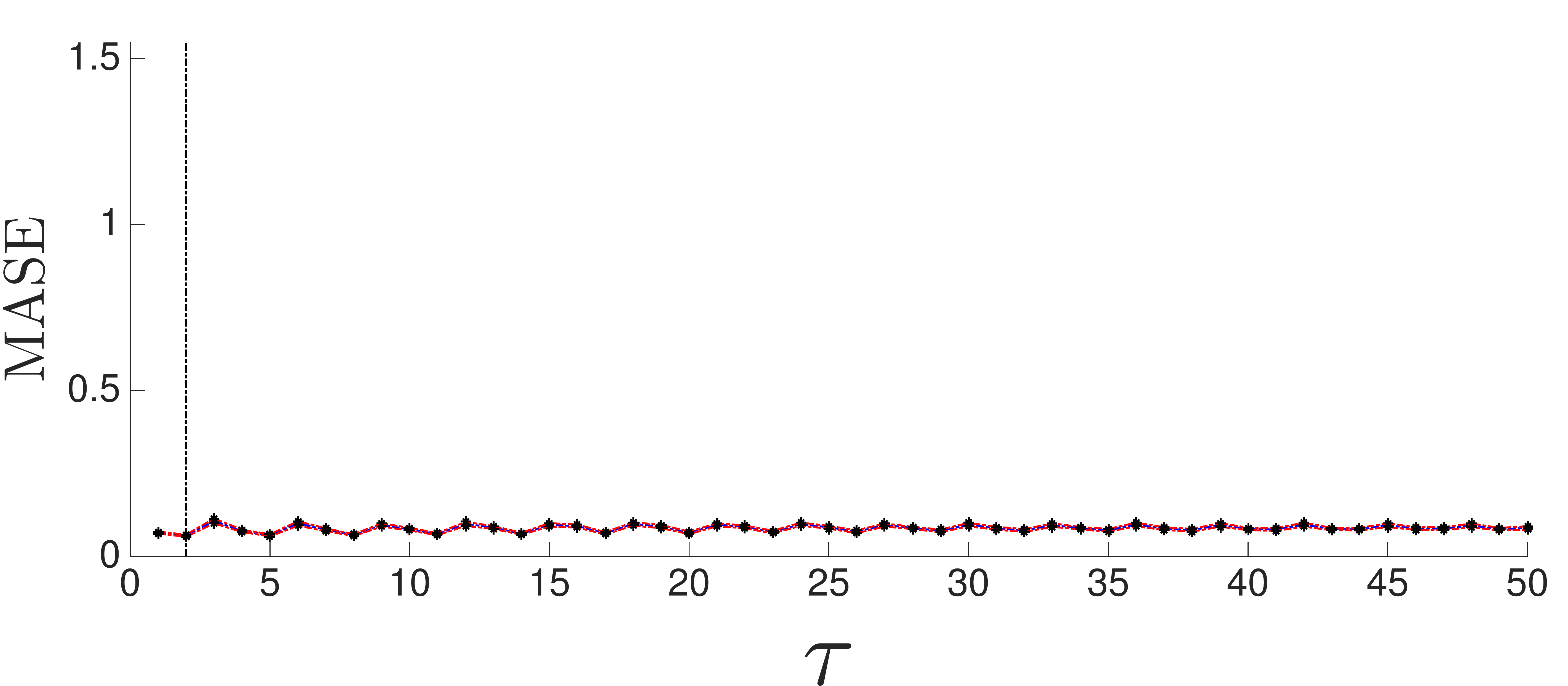}
                \caption{\roLMA on \col.}
        \end{subfigure}%
        \begin{subfigure}[b]{0.5\textwidth}
                \includegraphics[width=\textwidth]{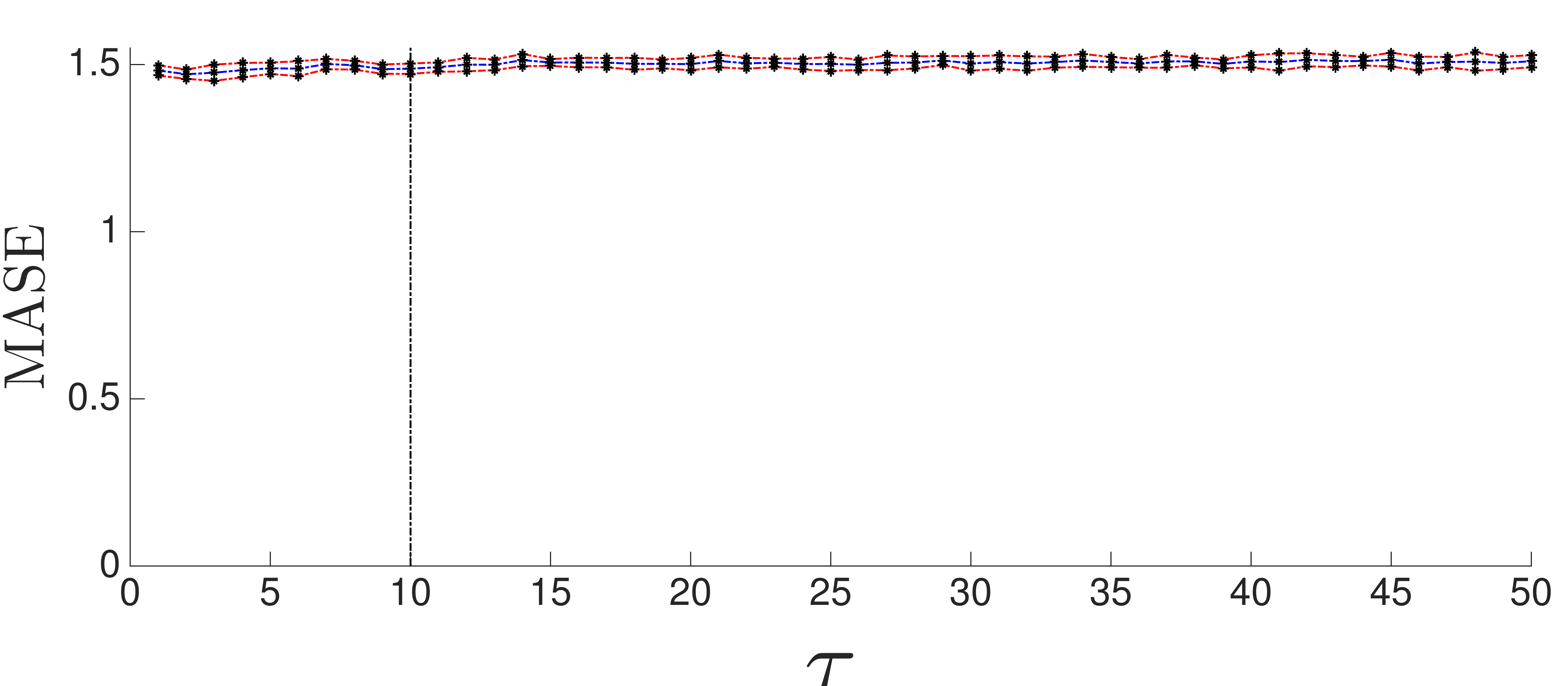}
                \caption{\roLMA on \gcc.}
        \end{subfigure}%
\caption{The effect of $\tau$ on \roLMA forecast accuracy.  The blue
dashed curves are the average 1-MASE of the \roLMA forecasts; the red
dotted lines show $\pm$ the standard deviation.  The black vertical
dashed lines mark the $\tau$ that is the first minimum of the mutual
information curve for each time
series.}
\label{fig:adapMASEvsTAU}
\end{figure}

Figure~\ref{fig:adapMASEvsTAU} explores this effect in more detail.
Across all $\tau$ values, the 1-MASE of \col was generally lower than
the other three experiments---again, simply because this time series
has more predictive structure.  The $K=22$ curve is generally lower
than the $K=47$ one for the same reason, as discussed at the end of
the previous section.  For both Lorenz-96 traces, prediction accuracy
increases monotonically with $\tau$.  It is known that increasing
$\tau$ can be beneficial for longer prediction horizons~\cite{kantz97}.
The situation in Figure~\ref{fig:adapMASEvsTAU} involves short prediction horizons, so it makes sense
that my observations are consistent with the contrapositive of that
result.

For the experimental traces, the relationship between $\tau$ and 1-MASE
score is less simple.  There is only a slight upward overall trend
(not visible at the scale of the Figure) and the curves are
nonmonotonic.  This latter effect is likely due to periodicities in
the dynamics, which are very strong in the \col signal ({\it viz.,} a
dominant unstable period-three orbit in the dynamics, which traces out
the top, bottom, and middle bands in Figure~\ref{fig:col12DEmbedding}).
Periodicities can cause obvious problems for delay
reconstructions---and forecast methods that employ them---if the delay
is a harmonic or subharmonic of their frequencies, simply because the
coordinates of the delay vector are not independent samples of the
dynamics.  It is for this reason that Takens mentions this condition
in his original paper.  Here, the effect of this is an oscillation in
the forecast accuracy vs. $\tau$ curve: low when it is an integer
multiple of the period of the dominant unstable periodic orbit in
the \col dynamics, for instance, then increasing with $\tau$ as more
independence is introduced into the coordinates, then falling again as
$\tau$ reaches the next integer multiple of the period, and so on.

This naturally leads to the issue of choosing a good value for the
delay parameter.  Recall that all of the experiments reported so far used a $\tau$ value chosen at the first
minimum of the mutual information curve for the corresponding time
series.  These values are indicated by the black vertical dashed lines
in Figure~\ref{fig:adapMASEvsTAU}.  This estimation strategy was
simply a starting point, chosen here because it is arguably the most
common heuristic used in the nonlinear time-series analysis community.
As is clear from Figure~\ref{fig:adapMASEvsTAU}, though, it
is {\it not} the best way to choose $\tau$ for reduced-order forecast
strategies.  Only in the case of \col is the $\tau$ value suggested by
the mutual-information calculation optimal for \roLMA---that is, does
it fall at the lowest point on the 1-MASE vs. $\tau$ curve.

\begin{figure}[bt!]
        \centering
        \vspace{0.5cm}
        \begin{subfigure}[b]{0.47\textwidth}
\includegraphics[width=\textwidth]{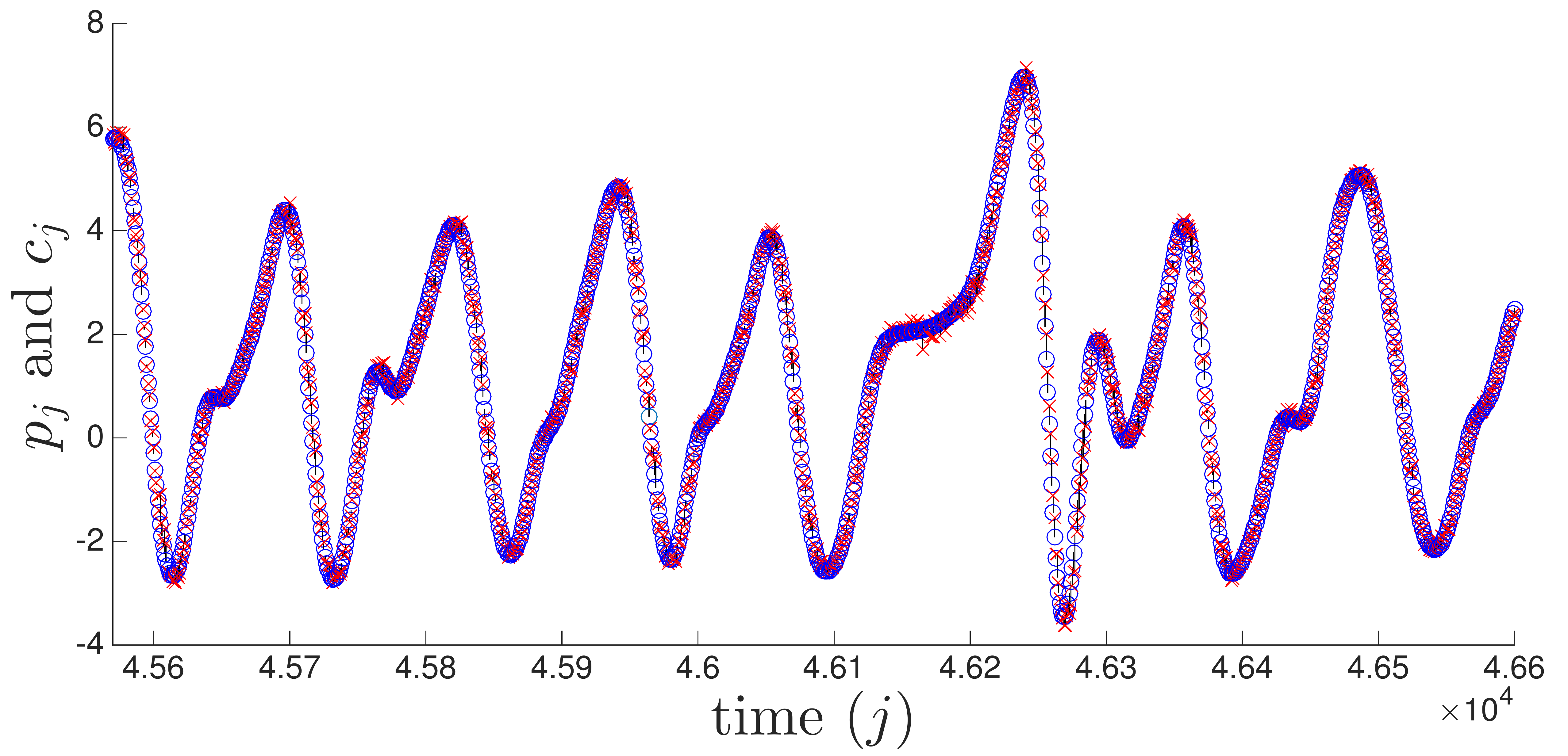}
                  \caption{ \roLMA forecast ($MASE=0.985$, default $\tau$) }
        \end{subfigure}
        \quad
        \begin{subfigure}[b]{0.47\textwidth}
  \includegraphics[width=\textwidth]{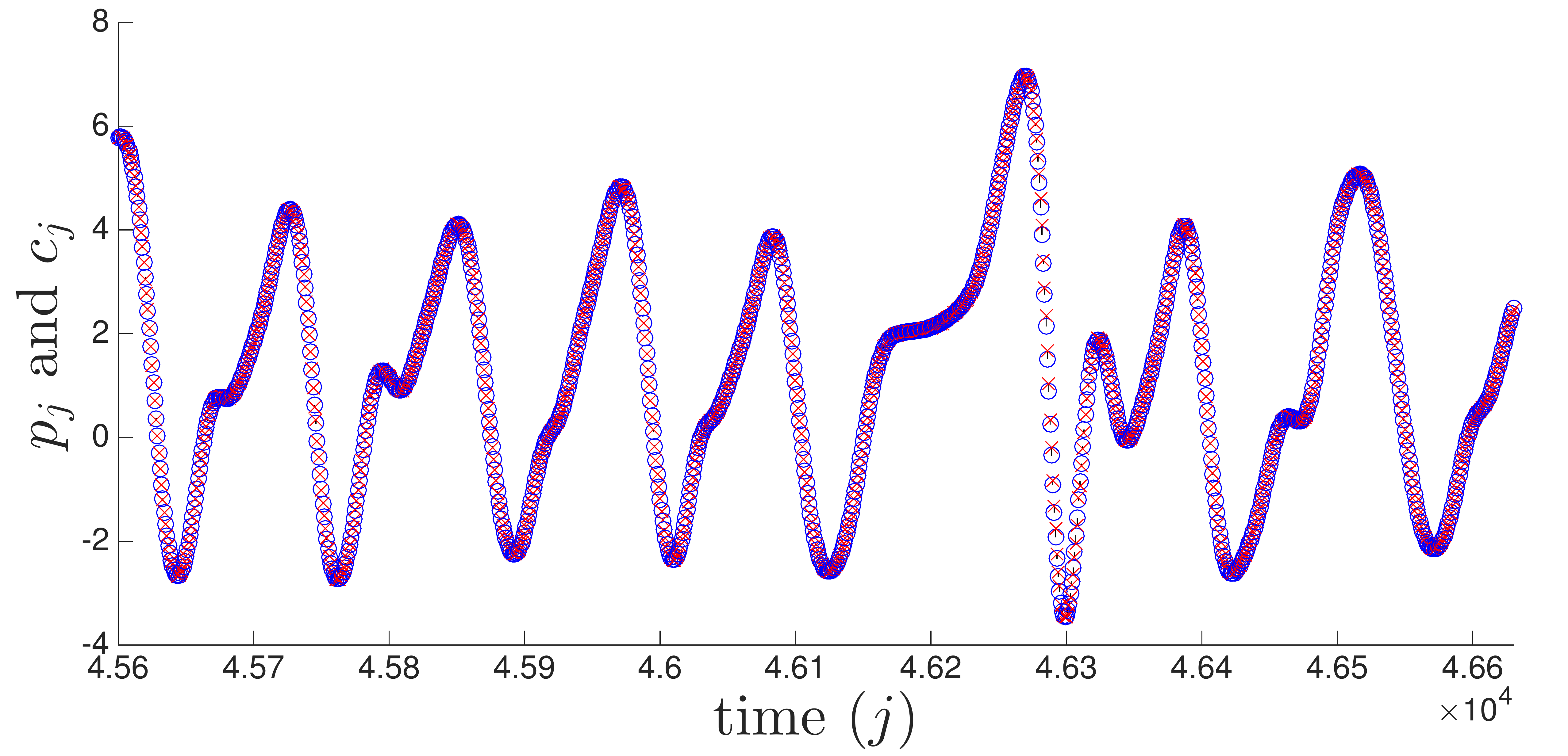}
                  \caption{ \roLMA forecast ($MASE=0.115$, best-case $\tau$) }
        \end{subfigure}
\caption{ Time-domain plots of \roLMA forecasts of  a $K=47$ Lorenz-96
trace with default and best-case $\tau$ values: (a) the first minimum of the time-delayed average mutual information and (b) the minimum of Figure~\ref{fig:adapMASEvsTAU}(b).}
\label{fig:K47BestCasepredictions}
\end{figure}

This suggests that one can often improve the performance of
\roLMA  simply by choosing a different $\tau$---i.e., by adjusting the
one free parameter of that reduced-order forecast method.  In all
cases (aside from \col, where the default $\tau$ was the optimal
value), adjusting $\tau$ allow \roLMA to outperform \fnnLMA. 
The improvement can be quite striking: for visual comparison,
Figure~\ref{fig:K47BestCasepredictions} shows
\roLMA forecasts of a representative 
$K=47$ Lorenz-96 trace using default and best-case values of $\tau$.
However, that comparison is not really fair.  Recall that the
embedding that is used by \fnnLMA, as defined so far, fixes $\tau$ at
the first minimum of the average mutual information for the
corresponding trace.  It may well be the case that that $\tau$ value
is suboptimal for {\it that} method as well---as it was for \roLMA.
To test this, I perform an additional set of experiments to find
the optimal $\tau$ for \fnnLMA.  Table~\ref{tab:best-case-error} shows
the numerical values of the 1-MASE scores, for forecasts made with
default and best-case $\tau$ values, for both methods and all traces.
In the two simulated examples, best-case \roLMA ~ significantly
outperforms best-case \fnnLMA; in the two experimental examples, 
best-case \fnnLMA is better, but not by a huge margin.  That is, even
if one optimizes $\tau$ individually for these two methods, \roLMA
keeps up with, and sometimes outperforms, \fnnLMA. Again, this supports the main point of this thesis: forecast methods
based on incomplete reconstructions of time-series data can be very
effective---and much less work than those that require a full
embedding.

\begin{table*}
\caption{The effects of the $\tau$ parameter.  The ``default'' value 
is fixed, for both \roLMA and \fnnLMA, at the first minimum of the
average mutual information for that trace; the ``best case'' value is
chosen individually, for each method and each trace, from plots like
the ones in Figure~\ref{fig:adapMASEvsTAU}.}

  \begin{center}
  \begin{tabular}{lcccc}
 \hline Signal & \roLMA              & \roLMA               & \fnnLMA             & \fnnLMA\\
  & (default $\tau$)  & (best-case $\tau$) & (default $\tau$)  & (best-case $\tau$) \\ 
  \hline
    Lorenz-96 $K=22$     & $0.391\pm0.016$   & $0.073 \pm 0.002$  & $0.441 \pm 0.033$ & $0.137\pm0.006$ \\
    Lorenz-96 $K=47$     & $0.985 \pm 0.047$ & $0.115 \pm 0.006$  & $1.007 \pm 0.043$ & $0.325\pm0.020$ \\
    \col          & $0.063\pm0.003$   & $0.063\pm 0.003$& $ 0.050 \pm0.002$ & $0.049\pm0.002$ \\
    \gcc          & $1.488\pm0.016$   & $1.471\pm0.014$   & $ 1.530\pm 0.021$ &$1.239\pm0.020$\\ \hline
  \end{tabular}
  \end{center}
 \label{tab:best-case-error}
  \end{table*}

In view of my claim that part of the advantage of \roLMA stems from the
natural noise mitigation effects of a low-dimensional {reconstruction}, it
may appear somewhat odd that \fnnLMA works better on the experimental
time-series data, which certainly contain noise.  Comparisons of large
1-MASE scores are somewhat problematic, however.  
Recall that $1$-MASE $>1$ means that the forecast is worse than an
in-sample random-walk forecast of the same trace.  The bottom row of
numbers in Table~\ref{tab:best-case-error}, then, indicate that both
LMA-based methods---no matter the $\tau$ values---generate poor
predictions for \gcc: 24--53\% worse, on the average, than simply
predicting that the next value will be equal to the previous value.
There could be a number of reasons for this poor performance.  This
signal has almost no predictive structure \cite{josh-pre}
and \fnnLMA's extra axes may add to its ability to capture that
structure---in a manner that outweighs the potential noise effects of
those extra axes.  The dynamics of \col, on the other hand, are fairly
low dimensional and dominated by a single unstable periodic orbit; it
could be that the embedding of these dynamics used in \fnnLMA captures
its structure so well that \fnnLMA is basically perfect and \roLMA
cannot do any better.


While the plots and 1-MASE scores in this section
suggest that \roLMA forecasts are quite good---better than traditional linear methods, and as good or better than LMA upon true embeddings---it is important to note
that both ``default'' and ``best-case'' $\tau$ values were chosen
after the fact in all of those experiments.  This is not useful in
practice.  A significant advantage of a reduced-order forecast
strategy like \roLMA is its ability to work `on the fly' in situations
where one may not have the leisure to run an average mutual
information calculation on a long segment of the trace and find a
clear minimum---let alone construct a plot like
Figure~\ref{fig:adapMASEvsTAU} and choose an optimal $\tau$ from it.
(Producing that plot required 3,000 runs involving a total of
22,010,700 forecasted points, which took approximately 44.5 hours on
an Intel Core i7.)

The results that I report in Section~\ref{sec:tdAIS} and in~\cite{josh-tdAIS}, however, suggest that 
it is possible to estimate optimal $\tau$ values for delay
reconstruction-based forecasting by calculating the value that
maximizes the information shared between each delay vector and the
future state of the system. 
%
%
For all of the examples in this thesis, that strategy produces the same
$\tau$ value as found with the exhaustive search mentioned above.
This is a fairly efficient calculation: $\mathcal{O}(n\log n)$ time
where $n$ is the length of the time series.  Even so, it can be
onerous if $n$ is very large.  However, this measure can be calculated
on very small subsets of the time series and still produce accurate
results, which could allow $\tau$ to be selected {\it adaptively} for
the purposes of forecasting nonstationary systems with \roLMA. 

\subsection{Prediction Horizon}\label{sec:hforecasts}

There are fundamental limits on the prediction of chaotic systems.
Positive Lyapunov exponents make long-term forecasts a difficult
prospect beyond a certain point for even the most sophisticated
methods~\cite{weigend-book, kantz97, josh-pre}.  Note that the
coordinates of points in higher-dimensional delay-reconstruction
spaces sample wider temporal spans of the time series.  In
theory, this means that one should be able to forecast further into
the future with a higher-dimensional reconstruction without losing
memory of the initial condition.  This raises an important concern
about \roLMA: whether its accuracy will degrade with increasing
prediction horizon more rapidly than that of \fnnLMA.

Recall that the formulations of both methods, as described and
deployed in the previous sections of this chapter, assume that
measurements of the target system are available in real time: they
``rebuild'' the LMA models after each step, adding new time-series
points to the embeddings or reconstructions as they arrive.
Both \roLMA and \fnnLMA can easily be modified to produce longer
forecasts, however---say, $h$ steps at a time, only updating the model
with new observations at $h$-step intervals.  Naturally, one would
expect forecast accuracy to suffer as $h$ increased for any
non-constant signal.  The question at issue in this section is whether
the greater temporal span of the data points used by \fnnLMA mitigates
that degradation, and to what extent.

\begin{table}[tb!]
\caption{The $h$-step mean
absolute scaled error ($h$-MASE) scores for different forecast
horizons ($h$).  As explained in Section~\ref{sec:accuracy}, $h$-MASE scores should
not be compared for different $h$ ({\it i.e.,} down the columns of this
table).}
  \begin{center}
  \begin{tabular}{ccccccc}
  \hline\hline Signal  & $h$ & \roLMA & \roLMA & \fnnLMA & \fnnLMA \\ 
& & (default $\tau$) & (best-case $\tau$) & (default $\tau$) &  (best-case $\tau$) \\
  \hline
Lorenz-96 $K=22$  &1& $0.391\pm0.016$ & $0.073\pm0.002$  & $0.441\pm0.003$ & $ 0.137\pm0.006$ \\
Lorenz-96 $K=22$ &10& $0.101\pm0.008$ & $0.066\pm0.003$   & $0.062\pm0.011$ & $0.033\pm0.002$ \\
Lorenz-96 $K=22$ &50& $0.084\pm0.007$ & $0.074\pm0.008$  &  $0.005\pm0.002$ & $0.004\pm0.001$ \\
Lorenz-96 $K=22$ &100& $0.057\pm0.005$& $0.050\pm0.004$ & $0.003\pm0.001$ & $0.003\pm0.001$ \\
\hline
Lorenz-96 $K=47$  &1& $0.985\pm0.047$ & $0.115\pm0.006$ & $0.995\pm0.053$ & $0.325\pm 0.020 $ \\
Lorenz-96 $K=47$  &10 &$0.223\pm0.011$  &  $0.116\pm 0.005$  & $ 0.488\pm0.042$& $0.218\pm0.012$ \\
Lorenz-96 $K=47$  &50 &$0.117\pm0.011$&$0.112\pm 0.010$& $0.127\pm0.011$& $0.119\pm0.010$\\

Lorenz-96 $K=47$  &100 &$0.075\pm0.006$ &$0.068\pm0.005$ & $0.079\pm0.005$& $0.075\pm0.004$ \\

\hline
\col  &1 &$0.063\pm 0.003$&$0.063\pm 0.003$ &$0.050\pm0.002$ & $0.049\pm0.002$  \\
\col  &10 & $0.054\pm0.006$& $0.046\pm0.003$& $0.021\pm0.001$& $0.018\pm0.001$ \\
\col  &50 & $0.059\pm0.009$& $0.037\pm0.003$& $0.012\pm0.003$& $0.009\pm0.001$ \\
\col  &100 & $ 0.044\pm0.004$& $0.028\pm0.006$& $0.010\pm0.003$& $0.007\pm0.001$ \\
\hline
\gcc  &1 &$1.488\pm0.016$&$1.471\pm0.014$&$1.530\pm0.021$&$1.239\pm0.020$\\
\gcc  &10 & $0.403\pm0.009$& $0.396\pm0.009$ & $0.384\pm0.007$ & $0.369\pm0.010$   \\
\gcc  &50 &$0.154\pm0.003$& $0.151\pm0.005$& $0.143\pm0.003$ & $0.141\pm0.003$ \\
\gcc  &100 & $0.101\pm0.002$& $0.101\pm0.003$  &$0.095\pm0.002$& $0.093\pm0.002$  \\
  \hline\hline
  \end{tabular}
  \end{center}
 \label{tab:hstep}
  \end{table}

In Table~\ref{tab:hstep}, I provide $h$-MASE scores---with $h\ge1$ to reflect the increased prediction horizons I am considering---for $h$-step
versions of the different forecast experiments\footnote{For an explanation of $h$-MASE see Section~\ref{sec:accuracy}.}from
Sections~\ref{sec:roLMALorenz96} and \ref{sec:compPerfProj}.   The important comparisons here are, as
mentioned above, across the rows of the table.  The different methods
``reach'' different distances back into the time series to build the
models that produce those forecasts, of course, depending on their
delay and dimension.  At first glance, this might appear to make it
hard to sensibly compare, say, default-$\tau$ \roLMA and
best-case-$\tau$ \fnnLMA, since they use different $\tau$s and
different values of the reconstruction dimension and thus are spanning
different ranges of the time series.  Because $h$ is measured in units
of the sample interval of the time series, however, comparing one
$h$-step forecast to another (for the same $h$) does make sense.

There are a number of interesting questions to ask about the patterns
in this table, beginning with the one that set off these experiments:
how do \fnnLMA and \roLMA compare if one individually optimizes $\tau$
for each method?  The numbers indicate that \roLMA beats \fnnLMA for
$h=1$ on the $K=22$ traces, but then loses progressively badly ({\it i.e.,} 
by more $\sigma$s) as $h$ grows.  \col follows the same pattern except
that \roLMA is worse even at $h=1$.  For \gcc, \fnnLMA performs better
at both $\tau$s and all values of $h$, but the disparity between the
accuracy of the two methods does not systematically worsen with
increasing $h$.  For $K=47$, \roLMA consistently beats \fnnLMA for
both $\tau$s for $h\leq 10$ but the accuracy of the two methods is
comparable for longer prediction horizons. These results suggest that 
optimizing $\tau$ can improve both \fnnLMA and \roLMA and that, depending 
on the signal, this optimization can change the relative accuracy of the two methods.  
This finding catalyzed the development of the forecast-specific 
parameter selection framework that is outlined in Section~\ref{sec:tdAIS} and \cite{josh-tdAIS}.

Another interesting question is whether the assertions in the previous
section stand up to increasing prediction horizon.  Those assertions
are based on the results that appear in the $h=1$ rows of
Table~\ref{tab:hstep}: \roLMA was better than \fnnLMA on the $K=22$
Lorenz-96 experiments, for instance, for both $\tau$ values.  This
pattern does not persist for longer prediction horizons:
rather, \fnnLMA generally outperforms \roLMA on the $K=22$ traces for
$h=10, 50,$ and 100.  The $h=1$ comparisons for $K=47$
and \col ~{\it do} generally persist for higher $h$, however.  As
mentioned before, \gcc is problematic because its 1-MASE scores are so
high, but the accuracies of the two methods are similar for all $h>1$.

\begin{figure*}[tb!]
        \centering
        \vspace{0.5cm}
\includegraphics[width=\textwidth]{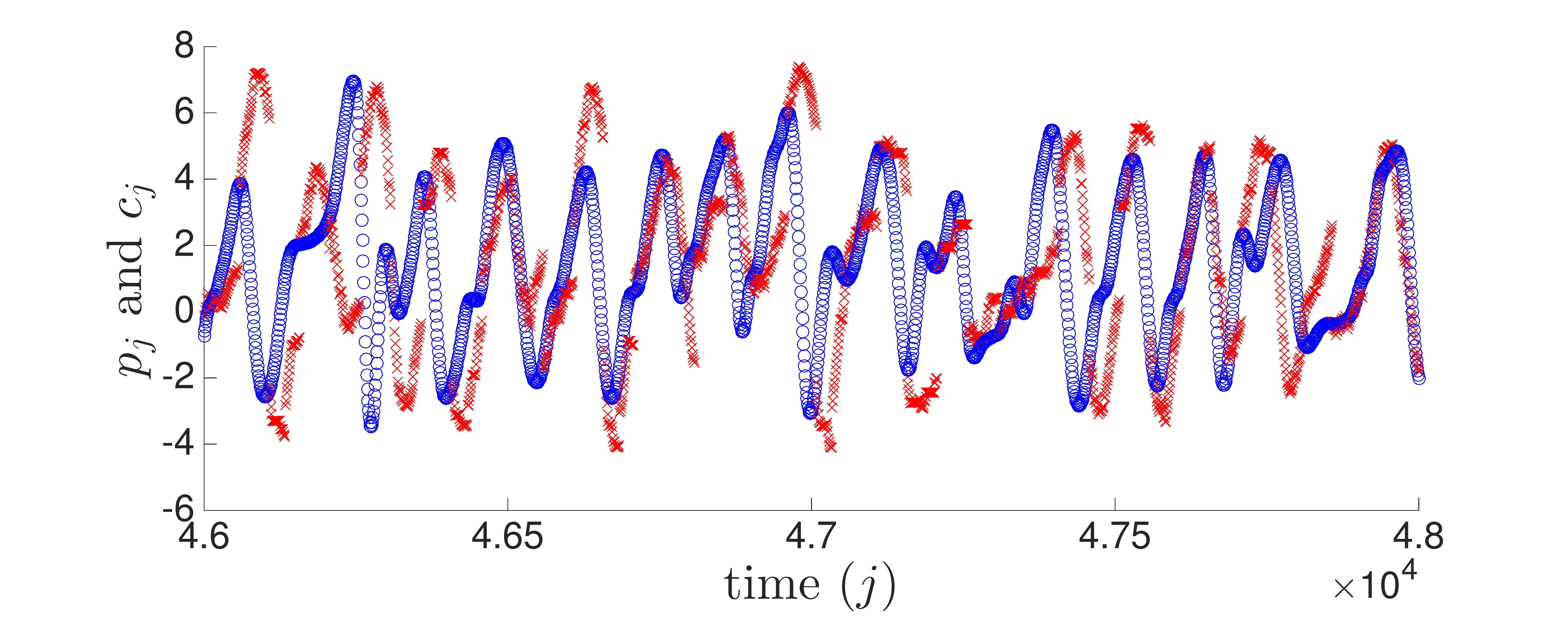}

\caption{A  best-case-$\tau$ \roLMA forecast of a $K=47$ Lorenz-96 
trace for $h=50$.  The forecast (red) follows the true trajectory
(blue) for a while, falls off onto a shadow trajectory, then gets
recorrected when a new set of observations are incorporated into the
model after $h$ time steps.}
\label{fig:shadowtrajectories}
\end{figure*}

The fact that \fnnLMA generally outperforms \roLMA for longer
prediction horizons makes sense simply because \roLMA samples less of
the time series and therefore has less `memory' about the dynamics.
This is a well-known effect~\cite{weigend-book}.  
In view of the fundamental limits on prediction of
chaotic dynamics, however, it is worth considering
whether {\it either} method is really making correct long-term
forecasts.  Indeed, time-domain plots of long-term forecasts (e.g.,
Figure~\ref{fig:shadowtrajectories}) reveal that both \fnnLMA
and \roLMA forecasts have fallen off the true trajectory and onto
shadow trajectories---another well-known phenomenon when forecasting
chaotic dynamics~\cite{sauer-delay}.

In other words, it appears that even a 50-step forecast of these
chaotic trajectories is a tall order: {\it i.e.,} that I am running up
against the fundamental bounds imposed by the Lyapunov exponents.  In
view of this, it is promising that \roLMA generally keeps up
with \fnnLMA in many cases---even when both methods are struggling
with the prediction horizon, and even though the model that \roLMA
uses has much less memory about the past history of the trajectory.
An important aspect of my future research on this topic will be
developing efficient methods for deriving bounds on reasonable
prediction horizons purely from the time series, {\it i.e.,} without using
traditional methods such as Lyapunov exponents, which are difficult to
estimate from experimental data.  See \cite{josh-tdAIS} for some of my
preliminary results on this line of research.

\subsection{Data Length}\label{sec:datalength}

Most real-world data sets are fixed in length and some are quite
short.  Moreover, many of the dynamical systems that one might like to
predict are nonstationary.  For these reasons, it is important to
understand the effects of data length upon forecast methods that employ delay reconstructions.  For the reconstruction to be an
actual embedding that supports accurate calculations of dynamical
invariants, the data requirements are fairly dire.  Traditional
estimates ({\it e.g.,} by Smith~\cite{smithdatabound} and by
Tsonis {\it et al.}~\cite{tsonisdatabound}) suggest that $\approx
10^{17}$ data points would be required to embed the Lorenz-96 $K=47$
data in Section~\ref{sec:roLMALorenz96}, where the known $d_{KY}$
values~\cite{KarimiL96} indicate that one might need at least $m=38$
dimensions to properly unfold the dynamics.  As described in Section~\ref{sec:dce}, however, that is only
truly necessary if one is interested in preserving the diffeomorphism
between true and reconstructed dynamics, down to the last detail.  For
the purposes of prediction, thankfully, one can make progress with far
less data.  For example, Sauer~\cite{sauer-delay} successfully
forecasted the continuation of a 16,000-point time series using a 16-dimensional embedding; Sugihara \& May~\cite{sugihara90} used
delay-coordinate embedding with $m$ as large as seven to successfully
forecast biological and epidemiological time-series data as short as
266 points.
%
%

While the results in the previous sections are based on far longer
traces than the examples mentioned at the end of the previous
paragraph, it is still worth exploring whether data-length issues are
affecting those results---and evaluating whether those effects
differentially impact \roLMA because of its lower-dimensional model.
This kind of test can be problematic in practice, of course, since it
requires varying the length of the data set.  In a synthetic example
like Lorenz-96, that is not a problem, since one can just run the ODE
solver for more steps.  

\begin{figure}[bt!]
\vspace{0.5cm}
        \centering
        \begin{subfigure}[b]{0.49\textwidth}
\includegraphics[width=\textwidth]{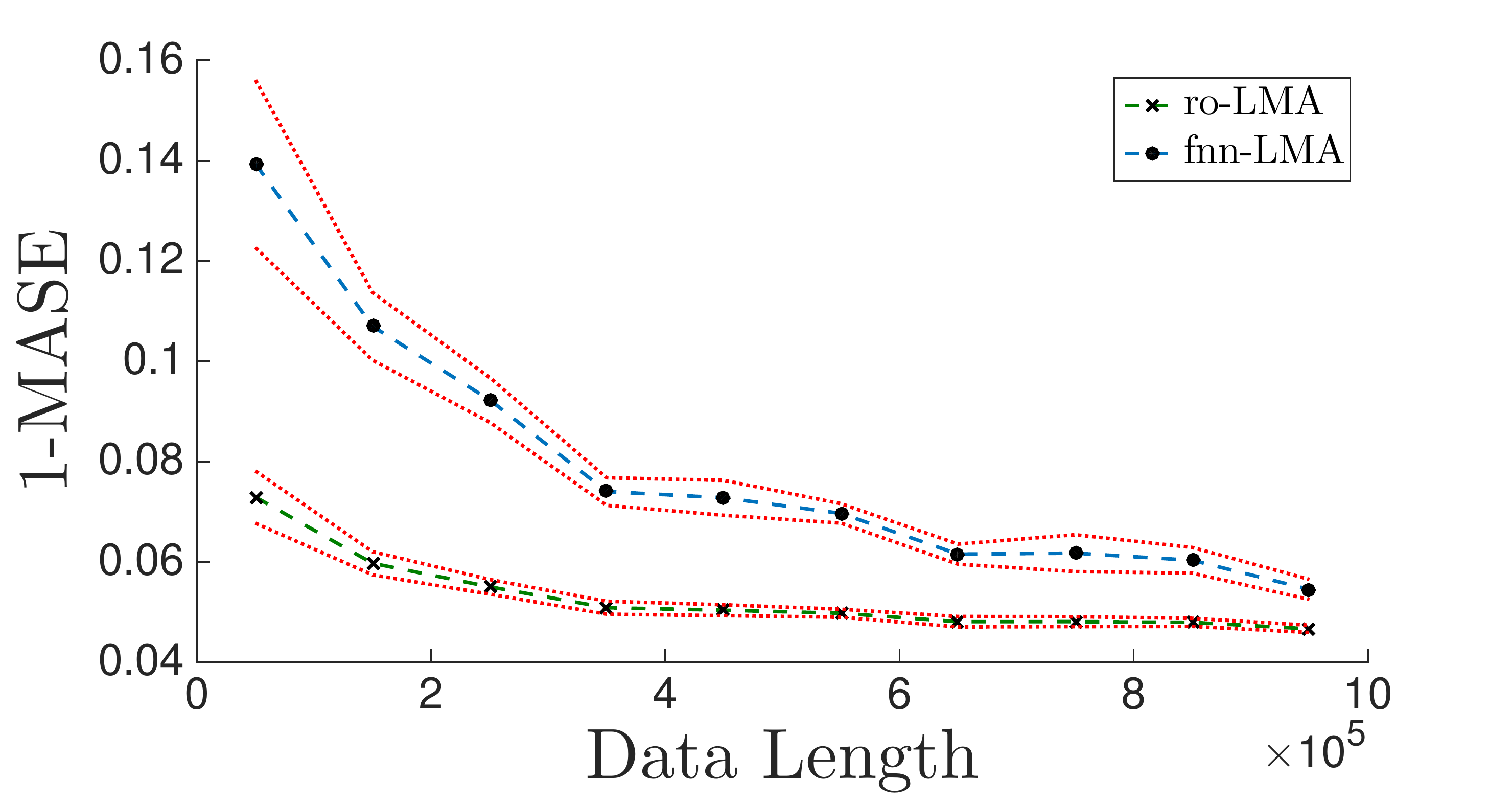}
                  \caption{1-MASE as a function of data length for predictions of
the Lorenz-96 $K=22,F=5$ traces.}
        \end{subfigure}
        \begin{subfigure}[b]{0.49\textwidth}
  \includegraphics[width=\textwidth]{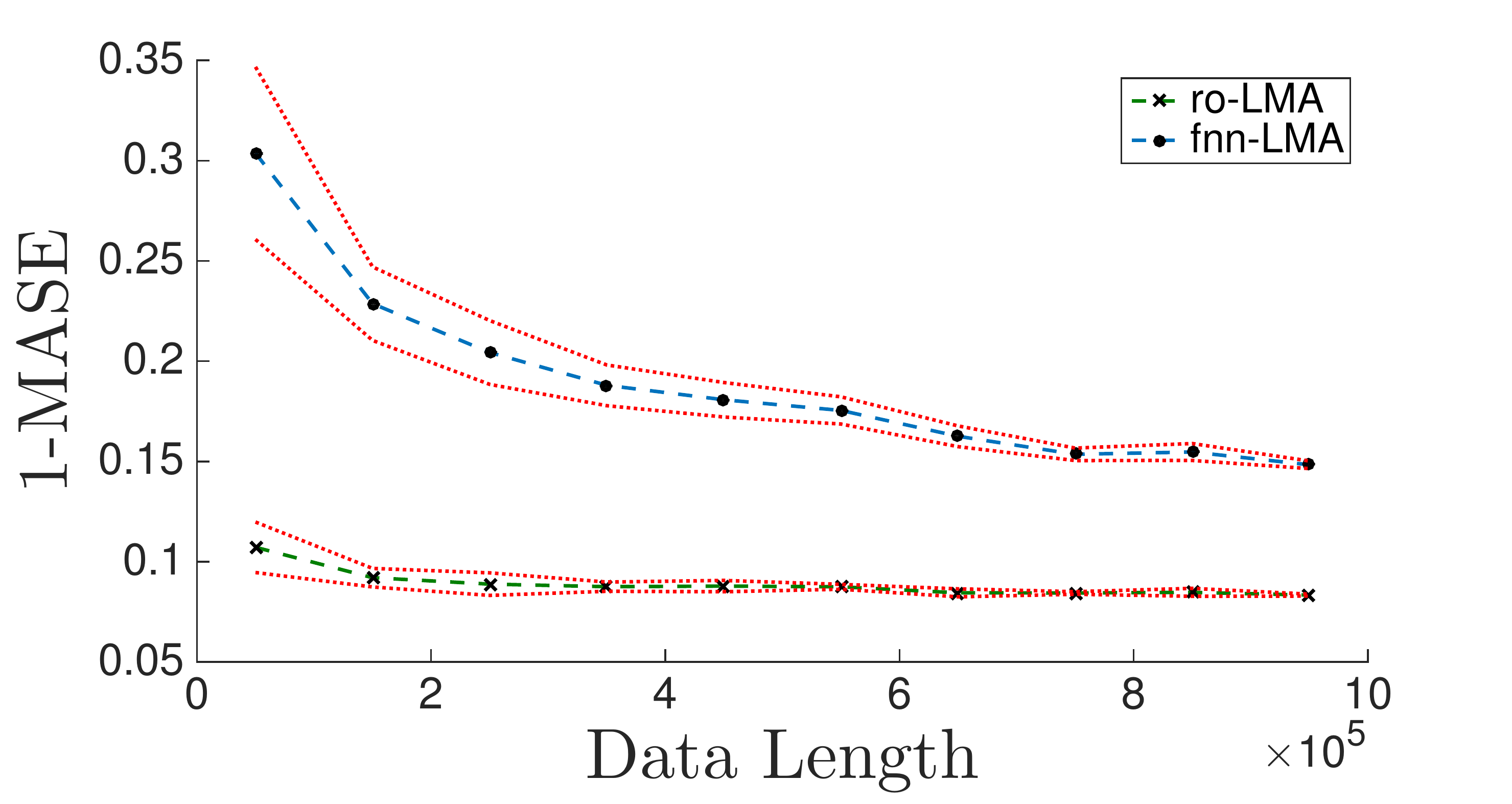}
                  \caption{1-MASE as a function of data length for predictions of
the Lorenz-96 $K=47,F=5$ traces.}
        \end{subfigure}
\caption{The effects of data length on forecast error for fixed prediction
horizon $h=1$. The dashed lines show the mean forecast error for each method; the
dotted lines indicate the range of the standard deviation.}
\label{fig:dataexperiments}
\end{figure}

Figure~\ref{fig:dataexperiments} shows 1-MASE scores for \roLMA
and \fnnLMA forecasts of the Lorenz-96 system as a function of data
length.
  For both
$K=22$ and $K=47$, the \fnnLMA error is higher than the \roLMA error,
corroborating the results in the previous sections.  Both generally
improve with data length, as one would expect---but only up to a
point.  The 1-MASE scores of the \roLMA forecasts, in particular,
reach a plateau at about 350,000 points in the $K=22$ case and 150,000
in the $K=47$ case.  \fnnLMA, on the other hand, keeps improving out
to the end of the Figure.  Eventually, there is a crossover for
$K=22$, but not until 1.6 million points.  For $K=47$, the
curves are still converging out to 4 million points. The
difference in crossover points is not surprising, given that the
dimension of the $K=47$ dynamics is so much higher: $d_{KY}<\approx3$,
versus $\approx 19$.  What {\it is} surprising, and useful, is that
for traces with 1 million points---far more data than a practitioner
can generally hope for---\roLMA is still outperforming \fnnLMA.
Moreover, it is doing so using fewer dimensions, which makes it
computationally efficient.

Plateaus in curves like the ones in Figure~\ref{fig:dataexperiments}
suggest that the corresponding forecast method has captured all of the
information that it can use, so that adding data does not improve the
forecast.  This effect, which is described at more length in Section~\ref{sec:tdAIS}, depends on dimension for the obvious reason that
filling out a higher-dimensional object requires more data.  This
suggests another piece of forecasting strategy: when one is data-rich,
it may be wise to choose \fnnLMA over \roLMA.  In making this choice,
though, one should also consider the added computational complexity,
which will be magnified by the longer data length.  When data are not
plentiful, though, or the system is nonstationary, my results suggest
that it is advantageous to ignore the theoretical bounds
of~\cite{sauer91} and use low-dimensional reconstructions in forecast
models.

\section{Summary}\label{sec:concl}

In summary,  it appears that incomplete embeddings---time-delay 
reconstructions that do not satisfy the formal conditions of the
embedding theorems---are indeed adequate for the purposes of
forecasting dynamical systems.  Indeed, they appear to offer simple
state-space based methods {\it even more} traction on this prediction
task than full embeddings, with greatly reduced computational
effort.

The study in this thesis specifically focuses on the $2D$
instantiation of the claim above, for two reasons.  First, that is the
extreme that, in a sense, most seriously violates the basic tenets of
the delay-coordinate embedding machinery; second, working in 2D enables the
largest reduction in the cost of the near-neighbor searches that are
the bulk of the computational effort in most state space-based
forecast methods.
%
%
A number of other issues arise when one considers increasing the
reconstruction dimension beyond two, besides raw computational
complexity.  Among other things, that would introduce another free parameter into the
method, thereby requiring some sort of principled strategy for
choosing its value (a choice that \roLMA completely avoids by fixing
$m=2$). In general, one would expect forecast accuracy to improve with the
number of dimensions in the reconstruction, but not without limit.  Among other things,
noise effects grow with that dimension (simply because every
noisy data point affects $m$ points in the reconstructed trajectory).
And from an information-theoretic standpoint, one would expect
diminishing returns when the span of the delay vector
($m \times \tau$) exceeds the ``memory'' of the system.  For all of
those reasons, it would seem that there should be a plateau beyond
which increasing the reconstruction dimension does not improve the
accuracy of a forecast methods that use the resulting models.  I
explore that issue further in \cite{josh-tdAIS} and synopsize the relevant material in Section~\ref{sec:tdAIS}. In that discussion, I use the measure
mentioned in the last paragraph of Section~\ref{sec:varyingproj} to
derive optimal reconstruction dimensions for near-neighbor forecasting
for a broad range of systems, noise levels, and forecast
horizons---all of which turn out to be $m=2$.  In the bulk of the
dynamical systems literature on forecasting, however, the optimal
reconstruction dimension for the purposes of forecasting was thought
to be near the value that provides a true embedding of the data.  The
results in this thesis suggest, again, that this is not the case.

I chose the classic Lorenz method of analogues as a good exemplar of
the class of state space-based forecast methods, but I believe that
my results will hold for other members of that class
({\it
e.g.},~\cite{casdagli-eubank92,weigend-book,Smith199250,sauer-delay,sugihara90,pikovsky86-sov}).
Working with a low-dimensional reconstruction could potentially reduce the
computational search and storage costs of {\it any} such method,
while also avoiding the so-called ``curse of dimensionality'' and
mitigating noise multipliers caused by extra embedding
dimensions~\cite{Casdagli:1991a}.  Reduced-order reconstructions also reduces data requirements,
since fewer points are required to fill out a lower-dimensional
object.  And when one fixes $m=2$, there is only a single free
parameter $\tau$ in the method---one that can be estimated effectively from a
short sample of the data set, allowing the reduced-order method to adapt to nonstationary
dynamics.  There may be some limitations on the class of methods for
which these claims hold, of course; matters may get more complicated,
and the results less clear, for forecast methods that perform other
kinds of projections.  On the flip side, however, my results can be
viewed as explaining why those methods work so well.

Again, no forecast model will be ideal for all noise-free
deterministic signals, let alone all real-world time-series data sets.
However, the proof of concept offered in this section is encouraging:
prediction in projection---a simple yet powerful reduction of a
time-tested method---appears to work remarkably well, even though
the models that it uses are not necessarily topologically faithful to
the true dynamics.  

%% file: explanation.tex

\chapter{Why it Works: A Deeper Understanding of Delay-Coordinate Reconstruction}\label{ch:explain}

The experimental validation of \roLMA provided in Chapter~\ref{ch:pnp} is promising but that analysis gave rise to several unanswered questions\label{page:questions}, {\it e.g.,} 
\begin{itemize}
\item {\it Why} does \roLMA work when it is effectively a heresy?,
\item If $m=2$ works so well, why not $m=3$?,
\item How much data is necessary before $m>2$ is the clear winner over higher-dimensional reconstructions?,
\item If one wants to forecast two or three steps into the future, is $m=2$ still efficient, or should $m$ be increased?,
\item Can $\tau$ be chosen {\it a priori} to optimize the accuracy of \roLMA?, and 
\item Can all of these questions be answered strictly by analyzing the data?
\end{itemize}
 This chapter provides a two-part analysis that answers many of these questions. The first part leverages information theory to select forecast-optimal parameters for delay-coordinate reconstruction. The second borrows methods from computational topology to gain new insight into the delay-coordinate embedding theory and machinery.

These two theoretical frameworks---information theory and computational topology---are mathematically disjoint but complementary in terms of developing a complete theory of reconstruction-based forecasting. In particular, the combination of these two tools allows, for the construction of a new paradigm in delay-coordinate reconstruction. 
Section~\ref{sec:tdAIS} offers a novel method, developed in collaboration with R.~G.~James, for leveraging the information that is stored in delay vectors to perform parameter selection that is tailored to the exact stipulations of the data set at hand ({\it e.g.,} data length, signal-to-noise ratio and desired forecast horizon). The traditional approach to this is based on the assumption that the diffeomorphism instantiated by the delay-coordinate map, which is essential for dynamical invariant calculations, is also optimal for forecasting. As the results in Chapter~\ref{ch:pnp} suggest, however, this may not be the best approach. Section~\ref{sec:tdAIS} further corroborates these findings and suggests a reason why. 
Section~\ref{sec:compTopo} provides a deeper theoretical understanding of delay-coordinate reconstruction through computational topology, offering yet another reason why \roLMA works. This exploration is based on the assumption that, when forecasting, one might only require knowledge of the topology of the invariant set; in collaboration with J.~D.~Meiss, I conjecture that the reconstructed
dynamics might be {\it homeomorphic} to the original dynamics at a
lower dimension than that needed for a diffeomorphically correct
embedding. This suggests why \roLMA gets traction despite its use of an incomplete reconstruction. 

The combination of these powerful mathematical tools---information theory and computational topology---allows me to construct a deeper and more complete story of reduced-order forecasting with delay-coordinate reconstruction.

\input{tdAIS} 
\input{compTopology} 

\section{Summary}

Chapter~\ref{ch:pnp} offered experimental validation of my proposed paradigm shift in the practice of delay-coordinate reconstruction, but left many unanswered questions about the theoretical underpinnings (and implications) of this approach. This chapter answered those questions by drawing on the complementary theoretical frameworks of information theory and computational topology. The novel  computationally-efficient metric (\mytau) explicitly leveraged information stored in a delay vector to select forecast-optimal parameters for delay coordinate reconstruction. This metric allowed for tailoring of reconstruction parameters to
available data length, the signal-to-noise ratio in the time series,  and the desired prediction horizon. In addition, this information-theoretic approach gave me the language and methods to answer difficult questions like the ones posed on page~\pageref{page:questions}. Perhaps most importantly, the results in Section~\ref{sec:tdAIS} validated that the state estimator of \roLMA often has more information about the near  future than a traditional embedding---a fact that is completely counter to all the current theory.

Section~\ref{sec:compTopo} deviated significantly from the tone of the rest of this thesis. Instead of discussing delay-coordinate reconstruction purely from a forecasting perspective, it turned a critical eye toward the theoretical foundation and assumptions of this powerful framework, which is the basis for all of nonlinear time-series analysis. The consistent story throughout the previous chapters is that the theoretical requirements of delay-coordinate reconstruction are not necessary---and can indeed be overkill---when one wishes to use it for the purposes of short-term forecasting. But, {\it why} is this the case? Is it simply because more information is present in lower-dimensional state estimators, as established in Section~\ref{sec:tdAIS}, or is there something deeper underpinning this method from a theoretical perspective? In Section~\ref{sec:compTopo}, I used the canonical Lorenz 63 system to argue that large-scale homology can be attainable at much lower dimensions than the theory might suggest. I believe this in turn suggests that a {\it homeomorphism}, in the form of the delay-coordinate map, can be achieved at lower dimensions than what is needed for a diffeomorphism.  This insight suggested an alternative explanation as to {\it why} \roLMA gets traction before it should: specifically, that a homeomorphic reconstruction of a dynamical system may be sufficient for short-term forecasting of dynamical systems. However, further work will be required to rigorously prove this broader claim.

%% file: tdAIS.tex

\section{Leveraging Information Storage to Select Reconstruction Parameters}\label{sec:tdAIS}

As has been discussed throughout this thesis, the task of choosing good values for the free parameters in delay-coordinate reconstruction has been the
subject of a large and active body of literature over the past few
decades.  The majority of these techniques focus
on the {\it geometry} of the reconstruction,  
which is
appropriate when one is interested in quantities like fractal
dimension and Lyapunov exponents. It is not necessarily the best
approach when one is building a delay reconstruction {\it for the
  purposes of prediction}, however, as I showed in Section~\ref{sec:varyingproj}.  That issue, which is the focus of this
section, has received comparatively little attention in the extensive
literature on delay reconstruction-based
prediction~\citeDCEFORECASTING.  

In this section, I propose a
robust, computationally efficient method that I call {\it time delayed active information storage}, \mytau, which can be
used to select parameter values that maximize the information shared 
between the past and the future---or, equivalently, that maximize the
reduction in uncertainty about the future given the current model of
the past~\cite{josh-tdAIS}.  The implementation details, and a complexity analysis of
the algorithm, are covered in Section~\ref{sec:sharedinfo}.  In
Section~\ref{sec:tdAISresults}, I show that simple prediction methods
working with \mytau-optimal reconstructions---{\it i.e.,} constructions using
parameter values that follow from the \mytau calculations---perform
better, on both real and synthetic examples, than those same forecast
methods working with reconstructions that are built using the
traditional parameter selection heuristics (time-delayed mutual information for $\tau$ and false-near neighbors for $m$). Finally, in
Section~\ref{sec:dataandhorizon} I explore the utility of \mytau in
the face of different data lengths and prediction horizons.

\subsection{Shared Information and Delay Reconstructions}
\label{sec:sharedinfo}

The information shared between the past and the future is known as the
excess entropy~\cite{crutchfield2003}.  I will denote it
here by $E = I[\overleftarrow{X}, \overrightarrow{X}]$, where $I$ is
the mutual information~\cite{yeung2012first} and $\overleftarrow{X}$
and $\overrightarrow{X}$ represent the infinite past and the infinite
future, respectively.  $E$ is often difficult to estimate from data
due to the need to calculate statistics over potentially infinite
random variables~\cite{james2014many}.  While this is possible in
principle, it is too difficult in practice for all but the simplest of
dynamics~\cite{strelioff2014bayesian}.  In any case, the excess entropy
is not exactly what one needs for the purposes of prediction, since it
is not realistic to expect to have the infinite past or to predict
infinitely far into the future.  For my purposes, it is more
productive to consider the information contained in the {\it recent}
past and determine how much that explains about the not-too-distant
future.  To that end, I define the {\it state active information storage}
\begin{equation}
  \smytau\equiv{I}[\mathcal{S}_j,~X_{j+h}]
\end{equation}
where $\mathcal{S}_j$ is an estimate of the state of the system at
time $j$ and $X_{j+h}$ is the state of the system $h$ steps in the
future. In the special case where the state estimate $\mathcal{S}$ takes the form of a standard $m$-dimensional delay vector, I will refer to \smytau as the {\it time delayed active information storage}
\begin{equation}
  \mytau\equiv{I}[[X_j,X_{j-\tau},\dots,X_{j-(m-1)\tau}],~X_{j+h}]
\end{equation}

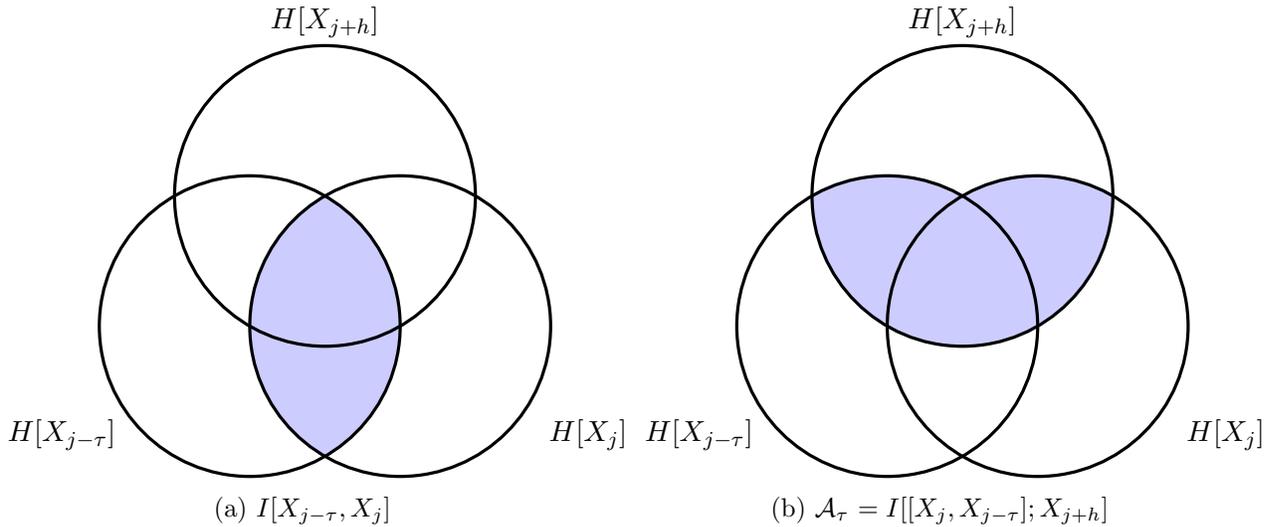
\begin{figure}[ht!]
  \begin{subfigure}[b]{0.49\textwidth}
    \begin{tikzpicture}[baseline=0,label=$m$]
      \setupdiagrams

    \begin{scope}
      \clip \Bcirc;
      \fill[filled] \Ccirc;
    \end{scope}
    \drawdiagram
    \end{tikzpicture}
    \caption{$I[X_{j-\tau}, X_{j}]$}
    \label{fig:MI-I-Diagram}
  \end{subfigure}%
  \quad
  \begin{subfigure}[b]{0.49\textwidth}
   \begin{tikzpicture}[baseline=0]
      \setupdiagrams

      \begin{scope}
        \clip \Bcirc \Ccirc;
        \fill[filled] \Acirc;
      \end{scope}

      \drawdiagram
    \end{tikzpicture}
  \caption{\mytau$ = I[[X_{j},X_{j-\tau}];X_{j+h}]$}
  \label{fig:mytau}      \end{subfigure}
\caption{ (a) An I-diagram of the time-delayed mutual information.  The
    circles represent uncertainties ($H$) in different variables; the
    shaded region represents $I[X_{j};X_{j-\tau}]$, the time-delayed
    mutual information between the current state $X_{j}$ and the state
    $\tau$ time units in the past, $X_{j-\tau}$.  Notice that the shaded
    region is indifferent to $H[X_{j+h}]$, the uncertainty about the future. (b) An I-diagram of \mytau, the quantity proposed in this
    section: $I[[X_{j},X_{j-\tau}];X_{j+h}]$.  This quantity captures
    the shared information between the past, present, and future
    independently, as well as the information that the past and
    present, together, share with the future.
}\label{fig:compare}
\end{figure}

\mytau can be nicely visualized---and compared to traditional methods
like time-delayed mutual information---using the
I-diagrams of
Yeung, introduced in Section~\ref{sec:multiinfoandidiagrams}. Figure~\ref{fig:compare}(a) shows an
$I$-diagram of time-delayed mutual information for a specific $\tau$.
Recall that in a diagram like this, each circle represents the uncertainty in a
particular variable.  The left circle in
Figure~\ref{fig:compare}(a), for instance, represents the average
uncertainty in observing $X_{j-\tau}$ ({\it i.e.,} $H[X_{j-\tau}]$); similarly, the top
circle represents $H[X_{j+h}]$, the uncertainty in the $h^{th}$ future
observation.  Also recall that each of the overlapping regions represents {\it shared}
uncertainty: {\it e.g.,} in Figure~\ref{fig:compare}(a), the shaded region
represents the shared uncertainty between $X_{j}$ and $X_{j-\tau}$---more precisely, the quantity $I[X_{j},~X_{j-\tau}]$.
Notice that minimizing the shaded region in Figure~\ref{fig:compare}(a)---that is, rendering $X_{j}$ and
$X_{j-\tau}$ as independent as possible---maximizes the total
uncertainty that is explained by the combined model
$[X_{j},X_{j-\tau}]$ (the sum of the area of the two circles).  This
is precisely the argument made by Fraser and Swinney in
~\cite{fraser-swinney}; see Section~\ref{sec:numericaltau} for a full explanation. However, it is easy to see from the $I$-diagram
that choosing $\tau$ in this way does not explicitly take into account
explanations of the {\it future}---that is, it does not reduce the
uncertainty about $X_{t+h}$. Moreover, this approach to selecting $\tau$ does not automatically extend
to higher dimensional embeddings, {\it e.g.,} minimizing $I[X_{j},X_{j-\tau}]$, may or may not minimize $I[X_{j},X_{j-\tau}, X_{j-2\tau}]$ and in fact this extension is non-trivial; see Section~\ref{sec:multiinfoandidiagrams} for a full discussion of why this is so.
The obvious next step would be to explicitly include the future in the estimation procedure. As I discussed in Section~\ref{sec:multiinfoandidiagrams}, however, explicitly including the future in the calculation {\it i.e.,} $I[X_{j},X_{j-\tau},X_{j+h}]$, is not straightforward. The rest of this section discusses some of the common interpretations of this quantity and why they are not appropriate for the task at hand.  

The interaction information~\cite{McGill-1954,Bell03theco-information} is one such interpretation of $I[X_{j},X_{j-\tau},X_{j+h}]$ depicted in Figure~\ref{fig:multi-information-panel}(b); this is the intersection of $H[X_{j}]$, $H[X_{j-\tau}]$ and $H[X_{j+h}]$. It describes the reduction in uncertainty that the {\it two} past states, together, provide regarding the future. While this is obviously an improvement over the time-delayed mutual information of Figure~\ref{fig:compare}(a), it does not take into account the information that is shared between $X_j$ and the future but {\it not shared with the past} ({\it i.e.,} $X_{j-\tau}$), and vice versa. The binding information and total correlation, 
depicted in Figures~\ref{fig:multi-information-panel}(c) and (d), address this shortcoming, but both also include information that is shared between the past and the present, but not with the future. This is not terribly useful for the purposes of prediction. Moreover, the total correlation overweights information that is shared between all three circles---past, present, and future---thereby artificially over-valuing information that is shared in all delay coordinates.  In the context of predicting $X_{t+h}$, the provenance of the information is irrelevant and so the total correlation also seems ill-suited to the task at hand.

Note that the total correlation has been used in a similar manner to the time-delayed mutual information method in estimating $\tau$~\cite{fraser-swinney}: {\it e.g.,} minimizing $\mathcal{M}[X_{j}; X_{j-\tau}; X_{j-2\tau}]$ for a three-dimensional embedding. Minimizing the total correlation is equivalent to maximizing the entropy, making the delay vectors maximally informative because dependencies among the dimensions have been minimized. While on the surface this may seem a boon to prediction, consider the issue of predicting the state of the system at time $j+\tau$: if the coordinates of the delay vector are maximally independent, they will also be independent {\it of the value being predicted}.  In light of this, the minimal total correlation approach is not well aligned with the goal of prediction.

Time-delayed active information storage addresses all of the issues raised in the previous
paragraphs. By treating the generic delay vector as a joint variable,
rather than a series of single variables, \mytau captures the shared
information between the past, present, and future independently---the
left and right colored wedges in Figure~\ref{fig:compare}(b)---as well as
the information that the past and present, together, share with the
future (the center wedge).  By choosing delay-reconstruction
parameters that maximize \mytau, then, one can explicitly maximize the
amount of information that each delay vector contains about the
future~\cite{josh-tdAIS}.

That property means that \mytau can be used to select $\tau$ for \roLMA. 
Specifically, to estimate a ``forecast-optimal" $\tau$ value for \roLMA using \mytau, one would simply calculate $\mathcal{A}_\tau =I[[X_{j},X_{j-\tau}],~X_{t+h}]$ for a range of $\tau$, choosing
the first maximum ({\it i.e.,} minimizing the uncertainty about the $h^{th}$
future observation).  
In Section~\ref{sec:tdAISresults}, I explore
that claim using \roLMA and \fnnLMA, but that exploration can be easily extended to any time-delayed state estimator---such as the methods used in~\cite{weigend-book,casdagli-eubank92,Smith199250,sugihara90}---by using the general form of \mytau, {\it viz.,} the state active information storage, \smytau. 
 For example, if the time series is pre-processed ({\it e.g.,} via a Kalman filter~\cite{sorenson}, a low-pass
filter and an inverse Fourier transform~\cite{sauer-delay}, or some
other local-linear
transformation~\cite{weigend-book,casdagli-eubank92,Smith199250,sugihara90,kantz97},) the state estimator simply becomes $\mathcal{S}_j= \hat{\vec{x}}_j$
where $\hat{\vec{x}}_j$ is the processed $m$-dimensional delay vector.

\subsection{Selecting ``Forecast-Optimal" Reconstruction Parameters}
\label{sec:tdAISresults}

This section demonstrates how to use \mytau to choose parameter
values for delay-coordinate reconstructions constructed specifically for the purposes of forecasting, using 
several of the case studies presented in Chapter~\ref{ch:systems}.  For the discussion that follows, the
term ``\mytau-optimal'' is used to refer to the parameter values ($m$ and
$\tau$) that maximize \mytau over a range of $m$ and $\tau$. The general parameter selection framework is presented at first, not assuming the use of either \fnnLMA or \roLMA, and then the \mytau-optimal reconstructions are compared to \fnnLMA and \roLMA. For simplicity, in this initial discussion, forecast horizons are fixed at $h=1$ for each experiment. For the \mytau calculations, this
means that $\mathcal{A}_{\tau} = I[\mathcal{S}_j,X_{j+1}]$, with
$\mathcal{S}_j= [X_{j},X_{j-\tau},\dots,X_{j-(m-1)\tau)}]^T$. Recall that with one-step forecasts, 1-MASE is the figure of merit. 
 Section~\ref{subsec:predictionhorizon} considers  
increasing the prediction horizon using $h$-MASE, with $h>1$, to assess accuracy. Section~\ref{subsec:datalength} considers the effects of the length of the traces.

\subsubsection{Synthetic Examples}
\label{subsec:synthetic}

The first step in this demonstration uses some standard synthetic examples,
both maps (H\'{e}non, logistic) and flows: the classic Lorenz 63
system~\cite{lorenz} and the Lorenz 96 atmospheric
model~\cite{lorenz96Model}.  
 The dynamics of each of these systems are reconstructed from the traces described in Chapter~\ref{ch:systems} using
different values $m$ and $\tau$. \mytau is computed for each of those reconstructed trajectories using a Kraskov-St\"ugbauer-Grassberger (KSG) estimator~\cite{KSG}, as described in Section~\ref{sec:KSG}. LMA is then used to
generate forecasts of every trace using each $\{m,\tau\}$ pair, their 1-MASE scores are computed as
described in Section~\ref{sec:accuracy}, and the relationship between the 1-MASE scores and the \mytau values
for the corresponding time series are discussed.  

\noindent\textbf{Flow Examples}
\label{subsubsec:lorenz96}

\begin{figure}[bt!]
  \centering
  \vspace{0.5cm}
  \begin{subfigure}[b]{0.49\textwidth}
    \includegraphics[width=\textwidth]{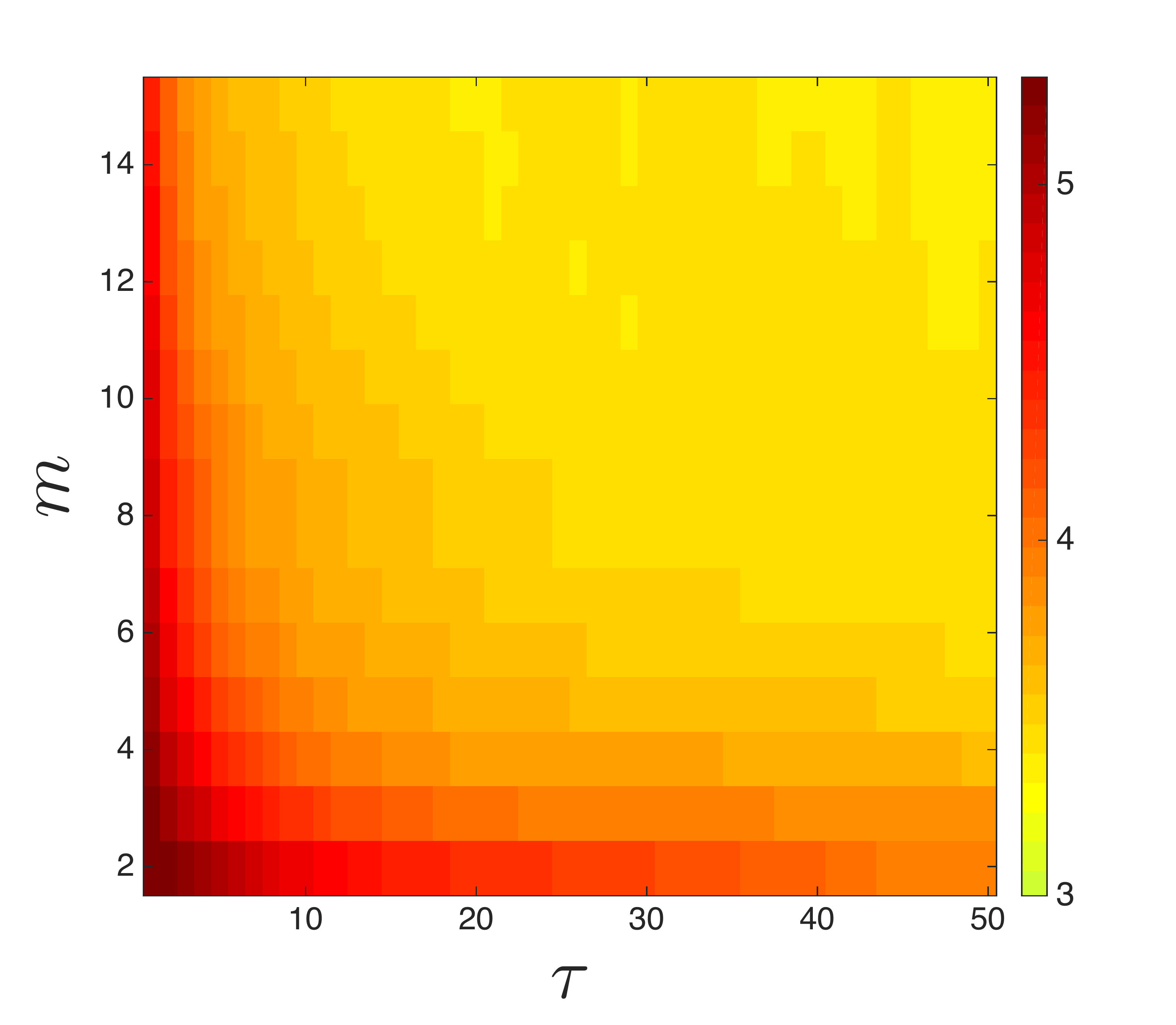}
    \caption{}
    \label{fig:L96N22F5SPI}
  \end{subfigure}
  \begin{subfigure}[b]{0.49\textwidth}
    \includegraphics[width=\textwidth]{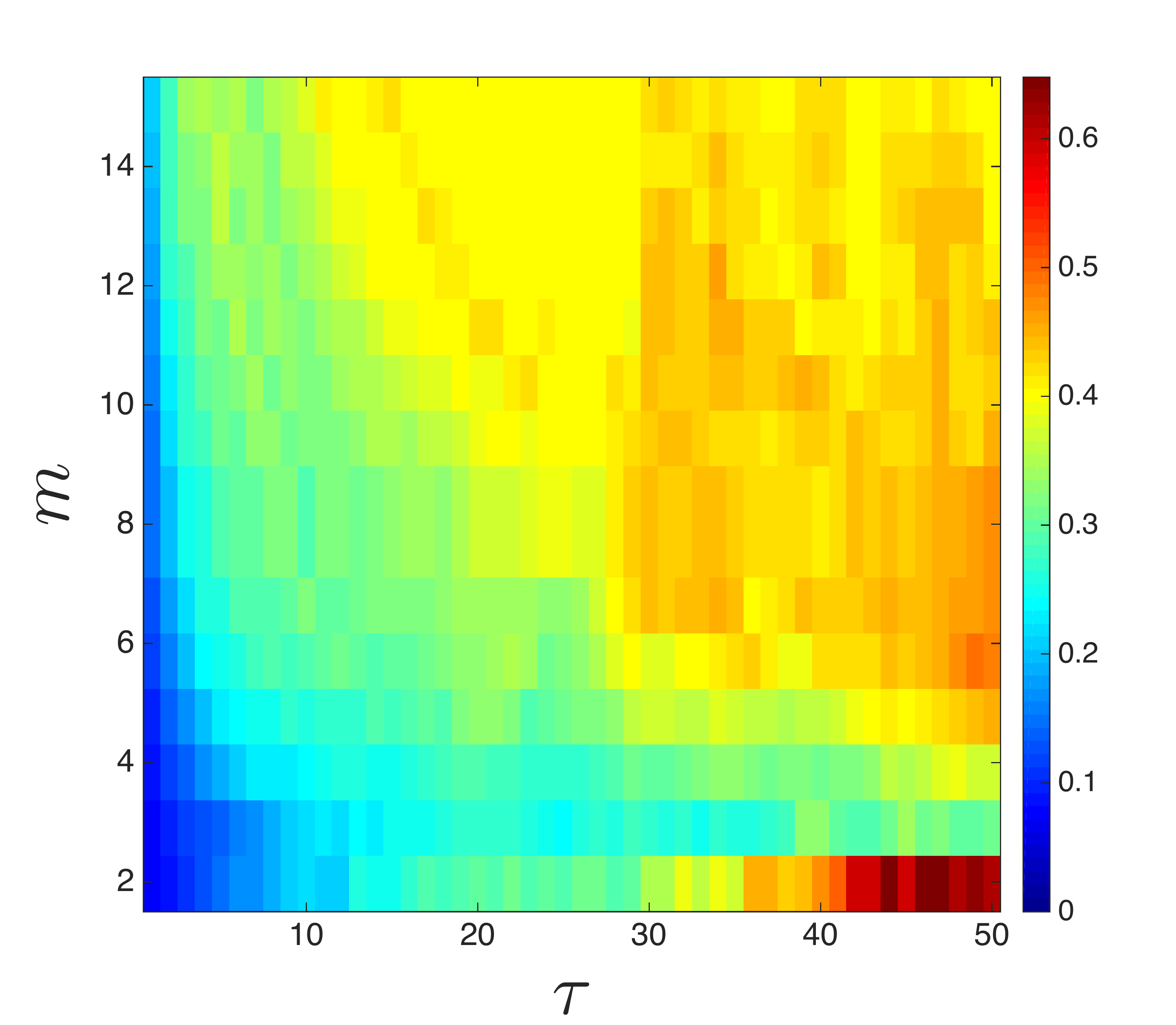}
    \caption{}
    \label{fig:L96N22F5MASE}
  \end{subfigure}
  \caption{The effects of reconstruction parameter values on \mytau and forecast accuracy for the Lorenz 96 system. (a) \mytau values for different delay reconstructions of a representative trace from that system with $\{K=22,F=5\}$. (b) 1-MASE scores for LMA forecasts on different delay reconstructions of that trace.}
  \label{fig:tauandmL96}
\end{figure}

Figure~\ref{fig:tauandmL96}(a) shows a heatmap of the \mytau values for
reconstructions of a representative trajectory from the Lorenz 96 system with
$\{K=22,F=5\}$, for a range of $m$ and $\tau$.
Not surprisingly, this image reveals a strong dependency between the
values of the reconstruction parameters and the reduction in
uncertainty about the near future that is provided by the
reconstruction.  Very low $\tau$ values, for instance, produce delay
vectors that have highly redundant coordinates, and so provide substantial
information about the immediate future. Again, the standard heuristics only focus on minimizing redundancy between coordinates, choosing the $\tau$ value that minimizes the mutual information between the first
two coordinates in the delay vector.  For this Lorenz 96 trajectory, that 
approach~\cite{fraser-swinney} yields $\tau=26$,
while standard dimension-estimation heuristics \cite{KBA92} suggest
$m=8$.  The \mytau value for a delay reconstruction built with those
parameter values is $3.471\pm0.051$.  This is {\it not}, however, the
\mytau-optimal reconstruction; choosing $m=2$ and $\tau=1$, for
instance, results in a higher value ($\mytau=5.301\pm0.019$)---{\it i.e.,}
significantly more reduction in uncertainty about the future.  This may
be somewhat counter-intuitive, since each of the delay vectors in the
\mytau-optimal reconstruction spans far less of the data set and thus
one would expect points in that space to contain {\it less}
information about the future.  Figure~\ref{fig:tauandmL96}(a) suggests,
however, that this in fact not the case; rather, the uncertainty
{\it increases} with both dimension and time delay.
\label{page:increase-with-tau}

The question at issue in this section is whether that reduction in
uncertainty about the future correlates with improved accuracy of an
LMA forecast built from that reconstruction.  Since the \mytau-optimal
choices
maximize the shared information between the state estimator and
$X_{j+1}$, one would expect a delay reconstruction model built with
those choices to afford LMA the most leverage.  To test that
conjecture, I perform an exhaustive search with $m=2,\dots,15$ and
$\tau=1,\dots,50$.  For each $\{m, \tau\}$ pair, I use LMA to
generate forecasts from the corresponding reconstruction, compute
their 1-MASE scores, and plot the results in a heatmap similar to
the one in Figure~\ref{fig:tauandmL96}(a).  As one would expect, the
1-MASE and \mytau heatmaps are generally antisymmetric.  This
antisymmetry breaks down somewhat for low $m$ and high $\tau$, where
the forecast accuracy is low even though the reconstruction contains a
lot of information about the future.

I suspect that this breakdown is due to a combination of overfolding (too-large values of $\tau$) and projection (low $m$).  Even though
each point in an overfolded reconstruction may contain a lot of information
about the future, the false crossings created by this combination of
effects pose problems for a near-neighbor forecast strategy like LMA.
The improvement that occurs if one adds another dimension is
consistent with this explanation.  Notice, too, that this effect only
occurs far from the maximum in the \mytau surface---the area that is
of interest if one is using \mytau to choose parameter values for
reconstruction models.

In general, though, maximizing the redundancy between the state
estimator and the future does appear to minimize the resulting
forecast error of LMA.  Indeed, the maximum on the surface of
Figure~\ref{fig:tauandmL96}(b) ($m=2,\tau=1$) is exactly the minimum on
the surface of Figure~\ref{fig:tauandmL96}(a).  The accuracy of this
forecast is almost six times higher (1-MASE $ = 0.074\pm0.002$) than that of a
forecast constructed with the parameter values suggested by the
standard heuristics ($0.441\pm0.033$).  Note that the optima of these
surfaces may be broad: {\it i.e.,} there may be {\it ranges} of $m$ and
$\tau$ for which \mytau and 1-MASE are optimal, and roughly constant.
In these cases, it makes sense to choose the lowest $m$ on the
plateau, since that minimizes computational effort, data requirements,
and noise effects. Notice that in this experiment, $m=2$ was 
actually the \mytau-optimal reconstruction dimension, and that correspondence let me calculate the forecast optimal $\tau$ for \roLMA without exhaustive search.  

\begin{figure}[tb!]
  \centering
  \vspace{0.5cm}
  \begin{subfigure}[b]{0.49\textwidth}
    \includegraphics[width=\textwidth]{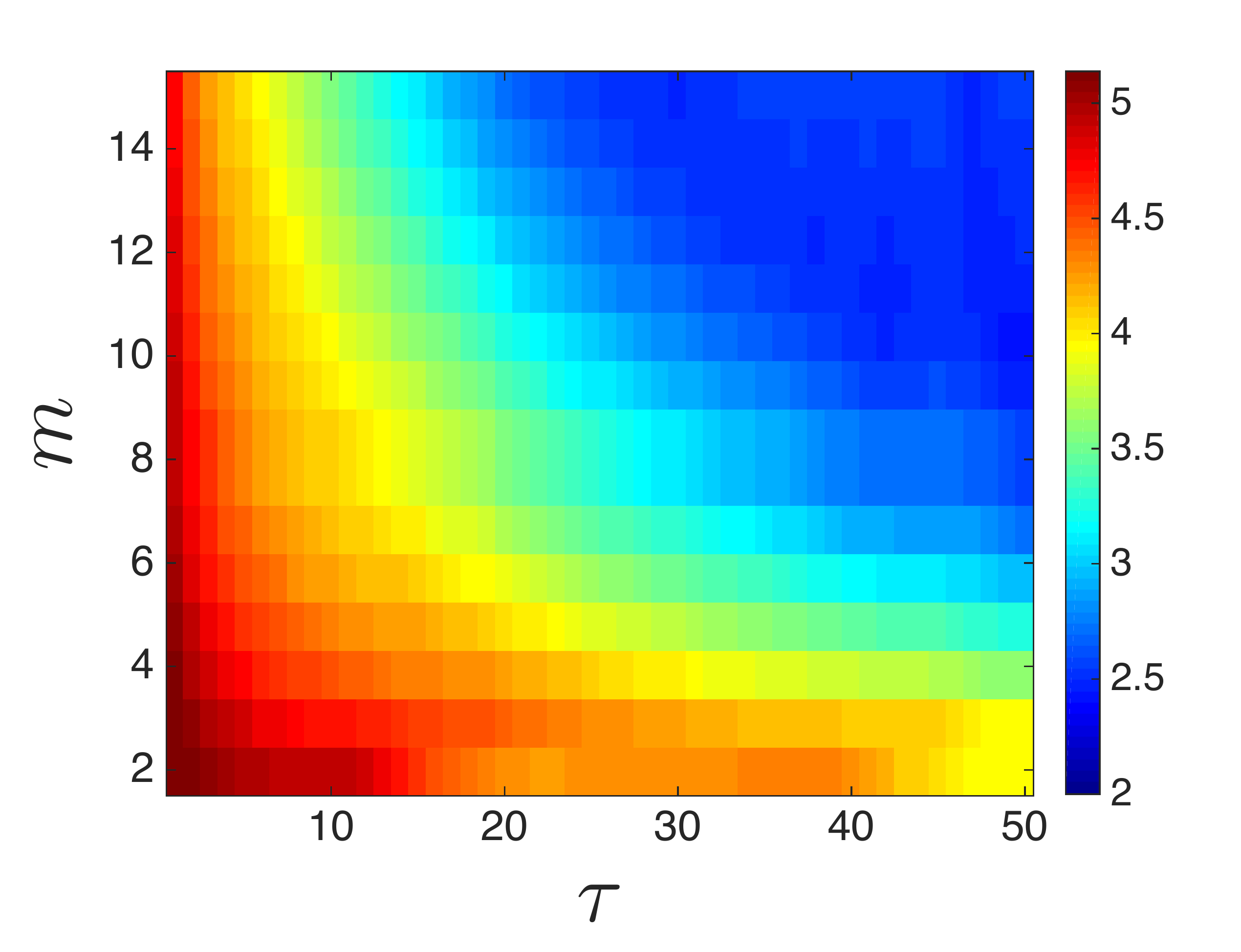}
    \caption{}
    \label{fig:L63SPI}
  \end{subfigure}
  \begin{subfigure}[b]{0.49\textwidth}
    \includegraphics[width=\textwidth]{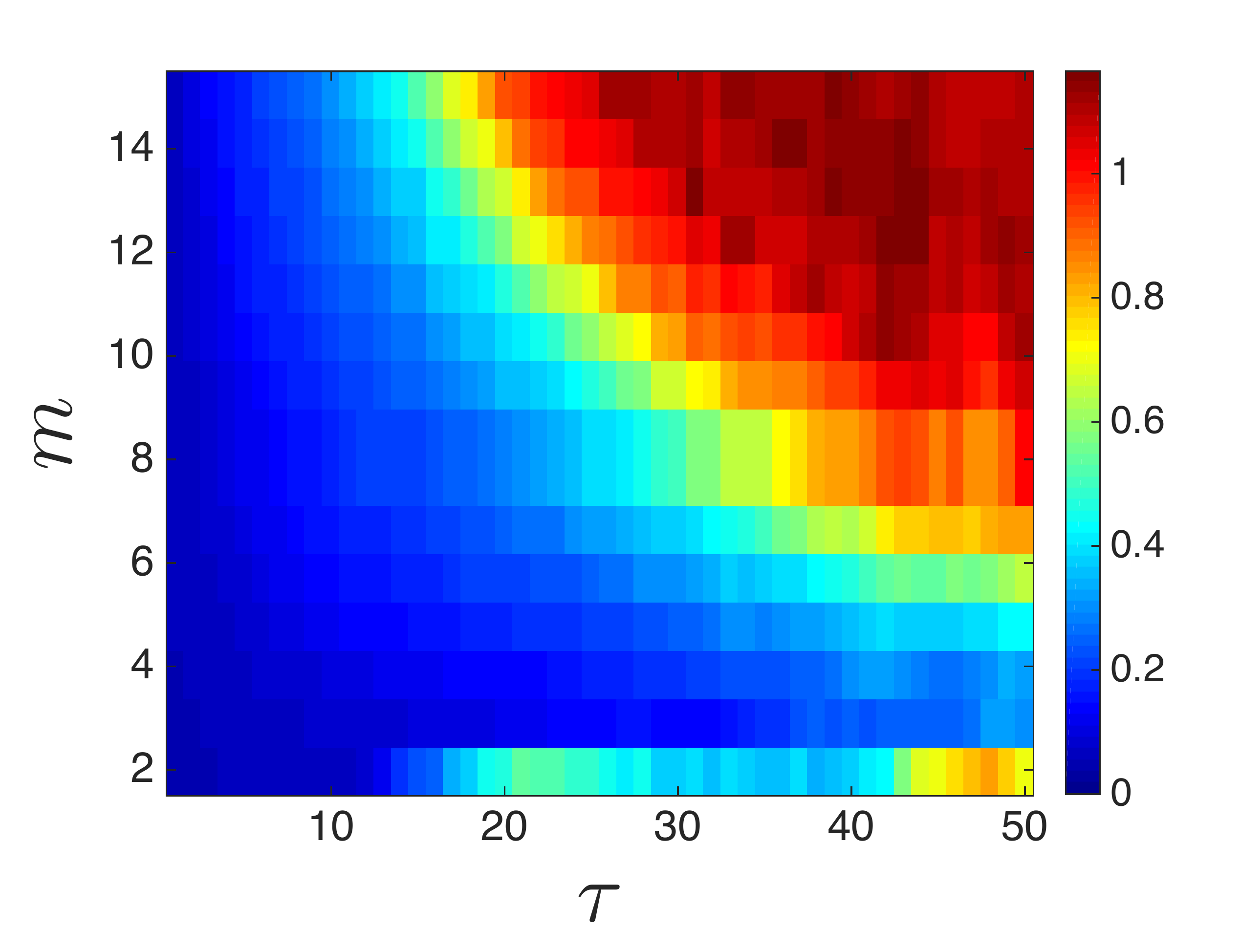}
    \caption{}
    \label{fig:L63MASE}
  \end{subfigure}
  \caption{The effects of reconstruction parameter values on \mytau and forecast accuracy for the Lorenz 63 system. (a) \mytau values for different delay reconstructions of a representative trace from that system. (b) 1-MASE scores for LMA forecasts on different delay reconstructions of that trace.}
  \label{fig:tauandmL63}
\end{figure}

While the results discussed in the previous paragraph do provide a
preliminary validation of the claim that one can use \mytau to select
good parameter values for delay reconstruction-based forecast
strategies, they only involve a single example system.  Similar
experiments on traces from the Lorenz 96 system with different
parameter values $\{K=47,F=5\}$ (not shown) demonstrate identical results---indeed, the
heatmaps are visually indistinguishable from the ones in
Figure~\ref{fig:tauandmL96}. Furthermore, for $\{K=47,F=5\}$, $m=2$ is 
again the \mytau-optimal reconstruction dimension, and \mytau again estimates the forecast optimal $\tau$ for \roLMA---quickly, without exhaustive search. Figure~\ref{fig:tauandmL63} shows
heatmaps of \mytau and 1-MASE for similar experiments on the canonical Lorenz 63 system~\cite{lorenz}.
As in the Lorenz 96 case, the heatmaps are generally antisymmetric,
confirming that maximizing \mytau is roughly equivalent to minimizing
1-MASE.  Again, though, the antisymmetry is not perfect; for high
$\tau$ and low $m$, the effects of projecting an overfolded attractor
cause false crossings that trip up LMA.  As before, adding a dimension
mitigates this effect by removing these false crossings.  Both the
Lorenz 63 and Lorenz 96 plots show a general decrease in
predictability for large $m$ and high $\tau$, with roughly hyperbolic
equipotentials dividing the colored\footnote{Note that the
  color map scales are not identical across all heatmap figures in
  this thesis; rather, they are chosen individually, to bring out the
  details of the structure of each experiment.}regions.  The locations and
heights of these equipotentials differ because the two signals are
not equally easy to predict.  This matter is discussed further at the
end of this section.

\begingroup
\renewcommand{\arraystretch}{1.3}
\begin{table*}[tb!]
  \caption{1-MASE values for various delay reconstructions of the
    different examples studied here.  1-MASE$_H$ is the representative
    accuracy of LMA forecasts that use delay reconstructions with
    parameter values ($m_{H}$ and $\tau_{H}$) chosen via
    standard heuristics for the corresponding traces.  Similarly, 1-MASE$_{\mytau}$
    is the accuracy of LMA forecasts that use reconstructions built
    with the $m$ and $\tau$ values that maximize \mytau, and 1-MASE$_E$
    is the error of the best forecasts for each case, found via
    exhaustive search over the $m,\tau$ parameter
    space. \alert{${**}$}: on these signals the standard heuristics
    failed.}
  \begin{center}
    \begin{tabular}{r|ccc}
      \hline\hline Signal & 1-MASE$_H$ & 1-MASE$_{\mytau}$ &  1-MASE$_E$ \\
      Parameters & $\{m_H,\tau_H\}$ &  $\{m_\mathcal{A_\tau},\tau_\mathcal{A_\tau}\}$ & $\{m_E,\tau_E\}$\\
      \hline
      \multirow{2}{*}{Lorenz-96 $K=22$} &$0.441 \pm 0.033$&$0.074\pm0.002$& $0.074\pm0.002$ \\
      				   &$\{8,26\}$&$\{2,1\}$& $\{2,1\}$ \\\hdashline

      \multirow{2}{*}{Lorenz-96 $K=47$} &$1.007 \pm 0.043$&$0.115 \pm 0.006$&$0.115 \pm 0.006$\\
      				   &$\{10,31\}$&$\{2,1\}$& $\{2,1\}$ \\\hdashline

      \multirow{2}{*}{Lorenz 63} &$0.144\pm0.008$& $0.062\pm0.006$& $0.058\pm0.005$\\
                       &$\{5,12\}$&$\{3,1\}$&$\{2,1\}$\\
      \hline
      \multirow{2}{*}{H\'enon Map} & \alert{$^{**}$}&$4.46\times 10^{-4}\pm 2.63\times 10^{-5}$&$4.46\times 10^{-4}\pm 2.63\times 10^{-5}$ \\ 
                             &$\{\alert{^{**}},\alert{^{**}}\}$&$\{2,1\}$&$\{2,1\}$\\\hdashline
     \multirow{2}{*}{Logistic Map} & \alert{$^{**}$}&$2.19\times10^{-5}\pm 2.72\times10^{-6}$&$2.19\times10^{-5}\pm 2.72\times10^{-6}$ \\ 
                             &$\{\alert{^{**}},\alert{^{**}}\}$&$\{1,1\}$&$\{1,1\}$\\
      \hline\hline
    \end{tabular}
  \end{center}
  \label{tab:myTauParams}
\end{table*}%
\endgroup
Numerical \mytau and 1-MASE values for LMA forecasts on different
reconstructions of both Lorenz systems are tabulated in the top three
rows of Table~\ref{tab:myTauParams}, along with the reconstruction
parameter values that produced those results. 
These results bring out two important points.  First, as
suggested by the heatmaps, the $m$ and $\tau$ values that maximize
\mytau (termed $m_{\mytau}$ and $\tau_{\mytau}$ in the table)
are close, or identical, to the values that minimize 1-MASE ($m_E$ and
$\tau_E$) for all three Lorenz systems.  This is notable because---as
discussed in Section~\ref{subsec:datalength}---the former can be
estimated quite reliably from a small sample of the trajectory in only
a few seconds of compute time, whereas the exhaustive search that is
involved in computing $m_E$ and $\tau_E$ for
Table~\ref{tab:myTauParams} required close to 30 hours of CPU time per system/parameter set ensemble.
A second important point that is apparent from Table~\ref{tab:myTauParams} is that
delay reconstructions built using the traditional heuristics---the
values with the $H$ subscript---are comparatively ineffective for the
purposes of LMA-based forecasting.  This is notable because that is
the default approach in the literature on state-space based
forecasting methods for dynamical systems. Moreover, in all cases $m_E$ and $m_{\mytau}$ are far 
lower than what the embedding theory would suggest, further corroborating the basic premise of this thesis. 

\begin{figure}[tb!]
  \centering
  \vspace{0.5cm}
  \begin{subfigure}[b]{0.49\textwidth}
    \includegraphics[width=\textwidth]{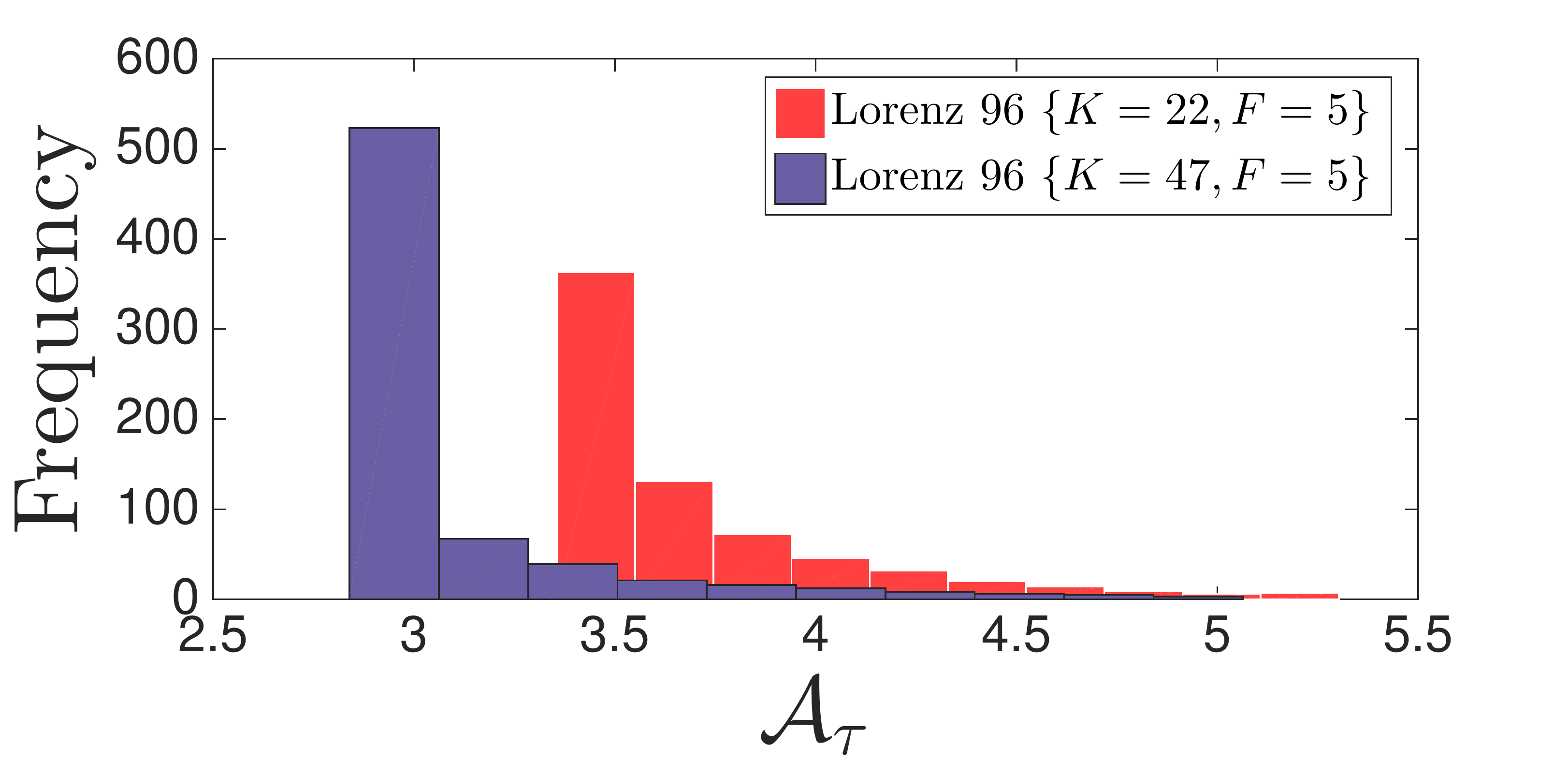}
    \caption{}
    \label{fig:L96SPIHIST}
  \end{subfigure}
  \begin{subfigure}[b]{0.49\textwidth}
    \includegraphics[width=\textwidth]{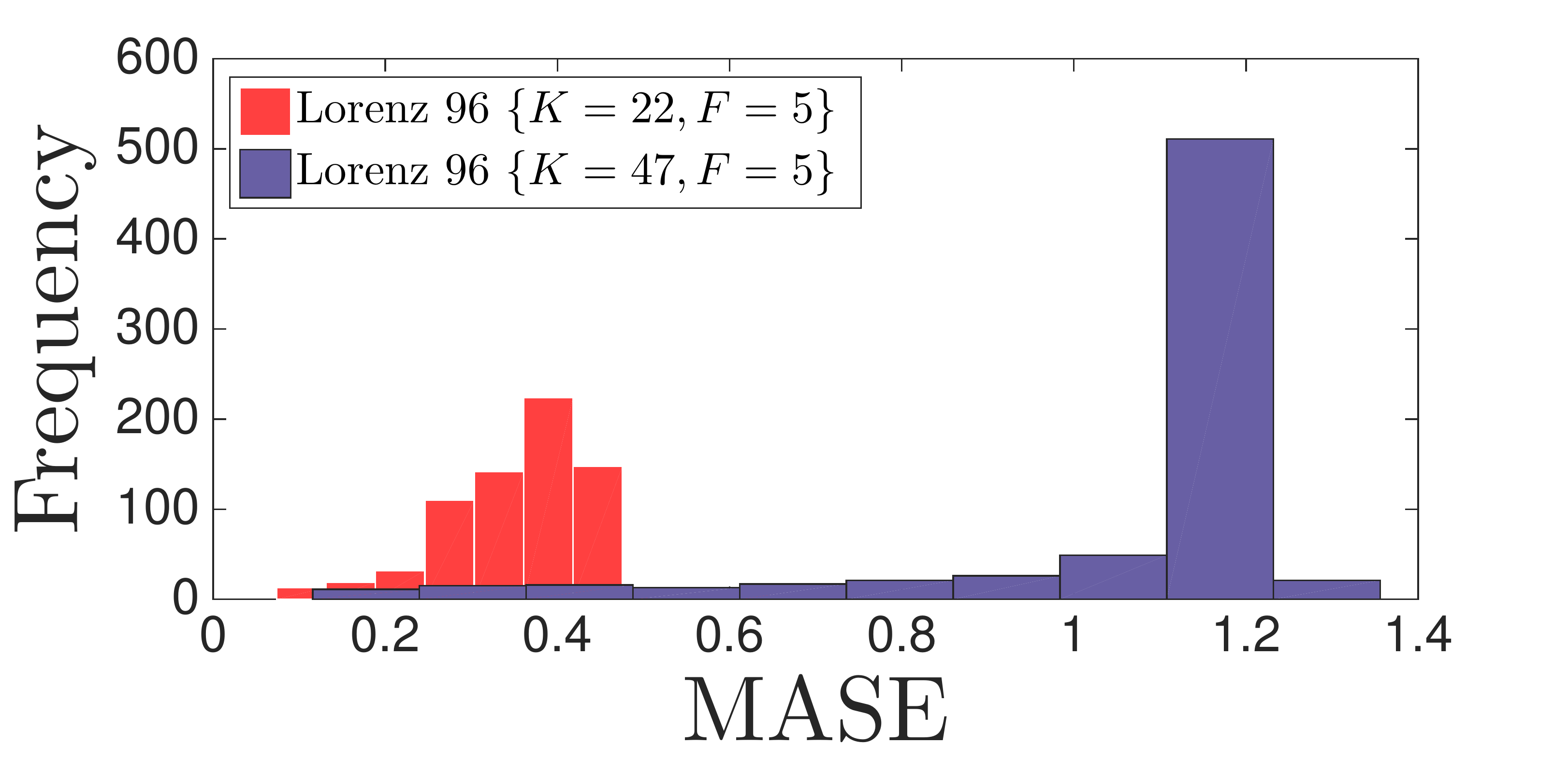}
    \caption{}
    \label{fig:L96MASEHIST}
  \end{subfigure}
  \caption{Histograms of \mytau and 1-MASE values for all traces from the Lorenz 96 $\{K=22,F=5\}$ and $\{K=47,F=5\}$ systems for all $\{m,\tau\}$ values in Figures~\ref{fig:tauandmL96} and~\ref{fig:tauandmL63}: (a) \mytau (b) 1-MASE.}
  \label{fig:L96HistCompares}
\end{figure}

A close comparison of Figures~\ref{fig:tauandmL96}
and~\ref{fig:tauandmL63} brings up another important point: some time
series are harder to forecast than others.
Figure~\ref{fig:L96HistCompares} breaks down the details of the two
suites of Lorenz 96 experiments, showing the distribution of \mytau
and 1-MASE values for all of the reconstructions.
Although there is some overlap in the $K=22$ and $K=47$
histograms---{\it i.e.,} best-case forecasts of the former are better than
most of the forecasts of the latter---the $K=47$ traces generally
contain less information about the future and thus are harder to
forecast accurately. As discussed in Section~\ref{sec:roLMALorenz96}, this is to be expected. 

\noindent\textbf{Map Examples}

Delay reconstruction of discrete-time dynamical systems, while
possible in theory, can be problematic in practice.  Although the
embedding theorems do apply in these cases, the heuristics for
estimating $m$ and $\tau$ often fail.  The time-delayed mutual
information of \cite{fraser-swinney}, for example, may decay
exponentially, without showing any clear minimum.  And the lack of
spatial continuity of the orbit of a map violates the underlying idea
behind the method of \cite{KBA92}.  State space-based forecasting
methods can, however, be very useful in generating predictions of
trajectories from systems like this---{\it if} one selects the two free 
parameters properly.

In view of this, it would be particularly useful if one could use
\mytau to choose embedding parameter values for maps.  This section
explores that notion using two canonical examples, shown in the bottom
two rows of Table~\ref{tab:myTauParams}.  For the H\'enon map
\begin{eqnarray}
  x_{n+1} &=& 1-ax^2_n +y_n  \\
  y_{n+1} &=& bx_n 
\end{eqnarray}
with $a=1.4$ and $b=0.3$, the \mytau-optimal parameter values,
$m=2$ and $\tau=1$, occur at \mytau$=6.617\pm0.011$, over the 15 trajectories generated from randomly-chosen initial conditions.  As in the flow examples, these are identical to
the values that minimized 1-MASE ($4.46 \times 10^{-4} \pm 2.63 \times 10^{-5}$).  These parameter values make sense, of course;
a first-return map of the $x$ coordinate is effectively the H\'enon
map, so $[x_j,x_{j-1}]$ is a perfect state estimator (up to a scaling
term).  But in practice, of course, one rarely knows the underlying
dynamics of the system that generated a time series, so the fact that
one can choose good reconstruction parameter values by maximizing
\mytau is notable---especially since standard heuristics for that
purpose fail for this system.

The same pattern holds for the logistic map, $x_{n+1} = rx_n(1-x_n)$,
with $r=3.65$. Again, for validation, I generate 15 trajectories from randomly-chosen initial conditions.  For this ensemble of experiments, the \mytau-optimal parameter values, which occur at \mytau$=9.057\pm0.001$, coincided
with the minimum of the 1-MASE surface ($2.19 \times10^{-5} \pm 2.72 \times10^{-6}$).  As in the H\'{e}non example,
these values ($m=1$ and $\tau=1$) make complete sense, given the form
of the map.  But again, one does not always know the form of the
system that generated a given time series.  In both of these map examples, the standard heuristics fail, but \mytau clearly
indicates that one does not actually need to reconstruct these
dynamics---rather, that near-neighbor forecasting {\it on the time series
  itself} is the best approach.

\subsubsection{Selecting Reconstruction Parameters of Experimental Time Series}
\label{subsec:experimental}

The previous section provided a preliminary verification
of the conjecture that parameters that maximize \mytau also maximize forecast accuracy  
for LMA, for both maps and flows.  While experiments with synthetic
examples are useful, it is important to show that \mytau is also a useful way to choose
parameter values for delay reconstruction-based forecasting of
real-world data, where the time series are noisy and perhaps short,
and one does not know the dimension of the underlying system---let
alone its governing equations.  This section extends the exploration in the previous section, using experimental data from two different dynamical
systems: a far-infrared laser and a laboratory computer-performance
experiment.

\noindent\textbf{A Far-Infrared Laser}
\label{subsubsec:laser}

\begin{figure}[tb!]
  \centering
  \vspace{0.5cm}
  \begin{subfigure}[b]{0.49\textwidth}
    \includegraphics[width=\textwidth]{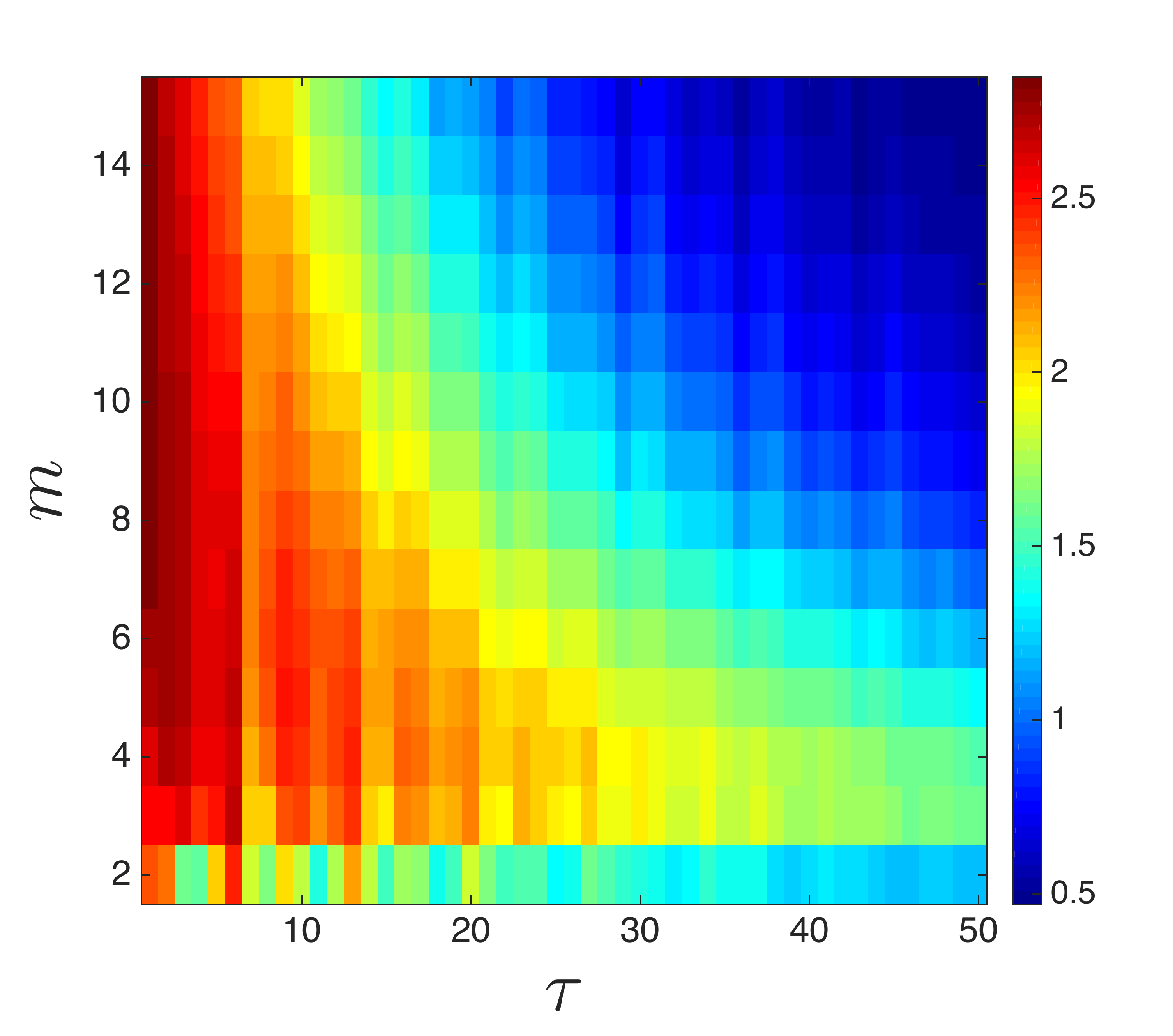}
    \caption{}
    \label{fig:LASERMASE}
  \end{subfigure}
  \begin{subfigure}[b]{0.49\textwidth}
    \includegraphics[width=\textwidth]{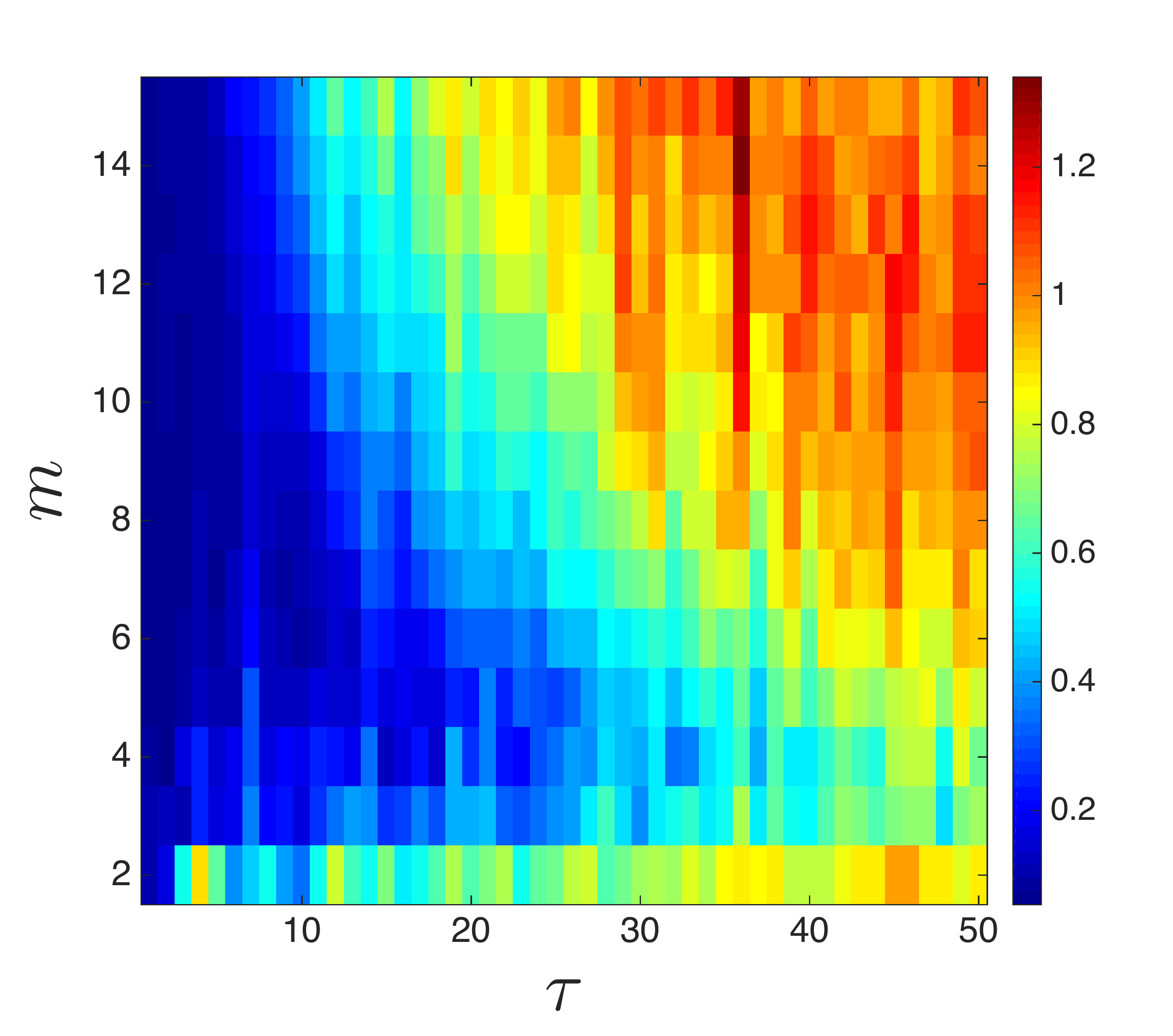}
    \caption{}
    \label{fig:LASERSPI}
  \end{subfigure}
  \caption{The effects of reconstruction parameter values on \mytau and forecast accuracy for SFI dataset A. (a) \mytau values for different delay reconstructions of that signal. (b) 1-MASE scores for LMA forecasts of those reconstructions.}
  \label{fig:MASEandmyTAULaser}
\end{figure}

I begin this discussion by returning to the canonical test case from Chapter~\ref{ch:overview}, SFI dataset A~\cite{weigend-book}, which was gathered from a far-infrared
laser.  As in the synthetic examples in Section~\ref{subsec:synthetic}, the
\mytau and 1-MASE heatmaps (Figure~\ref{fig:MASEandmyTAULaser}) are
largely antisymmetric for this signal.
Again, there is a band across the bottom of each image because of the
combined effects of overfolding and projection.  Note the similarity
between Figures ~\ref{fig:MASEandmyTAULaser} and
~\ref{fig:tauandmL63}: the latter resemble ``smoothed'' versions of
the former.  It is well known~\cite{weigend-book} that the SFI dataset A is well described by the Lorenz 63 system with some added
noise, so this similarity is reassuring.  An LMA
forecast using the \mytau-optimal reconstruction of this trace was
more accurate\footnote{Note that the SFI dataset A 1-MASE values are not averages as there is only one trace.}than similar forecasts using  reconstructions built
using traditional heuristics---1-MASE$_{\mytau} = 0.0592$ versus
1-MASE$_H=0.0733$---and only slightly worse than the optimal value 1-MASE$_E=0.0538$.  
However, the values of $\{m_{\mytau},\tau_{\mytau}\}$ and
$\{m_E,\tau_E\}$ are not identical for this signal.  This is because
the optima in the heatmaps in Figure~\ref{fig:MASEandmyTAULaser} are
bands, rather than unique points---as was the case in the synthetic
examples in Section~\ref{subsec:synthetic}.  In a situation like this,
a range of $\{m,\tau\}$ values are statistically indistinguishable,
from the standpoint of the forecast accuracy afforded by the
corresponding reconstruction.  The values suggested by the \mytau
calculation ($m_{\mytau}=9$ and $\tau_{\mytau}=1$) and by the
exhaustive search ($m_E=7$, $\tau_E=1$) are all on this
plateau, those suggested by the traditional heuristics ($m_H=7$, $\tau_H=3$) however are not.
 Again, these results suggest that one can use \mytau to choose good parameter
values for delay reconstruction-based forecasting, but SFI dataset A is only a
single trace from a fairly simple system.

\noindent\textbf{Computer Performance Dynamics}
\label{subsubsec:computer}

\begin{figure}[ht!]
  \centering
  \vspace{0.5cm}
  \begin{subfigure}[b]{0.49\textwidth}
    \includegraphics[width=\textwidth]{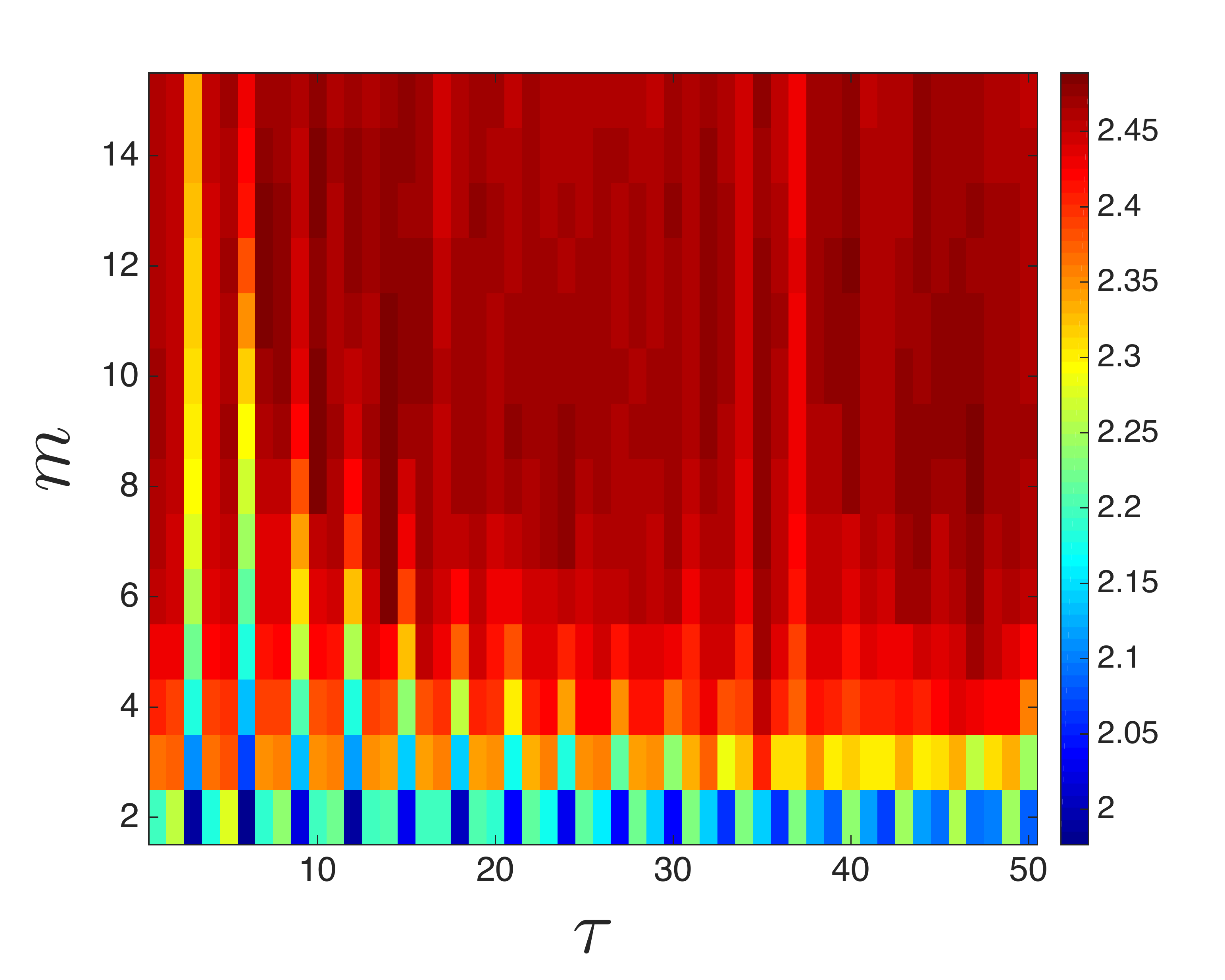}
    \caption{}
    \label{fig:colMASE}
  \end{subfigure}
  \begin{subfigure}[b]{0.49\textwidth}
    \includegraphics[width=\textwidth]{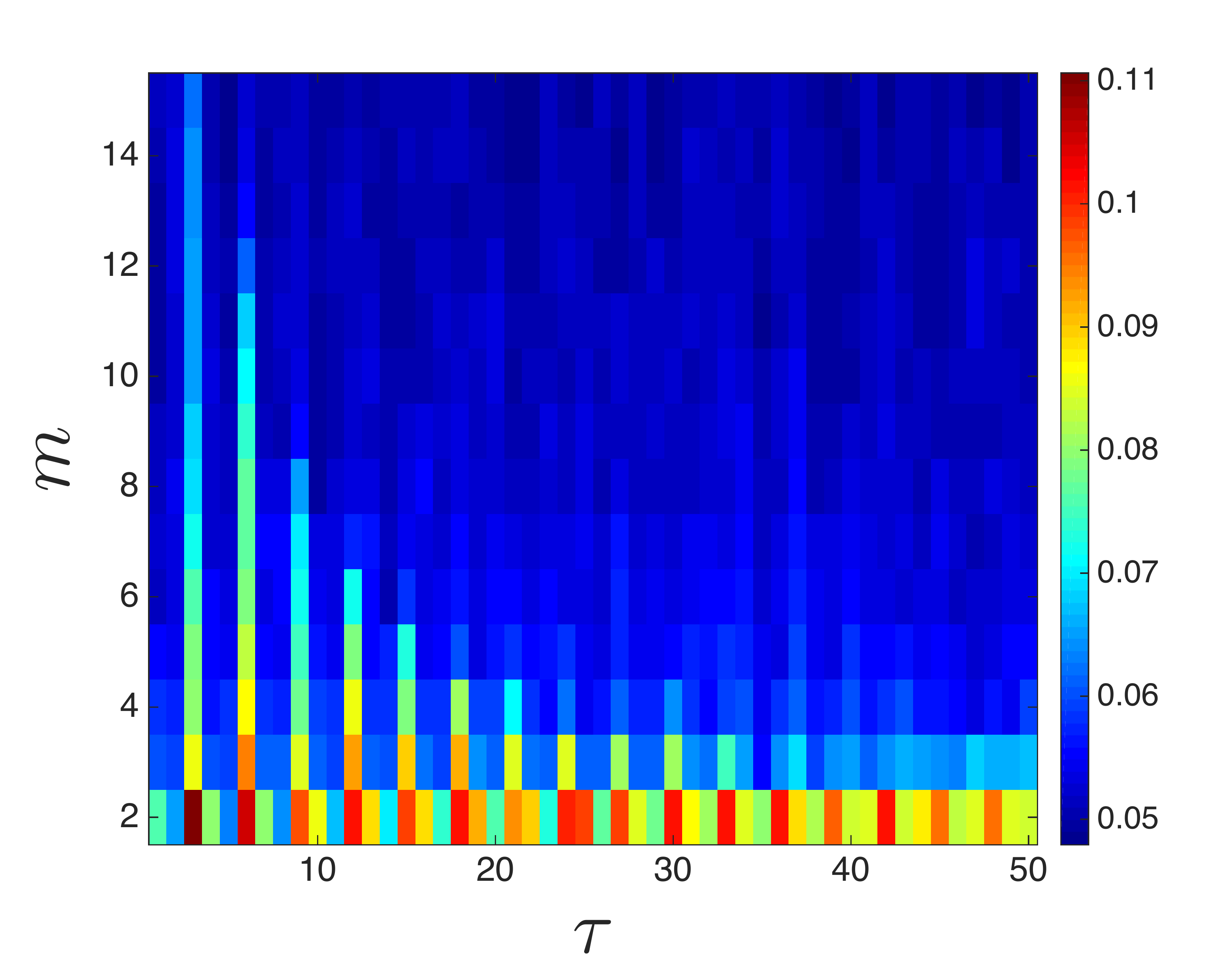}
    \caption{}
    \label{fig:colmytau}
  \end{subfigure}
  \caption{The effects of reconstruction parameter values on \mytau and forecast accuracy for a representative trace of \col. (a) \mytau values for different delay reconstructions of that trace. (b) 1-MASE scores for LMA forecasts on those reconstructions.}
  \label{fig:tauandmcol}
\end{figure}

Finally, I will return to the computer performance dynamics of \col and \gcc: experiments that involve multiple traces from each system, which allows for statistical analysis. 
As in the previous examples (Lorenz 63, Lorenz 96, H\'enon Map, Logistic Map, SFI dataset A), heatmaps of 1-MASE and \mytau for a representative 
\col time series---Figure~\ref{fig:tauandmcol}(b)---are largely
antisymmetric.
And again, reconstructions using the \mytau-optimal parameter values
allowed LMA to produce highly accurate forecasts of this signal:
1-MASE$_{\mytau}=0.050\pm0.002$, compared to the optimal 1-MASE$_E=0.049\pm0.002$.  There
are several major differences between these plots and the previous ones, though, beginning with the vertical stripes.
These are due to the dominant unstable periodic orbit of period 3 in
the chaotic attractor in the \col dynamics.  When $\tau$ is a multiple
of this period ($\tau=3\kappa$),
\label{page:3k} the coordinates of the delay
vector are not independent, which lowers \mytau and makes forecasting
more difficult.  (There is a nice theoretical discussion of this
effect in \cite{sauer91}.)
Conversely, \mytau spikes and 1-MASE plummets when $\tau=3\kappa-1$, since
the coordinates in such a delay vector cannot share any prime factors
with the period of the orbit.  The band along the bottom of both
images is, again, due to a combination of overfolding and projection. The other 14 traces in this experiment yield structurally identical heatmaps and the variance between these trials were only $\pm0.037$ on average.

\begin{figure}[tb!]
  \centering
  \vspace{0.5cm}
  \begin{subfigure}[b]{\columnwidth}
    \includegraphics[width=\columnwidth]{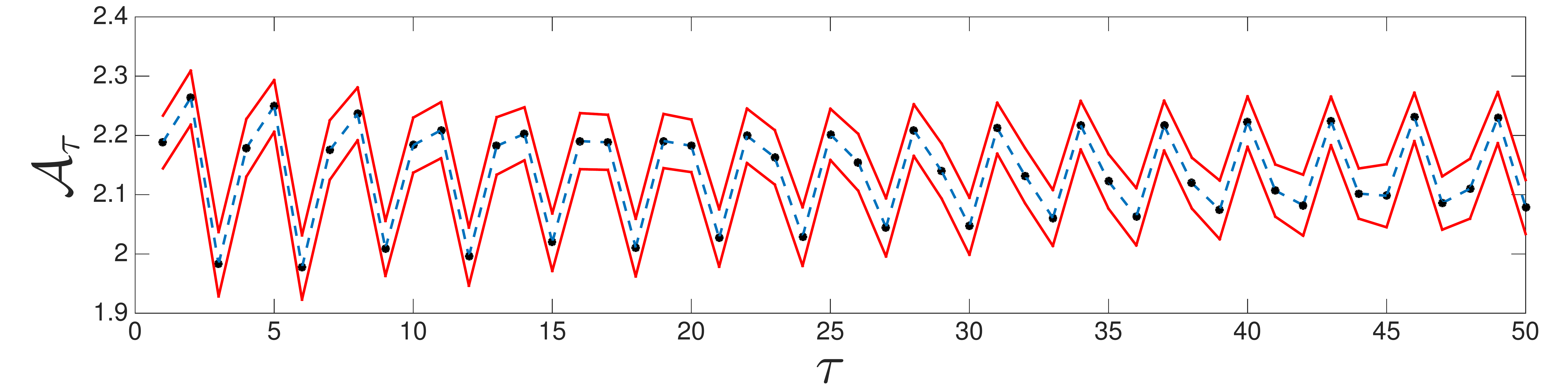}
    \caption{}
    \label{fig:colM2mytau}
  \end{subfigure}
  \begin{subfigure}[b]{\columnwidth}
    \includegraphics[width=\columnwidth]{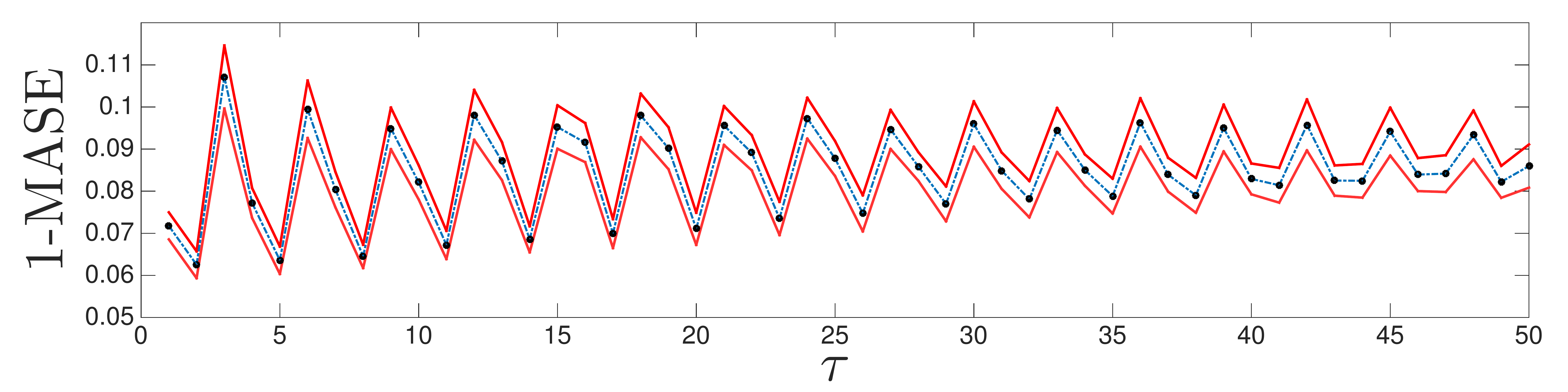}
    \caption{}
    \label{fig:L96N47F5MASE}
  \end{subfigure}%
   \caption{1-MASE and \mytau for \roLMA forecasts of all 15 \col traces, plotted as a function of $\tau$.  The blue dashed curves show the averages across all trials; the red dotted lines are that average $\pm$ the standard deviation. (a) \mytau values for delay reconstructions of these traces with $m=2$ and a range of values of $\tau$. (b) 1-MASE scores for \roLMA forecasts of those reconstructions.}
  \label{fig:MASEandmyTAUcol}
\end{figure}

Another difference between the \col heatmaps and the ones in
Figures~\ref{fig:tauandmL96}, \ref{fig:tauandmL63},
and~\ref{fig:MASEandmyTAULaser} is the apparent overall trend: the
``good'' regions (low 1-MASE and high \mytau) are in the lower-left
quadrants of those heatmaps, but in the upper-right quadrants of
Figure~\ref{fig:tauandmcol}.  This is partly an artifact of the
difference in the color-map scale, which is chosen here to bring out
some important details of the structure, and partly due to that
structure itself.  Specifically, the optima of the \col heatmaps---the
large dark red and blue regions in Figures~\ref{fig:MASEandmyTAUcol}(a) 
and (b), respectively---are much broader than the ones
in the earlier discussion of this section, perhaps because the signal is
so close to periodic.  (This is also the case to some extent in the
SFI Dataset A example, for the same reason.)  This geometry makes precise
comparisons of \mytau-optimal and 1-MASE-optimal parameter values
somewhat problematic, as the exact optima on two almost-flat but
slightly noisy landscapes may not be in the same place.  Indeed, the
\mytau values at $\{m_{\mytau},\tau_{\mytau}\}$ and $\{m_E,\tau_E\}$
are within a standard error across all 15 traces of \col.

And that brings up an interesting tradeoff.  For practical purposes,
what one wants is $\{m_{\mytau},\tau_{\mytau}\}$ values that produce a
1-MASE value that is {\it close to} the optimum 1-MASE$_E$.  However,
the algorithmic complexity of most nonlinear time-series analysis and
prediction methods scales badly with $m$.  In cases where the \mytau
maximum is broad, then, one might want to choose the lowest value of
$m$ on that plateau---or even a value that is on the {\it shoulder}
of that plateau, if one needs to balance efficiency over accuracy.
As I showed in Chapter~\ref{ch:pnp}, using \roLMA does just that, and that appears to work quite well for the \col data.
This amounts to marginalizing the heatmaps in
Figure~\ref{fig:tauandmcol} with $m=2$, which produces cross sections like the
ones shown in Figure~\ref{fig:MASEandmyTAUcol}.
The antisymmetry between \mytau and 1-MASE is quite apparent in these
plots; the global maximum of the former coincides with the global
minimum ($0.0649 \pm0.003$) of the latter, at $\tau=2$.  
 This is not much higher than the overall optimum of $0.0496 \pm0.002$---a value
from forecasts whose free parameters requires almost six orders of
magnitude more CPU time to compute. This result not only corroborates the main premise of this thesis, but also suggests a more effective way to calculate \mytau one simply fixes $m=2$, as is done with \roLMA, then selects $\tau$ by calculating \mytau across a
range of $\tau$s, rather than across a 2D $\{m,\tau\}$ space.  

The correspondence between 1-MASE and \mytau also holds true for other
marginalizations: {\it i.e.,} the minimum 1-MASE and the maximum \mytau
occur at the same $\tau$ value for all $m$-wise slices of the \col
heatmaps, to within statistical fluctuations.  The methods of~\cite{fraser-swinney} and \cite{KBA92}, incidentally, suggest
$\tau_H=2$ and $m_H=12$ for these traces; the average 1-MASE of an LMA
forecast on such a reconstruction is $0.0530 \pm 0.002$, which is somewhat better
than the best result from the $m=2$ marginalization, although still
short of the overall optimum.  The correspondence between $\tau_H$ and
$\tau_{\mytau}$ is coincidence; for this particular signal, maximizing
the independence of the coordinates happens to maximize the
information about the future that is contained in each delay vector. This is most likely due to the strength of the unstable three cycle present in these dynamics. In this case, the coordinates would be maximally independent {\it and} contain the most information about the future when $\tau=\rho-1$, where $\rho$ is the period of the dynamics.  The
$m=12$ result is not coincidence---and quite interesting, in view of
the fact that the $m=2$ forecast is so good.  It is also surprising in
view of the huge number of transistors---potential state
variables---in a modern computer.  As described in \cite{mytkowicz09},
however, the hardware and software constraints in these systems
confine the dynamics to a much lower-dimensional manifold. 

The \col program is what is known in the computer-performance
literature as a ``micro-kernel''---a extremely simple example that is
used in proof-of-concept testing.  The fact that its dynamics are so
rich speaks to the complexity of the hardware (and the
hardware-software interactions) in modern computers; again, see~\cite{todd-phd,mytkowicz09} for a much deeper discussion of these
issues.  Modern computer programs are far more complex than this
simple micro-kernel, of course, which begs the question: what does
\mytau tell us about the dynamics of truly complex systems like
the memory or processor usage patterns of sophisticated programs---which the computer performance community models as
stochastic systems?

\begin{figure}[tb!]
  \centering
           \vspace*{0.5cm}
  \begin{subfigure}[b]{0.49\textwidth}
    \includegraphics[width=\textwidth]{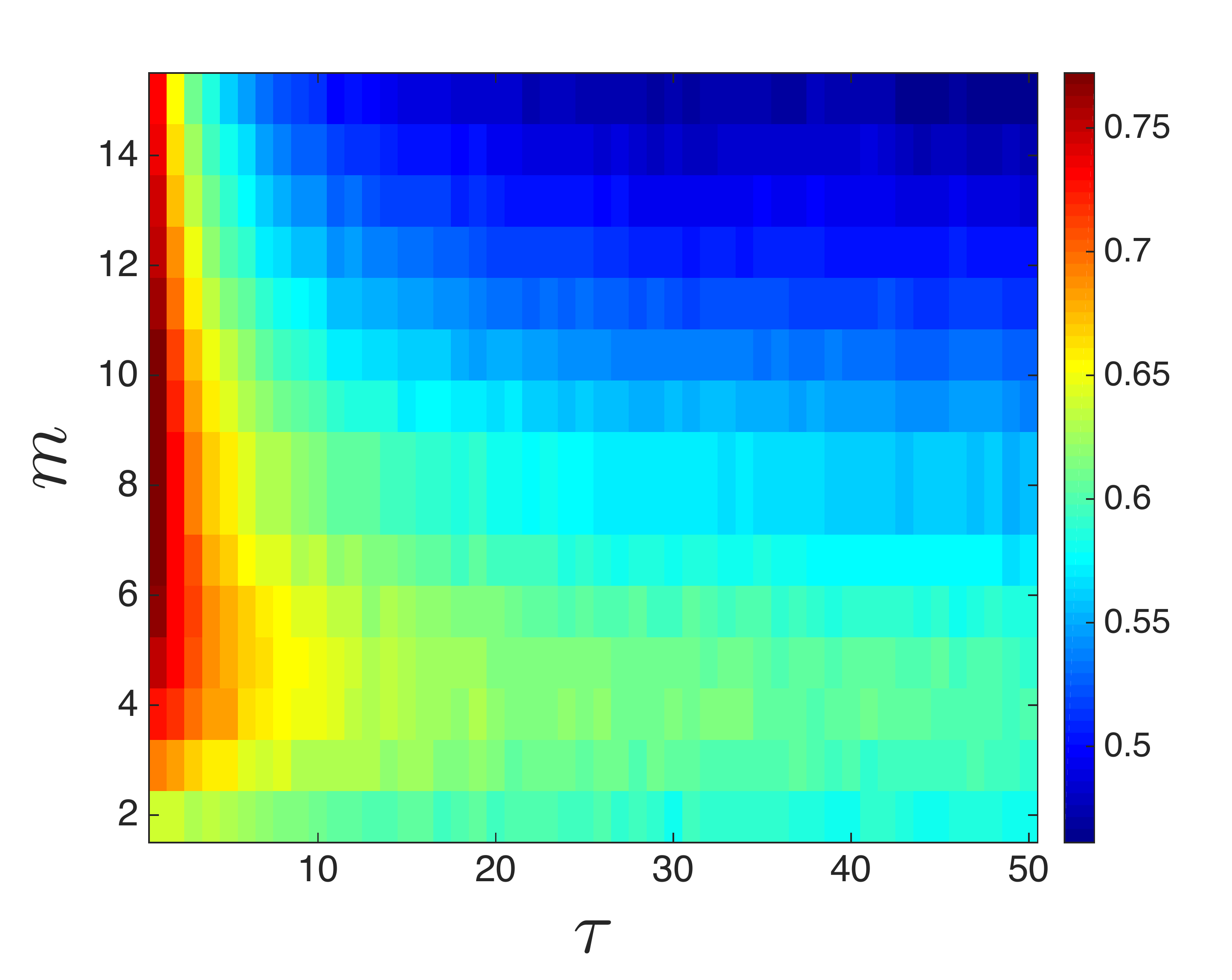}
    \caption{}
    \label{fig:gccmytau}
  \end{subfigure}
  \begin{subfigure}[b]{0.49\textwidth}
    \includegraphics[width=\textwidth]{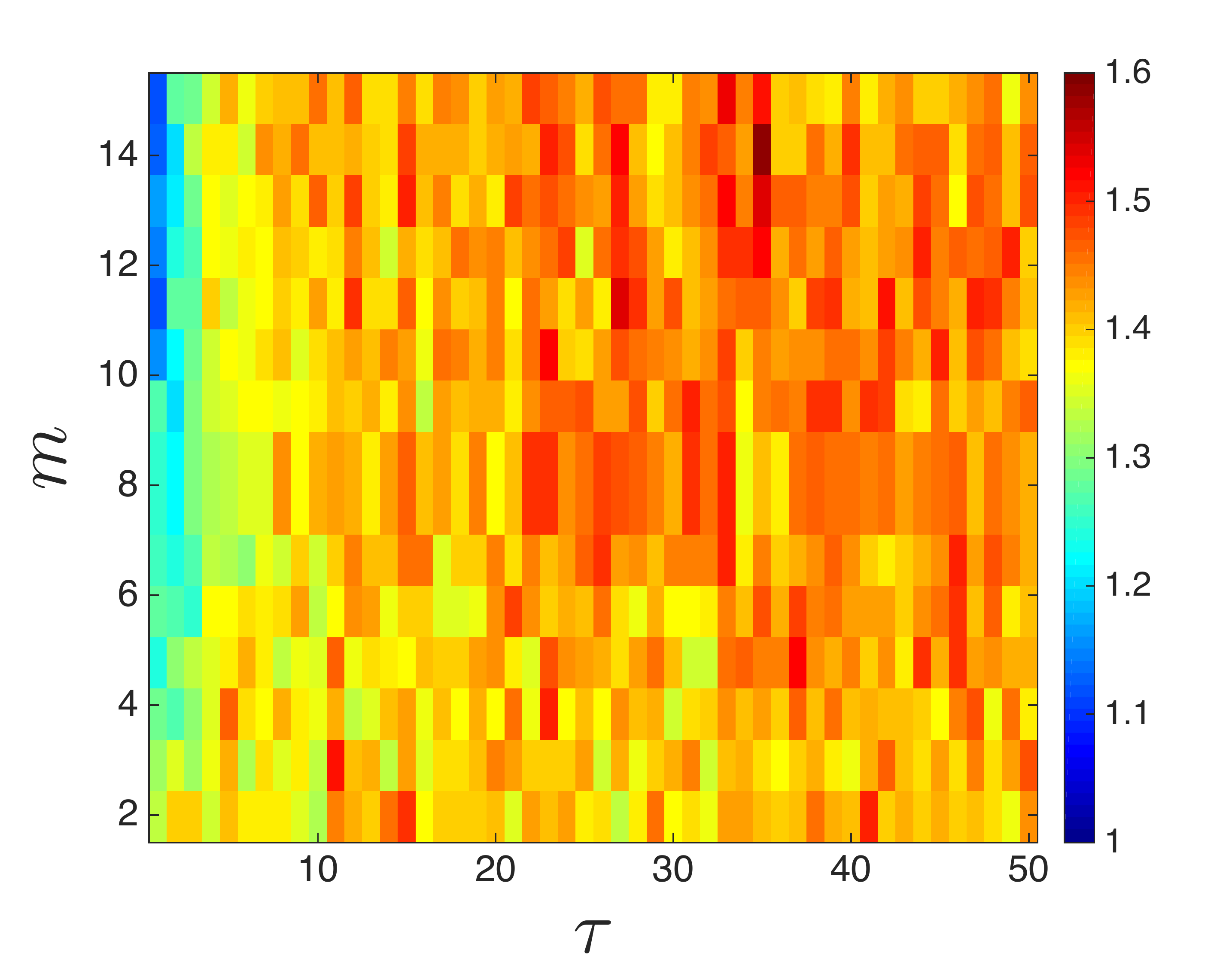}
    \caption{}
    \label{fig:gccMASE}
  \end{subfigure}
  \caption{The effects of reconstruction parameter values on \mytau and forecast accuracy for a representative trace from a computer-performance dynamics experiment using the \gcc benchmark. (a) \mytau values for different delay reconstructions of this trace. (b) 1-MASE scores for LMA forecasts on those reconstructions.}
  \label{fig:tauandmgcc}
\end{figure}

For \gcc, the answer is, again, that \mytau appears to be an effective
and efficient way to assess predictability. As shown in \cite{josh-pre} and synopsized in Chapter~\ref{ch:wpe}, this time series shares little to no
information with the future: {\it i.e.,} it {\it cannot} be predicted
using delay reconstruction-based forecasting methods, regardless of
$\tau$ and $m$ values.  The experiments in \cite{josh-pre} required
dozens of hours of CPU time to establish that conclusion; \mytau gives
the same results in a few seconds, using much less data.  The
structure of the heatmaps for this experiment, as shown in
Figure~\ref{fig:tauandmgcc}, is radically different.
The patterns visible in the previous 1-MASE plots, and the
antisymmetry between \mytau and 1-MASE plots, are absent from this pair of images, reflecting the lack of predictive content
in this signal.  Note, too, that the color map scales are different in this
figure.  This reflects the much-lower values of \mytau for this 
signal: over all 15 experiments of \gcc, for the parameter range in Figure~\ref{fig:tauandmgcc}, \mytau reached an absolute maximum of 0.7722, compared to the absolute maximum of 5.3026 for all experiments of the 
Lorenz 96 with $K=22$.
%
%
Indeed, the 1-MASE surface in Figure~\ref{fig:tauandmgcc}(b) never dips
below 1.0, Figure~\ref{fig:tauandmL96}, in contrast, never
  exceeds $\approx 0.6$ and generally stays below 0.3. These results are consistent across all traces in these experiments, {\it i.e.,} for all 15 traces of \gcc, 1-MASE never drops below 1.0. That is,
regardless of parameter choice, LMA forecasts of \gcc are no better
than simply using the prior value of this scalar time series as the
prediction. In comparison, with every experiment with Lorenz 96 $K=22$---regardless of parameter choice---the 1-MASE for LMA generally stays below 0.3---more than twice as good as a random walk.  The uniformly low \mytau values in
Figure~\ref{fig:tauandmgcc}(a) are an effective indicator of this---and,
again, they can be calculated quickly, from a relatively small sample
of the data.  It is to that issue that I turn next.

\subsection{Data Requirements and Prediction Horizons}
\label{sec:dataandhorizon}

In some real-world situations, it may be impractical to rebuild
forecast models at every step, as I have done in the previous
sections of this thesis---because of computational expense, for
instance, or because the data rate is very high.  In these situations,
one may wish to predict $h$ time steps into the future, then stop and
rebuild the model to incorporate the $h$ points that have arrived
during that period, and repeat.  

In chaotic systems, of course, there
are fundamental limits on prediction horizon even if one is working
with infinitely long traces of all state variables.  A key question at
issue in this section is how that effect plays out in forecast models
that use delay reconstructions from scalar time-series data. I explore that issue in Section~\ref{subsec:predictionhorizon}.  And
since real-world data sets are not infinitely long, it is also important to
understand the effects of data length on the estimation of \mytau. I explore this question in the following section, using one-step-ahead forecasts so that I can compare the results to those in the previous sections.

\subsubsection{Data Requirements for \mytau Estimation}
\label{subsec:datalength}


The quantity of data used in a delay reconstruction directly impacts
the usefulness of that reconstruction.  If one is interested in
approximating the correlation dimension via the Grassberger-Procaccia
algorithm, for instance, it has been shown that one needs
$10^{(2+0.4m)}$ data points~\cite{tsonisdatabound,smithdatabound}.
Those bounds are overly pessimistic for forecasting, however, as mentioned in Section~\ref{sec:datalength}.  A key challenge, then, is to determine whether
one's time series {\it really} calls for as many dimensions and data
points as the theoretical results require, or whether one can get away
with fewer dimensions---and how much data one needs in order to figure
all of that out.

\mytau is a useful solution to those challenges.  As
established in the previous sections, calculations of this quantity
can reveal what dimension is required for delay reconstruction-based forecasting of dynamical systems.  And, as alluded to there, \mytau can be
estimated accurately from a surprisingly small number of points.  The
experiments in this section explore that intertwined pair of claims in
more depth by increasing the length of the Lorenz 96 traces and
testing whether the information content of the state estimator derived
from standard heuristics converges to the \mytau-optimal
estimator. (This kind of experiment is not possible in
practice, of course, when the time series is fixed, but can be done
in the context of this synthetic example.)

\begin{figure}[b!]
  \centering
  \begin{subfigure}[b]{\columnwidth}
    \includegraphics[width=\columnwidth]{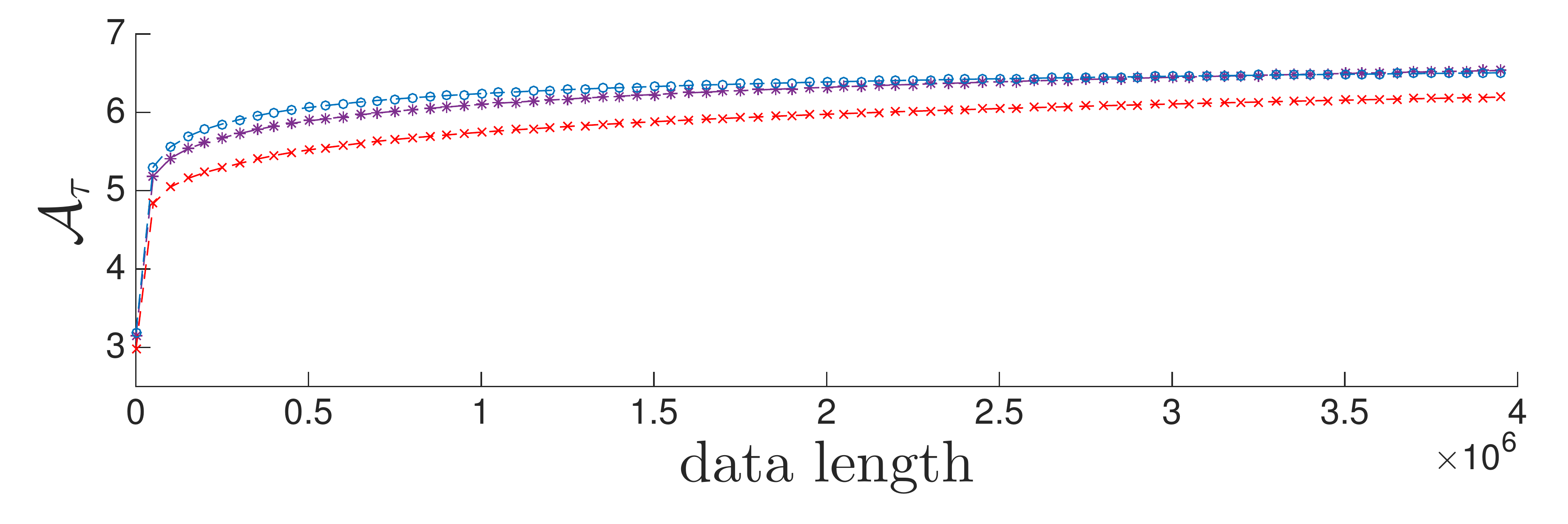}
    \caption{Lorenz-96 $\{K=22,F=5\}$ system}
    \label{fig:data22}
  \end{subfigure}
  ~
  \begin{subfigure}[b]{\columnwidth}
    \includegraphics[width=\columnwidth]{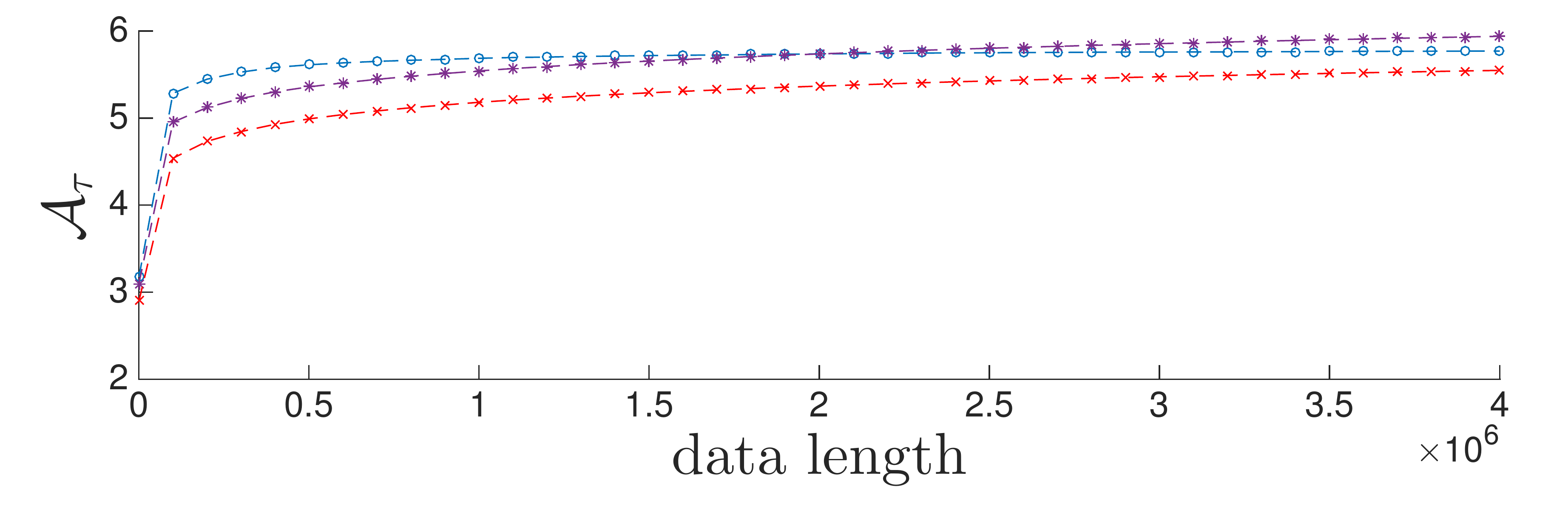}
    \caption{Lorenz-96 $\{K=47,F=5\}$ system}
  \end{subfigure}%
  \caption{Average optimal \mytau versus data length for all 15 traces from the Lorenz-96 system using $\tau=1$ in all cases. Blue circles corresponds to $m=2$, purple diamonds to $m=4$, and red xs to $m=8$. (a) Optimal \mytau for traces from the $\{K=22,F=5\}$ system. (b) Optimal \mytau for traces from the $\{K=47,F=5\}$ system.}
  \label{fig:mytaudata}
\end{figure}

Figure~\ref{fig:mytaudata} shows the results.  
When the data length is short, the $m=2$ state estimator has the most
information about the future.  This makes perfect sense; a short time
series cannot fully sample a complicated object, and when an
ill-sampled high-dimensional manifold is projected into a low-dimensional space, infrequently visited regions of that manifold can
act effectively like noise.  From an information-theoretic standpoint,
this would increase the effective Shannon entropy rate of each of the
variables in the delay vector.  In the $I$-diagram in
Figure~\ref{fig:compare}(b), this would manifest as drifting apart of the
two circles, decreasing the size of the shaded region that one needs to maximize
for effective forecasting.

If that reasoning is correct, longer data lengths should fill out the
attractor, thereby mitigating the spurious increase in the Shannon
entropy rate and allowing higher-dimensional reconstructions to
outperform lower-dimensional ones.  This is indeed what happens, as shown in
Figure~\ref{fig:mytaudata}.  For both the $K=22$ and $K=47$ traces,
once the signal is 2 million points long, the four-dimensional
estimator stores more information about the future than the two-dimensional case.
Note, though, that the optimal \mytau of
the $m=8$ reconstruction model is still lower than in the $m=2$ or $m=4$
cases, even at the right-hand limit of the plots in
Figure~\ref{fig:mytaudata}.  That is, even with a time series that
contains $4 \times 10^6$ points, it is more effective to use a lower-dimensional reconstruction to make an LMA forecast.  But the really
important message here is that \mytau allows one to determine the best
reconstruction parameters {\it for the available data}, which is an
important part of the answer to the challenges outlined at the
beginning of this chapter.

Something very interesting happens in the $m=2$ results for Lorenz 96
model with $K=47$: the \mytau curve reaches a maximum value around
100,000 points and stops increasing, regardless of data length.  What
this means is that this two-dimensional reconstruction contains as
much information about the future as can be ascertained from the \roLMA state estimator, suggesting that increasing the length of the training set would
not improve forecast accuracy.  To explore this, I construct LMA
forecasts of different-length traces (100,000--2.2 million points)
from this system, then reconstruct their dynamics with different $m$
values and the appropriate $\tau_{\mytau}$ for each case, and---again---repeat this full experiment 15 times for statistical validation.  For $m=2$,
both \mytau and 1-MASE results did indeed plateau at 200,000
points---at $5.736\pm0.0156$ and $0.0809\pm 0.0016$, respectively.  As before, more data does
afford higher-dimensional reconstructions more traction on the
prediction problem: the $m=4$ forecast accuracy surpassed $m=2$ at
around 2 million points. 
In neither case, by the way,
did $m=8$ catch up to either $m=2$ or $m=4$, even at 4 million data
points.  Of course, one must consider the cost of storing the
additional variables in a higher-dimensional model,
particularly in data sets this long, so it may be worthwhile in
practice to settle for the $m=2$ forecast---which is only slightly
less accurate and requires only 200,000 points.  This has another
major advantage as well.  If the time series is non-stationary, a
forecast strategy that requires fewer points is particularly useful because it can adapt more quickly.

\subsubsection{Choosing reconstruction parameters for increased prediction horizons.}
\label{subsec:predictionhorizon}

So far in this section, I have considered forecasts that were
constructed one step at a time and studied the correspondence of their
accuracy with one-step-ahead calculations of \mytau.  Here,
I consider longer prediction horizons ($h$) and explore whether one
can use a $h$-step-ahead version of \mytau---i.e.,
$I[\mathcal{S}_j,X_{j+h}]$, with $h>1$---to choose parameter values
that maximize the information that each delay vector contains
 about the
value of the time series $h$ steps in the future. 

\begin{figure}[tb!]
  \centering
            \vspace*{1mm}
  \includegraphics[width=\columnwidth]{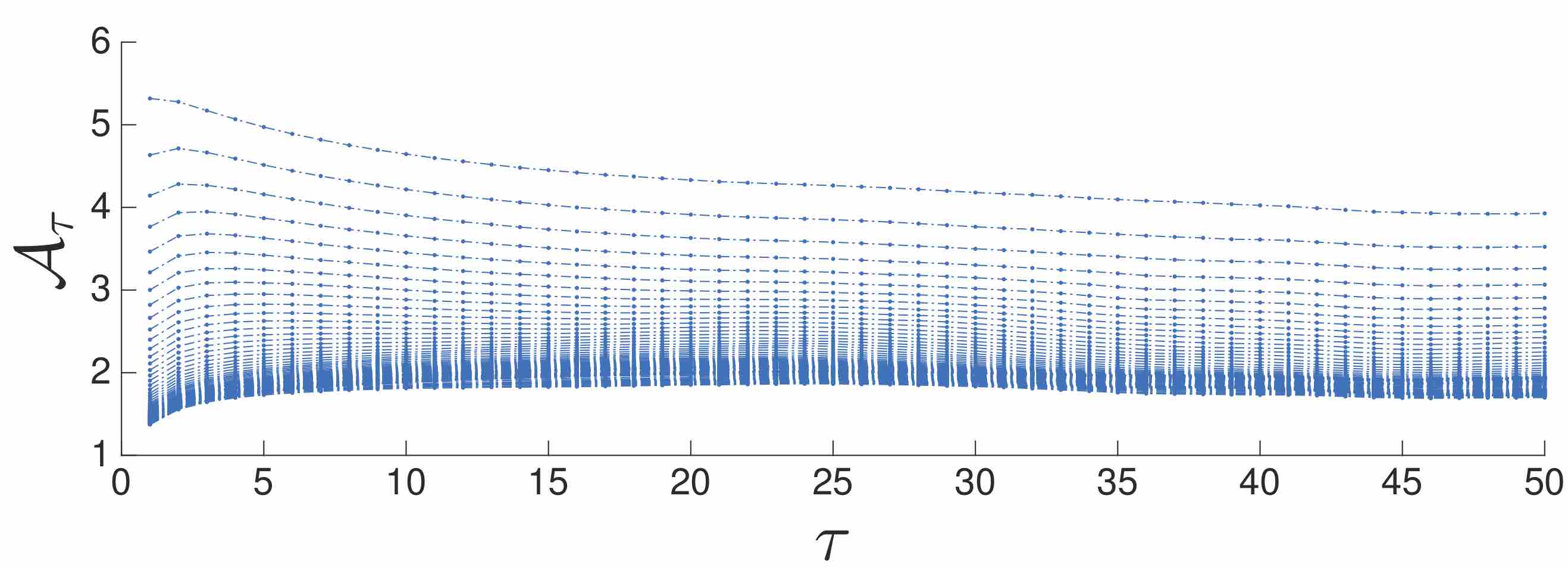}
  \caption{The effects of prediction horizon ($h$) on \mytau for a representative time series of the $K=22$ Lorenz 96 system for a fixed reconstruction dimension ($m=2$).  The curves in the image, from top to bottom, correspond to prediction horizons of $h=1$ to $h=100$.}
  \label{fig:L22pvsTau}
\end{figure}

Of course, one would expect the \mytau-optimal $\{m,\tau\}$ values for a given
time series to depend on the prediction horizon.  It has been shown,
for instance, that longer-term forecasts generally do better with
larger $\tau$~\cite{kantz97}, and conversely~\cite{joshua-pnp}.  It also
makes sense that one might need to reach different distances into the
past (via the span of the delay vector) in order to reduce the
uncertainty about events that are further into the
future~\cite{weigend-book}.  All of these effects are corroborated by \mytau.
Figure~\ref{fig:L22pvsTau} demonstrates this with a representative trace of the
$K=22$ Lorenz 96 system, focusing on $m=2$ for simplicity.
The topmost dashed curve in this figure is for the $h=1$ case---{\it i.e.,} a
horizontal slice of Figure~\ref{fig:tauandmL96}(a) made at $m=2$.  The
maximum of this curve is the optimal $\tau$ value ($\tau_{\mytau}$)
for this reconstruction.  The overall shape of this curve reflects the
monotonic increase in the uncertainty about the future with $\tau$
that is noted on page~\pageref{page:increase-with-tau}.  The other
curves in Figure~\ref{fig:L22pvsTau} show \mytau as a function of
$\tau$ for $h=2, 3, \dots$, down to $h=100$.  Note that the lower curves do not decrease monotonically; rather, there
is a slight initial rise.  This is due to the issue raised above about
the span of the delay vector: if one is predicting further into the
future, it may be useful to reach further into the past.  In general,
this causes the optimal $\tau$ to shift to the right as prediction
horizon increases, going down the plot---{\it i.e.,} longer prediction
horizons require larger $\tau$s ({\it cf.} \cite{kantz97}).  For very long
horizons, the choice of $\tau$ appears to matter very little.  In
particular, \mytau is fairly constant (and quite low) for $5<\tau<50$
when $h>30$---{\it i.e.,} regardless of the choice of $\tau$, there is very
little information about the $h$-distant future in any delay
reconstruction of this signal for $h>30$.  This effect should not be
surprising, and is well corroborated in the literature.  However,
it can be hard to know {\sl a priori}, when one is confronted with a
data set from an unknown system, what prediction horizon makes
sense.  \mytau offers a computationally efficient way to answer that
question from not very much data.

\begin{figure}[tb!]
  \centering
          \vspace*{1mm}
  \includegraphics[width=\textwidth]{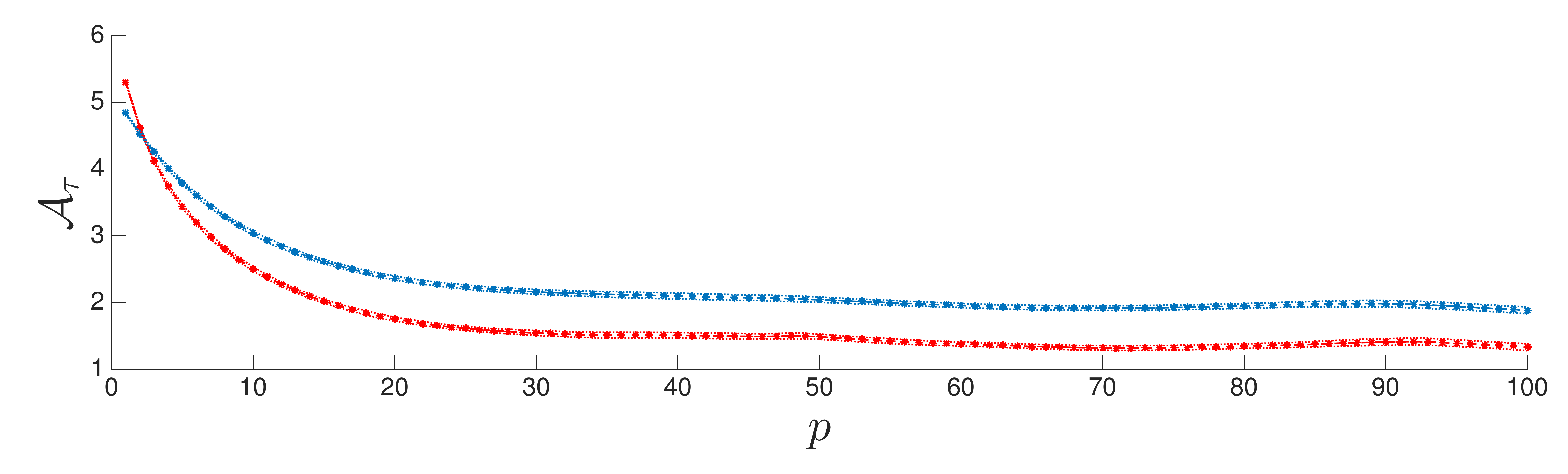}
  \caption{The effects of prediction horizon ($h$) on average \mytau (over 22 trials) of the $K=22$ Lorenz 96 system for a fixed time delay ($\tau=1$) and two different reconstructions of the system. The red line represents $m=2$; the blue represents $m_H=8$, the value suggested for this signal by the technique of false neighbors.}
  \label{fig:L22mytauhorizon}
\end{figure}

Figure~\ref{fig:L22mytauhorizon} shows a similar exploration but considers the effects of the reconstruction
dimension on \mytau as forecast horizon increases, this time fixing $\tau=1$ for simplicity.
These results indicate that the $m=2$ state estimator contains more information about the future
for short prediction horizons.  This ties back to a central theme of this thesis: low-dimensional
reconstructions can work quite well. Unsurprisingly, that does not always hold for arbitrary prediction horizons, Figure~\ref{fig:L22mytauhorizon} shows that the full
reconstruction is better for longer horizons.  This is to be expected, since a higher reconstruction dimension allows the state
estimator to capture more information about the past.  Finally, note
that \mytau decreases monotonically with prediction horizon for both
$m=2$ and $m_H$.  This, too, is unsurprising.  Pesin's
relation~\cite{pesin1977characteristic} says that the sum of the
positive Lyapunov exponents is equal to the entropy rate, and if there
is a non-zero entropy rate, then generically observations will become
increasingly independent the further apart they are.  This explanation
also applies to Figure~\ref{fig:L22pvsTau}, of course, but it does
{\it not} hold for signals that are wholly (or nearly) periodic.

\begin{figure}[tb!]
  \centering
          \vspace*{1mm}
  \includegraphics[width=\textwidth]{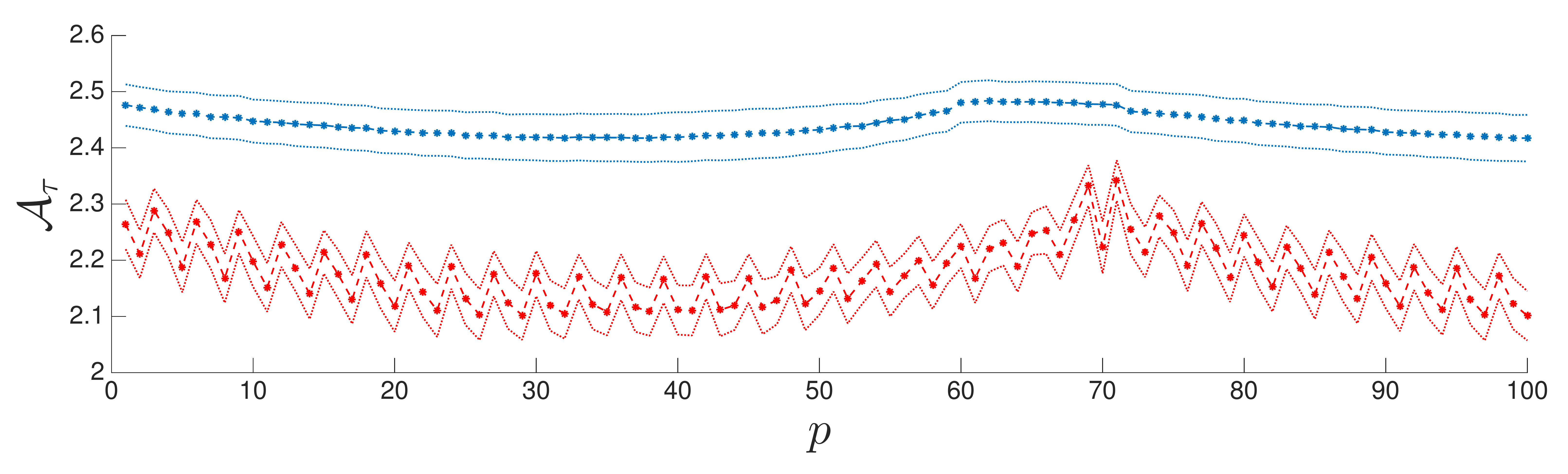}
  \caption{The effects of prediction horizon ($h$) on average \mytau over the 15 trials of \col for a fixed time delay ($\tau=1$) and two different reconstruction dimensions.  The red line represents $m=2$; the blue represents $m_H=12$, the value suggested for this signal by the technique of false neighbors.}
  \label{fig:colmytauhorizon}
\end{figure}

Recall that the \col dynamics in Section~\ref{subsubsec:computer} are
chaotic, but with a dominant unstable periodic orbit---which have a
variety of interesting effects on the results.
Figure~\ref{fig:colmytauhorizon} explores the effects of prediction
horizon on those results.
Not surprisingly, there is some periodicity in the \mytau versus $h$
relationships, but not for the same reasons that caused the stripes in
Figure~\ref{fig:tauandmcol}(b).  Here, the {\it peaks} in \mytau do occur at
multiples of the period.  That is, the $m=2$ state estimator can
forecast with the most success when the value being predicted is in
phase with the delay vector.  Note that this effect is far stronger
for $m=2$ than $m_H$, simply because of the instability of that
periodic orbit; the visits made to it by the chaotic trajectory are more likely to be short than long.  As expected, \mytau
decays with prediction horizon---but only at first, after which it
begins to rise again, peaking at $h=69$ and $h=71$.  This may be due
to a second higher-order unstable periodic orbit in the \col dynamics.

In theory, one can derive rigorous bounds on prediction horizon.  The
time at which $\mathcal{S}_j$ will no longer have any information
about the future can be determined by considering:
\begin{equation}
  R(h) = \frac{I[\mathcal{S}_j,X_{j+h}]}{H[X_{j+h}]}
\end{equation}
{\it i.e.,} the percentage of the uncertainty in $X_{j+h}$ that can be
reduced by the delay vector.  Generically, this will limit to some
small value equal to the amount of information that the delay vector
contains about any arbitrary point on the attractor.  Given some
criteria regarding how much information above the ``background'' is
required of the state estimator, one could use an $R(h)$ versus $h$ curve
to determine the maximum practical horizon.

In practice, one can select parameters for delay reconstruction-based
forecasting by explicitly including the prediction horizon in the
\mytau function, fixing the horizon $h$ at the required value, performing
the same search as I did in earlier sections over a range of $m$ and
$\tau$, and then choosing a point on (or near) the optimum of that
\mytau surface.  The computational and data requirements of this
calculation, as shown in Section~\ref{subsec:datalength}, are far
superior to those of the standard heuristics used in delay
reconstructions.

\subsection{Summary}
\label{sec:conclusion}

\mytau is a novel metric for quantifying how much
information about the future is contained in a delay reconstruction.
Using a number of different dynamical systems, I demonstrated a
direct correspondence between the \mytau value for different delay
reconstructions and the accuracy of forecasts made with Lorenz's
method of analogues on those reconstructions.  Since \mytau can be
calculated quickly and reliably from a relatively small amount of
data, without any knowledge about the governing equations
or the state space dynamics of the system, that correspondence is a
major advantage, in that it allows one to choose parameter values for
delay reconstruction-based forecast models without doing an exhaustive
search on the parameter space.  Significantly, \mytau-optimal
reconstructions are better, for the purposes of forecasting, than
reconstructions built using standard heuristics like mutual
information and the method of false neighbors, which can require large
amounts of data, significant computational effort, and expert human
interpretation.  Perhaps, most importantly \mytau allows one to answer other questions regarding
forecasting with theoretically unsound models---{\it e.g.,}
why it is possible to obtain a better forecast using a low-dimensional
reconstruction than with a true embedding.  It also allows one to
understand bounds on prediction horizon without having to estimate
Lyapunov spectra or Shannon entropy rates, both of which are difficult to
obtain for arbitrary real-valued time series.  That, in turn, allows
one to tailor one's reconstruction parameters to the amount of
available data and the desired prediction horizon---and to know if a
given prediction task is just not possible.

The experiments reported in this section involved a simple near-neighbor forecast
strategy and state estimators that are basic delay reconstructions of
raw time-series data.  The definition and calculation of \mytau do not
involve any assumptions about the state estimator, though, so the
results presented here should also hold for other state estimators.
For example, it is common in forecasting applications to pre-process
the time series: for example, low-pass filtering or interpolating to
produce additional points.  Calculating \mytau after performing such
an operation will accurately reflect the amount of information in that
new time series---indeed, it would reveal if that pre-processing step
{\it destroyed} information.  And I believe that the basic
conclusions in this section extend to other state-space based forecast
schemas besides LMA, such as those used in \cite{weigend-book,casdagli-eubank92,Smith199250,sugihara90,sauer-delay}---although
\mytau may not accurately select optimal parameter values for
strategies that involve {\it post}-processing the data (e.g.,
GHKSS~\cite{ghkss}).  

There are many other interesting potential ways to leverage \mytau in the practice of forecasting. If the \mytau-optimal $\tau =1$, that may be a signal that the time series is undersampling the dynamics and that one should increase the sample rate. One could use the more general form \smytau at a finer grain to optimizing $\tau$ individually for each dimension, as suggested in~\cite{pecoraUnified,Small2004283,PhysRevE.87.022905}, where optimal values are selected based on criteria that are not directly related to prediction. To do this, one could define $\mathcal{S}_j= [X_{j}, X_{j-\tau_1}, X_{j-\tau_2}, \dots, X_{j-\tau_{m-1}}]$ and then simply maximize \smytau using that state estimator constrained over $\{\tau_i\}_{i=1}^{m-1}$.

%% file: compTopology.tex


\newtheorem{lem}[thm]{Lemma}
\newtheorem{cor}[thm]{Corollary}
\newtheorem{con}{Conjecture}

\newcommand{\bC}{{\mathbb{ C}}}
\newcommand{\bN}{{\mathbb{ N}}}
\newcommand{\bQ}{{\mathbb{ Q}}}
\newcommand{\bR}{{\mathbb{ R}}}
\newcommand{\bS}{{\mathbb{ S}}}
\newcommand{\bT}{{\mathbb{ T}}}
\newcommand{\bZ}{{\mathbb{ Z}}}

\newcommand{\cC}{{\cal C}}
\newcommand{\cF}{{\cal F}}
\newcommand{\cK}{{\cal K}}
\newcommand{\cL}{{\cal L}}
\newcommand{\cN}{{\cal N}}
\newcommand{\cO}{{\cal O}}
\newcommand{\cR}{{\cal R}}
\newcommand{\cS}{{\cal S}}
\newcommand{\cU}{{\cal U}}
\newcommand{\cW}{{\cal W}}
\newtheorem*{defn}{Definition}

\newcommand{\eps}{\varepsilon}
\newcommand{\diam}[1]{\mathrm{diam}(#1)}

\newcommand{\Eq}[1]{Equation~(\ref{eqn:#1})}
\newcommand{\Fig}[1]{Figure~\ref{fig:#1}}
\newcommand{\Sec}[1]{\S\ref{sec:#1}}
\newcommand{\Def}[1]{Definition~\ref{def:#1}}

\newcommand{\beq}[1]{\begin{equation}\label{eqn:#1}}
\newcommand{\eeq}{\end{equation}}

\section{Exploring the Topology of Dynamical Reconstructions}\label{sec:compTopo}

Topology is of particular interest in forecasting dynamics, since many
properties---the existence of periodic orbits, recurrence, entropy, etc.---depend only upon topology.  However, computing topology from time
series can be a real challenge---even with the aid of delay-coordinate reconstruction.   
As I have mentioned repeatedly throughout this thesis, success of this reconstruction
procedure depends heavily on the choice of the two free parameters, but the
embedding theorems provide little guidance regarding how to choose good values for 
these parameters.  
 The delay-coordinate reconstruction machinery
(both theorems and heuristics) targets the computation of dynamical
invariants like the correlation dimension and the Lyapunov exponent.
However, if one just wants to extract the topological structure of an invariant set---as is the case with forecasting---a scaled-back version of that machinery
may be sufficient. In the following discussion, I adopt the philosophy that one might only desire
knowledge of the topology of the invariant set. In collaboration with J.~D.~Meiss, I conjecture that
this might be possible with a lower reconstruction dimension than that
needed to obtain a true embedding.  That is, the reconstructed
dynamics might be {\it homeomorphic} to the original dynamics at a
lower dimension than that needed for a diffeomorphically correct
embedding~\cite{josh-physicaD}. This is an alternative validation of the central premise of my thesis.

To compute topology from data, one can use a
simplicial complex--{\it e.g.,} the {\it witness complex} of \cite{deSilva04}.  To
construct such a complex, one chooses a set of ``landmarks,''
typically a subset of the data, that become the vertices of the
complex.  The connections between the landmarks are determined by
their nearness to the rest of the data---the ``witnesses.''  Two
landmarks in the complex are joined by an edge, for instance, if they
share at least one witness. 

My initial work on this approach~\cite{josh-physicaD} suggests that
 the witness complex correctly resolves the homology of the underlying 
invariant set---{\it viz.,} its Betti numbers---even if the reconstruction dimension is well below the thresholds for
which the embedding theorems assure smooth conjugacy between the true
and reconstructed dynamics. This means that some features of the large-scale topology are present even if the reconstruction dimension 
does not satisfy the associated theorems.  I conjecture that this structure affords 
 \roLMA the means necessary to generate accurate forecasts.  
%
The witness complex is covered in more depth in Section~\ref{sec:wc}, which also
describes the notion of {\it persistence} and demonstrates how that idea is
used to choose scale parameters for a complex built from reconstructed
time-series data. In Section~\ref{sec:meat}, I explore how the homology of such a
complex changes with reconstruction dimension. 

\subsection{Witness Complexes for Dynamical Systems}
\label{sec:wc}

To compute the topology of data that sample an invariant set of
a dynamical system, one needs a complex that captures the shape of the
data but is robust with respect to noise and other sampling issues.
A witness complex is an ideal choice for these purposes.  Such a complex
is determined by the reconstructed time-series data, $W \subset
\bR^m$---the {\it witnesses}---and an associated set $L \subset
\bR^m$, the {\it landmarks}, which can (but need not) be chosen from
among the witnesses.
%
%
The landmarks form the vertex set of the complex; the connections
between them are dictated by the geometric relationships between $W$
and $L$.  In a general sense, a witness complex can be defined through
a relation $R(W,L) \subset W \times L$. As Dowker noted~\cite{Dowker52}, any relation gives rise to a pair of simplicial
complexes.  In the one used here, a point $w \in W$ is a witness to an
abstract $k$-dimensional simplex $\sigma = \langle l_{i_1}, l_{i_2}, \ldots
l_{i_{k+1}}\rangle \subset L$ whenever $\{w\} \times \sigma \subset
R(W,L)$. The collection of simplices that have witnesses is a complex
relative to the relation $R$. For example, two landmarks are connected
if they have a common witness---this is a one-simplex.  Similarly, if
three landmarks have a common witness, they form a two-simplex, and so
on.

\begin{figure}[tb!]
        \centering
                \includegraphics[width=0.4\textwidth]{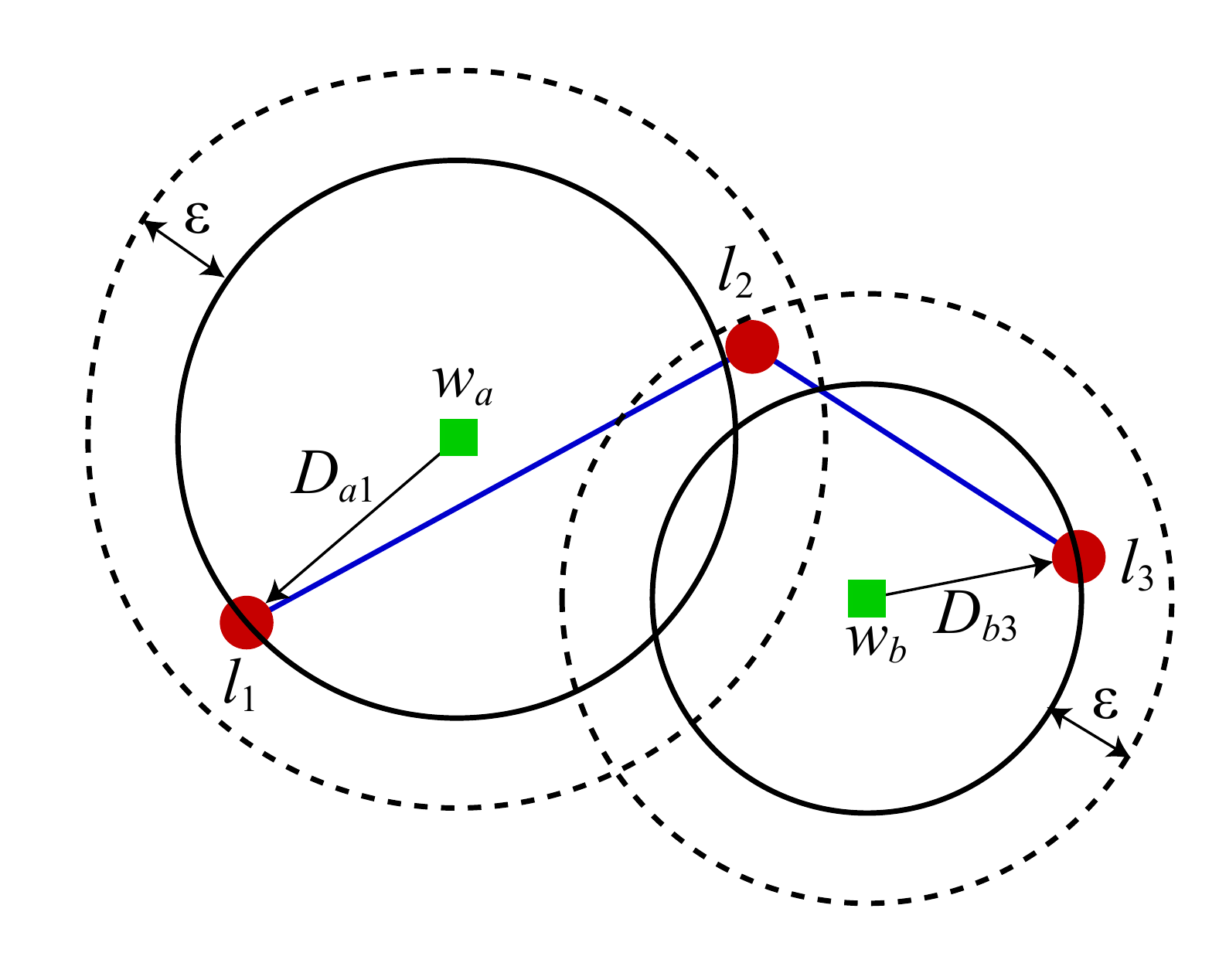}
 \caption{Illustration of the fuzzy witness relation \Eq{witnessSet}. The closest landmark to witness $w_a$ is $l_1$, and since $\|w_a-l_2\| < D_{a1} + \eps$, the simplex $\langle l_1,l_2 \rangle$ is in the complex. Similarly $w_b$ witnesses the edge $\langle l_2, l_3 \rangle$.}
 \label{fig:WitnessRelation}
\end{figure}

There are many possible definitions for a witness relation $R$; see~\cite{josh-physicaD} for a discussion.  
A relation that is particularly useful for analyzing noisy real data~\cite{josh-physicaD,Alexander15} is the $\eps$-weak witness~\cite{Carlsson14},
or what is called a ``fuzzy" witness \cite{Alexander15}: a point witnesses a simplex 
if all the landmarks in that simplex are within $\eps$ of the closest landmark to
the witness:
\begin{defn}[Fuzzy Witness]
The fuzzy witness set for a point $l \in L$ is the set
of witnesses
\beq{witnessSet}
	\cW_\eps(l) = \{w \in W : \|w-l\| \le \min_{l' \in L} \|w-l'\| + \eps \}
\eeq
\end{defn}
\noindent In this case, the relation consists of the collections
$R = \cup_{l \in L} (\cW_\eps(l) \times \{ l \})$ and a simplex $\sigma$ is
in the complex whenever $\cap_{l \in \sigma} \cW_\eps(l) \neq \emptyset$---that
is, when all of its vertices share a witness. This relation is illustrated in Figure~\ref{fig:WitnessRelation}.

The fuzzy witness complex reduces to the ``strong witness complex'' of de
Silva and Carlsson \cite{deSilva04} when $\eps =0$.  In such a
complex, an edge exists between two landmarks {\it iff} there exists
a witness that is exactly equidistant from those landmarks.  This is
not a practical notion of shared closeness for finite noisy data sets. In this case, $\eps$ in \Eq{witnessSet} allows for some amount of immunity to finite data and noise effects, but must be chosen correctly as I discuss in the next paragraph.  A simpler implementation of the fuzzy witness complex consists of simplices whose pairs of vertices have a common
witness; this implementation gives a ``clique"
or ``flag" complex, analogous to the Rips complex \cite{Ghrist08}.  This is called a ``lazy" complex in \cite{deSilva04} and
instantiated as the \begin{tt}LazyWitnessStream\end{tt} class in the
  {\tt javaPlex} \cite{tausz2012javaplex} software.  In the notation introduced above,
  the complex is \vspace{-0.75cm}
  \beq{fuzzyComplex} \cK_\eps(W,L) = \{ \sigma \subset
  L: \cW_\eps(l) \cap \cW_\eps(l') \neq \emptyset, \;\; \forall l,l'
  \in \sigma\}   \eeq
   Following \cite{Alexander15}, I will use this
  particular construction in the following discussion. My goal is to study the topology of witness complexes of delay-coordinate reconstructions and determine whether the topology is resolved correctly when the reconstruction dimension is low.
  
Figure~\ref{fig:witness-complexes-different-eps} shows four witness
complexes built from the 100,000-point trajectory of
the Lorenz 63 system that is shown in Figure~\ref{fig:Lorenz63Example}(a) for varying values of the fuzziness parameter, $\epsilon$.  The
landmarks (red dots) consist of $\ell = 201$ points equally spaced
along the trajectory, {\it i.e.,} every $\Delta t = 500^{th}$ point of the
time series. There are many ways to choose\footnote{For a deeper discussion of the number of landmarks to use and landmark selection choice see~\cite{josh-physicaD}.}landmarks; this particular strategy distributes them according to the invariant measure of the attractor.  
{One could also choose landmarks randomly from the trajectory or using 
the ``max-min" selector\footnote{Choose the first landmark at random,  
and given a set of landmarks, choose the next to be the data point farthest away from the current set.}of \cite{deSilva04}; each of these gives results similar to those shown.
 When $\eps$ is small, very few
witnesses fall in the thin regions required by \Eq{witnessSet}, so the
resulting complex does not have many edges and is thus not a good
representation of the shape of the data. As $\eps$ grows, more
witnesses fall in the ``shared'' regions and the complex fills in,
revealing the basic homology of the attractor of which the trajectory
is a sample.  There is an obvious limit to this, however: when $\eps$
is very large, even the largest holes in the complex are obscured.

\begin{figure}[tb!]
        \begin{subfigure}[b]{0.45\textwidth}
                \includegraphics[width=\textwidth]{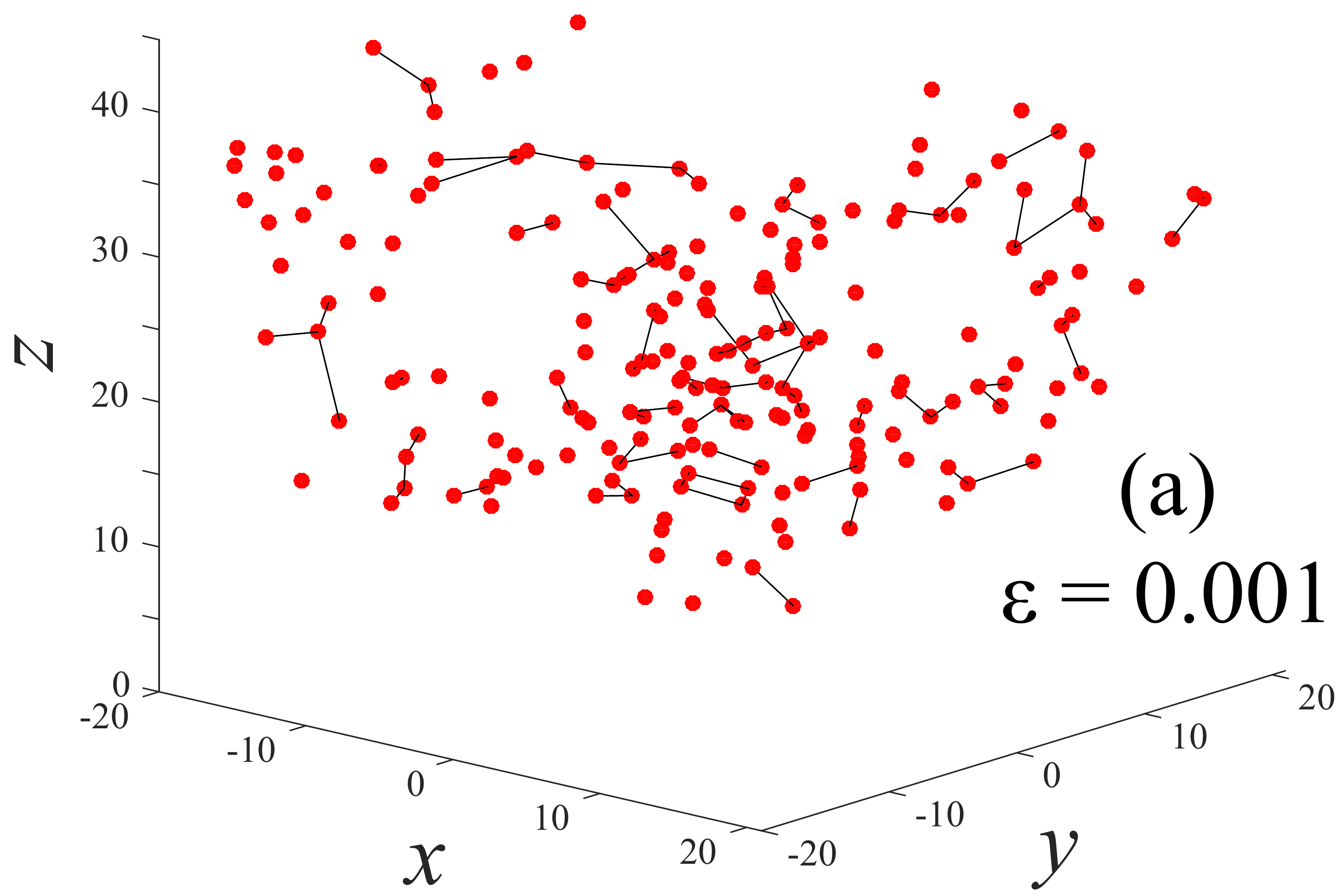}
        \end{subfigure}
        \begin{subfigure}[b]{0.45\textwidth}
                \includegraphics[width=\textwidth]{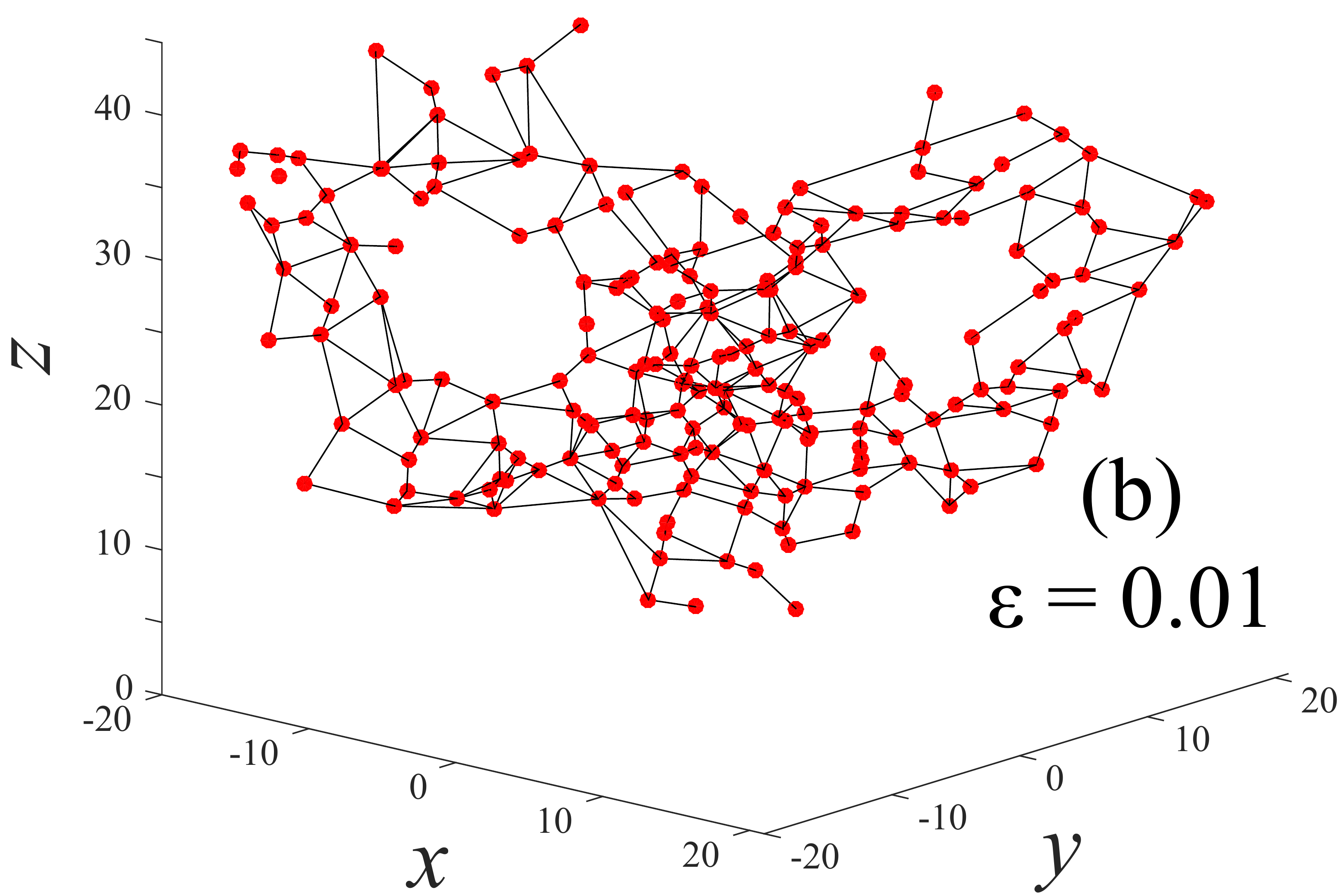}
        \end{subfigure}
        \begin{subfigure}[b]{0.45\textwidth}
                \includegraphics[width=\textwidth]{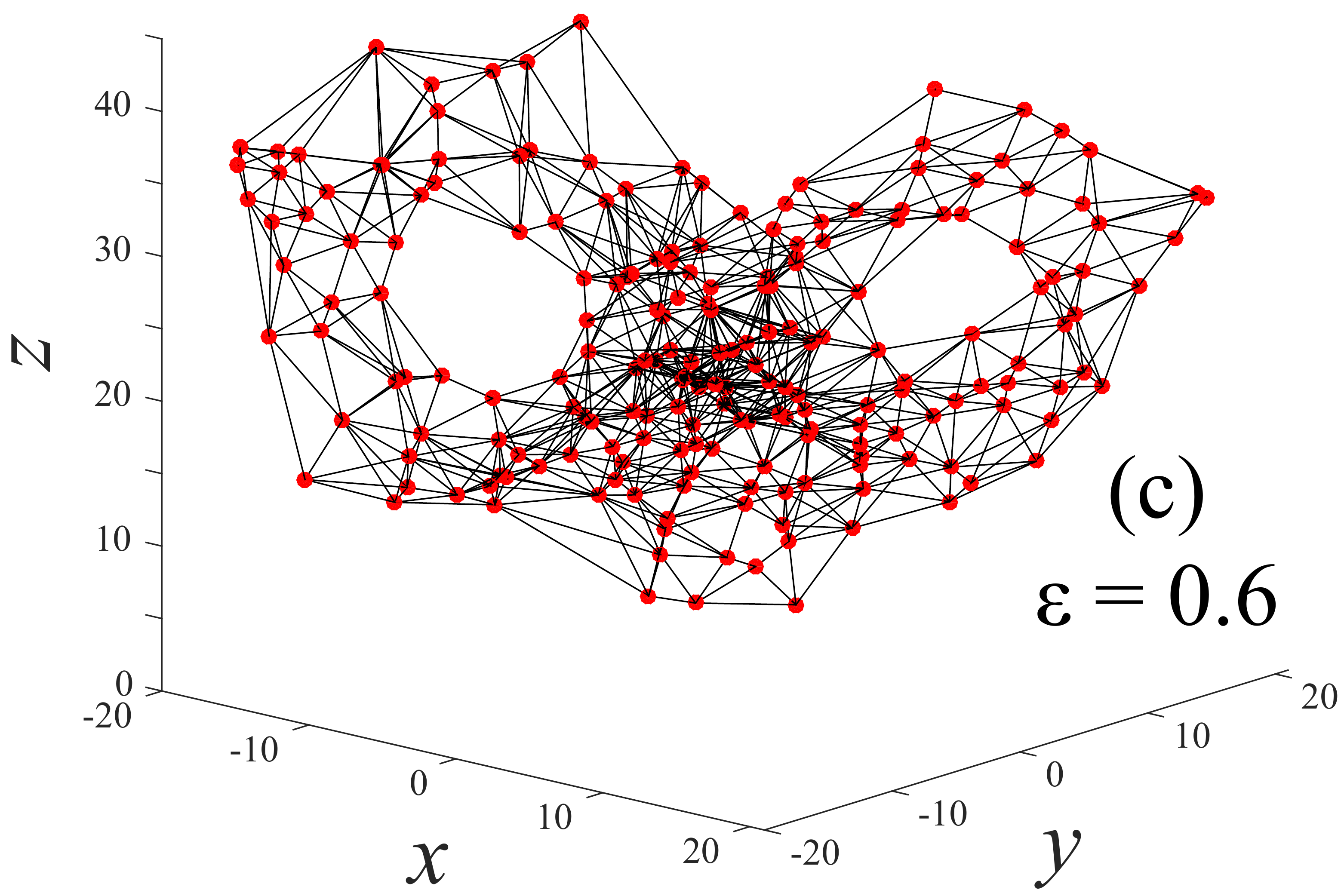}
        \end{subfigure} \quad \quad \quad \quad
        \begin{subfigure}[b]{0.45\textwidth}
                \includegraphics[width=\textwidth]{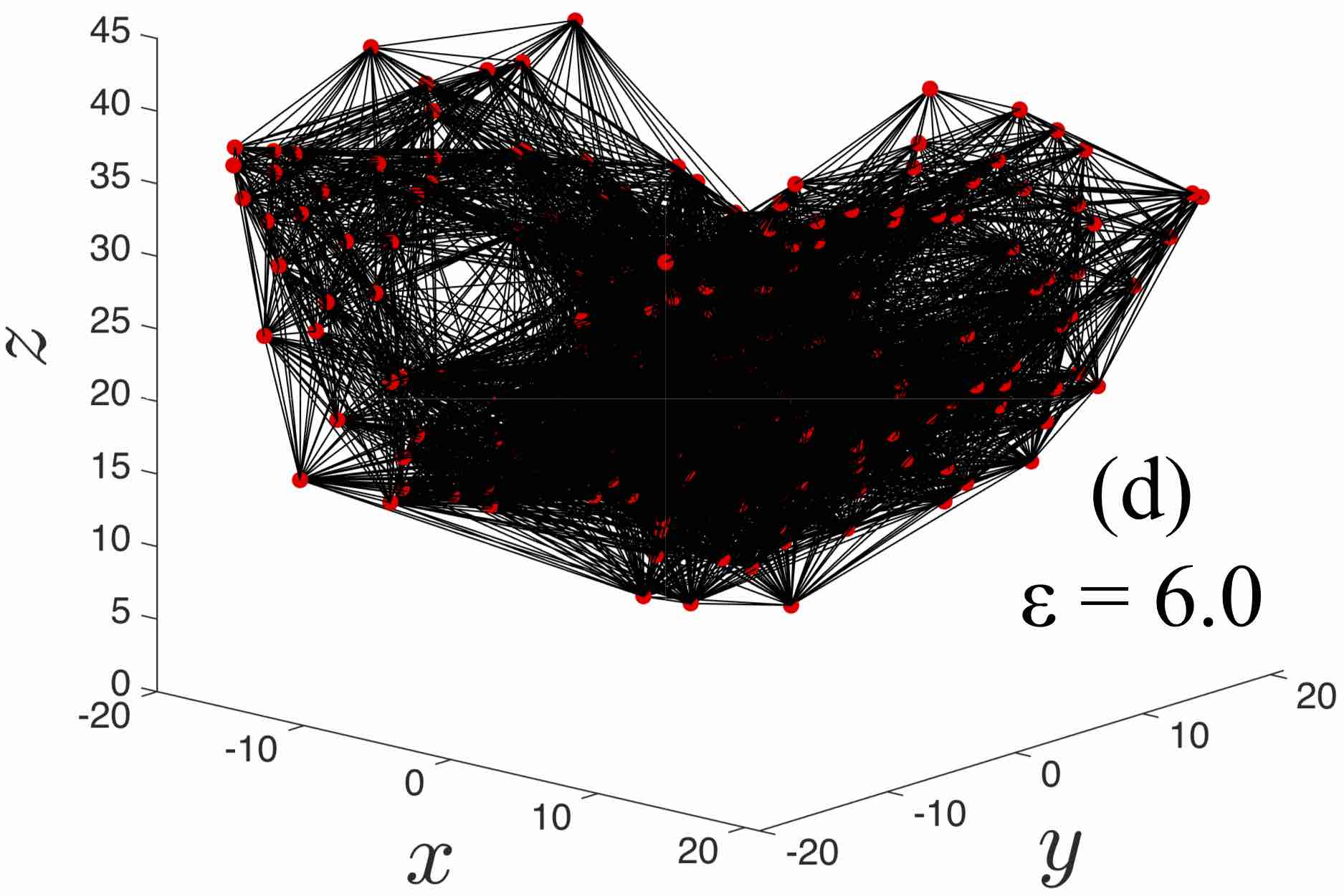}
        \end{subfigure}
\caption{
 %
 ``One-skeletons" of witness complexes constructed from the trajectory of Figure~\ref{fig:Lorenz63Example}(a) using the fuzzy witness relation depicted in Figure~\ref{fig:WitnessRelation}. In each one-skeleton, the red dots are the $\ell = 201$ equally-spaced landmarks. A black edge between two landmarks $l_i$ and $l_j$ signifies the existence of a one simplex $\langle l_i,l_j \rangle$ in the complex, {\it i.e.,}  $l_i$ and $l_j$ shared at least one witness. As $\epsilon$ increases, more landmarks will satisfy $\|w_a-l_k\| < D_{a1} + \eps$ for each $w_a$ and the one complex will fill in. 
}
\label{fig:witness-complexes-different-eps}
\end{figure}

In order to evaluate the topology of incomplete reconstructions, one needs to ensure the correctness of the topology. However, as Figure~\ref{fig:witness-complexes-different-eps} illustrates, the simplex, from which one estimates the topology, depends on the choice of $\epsilon$, and choosing the right $\epsilon$ for that job is non-trivial.  One can do so using the progression of images in Figure~\ref{fig:witness-complexes-different-eps} and the notion of {\it persistence}. Studying the change in homology under changing scale parameters is a
well-established notion in computational topology.  The underlying
idea of { persistence} \cite{ELZ01,Robins02,Zomorodian05} is that
any topological property of physical interest should be (relatively)
independent of parameter choices in the associated algorithms.  

One useful way to represent information about the changing topology of a
complex is the {\it barcode persistence diagram} \cite{Ghrist08}.
Figure~\ref{fig:Lorenz63-barcodes} shows barcodes of the first two
Betti numbers for the witness complexes of
Figure~\ref{fig:witness-complexes-different-eps}.  Each horizontal line in the
barcode is the interval in $\eps$ for which there exists a particular non-bounding
cycle, thus the number of such lines is the rank of the homology
group---a Betti number.  The values for $\beta_0$ and
$\beta_1$ are computed using {\tt javaPlex} \cite{tausz2012javaplex} over the range $0.017\le \eps \le
1.7$, using the \begin{tt}ExplicitMetricSpace\end{tt} to choose the equally
spaced points and the \begin{tt}LazyWitnessStream\end{tt} to obtain a clique
complex from the $\ell=201$ landmarks.  There are no three-dimensional voids
in the results, {\it i.e.,} $\beta_2$ was 
always zero for this range of $\eps$---a reasonable implication for 
this $2.06$-dimensional attractor.  When $\eps$ is very small,
as in Figure~\ref{fig:witness-complexes-different-eps}(a), the witness complex
has many components and the $\beta_0$ barcode shows a large number of
entries. As $\eps$ grows, the spurious gaps between these components
disappear, leaving a single component that persists above $\eps
\approx 0.014$.  That is, witness complexes constructed with $\eps >
0.014$ correctly capture the connectedness of the underlying
attractor.  The $\beta_1$ barcode plot shows a similar pattern: there
are many holes for small $\eps$ that are successively filled in as
that parameter grows, leaving the two main holes ({\it i.e.,} $\beta_1=2$)
for $\eps > 1.01$.  Above $\eps > 3.2$ (not shown in
Figure~\ref{fig:Lorenz63-barcodes}), one of those holes disappears; eventually,
for $\eps> 4.05$, the complex becomes topologically trivial.  Above
this value, the resulting complexes---recall
Figure~\ref{fig:witness-complexes-different-eps}(d)---have no non-contractible
loops and are homologous to a point (acyclic).

\begin{figure}[tb]
\centering
        \begin{subfigure}[b]{0.45\textwidth}
                \includegraphics[width=\textwidth]{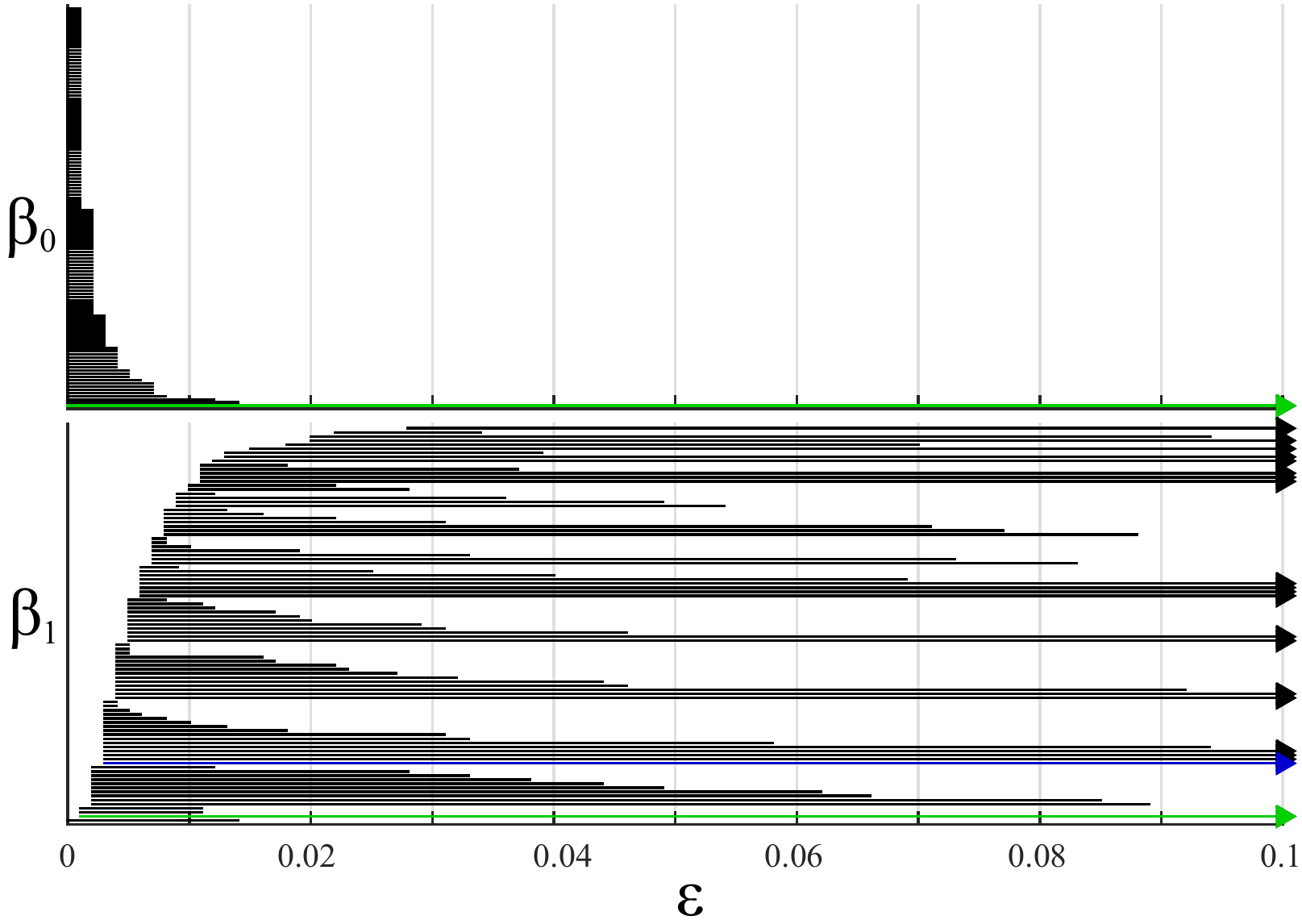}
        \end{subfigure}
        \begin{subfigure}[b]{0.45\textwidth}
			\includegraphics[width=\textwidth]{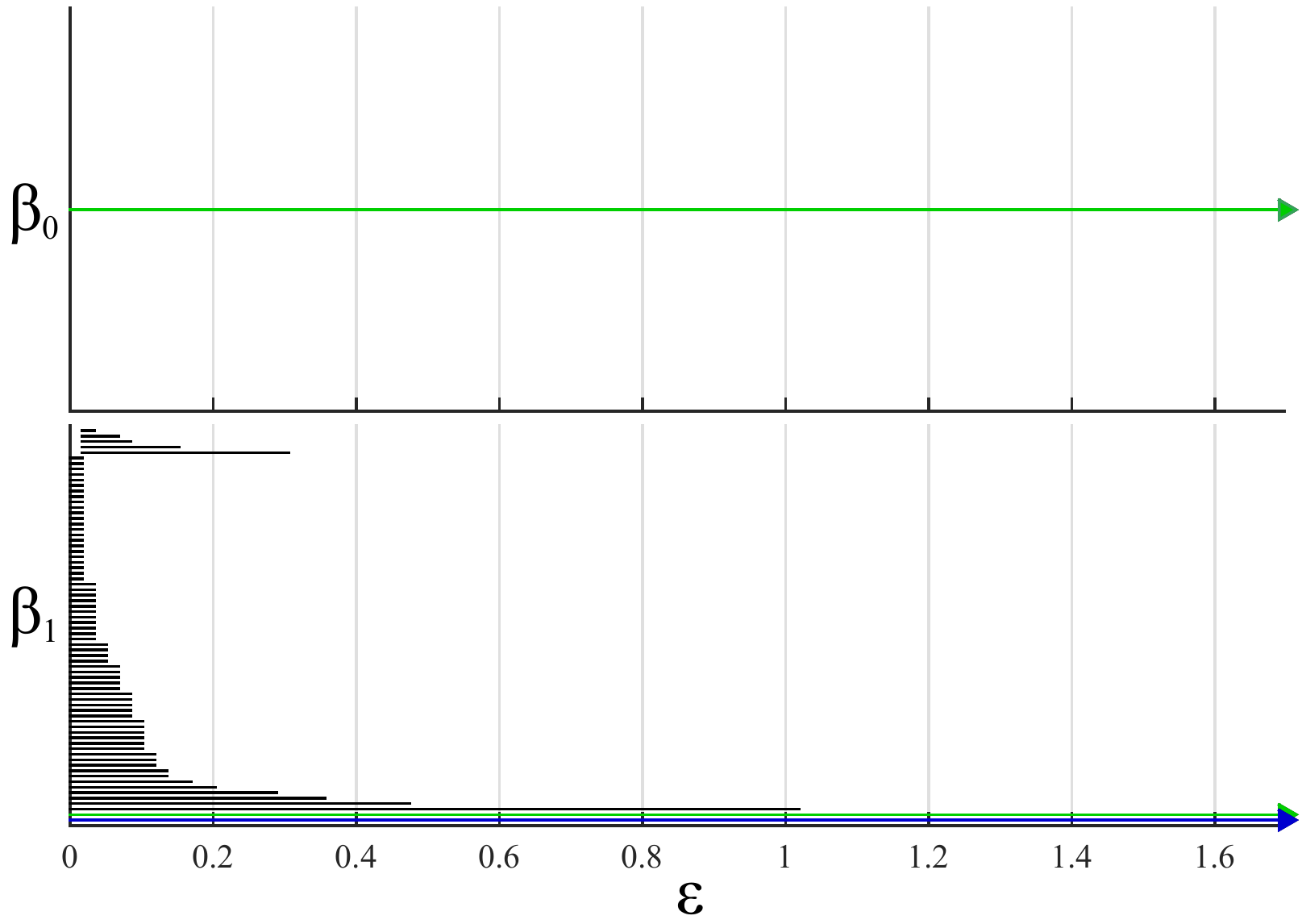}
        \end{subfigure}		
\caption{Persistence barcodes computed using {\tt javaPlex} for a
  $\ell=201$ witness complex of the trajectory of
  Figure~\ref{fig:Lorenz63Example}(a).  Each plot tabulates the two lowest Betti
  numbers of the complex for $100$ values of the scale parameter
  $\eps$.  The left panel shows the behavior when $0.001\le \eps\le
  0.1$, the right panel zooms out to the range $0.017 \le \eps \le 1.7$.}
 
\label{fig:Lorenz63-barcodes}
\end{figure}

As alluded to above, this notion of persistence can be turned around and used to select
good values for the parameters that play a role in topological data
analysis---{\it e.g.,} looking for the $\eps$ value at which the homology
stabilizes or selecting the number of landmarks that are necessary to construct a topologically faithful complex.  
However, definitions of what constitutes stabilization
are subjective and can be problematic.  
Even so, persistence is a powerful technique and I make
use of it in a number of ways in the rest of this section.

The examples in Figures~\ref{fig:witness-complexes-different-eps} and~\ref{fig:Lorenz63-barcodes} involve a full trajectory from a
dynamical system. This thesis focuses on reconstructions of scalar time-series data---structures
whose topology is guaranteed to be identical to that of the underlying
dynamics if the reconstruction process is carried out properly.  But
what if the dimension $m$ does {\it not} satisfy the requirements off
the theorems?  Can one still obtain useful results about the {\it topology} of
that underlying system, even if those dynamics are not
properly unfolded in the sense of \cite{packard80,takens,sauer91}? Throughout this thesis, I have argued that in the context of forecasting, the answer to that question is yes. In the next section, I take a step away from forecasting and examine 
whether the answer is also yes in the case of {\it topology}---specifically homology---and discuss the implications of that answer for the central theme of this thesis.

\subsection{Topologies of Reconstructions}
\label{sec:meat}
As discussed in Section~\ref{sec:dce}, a scalar time series of a dynamical system is a projection of the
$d$-dimensional dynamics onto $\bR^1$---an action that does not
automatically preserve the topology of the object. Delay-coordinate embedding allows one to reconstruct the underlying dynamics, up to diffeomorphism, if the reconstruction dimension is large enough.  
The question at issue in this section is whether one can use the
witness complex to obtain a useful, coarse-grained description of the
topology from lower-dimensional reconstructions, namely the homology---{e.g.,} the basic
connectivity of the invariant set, or the number of holes in it that
are larger than a certain scale. The answer to this question can provide a deeper understanding of 
the mechanics of \roLMA.

The short answer is yes.
Figure~\ref{fig:homology-figure1} shows a side-by-side comparison of
witness complexes and barcode diagrams for the Lorenz 63 trajectory of Figure~\ref{fig:Lorenz63Example}(a) and a two-dimensional reconstruction ($m=2$)
using the $x$ coordinate of that trajectory.  The full $3D$
trajectory on the left and the
$2D$ reconstruction on the right have the same
homology. {\it In other words, the correct large-scale homology is accessible from a witness 
complex of a $2D$ reconstruction, even though the reconstruction does not satisfy the conditions of the
associated theorems.}

\begin{figure}[tb!]
 \centering
 
         \begin{subfigure}[b]{0.45\textwidth}
                \includegraphics[width=\textwidth]{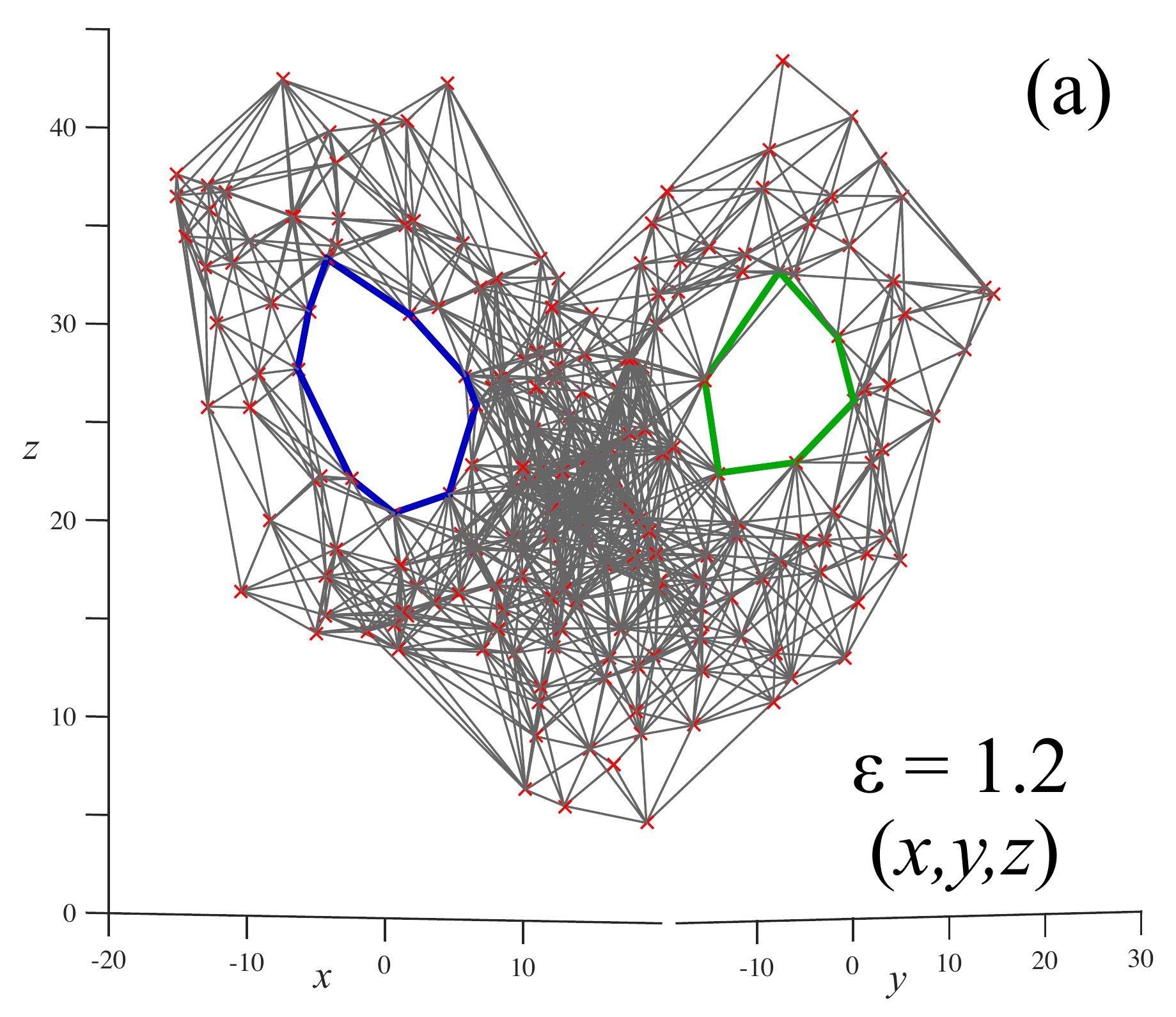}
        \end{subfigure}%
                 \begin{subfigure}[b]{0.45\textwidth}
                \includegraphics[width=\textwidth]{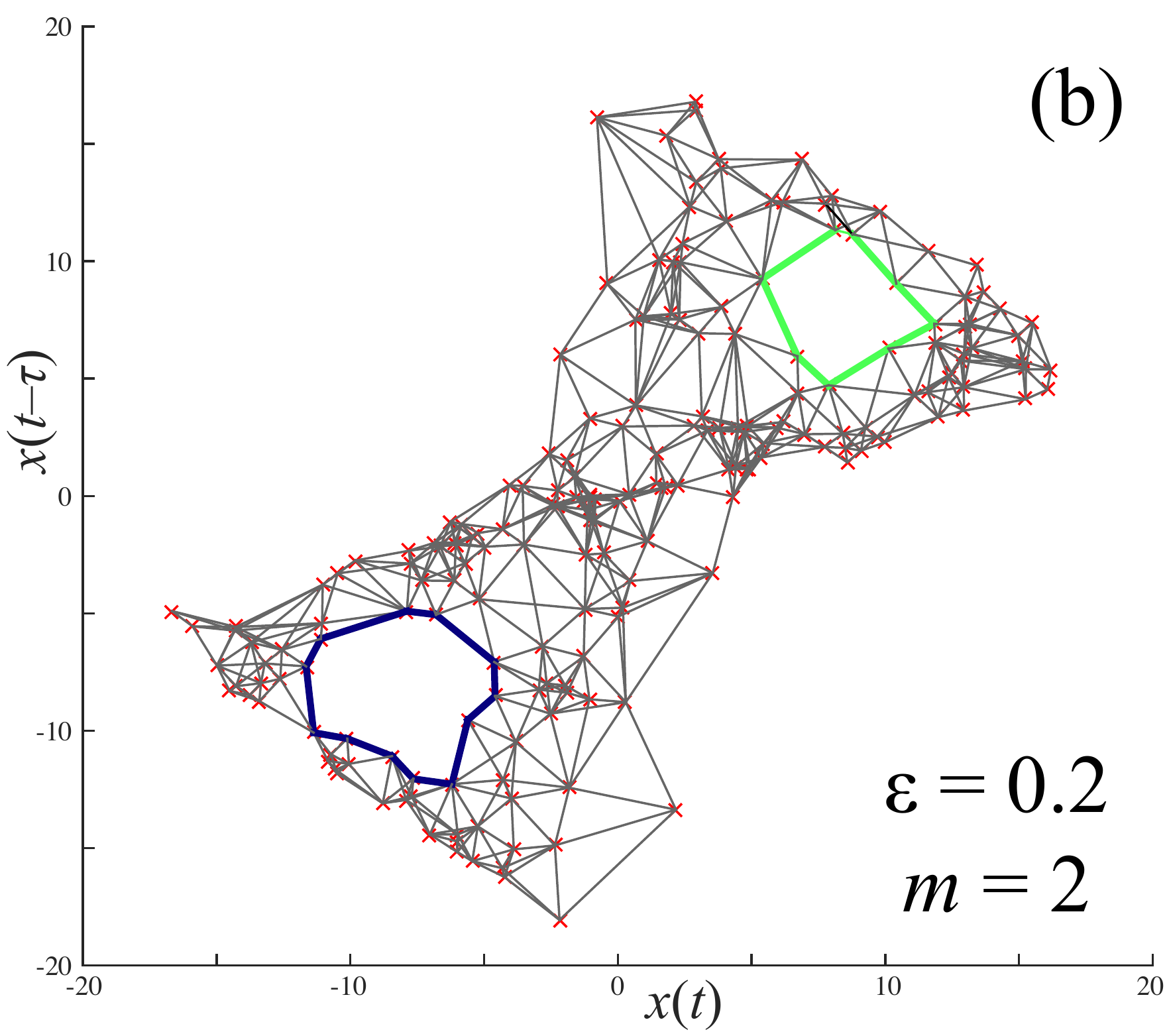}
        \end{subfigure}%
        \vspace*{3mm}
 
          \begin{subfigure}[b]{0.45\textwidth}
                \includegraphics[width=\textwidth]{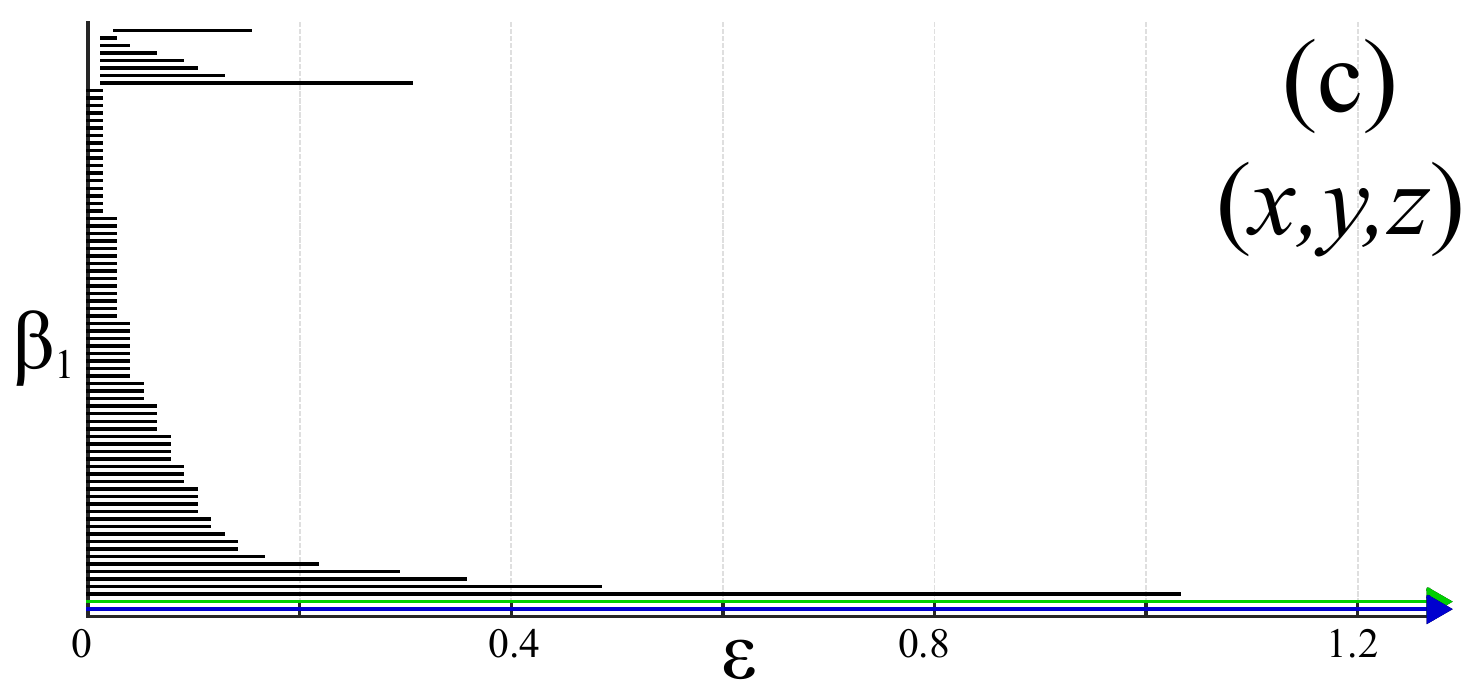}
        \end{subfigure}%
                 \begin{subfigure}[b]{0.45\textwidth}
                \includegraphics[width=\textwidth]{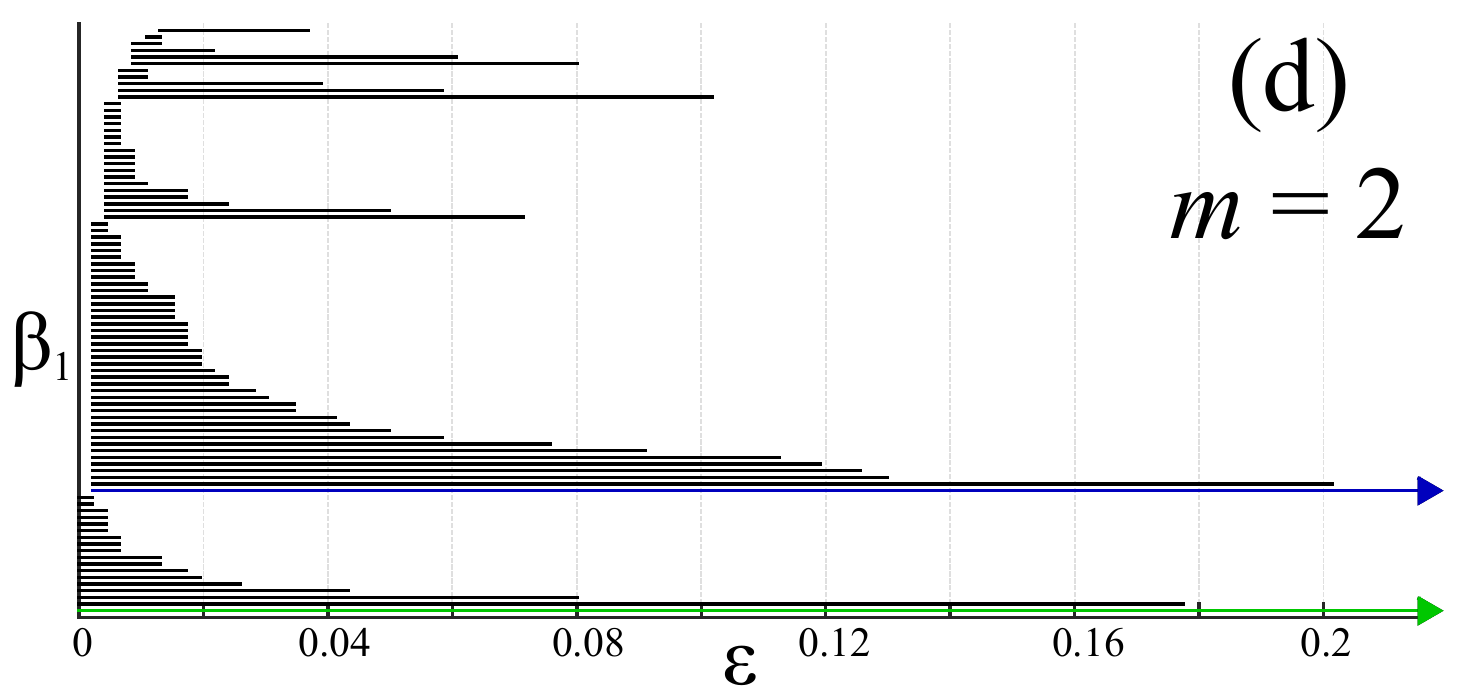}
        \end{subfigure}%
        
\caption{One-skeletons of the witness complexes (top row) and barcode
  diagrams for $\beta_1$ (bottom row) of the Lorenz system.  The plots
  in the left-hand column are computed from the three-dimensional
  $(x,y,z)$ trajectory of Figure~\ref{fig:Lorenz63Example}(a); those in the
  right-hand column are computed from a two-dimensional ($m=2$)
  delay-coordinate reconstruction from the $x$ coordinate of that
  trajectory with $\tau = 174 $.  In both cases, $\ell = 201$ equally
  spaced landmarks (red {\color{red}{$\times$}}s) are used.  Both
  complexes have two persistent nonbounding cycles (green and blue
  edges) but the $2D$ reconstruction requires only $\approx1900$
  simplices to resolve those cycles (at $\eps = 0.2$), while the full
  $3D$ trajectory requires $\approx 7000$ simplices (at $\eps = 1.2$)
  to eliminate spurious loops.}
\label{fig:homology-figure1}
\end{figure}

And that leads to a fundamental question for this thesis: how does the
homology of a delay-coordinate reconstruction change with the 
dimension $m$?  The standard answer to this in the delay-coordinate embedding literature is that the topology should change at
first, then stabilize when $m$ became large enough\footnote{{\it For example,} recall the method of dynamical invariants.}to correctly unfold
the topology of the underlying attractor. In practice, however, a too-large $m$ will invoke the curse of dimensionality and destroy the
fidelity of the reconstruction.  Moreover, increasing $m$ exacerbates both
noise effects and computational expense.  For all of these reasons, it
would be a major advantage if one could obtain useful information
about the homology of the underlying attractor---even if not the full topology---from a low-dimensional
delay-coordinate reconstruction.

Again, it appears that this is possible.
Figure~\ref{fig:LorenzEmbedSkeletons} shows witness complexes for
$m=2$ and $m=3$ reconstructions of the Lorenz time series of
Figure~\ref{fig:Lorenz63Example}(b).
The barcodes for the first two Betti numbers of these two complexes, as computed using
{\tt javaPlex}, have similar structure: the complexes become connected ($\beta_0 = 1$) at a small value of $\eps$, and the dominant, persistent homology corresponds to the two primary holes ($\beta_1 = 2$) in the attractor. 
Note, by the way, that Figure~\ref{fig:LorenzEmbedSkeletons}(a) is not simply a
$2D$ projection of Figure~\ref{fig:LorenzEmbedSkeletons}(b); the edges in each
complex reflect the geometry of the witness relationships in different
spaces, and so may differ.
Higher-dimensional reconstructions---not
easily displayed---have the same homology for suitable choices of
$\eps$, though for $m >5$, it is necessary to increase the number of
landmarks to obtain a persistent $\beta_1 = 2$.

\begin{figure}[htb]
        \centering
        \begin{subfigure}[b]{0.45\textwidth}
                \includegraphics[width=\textwidth]{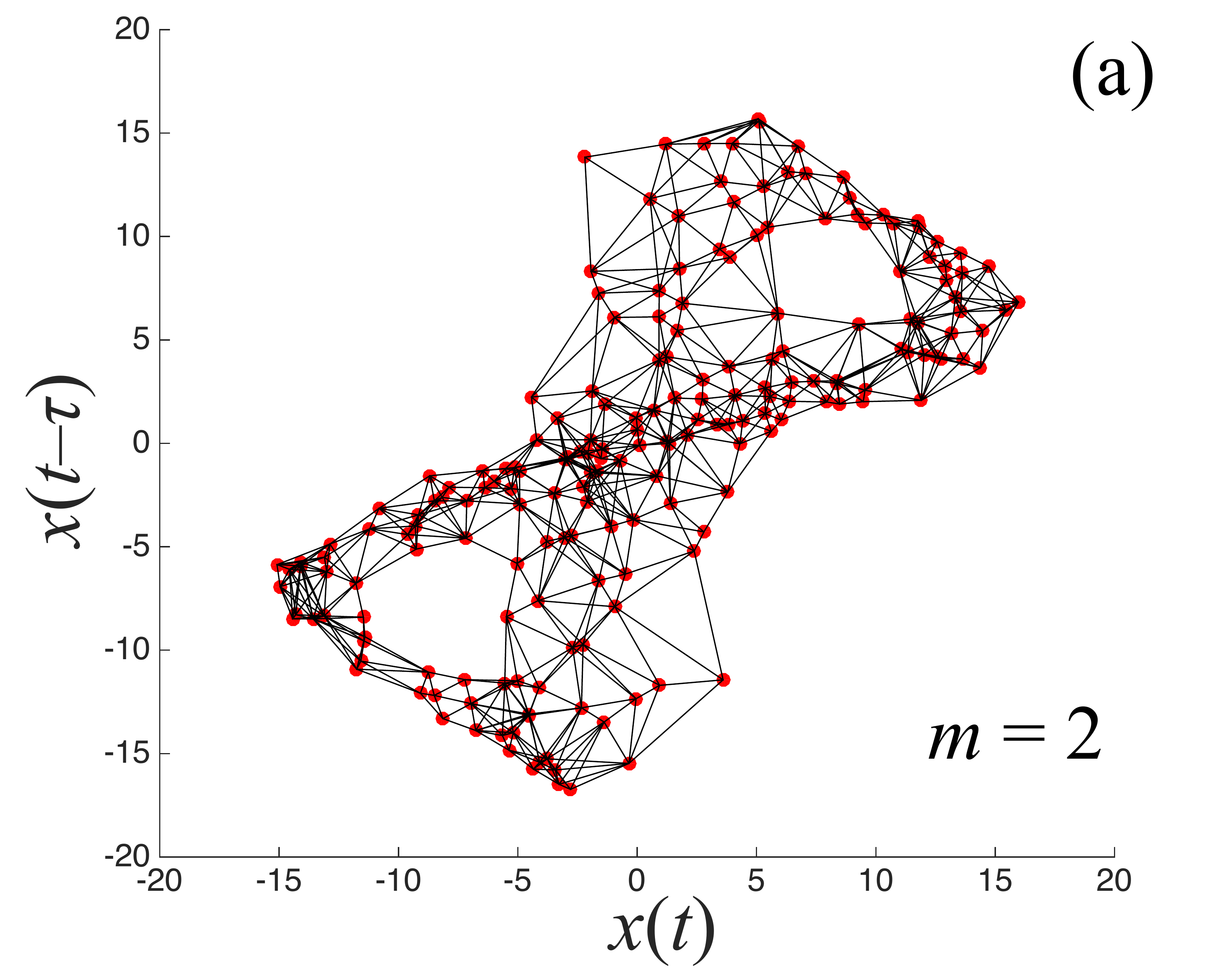}
        \end{subfigure}%
        \begin{subfigure}[b]{0.45\textwidth}
                \includegraphics[width=\textwidth]{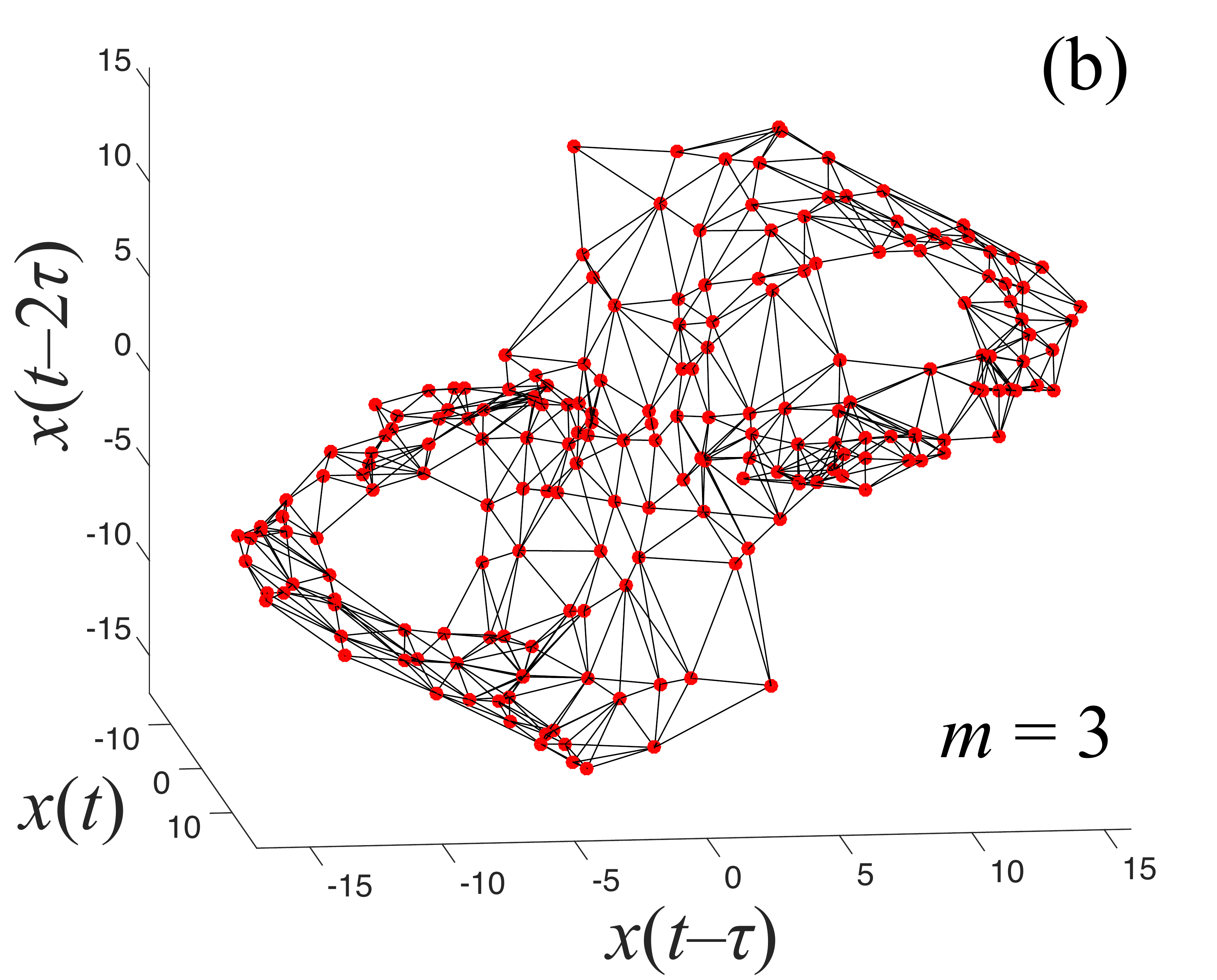}
        \end{subfigure}
\caption{The effects of reconstruction dimension: One-skeletons of
  witness complexes of different reconstructions of the scalar time
  series of Figure~\ref{fig:Lorenz63Example}(b).  Both reconstructions use $\tau
  =174$, the first minimum of the average time-delayed mutual
  information, $\ell=198$ equally spaced
  landmarks (red dots), and $\xi = 0.54\%$, as defined in \Eq{scaledXiDef}.}
\label{fig:LorenzEmbedSkeletons}
\end{figure}

That brings up an important point: if one wants to sensibly compare
witness complexes constructed from different reconstructions of a
single data set, one has to think carefully about the $\ell$ and
$\eps$ parameters.  Here, I use persistence to choose a good value
of $\ell$.  I find that the results are robust with respect to
changes in that value, across all reconstruction dimension values in
this study, so I fix $\ell \approx 200$ for all the experiments
reported in this section.\footnote
{The precise value varies slightly because the length of a
  trajectory reconstructed from a fixed-length data set decreases with
  increasing $m$ (since one needs a full span of $m \times (\tau)$ data
  points to construct a point in the reconstruction space).}

In the experiments in the previous section, the scale parameter $\eps$
was given in absolute units.  To generalize this approach across
different examples and different reconstruction dimensions, it makes sense to compare reconstructions with
$\eps$ chosen to be a fixed fraction of the diameter, $\diam{W}$, of the set $W$
\vspace{-0.75cm}
\beq{scaledXiDef} \eps =  \xi \, \diam{W} \eeq
  For example, for the
full $3D$ attractor in Figure~\ref{fig:Lorenz63Example}(a)
\begin{equation}
	\diam{W_{xyz}} = \sqrt{(x_{max} - x_{min})^2 +
  (y_{max} - y_{min})^2 + (z_{max} - z_{min})^2} = 75.3 
\end{equation}
so the $\eps$ values used in Figure~\ref{fig:Lorenz63-barcodes}---$0.017 \le \eps
\le 1.7$ in absolute units---translate to $2.3 \times 10^{-4} \le
\xi \le 0.023$ in this diameter-scaled measure.

The diameter of the reconstruction varies in a natural way with the dimension $m$.
Since delay-coordinate reconstruction of scalar data unfolds
the full range of those data along every added dimension, the diameter
of an $m$-dimensional reconstruction will be
\begin{equation}
	\diam{W_m} = \sqrt{m(x_{max}-x_{min})^2} = 37.0 \sqrt{m} 
\end{equation}
 for this dataset, where $x$ represents the scalar time-series data.
Since this unfolding will change the geometry of the reconstruction,
I need to scale $\eps$ accordingly.  The witness
complexes in Figure~\ref{fig:LorenzEmbedSkeletons} are constructed with a fixed
value of $\xi = 0.54\%$.
%
%
That is, for Figure~\ref{fig:LorenzEmbedSkeletons}(a), $\eps = 37.0\sqrt{2}(0.0054)
= 0.283$ in absolute units, while for Figure~\ref{fig:Lorenz63Example}(b),
$\diam{W_3}=37.0\sqrt{3}$ and $\eps= 0.346$.
This scaling of
$\eps$---which is used throughout the rest of this section---should
allow the witness complex to adapt appropriately to the effects of
changing reconstruction dimension and finite data.

To formalize the exploration of the reconstruction homology and extend
that study across multiple dimensions, one can use a variant of the
classic barcode diagram that shows, for each simplex, the
reconstruction dimension values at which it appears in and vanishes
from the complex.  Figure~\ref{fig:EdgeStability}(a) shows such a plot for edges
that involve $l_0$, the first landmark on the reconstructed
trajectory.  A number of interesting features are apparent in this
image.  Unsurprisingly, most of the one-simplices that exist in the
$m=1$ witness complex---many of which are likely due to the strong
effects of the projection of the underlying $\bR^d$ trajectory onto
$\bR^1$---vanish when one moves to $m=2$.  There are other short-lived
edges in the complex as well: e.g., the edge from $l_0$ to $l_{120}$
that is born at $m=2$ and dies at $m=3$.  The sketch in
Figure~\ref{fig:EdgeStability}(b) demonstrates how edges can be born as the
dimension increases: in the $m=2$ reconstruction, $\ell_1$ and
$\ell_3$ share a witness (the green square); when one moves to $m=3$,
spreading all of the points out along the added dimension, that
witness is moved far from $\ell_3$---and into the shared region
between $\ell_1$ and $\ell_2$.  There are also long-lived edges in the
complex of Figure~\ref{fig:EdgeStability}(a).  The one between $l_0$ and
$l_{140}$ that persists from $m=1$ to $m=8$ is particularly
interesting: this pair of landmarks has shared witnesses in the scalar
data {\it and in all reconstructions}.  Possible causes for this are
explored in more depth below.  All of these effects depend on $\xi$,
of course; lowering $\xi$ will decrease both the number and average
length of the edge persistence bars.

\begin{figure}[tb!]
\centering
        \centering
        \begin{subfigure}[b]{0.45\textwidth}
		\includegraphics[height=1.8truein]{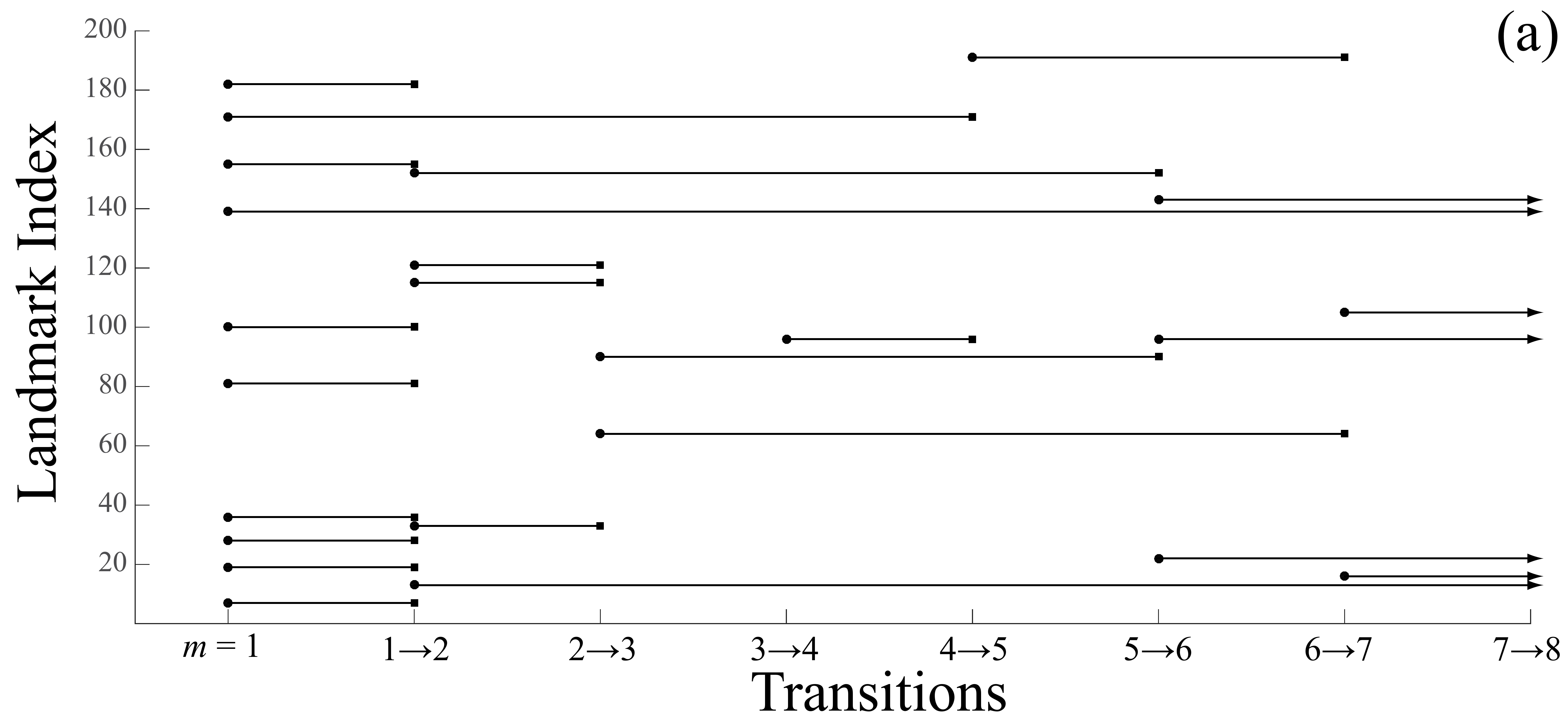}

        \end{subfigure}%
\hspace*{25mm}
        \begin{subfigure}[b]{0.35\textwidth}
		\includegraphics[height=1.8truein]{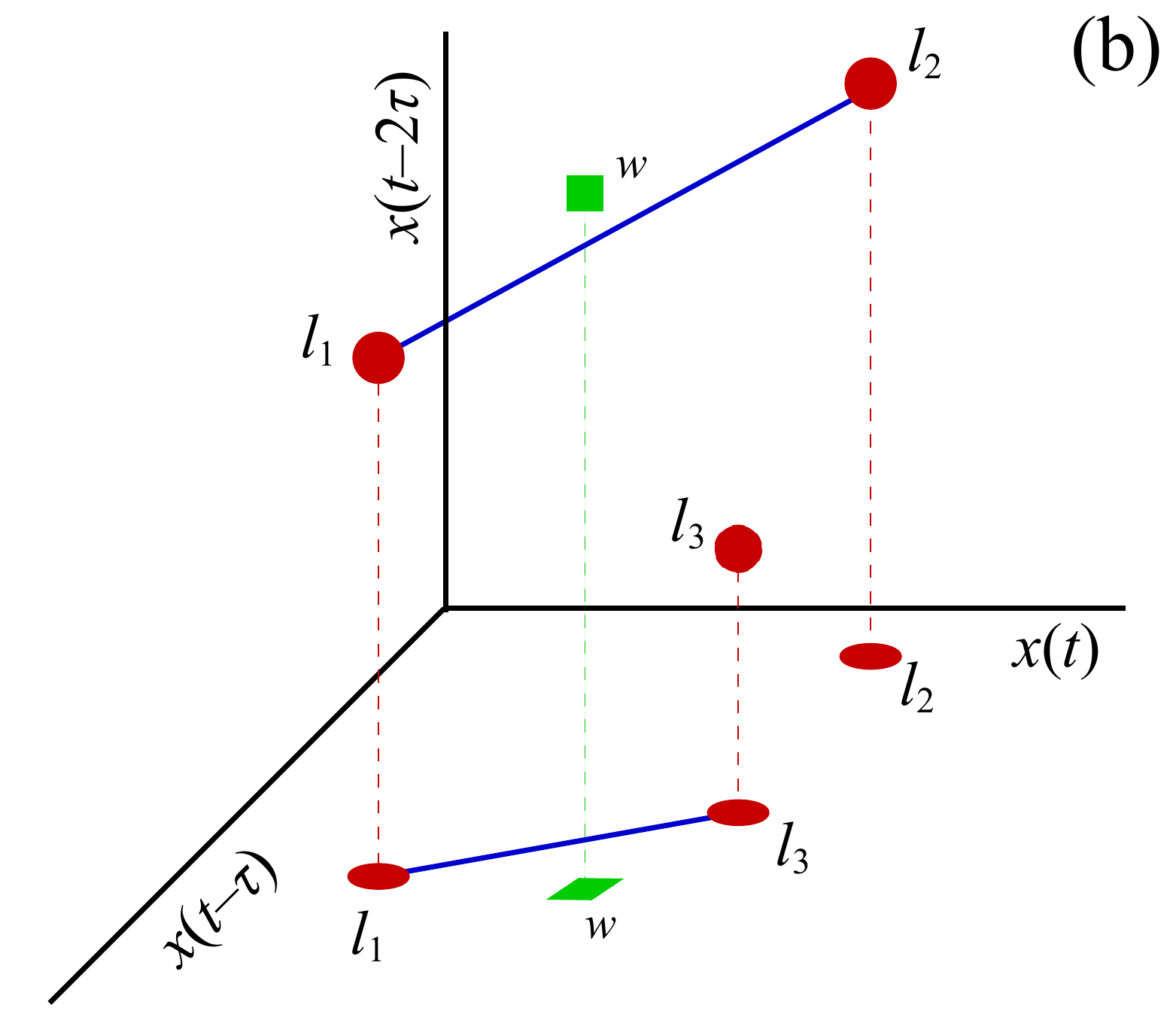}

        \end{subfigure}
        \caption{(a) Dimension barcode for edges in the witness
          complex of the reconstructed scalar time series of
          Figure~\ref{fig:Lorenz63Example}(b) that involve $l_0$, the first
          landmark, for reconstructions with $m=1,\ldots, 8$.  The
          vertical axis is labeled with the indices of the remaining
          $197$ landmarks in the complex; a circle at the
          $m-1\rightarrow m$ tickmark on the horizontal axis indicates
          the transition at which an edge between $l_0$ and $l_i$ is
          born; a square indicates the transition at which that
          edge vanishes from the complex.  An arrow at the
          right-hand edge of the plot indicates an edge that was still
          stable when the algorithm completed. For all
          reconstructions, $\tau =174$, $\ell=198$, and
          $\xi=0.54\%$.  (b) Sketch of the birth and death of edges at
          the $m=2\rightarrow3$ transition. }
\label{fig:EdgeStability}
\end{figure}

While this $\Delta m$ barcode image is interesting, the amount of
detail that it contains makes it somewhat unwieldy.  To study the
$m$-persistence of all of $\ell \times \ell$ edges in a witness
complex, one would need to examine $\ell$ of these plots---or condense them into a
single plot with $\ell^2$ entries on the vertical axis.  Instead, one can plot
what I call an {\it edge lifespan diagram}: an $\ell \times \ell$
matrix whose $(i,j)^{th}$ pixel is colored according to the maximum
range of $m$ for which an edge exists in the complex between the $i^{th}$
and $j^{th}$ landmarks; see  Figure~\ref{fig:BirthDeathLorenz}.
If the edge $\{l_i,l_j\}$ existed in the complex for
$2 \le m<3$ and $5 \le m<8$, for instance, $\Delta m$ would be three and the
$i,j^{th}$ pixel would be coded in cyan. Edges that do not exist for
any dimension are coded white.  

\begin{figure}[tb!]
        \centering
        \begin{subfigure}[b]{0.6\textwidth}
                \includegraphics[width=\textwidth]{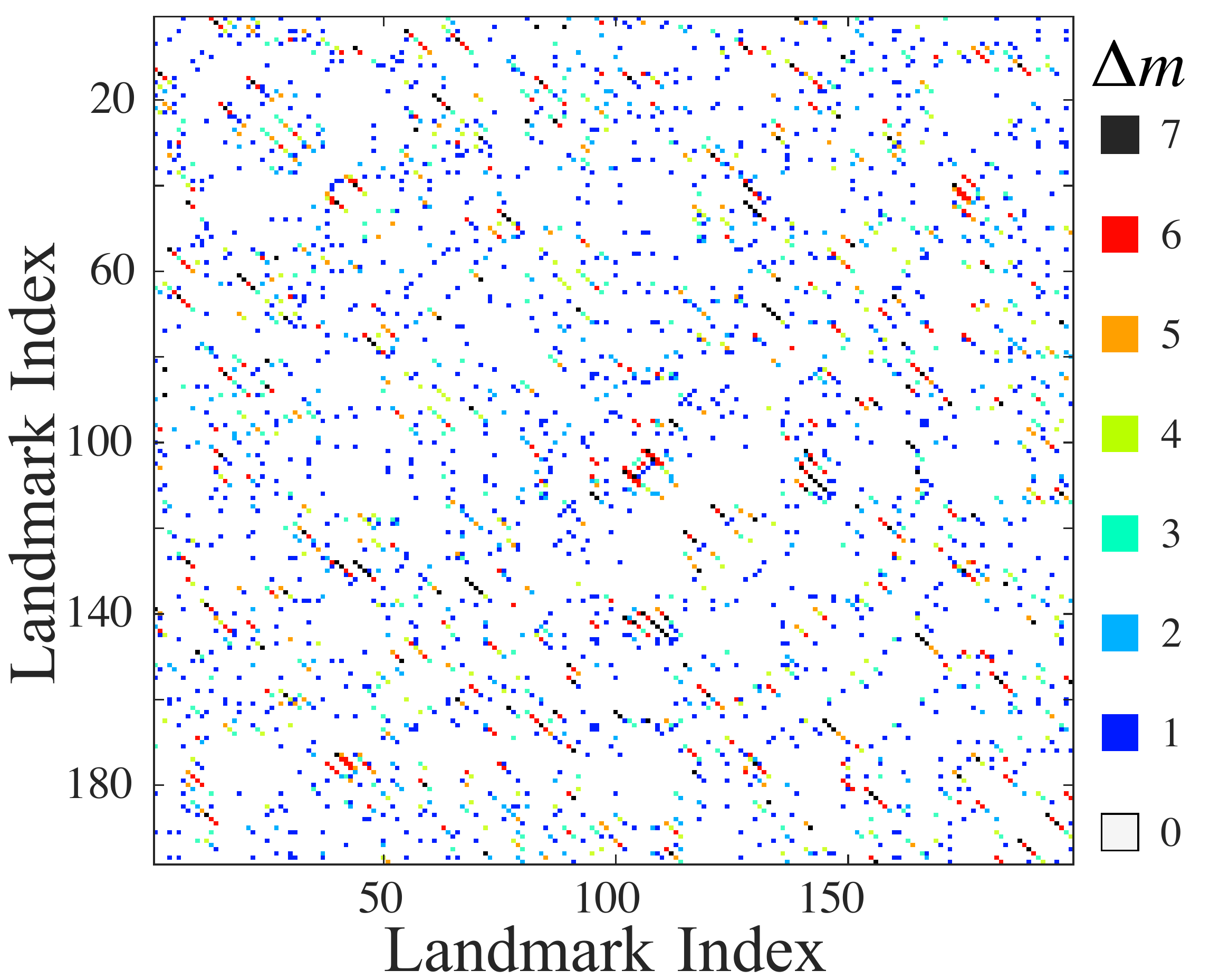}
        \end{subfigure}
%
\caption{Edge lifespan diagram: pixel $i,j$ on this image is
  color-coded according to the maximum range $\Delta m$ of dimension
  for which an edge exists between landmarks $l_i$ and $l_j$ in the
  witness complex of the reconstructed scalar time series of
  Figure~\ref{fig:Lorenz63Example}(b) for $m=1,\ldots, 8$.
%
%
For all reconstructions, $\tau =174$, $\ell=198$, and $\xi=0.54\%$.
}
\label{fig:BirthDeathLorenz}
\end{figure}

A prominent feature of Figure~\ref{fig:BirthDeathLorenz} is a large number
($683$) of edges with a lifespan $1$ (blue). Of these edges, $463$
exist for $m=1$, but not for $m=2$, and thus reflect the anomalous
behavior of projecting a $2.06$ dimensional object onto a line.  This
is also seen, as described above, in the barcode of Figure~\ref{fig:EdgeStability}.

Another interesting set of features in the lifespan diagram is the
diagonal line segments.  Note that the {\it color} of the pixels in
these segments varies, though most of them correspond to edges with
longer lifespans.  These segments indicate the existence of $\Delta
m$-persistent edges $\{l_i,l_j\}, \{l_{i+1},l_{j+1}\}, \{l_{i+2},
l_{j+2}\} \ldots$.  This is likely due to the continuity of the
dynamics \cite{Alexander12}.  Recall that the landmarks are evenly
spaced in time, so $l_{i+1}$ is the $\Delta t$-forward image of $l_i$.
Thus a diagonal segment may indicate that the $\Delta t$-forward
images of (at least one) witness that is shared between $l_i$ and
$l_j$ is shared between $l_{i+1}$ and $l_{j+1}$, and so on.  The
lengths of the longer line segments suggest that that continuity fails
after $5$-$10$ $\Delta t$ steps, probably because of the positive
Lyapunov exponents on the attractor.  As a simple check on this
reasoning, one can compute an edge lifespan diagram for a dynamical
system with a limit cycle.
The structure of such a plot (not shown) is dominated by diagonal
lines of high $\Delta m$-persistence, with a few other scattered
one-persistent edges.  

The rationale behind studying the {\it maximal} $m$-lifespan goes
back to one of the basic premises of persistence: that features that
persist for a wide range of parameter values are in some sense
meaningful.  To explore this, Figure~\ref{fig:delta-m} shows the witness complex
of Figure~\ref{fig:LorenzEmbedSkeletons}(a), highlighting the $\Delta m \ge
2$-persistent edges: those that exist at $m=2$ and persist at least to
$m=4$.  There exists a fundamental core to the complex that persists
as the dimension grows and thus is robust to geometric distortion, but
there are also short-lived edges that fill in the complex in accord
with the local geometric structure of the reconstruction.  Indeed,
when $m=2$, the projection artificially compresses near the origin;
small simplicies fill in this region due to the landmark clustering
there. However, in the transition to $m=3$---{\it viz.,}
Figure~\ref{fig:LorenzEmbedSkeletons}(b)---this region stretches away from the
origin, spreading the landmarks out.  There is a similar cluster of
``fragile" edges near the lower left corner of the complex.

Even though geometric evolution with increasing reconstruction dimension
leads to the death of many local edges, the large-scale
homology is correct in both complexes of Figure~\ref{fig:LorenzEmbedSkeletons},
although the fine-scale topology is resolved differently by the
dimension-dependent geometry.  So while the edges with longer lifespan
are indeed more important to the core structure, the short-lived edges
are also important because they allow the complex to adapt to the
geometric evolution of the attractor and fill in the details of the
skeleton that are necessary and meaningful in that dimension.

In the spirit of the false near-neighbor method \cite{KBA92}, one
might be tempted to take the short-lived edges as an indication that
the reconstruction dimension is inadequate.  However, one computes
homology {\it from the overall complex}.  As the example above shows,
homology is relatively robust with respect to individual edges.  The
moral of this story is that the lifespan of an edge is not necessarily
an obvious indication of its importance to the homology of the
complex; $\Delta m$-persistence plays a different role here than the
abscissa of traditional barcode persistence plots.

\begin{figure}[tb!]
        \centering
                \includegraphics[width=0.95\textwidth]{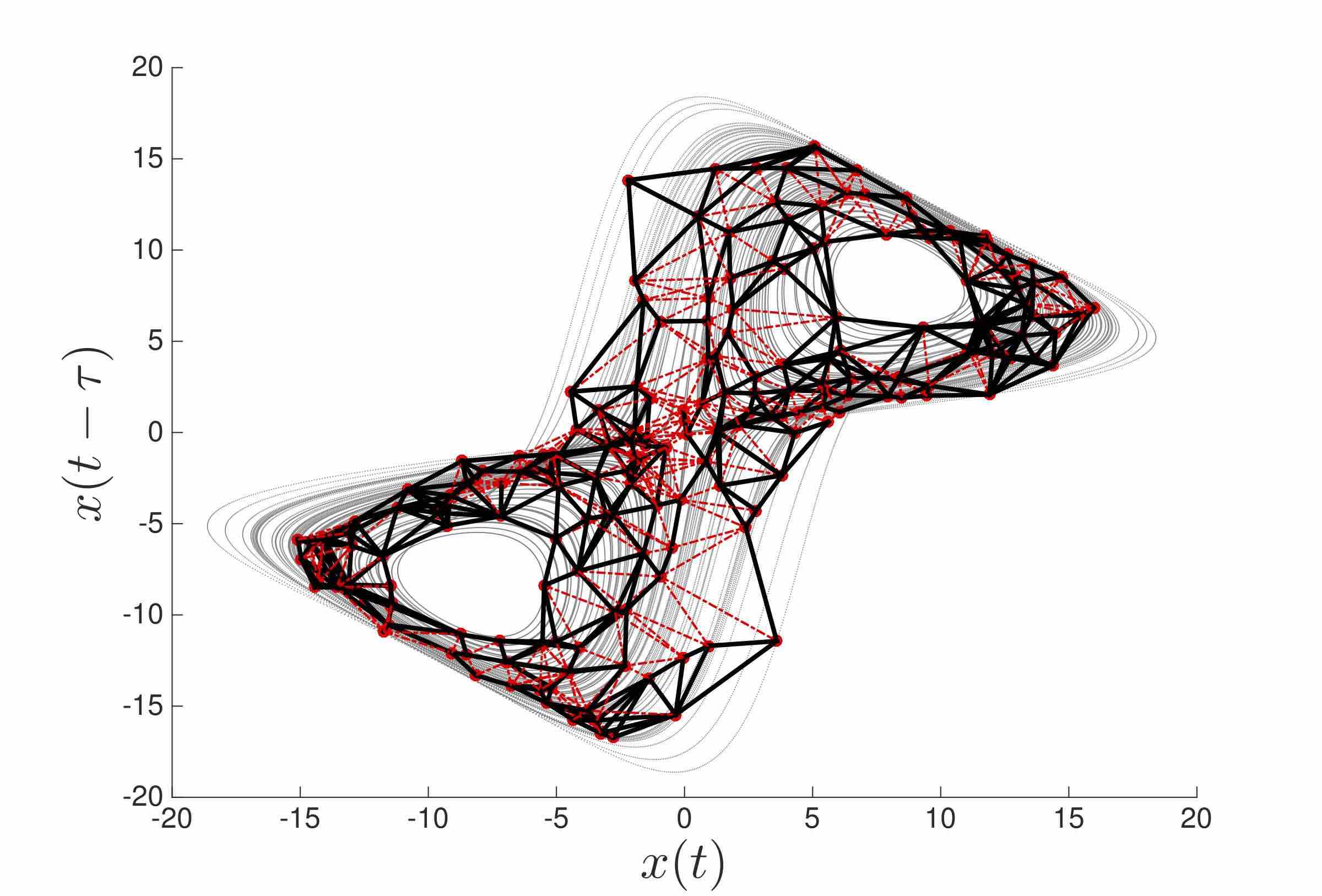}
 \caption{Witness complex of Figure~\ref{fig:LorenzEmbedSkeletons}(a) with $\Delta m \ge 2$-persistent edges shown as thick (black) lines, and the $\Delta m = 1$ edges as (red) dashed lines.}
 \label{fig:delta-m}
 \end{figure}

The results in this section show that it is possible to compute the homology of an invariant set
of a dynamical system using a simplicial complex built 
 from a low-dimensional reconstruction
of a scalar time series.  These results have a
number of interesting implications.  Among other things, they suggest
that the traditional delay-coordinate reconstruction process may be
excessive if one is only interested in large-scale topological structure such as the homology.  
This is directly apropos of the central claim of this thesis, as it explains why it is possible to construct
accurate predictions of the future state of a high-dimensional
dynamical system using a two-dimensional delay-coordinate
reconstruction.  The delay-coordinate machinery
strives to obtain a diffeomorphism---not a homeomorphism---between the true and
reconstructed attractors.  However, many of the properties of
attractors that are important for forecasting (continuity, recurrence, entropy, etc.) are
topological, so requiring only a homeomorphism is not only natural, but also more efficient \cite{mischaikow99}.

This section uses a single example---the Lorenz 63 system---but I
believe that the approach will work on other dynamical systems and I
plan in the future to do a careful exploration of additional systems,
both maps and flows.

%% file: pecomplexity.tex

\chapter{Model-Free Quantification of Time-Series Predictability}\label{ch:wpe}

Time-series data can span a wide range of complexities and no {\it single} forecast algorithm can be expected to handle this entire spectrum effectively. This poses an interesting problem when one is developing any new forecasting technology, such as the reduced order framework outlined in this thesis. In particular, given an arbitrary time series---with an undefined level of complexity---can one expect that method to be effective? The first step to answering this question is to define a spectrum of predictive complexity~\cite{josh-pre}.

On the low end of this spectrum are time series that exhibit perfect
predictive structure, {\it i.e,} signals whose future values can be
perfectly predicted from past values.  Signals like this can be
viewed as the product of an underlying process that generates
information and/or transmits it from the past to the future in a
perfectly predictable fashion.  Constant or periodic signals, for
example, fall in this class.  On the opposite end of this spectrum are
signals that are---from a forecasting perspective---{\it fully complex}, where the
underlying generating process transmits no information at all from the
past to the future.  White noise processes fall in this class.  In
fully complex signals, knowledge of the past gives no insight into the
future, regardless of what model one chooses to use. Signals in the
midrange of this spectrum, {\it e.g.,} deterministic chaos, pose interesting
challenges from a modeling perspective.  In these signals, enough
information is being transmitted from the past to the future that an
{\it ideal} model---one that captures the generating process---can
forecast the future behavior of the observed system with high
accuracy. \label{page:spectrum}

This leads naturally to an important and challenging question to which I alluded in the first paragraph of this chapter: given a
noisy real-valued time series from an unknown system, does there exist
any forecast model that can leverage the information (if any) that is
being transmitted forward in time by the underlying generating
process?  A first step in answering this question is to reliably
quantify where on the complexity spectrum a given time series falls; a
second step is to determine how complexity and predictability are
related in these kinds of data sets.  With these answers in hand, one
can develop a practical strategy for assessing appropriateness of
forecast methods for a given time series.  If the forecast produced by
\roLMA is poor, for example, but the time series
contains a significant amount of predictive structure, one can
reasonably conclude that \roLMA is inadequate to the task and
that one should seek another method.  

The goal of this chapter is to
develop effective heuristics to put that strategy into practice. Recall that Chapter~\ref{ch:pnp} of this thesis demonstrated that \roLMA~{\it can} be effective, and Chapter~\ref{ch:explain} provided reasons {\it why} and {\it how} this is the case. The heuristic proposed in this chapter goes one step further and addresses {\it when} \roLMA---and indeed any forecast algorithm---can be expected to be effective.

The information in an observation can be partitioned into two pieces:
redundancy and entropy generation~\cite{crutchfield2003}.
\label{page:redundancy}
 The approach exploits this decomposition in order to assess how much
predictive structure is present in a signal---{\it i.e.,} where it falls on
the complexity spectrum mentioned above.  I define {\it complexity}
as a particular approximation of Kolmogorov-Sinai
entropy~\cite{lind95}.  That is, I view a random-walk time series
(which exhibits high entropy) as purely complex, whereas a low-entropy
periodic signal is on the low end of the complexity spectrum.  This
differs from the notion of complexity used by {\it e.g.,} \cite{Shalizi2008},
which would consider a time series without any statistical
regularities to be non-complex.  In collaboration with R.~G.~James, I argue that an extension of
{\it permutation entropy}~\cite{bandt2002per}---a method for
approximating the entropy through ordinal analysis---is an effective
way to assess the complexity of a given time series.  Permutation
entropy, introduced in Section~\ref{sec:meaComplex-intro}, is ideal for this purpose because it works with real-valued
data and is known to converge to the true entropy value. Other
existing techniques either require specific knowledge of the
generating process or produce biased values of the
entropy~\cite{bollt2001}.

I focus on real-valued, scalar, time-series data from physical experiments.
I do not assume any knowledge of the generating process or its
properties: whether it is linear, nonlinear, deterministic,
stochastic, etc.  To explore the relationship between complexity,
predictive structure, and actual predictability, I generate forecasts
for several experimental computer performance time-series datasets using the five
different prediction strategies discussed in this thesis, then compare the accuracy of those
predictions to the permutation entropy of the associated signals. This results in two primary findings: \label{page:wpefindings}
\begin{enumerate}

\item The permutation entropy of a noisy real-valued time series from
  an unknown system is correlated with the accuracy of an appropriate
  predictor.

\item The relationship between permutation entropy and prediction
  accuracy is a useful empirical heuristic for identifying mismatches
  between prediction models and time-series data.

\end{enumerate}
There has, of course, been a great deal of good work on different ways
to measure the complexity of data, and previous explorations have
confirmed repeatedly that complexity is a challenge to prediction.  It
is well known that the way information is generated and processed
internally by a system plays a critical role in the success of
different forecasting methods---and in the choice of which method is
appropriate for a given time series.  This constellation of issues has
not been properly explored, however, in the context of noisy, poorly
sampled, real-world data from unknown systems.  That exploration, and
the development of strategies for putting its results into effective
practice, is the primary contribution of this chapter.  The empirical
results in Section~\ref{sec:wperesults} not only elucidate the
relationship between complexity and predictability, but also provide a practical strategy to aid practitioners in assessing the
appropriateness of a prediction model for a given real-world noisy time series
from an unknown system---a challenging task for which little guidance
is currently available.  In the context of this thesis, the value of this is that it provides a general framework for assessing whether or not \roLMA is an appropriate choice for a given time series, or if a more sophisticated or even a {\it simpler} strategy is required.

The rest of this chapter is organized as follows.
Section~\ref{sec:related} discusses previous results on generating
partitions, local modeling, and error distribution analysis, and
situates this work in that context.  In
Section~\ref{sec:wperesults}, I estimate the complexity of a number of time-series traces and compare that complexity to the accuracy of
various predictions models 
operating on that time series.  In Section~\ref{sec:conc}, I discuss
these results and their implications, and consider future areas of
research.

\section{Traditional Methods for Predicting Predictability}\label{sec:related}

Hundreds, if not thousands, of strategies have been developed
for a wide variety of prediction tasks.  The purpose of this chapter is
not to add a new weapon to this arsenal, nor to do any sort of
theoretical assessment or comparison of existing methods.  In the spirit of this thesis, the goals here 
are focused more on the {\it practice} of prediction: {\it (i)} to
empirically quantify the predictive structure that is present in a
real-valued scalar time series and {\it (ii)} to explore how the
performance of prediction methods is related to that inherent
complexity.  It would, of course, be neither practical nor interesting
to report results for every existing forecast strategy; instead, I use the same representative set of methods that appear throughout this thesis, as described in Section~\ref{sec:forecastmodels}.

Quantifying predictability, which is sometimes called ``predicting
predictability,'' is not a new problem.  Most of the corresponding
solutions fall into two categories that I call model-based error
analysis and model-free information analysis.
The first class focuses on errors produced by a specific forecasting
schema.  This analysis can proceed locally or globally.  The local
version approximates error distributions for different regions of a
time-series model using local ensemble in-sample\footnote{The terms ``in sample'' and ``out of sample'' are
  used in different ways in the forecasting community.  Here, I
  distinguish those terms by the part of the time series that is the
  focus of the prediction: the observed data for the former and the
  unknown future for the latter.  In-sample forecasts---comparisons of
  predictions generated from {\it part} of the observed time
  series---are useful for assessing model error and prediction
  horizons, among other things.}forecasting.
These distributions are then used as estimates of out-of-sample
forecast errors in those regions.  For example, Smith {\it et al.}
make in-sample forecasts using ensembles around selected points in
order to predict the local predictability of a time
series~\cite{Smith199250}.  This approach can be used to show that
different portions of a time series exhibit varying levels of
local predictive uncertainty.  

Local model-based error analysis works quite well, but it only
approximates the {\it local} predictive uncertainty {\it in relation
  to a fixed model}.  It cannot quantify the {\it inherent} predictability
of a time series and thus cannot be used to draw conclusions about
predictive structure that other forecast methods may be able to leverage.
Global model-based error analysis moves in this direction.  It uses
out-of-sample error distributions, computed {\it post facto} from a
class of models, to determine which of those models was best.  After
building an autoregressive model, for example, it is common to
calculate forecast errors and verify that they are normally
distributed.  If they are not, that suggests that there is structure
in the time series that the model-building process was unable to
recognize, capture, and exploit.  The problem with this approach is
lack of generality.
\label{page:normal-errors}
Normally distributed errors indicate that a model has captured the
structure in the data insofar as is possible, {\it given the
  formulation of that particular model} ({\it viz.,} the best possible
linear fit to a nonlinear dataset).  This gives no indication as to
whether another modeling strategy might do better.

A practice known as deterministic vs. stochastic
modeling~\cite{weigend93, Casdagli92dvsplots} bridges the gap
between local and global approaches to model-based error analysis.
The basic idea is to construct a series of local linear fits,
beginning with a few points and working up to a global linear fit that
includes all known points, and then analyze how the average
out-of-sample forecast error changes as a function of number of points
in the fit. The shape of such a ``DVS" graph indicates the amounts of
determinism and stochasticity present in a time series.

The model-based error analysis methods described in the previous three
paragraphs are based on specific assumptions about the underlying
generating process and knowledge about what will happen to the error
if those assumptions hold or fail.  Model-{\it free} information
analysis moves away from those restrictions.  My approach falls into
this class: I wish to measure the inherent complexity of an arbitrary empirical
time series, then study the correlation of that complexity with the
predictive accuracy of forecasts made using a number of different
methods.

I build on the notion of {\it redundancy} that was introduced on
page~\pageref{page:redundancy}, which formally quantifies how
information propagates forward through a time series:
{\it i.e.,} the mutual information between the past $n$ observations and the
current one.
The redundancy of i.i.d. random processes, for instance, is zero,
since all observations in such a process are independent of one
another.  On the other hand, deterministic systems, including chaotic
ones, have high redundancy---in fact, {\it maximal} redundancy in the infinite limit---and
thus they can be perfectly predicted if observed for long
enough~\cite{weigend93}.  In practice, it is quite difficult to
estimate the redundancy of an arbitrary, real-valued time series.
Doing so requires knowing either the Kolmogorov-Sinai entropy or the
values of all positive Lyapunov exponents of the system.  Both of
these calculations are difficult, the latter particularly so if the
data are very noisy or the generating system is stochastic.

Using entropy and redundancy to quantify the inherent predictability
of a time series is not a new idea.  Past methods for this, however,
({\it e.g.,}~\cite{Shannon1951, mantegna1994linguistic}) have hinged on
knowledge of the {\it generating partition} of the underlying
process, which lets one transform real-valued observations into
symbols in a way that preserves the underlying dynamics~\cite{lind95}.
Using a partition that is not a generating partition---{\it e.g.,} simply
binning the data---can introduce spurious complexity into the
resulting symbolic sequence and thus misrepresent the entropy of the
underlying system~\cite{bollt2001}.  Generating partitions are
luxuries that are rarely, if ever, afforded to an analyst, since one
needs to know the underlying dynamics in order to construct one.  And
even if the dynamics are known, these partitions are difficult to
compute and often have fractal boundaries~\cite{eisele1999}. (See Section~\ref{sec:binning} for a review of these issues.)

In the development described in the following section, I sidestep these issues by using a variant of the {\it permutation
  entropy} of Bandt and Pompe~\cite{bandt2002per} to estimate the
value of the Kolmogorov-Sinai entropy of a real-valued time
series---and thus the redundancy in that data, which my results
confirm to be an effective proxy for predictability.  This differs
from existing approaches in a number of ways.  It does not rely on
generating partitions---and thus does not introduce bias into the
results if one does not know the dynamics or cannot compute the
partition.  Permutation entropy makes no assumptions about, and
requires no knowledge of, the underlying generating process: whether it is linear or nonlinear, what its Lyapunov spectrum is, etc.  These features make my
approach applicable to noisy real-valued time series from all classes
of systems, deterministic and stochastic.

%

 \section{Predictability, Complexity, and Permutation Entropy 
 } 
 \label{sec:wperesults}

In this section, I offer an empirical validation of the two findings
introduced on page~\pageref{page:wpefindings}, namely:

\begin{enumerate}

\item The weighted permutation entropy (WPE) of a noisy real-valued
  time series from an unknown system is correlated with prediction
  accuracy---{\it i.e.,} the predictable structure in an empirical
  time-series data set can be quantified by its WPE.

\item The relationship between WPE and mean absolute scaled error
  (1-MASE) is a useful empirical heuristic for identifying mismatches
  between prediction models and time-series data---{\it i.e.,} when there is
  structure in the data that the model is unable to exploit.

\end{enumerate}

The experiments below involve four different prediction methods: \fnnLMA, \naive, ARIMA and random walk,
applied to time-series data from eight different experimental systems: {\tt
  col\_major}, {\tt 403.gcc}, and six different segments of a computer performance experiment that I have not yet discussed in this thesis called \svd, a Fortran program from the LAPACK linear algebra package~\cite{lapack} that calculates the singular value decomposition of a
rectangular $M$ by $N$ matrix with real-valued entries.  For my
experiments, I choose $M=750$ and $N=1000$ and generate the matrix
entries randomly.

\begin{figure}[tb!]
    \centering
    \vspace*{0.5cm}
    \includegraphics[width=\columnwidth]{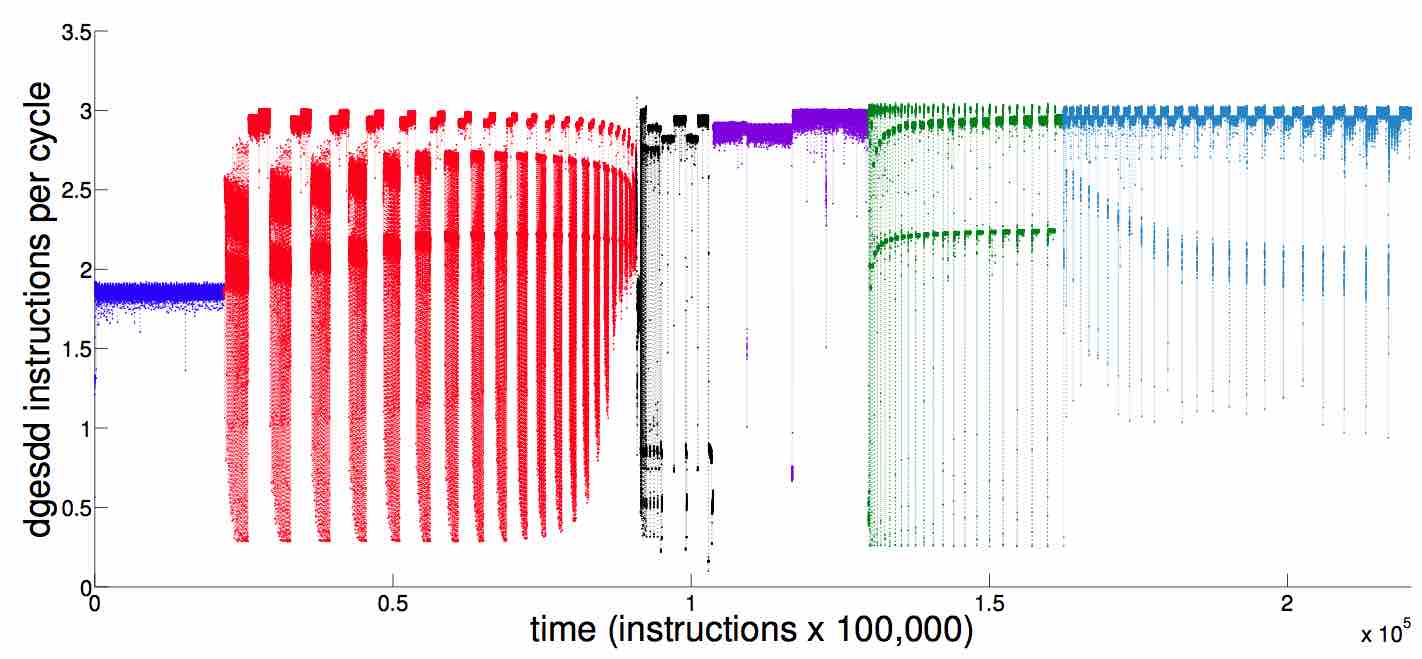}
    \caption{A processor performance trace of instructions per cycle (IPC) during the execution of \svd.  The colors (also separated by vertical dashed lines)
      identify the different {\it segments} of the signal that are
      discussed in the text.}
    \label{fig:svd-ts-colored}
  \end{figure}

 The behavior of this program as it computes the singular values of
this matrix is complex and interesting, as is clearly visible in Figure
\ref{fig:svd-ts-colored}.  
As the code moves though its different phases---diagonalizing the
matrix, computing its transpose, multiplying, etc.---the processor
utilization patterns change quite radically.  For the first
$\sim$21,000 (in units of 100,000 instructions) measurements,
roughly 1.8 instructions are executed per cycle, on the average, by
the eight processing units on this chip.  After that, the IPC moves
through a number of different oscillatory regimes, which I have
color-coded in the figure in order to make textual cross-references
easy to track.

The wide range of behaviors in Figure~\ref{fig:svd-ts-colored}
provides a distinct advantage, for the purposes of this chapter, in that
a number of different generating processes---with a wide range of
complexities---are at work in different phases of a single time
series.  The \col and \gcc traces in Figures~\ref{fig:col12DEmbedding} and
\ref{fig:gccipc}, in contrast, appear to be far more consistent over time---probably
the result of a single generating process with consistent complexity.
\svd, has multiple regimes, each probably the result of different
generating processes.  To take advantage of this rich experimental data set, I split the signal
into six different segments, thereby obtaining an array of examples
for the analyses in the following sections.  For notational
convenience, I refer to these 90 time-series data sets\footnote{15
  runs, each with six regimes} as {\tt dgesdd$_i$}, with $i \in
\{1\dots6\}$ where $i$ corresponds to one of the six segments of the
signal, ordered from left to right.  These segments, which were
determined visually, are shown in different colors in
Figure~\ref{fig:svd-ts-colored}.  Visual decomposition is subjective,
of course, particularly since the regimes exhibit some fractal
structure.  Thus, it may be the case that more than one generating
process is at work in each of our segments.  This is a factor in the
discussion that follows.

The objective
of these experiments is to explore how prediction accuracy is related
to WPE.
Working from the first 90\% of each signal, I generate a prediction
of the last 10\% using all four 
prediction methods, then
calculate the 1-MASE value of those predictions.  I also calculate
the WPE of each time series using a wordlength chosen via the
procedure described at the end of Section~\ref{sec:meaComplex-intro}.  In
order to assess the run-to-run variability of these results, I
repeat all of these calculations on 15 separate trials: {\it i.e.,} 15
different runs of each program.

\begin{figure}[tb!]
  \centering
  \vspace{0.5cm}
  \includegraphics[width=\textwidth]{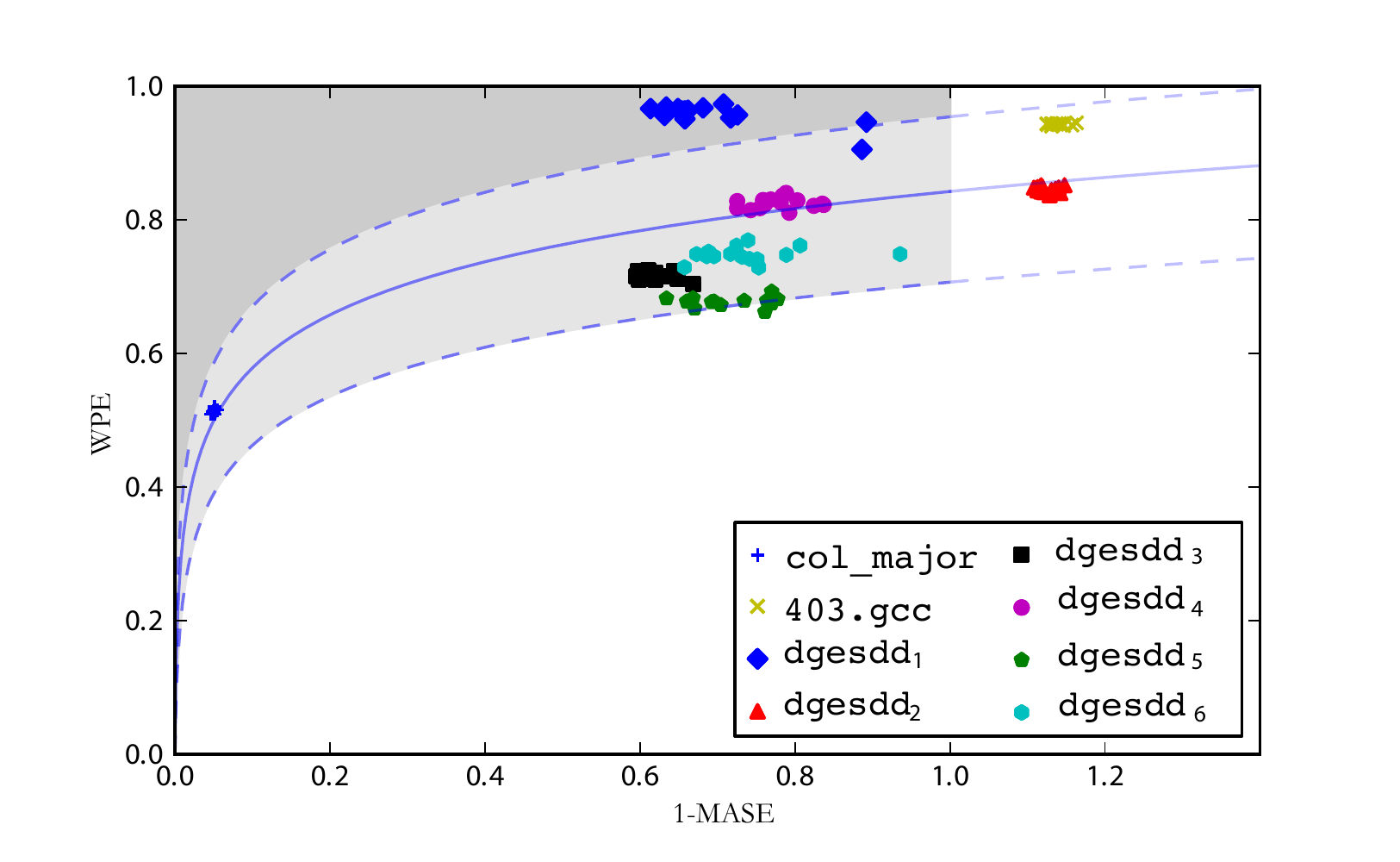}
  \caption{ Weighted permutation entropy versus 1-MASE of the best prediction of a number of different time series.  The solid
    curve is a least-squares log fit of these points.
%
%
The dashed curves reflect the standard deviation of the model in its
parameter space.  Points that lie below and to the right of the shaded
region indicate that the time series has more predictive structure
than the forecast strategy is able to utilize.}
  \label{fig:wpe_vs_mase_best}
\end{figure}

Figure~\ref{fig:wpe_vs_mase_best} plots the WPE values versus the
corresponding 1-MASE values of the {\it best} prediction for each of
the 120 time series in this study.
%
%
The obvious upward trend is consistent with the notion
that there is a pattern in the WPE-MASE relationship.  However, a
simple linear fit is a bad idea here.  First, any signal with zero
entropy should be perfectly predictable ({\it i.e.,} 1-MASE $\approx 0$), so
any curve fitted to these data should pass through the origin.
Moreover, WPE does not grow without bound, so one would expect the
patterns in the WPE-MASE pairs to reach some sort of asymptote.  For
these reasons, I choose to fit a function\footnote{The specific values of the coefficients
  are $a=7.97 \times 10^{-2}$ and $b=1.52 \times 10^3$.}of the form $y = a \log(b x
+ 1)$ to these points, with $y =$
WPE and $x=$ 1-MASE.  The solid curve in the figure shows this fit; the
dashed curves show the standard deviation of this model in its
parameter space: {\it i.e.,} $y = a \log(b x + 1)$ with $\pm$ one standard
deviation on each of the two parameters.  Points that fall within this
deviation volume (light grey) correspond to predictions that are
comparable to the best ones found in this study; points that fall
{\it above} that volume (dark grey) are better still.  I choose to
truncate the shaded region because of a subtle point regarding the
1-MASE of an ideal predictor, which should not be larger than 1 unless
the training and test signals are different.  This is discussed at
more length below.

The curves and regions in Figure~\ref{fig:wpe_vs_mase_best} are a
graphical representation of the first finding introduced on page~\pageref{page:wpefindings}. This representation
is, I believe, a useful heuristic for determining whether a given
prediction method is well matched to a particular time series. If your point is outside the grey region, then that particular model is not capturing all the available structure of the time series.  This is
not, of course, a formal result.  The forecast methods and data sets
used here were chosen to span the space of standard prediction
strategies and the range of dynamical behaviors, but they do not cover
those spaces exhaustively.  My goal here is an {\it empirical}
assessment of the relationship between predictability and complexity,
not formal results about a ``best'' predictor for a given time series.
There may be other methods that produce lower 1-MASE values than those
in Figure~\ref{fig:wpe_vs_mase_best}, but the sparseness of the points
above and below the one-$\sigma$ region about the dashed curve in this
plot strongly suggests a pattern of correlation between the underlying
predictability of a time series and its WPE.  The rest of this section
describes these results and claims in more detail---including the
measures taken to assure meaningful comparisons across methods,
trials, and programs---elaborates on the meaning of the different
curves and limits in the figure, and ties these results into the overall thesis goals.
 \begin{table*}[tb!]
\caption{1-MASE scores and weighted
  permutation entropies for all eight examples studied in this chapter.
  LMA $=$ Lorenz method of analogues; RW $=$ random-walk prediction.
  }
  \begin{center}
  \begin{tabular*}{\textwidth}{@{\extracolsep{\fill} } cccccc}
  \hline\hline 
Signal & RW 1-MASE & na\"{i}ve 1-MASE & \arima 1-MASE & \fnnLMA 1-MASE & WPE \\
\hline
  \col       & $1.001 \pm 0.002$ & $0.571 \pm 0.002$  & $0.599 \pm 0.211$ & $0.050 \pm 0.002$ & $0.513$\\
  \gcc       & $1.138 \pm 0.011$ & $1.797 \pm 0.010$  & $1.837 \pm 0.016$ & $1.530 \pm 0.021$ & $0.943$\\
  \svdone    & $0.933 \pm 0.095$ & $2.676 \pm 4.328$  & $0.714 \pm 0.075$ & $0.827 \pm 0.076$ & $0.957$\\
  \svdtwo    & $1.125 \pm 0.012$ & $3.054 \pm 0.040$  & $2.163 \pm 0.027$ & $1.279 \pm 0.020$ & $0.846$\\
  \svdthree  & $0.707 \pm 0.009$ & $31.386 \pm 0.282$ & $0.713 \pm 0.010$ & $0.619 \pm 0.021$ & $0.716$\\
  \svdfour   & $1.034 \pm 0.035$ & $2.661 \pm 0.074$  & $0.979 \pm 0.032$ & $0.779 \pm 0.036$ & $0.825$\\
  \svdfive   & $1.001 \pm 0.047$ & $20.870 \pm 0.192$ & $2.370 \pm 0.051$ & $0.718 \pm 0.048$ & $0.678$\\
  \svdsix    & $1.060 \pm 0.055$ & $2.197 \pm 0.083$  & $1.438 \pm 0.061$ & $0.739 \pm 0.068$ & $0.748$\\
    \hline\hline
  \end{tabular*}
  \end{center}
 \label{tab:wpeerror}
  \end{table*}%

\begin{figure}[bt!]
  \centering
  \vspace{0.5cm}
  \includegraphics[width=\textwidth]{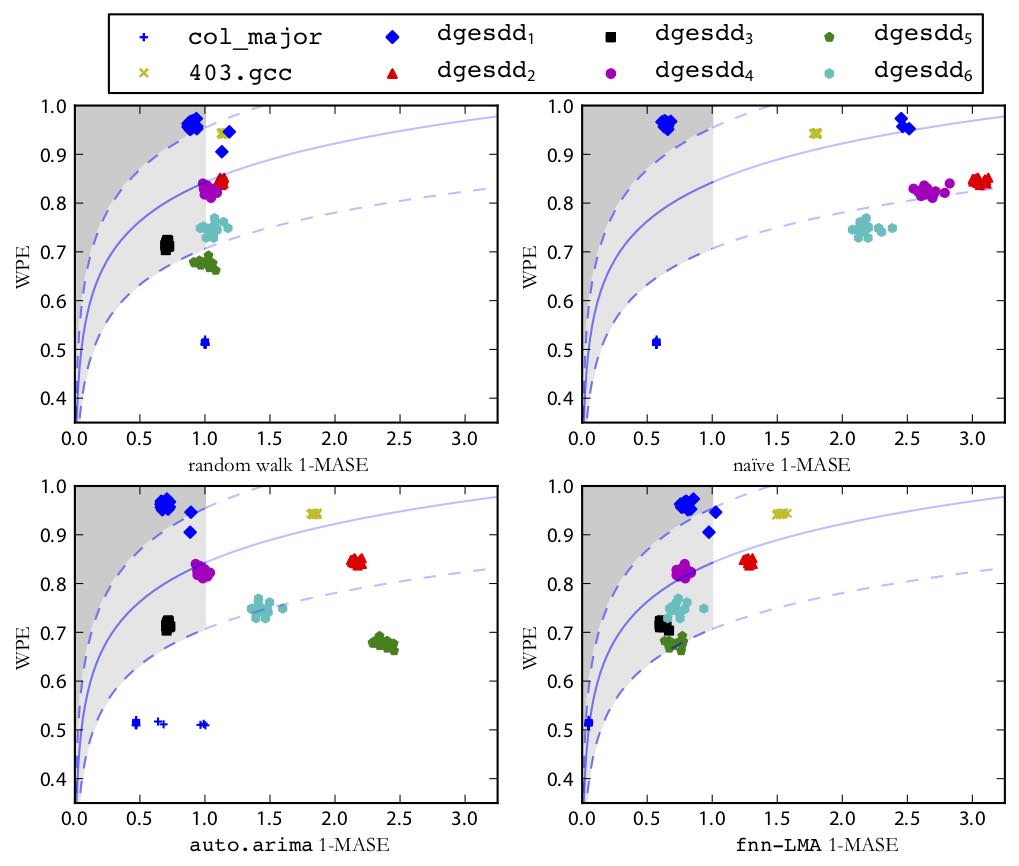}
\caption{ WPE vs. 1-MASE for all trials, methods, and systems---with the
  exception of \svdone, \svdthree, and \svdfive, which are omitted from the
  top-right plot for scale reasons, as described in the text.
%
%
Numerical values, including means and standard deviations of the
errors, can be found in Table~\ref{tab:wpeerror}.  The curves and shaded
regions are the same as in the previous figure.  }
    \label{fig:wpe_vs_mase_all}
\end{figure}

Figure~\ref{fig:wpe_vs_mase_all} shows WPE vs. 1-MASE plots for the full
set of experiments; Table~\ref{tab:wpeerror} contains all the associated numerical values.
There are 15 points in each cluster, one for each trial.  (The points
in Figure~\ref{fig:wpe_vs_mase_best} are the leftmost of the points
for the corresponding trace in any of the four plots in
Figure~\ref{fig:wpe_vs_mase_all}.)  The WPE values do not vary very
much across trials.  For most traces, the variance in 1-MASE scores is
low as well, resulting in small, tight clusters.  In some
cases---\arima predictions of {\tt col\_major}, for instance---the
1-MASE variance is larger, which spreads out the clusters horizontally.
The mean 1-MASE scores of predictions generated with \fnnLMA are generally closer to the dashed curve; the \arima method
clusters are more widely spread, the \naive ~clusters even more so.  A
few of the clusters have very high variance; these are discussed later
in this section.

The main thing to note here, however, is not the details of the shapes
of the clusters, but rather their positions in the four plots:
specifically, the fact that many of them are to the right of and/or
below the dashed curve that identifies the boundary of the shaded
region.  {\it These predictions are not as good as my heuristic
  suggests they could be.}  Focusing in on any single signal makes
this clear: \fnnLMA works best for \svdsix, for instance, followed by the
random-walk prediction method, then \arima and \naive.  Again, this provides
some practical leverage: if one calculates an WPE vs. 1-MASE value that
is outside the shaded region, that suggests that the prediction
method is not well matched to the task at hand---that is, the time
series has more predictive structure than the method is able to use.
%
%
The results of \arima on \svdsix, for instance, suggest that one should try a
different method.  The position of the \fnnLMA cluster for \svdsix, on the
other hand, reflects the ability of that method to capture and exploit the
structure that is present in this signal.  WPE vs. 1-MASE values like
this, which fall in the shaded region, indicate to the
practitioner that the prediction method is well-suited to the task.
The following discussion uses a number of examples to lay out the details that underlie these
claims.
%
%

Though \col is a very simple program, its dynamics are actually quite
complicated, as discussed in Section~\ref{sec:compPerfExperiments}.  Recall from
Figure~\ref{fig:forecast-example} and Table~\ref{tab:comp-perf-error} that the
\naive, \arima, and (especially) random-walk prediction methods do not
perform very well on this signal.  The 1-MASE scores of these
predictions are $0.571 \pm 0.002$, $1.001 \pm 0.002$, and $0.599 \pm
0.211$, respectively, across all 15 trials.  That is, \naive ~and
\arima perform only $\approx 1.7$ times better than the random-walk
method, a primitive strategy that simply uses the current value as the
prediction.  However, the WPE value for the \col trials is $0.513 \pm
0.003$, which is in the center of the complexity spectrum
described on page~\pageref{page:spectrum}.

This disparity---WPE values that suggest a high rate of forward
information transfer in the signal, but predictions with comparatively
poor 1-MASE scores---is obvious in the geometry of three of the four
images in Figure~\ref{fig:wpe_vs_mase_all}, where the \col clusters
are far to the right of and/or below the dashed curve.  Again, this
indicates that these methods are not leveraging the available
information in the signal.  The dynamics of \col may be complicated,
but they are not unstructured.  This signal is nonlinear and
deterministic~\cite{mytkowicz09}, and if one uses a prediction
technique that is based a nonlinear model (\fnnLMA)---rather than a method
that simply predicts the running mean (\naive) or the previous value
(random walk), or one that uses a linear model (\arima)---the 1-MASE
score is much improved: $0.050 \pm 0.001$.  This prediction is 20
times more accurate than a random-walk forecast, which is more in line
with the level of predictive structure that the low WPE value suggests
is present in the \col signal.  The 1-MASE scores of random-walk
predictions of this signal are all $\approx 1$---as one would
expect---pushing those points well below the shaded region.  Clearly
the stationarity assumption on which that method is based does not
hold for this signal.

The \col example also brings out some of the shortcomings of automated
model-building processes.  Note that the {\color{blue}$+$} points are
clustered very tightly in the lower left quadrant of the \naive,
random-walk, and \fnnLMA plots in Figure~\ref{fig:wpe_vs_mase_all}, but
spread out horizontally in the \arima plot.  This is because of the
way the {\tt auto.arima} process---the fitting procedure that I use for ARIMA models---works \cite{autoARIMA}.  If a
KPSS test\footnote{A Kwiatkowski-Phillips-Schmidt-Shin (KPSS) test~\cite{KPSSunit}---one of many tests performed by {\tt auto.arima} to chose the ARIMA parameters---is used for testing that an observable time series is stationary around a deterministic trend.}of the time series in question indicates that it is
nonstationary, the \arima recipe adds an integration term to the
model.  This test gives mixed results in the case of the \col process,
flagging five of the 15 trials as stationary and ten as nonstationary.
I conjectured that ARIMA models without an integration term perform
more poorly on these five signals, which increases the error and
thereby spreads out the points.  I tested this hypothesis by forcing
the inclusion of an integration term in the five cases where a KPSS
test indicated that such a term was not needed.  This action removes 
the spread, pushing all 15 of the \col ~ \arima points in
Figure~\ref{fig:wpe_vs_mase_all} into a tight cluster.

The discussion in the previous paragraph highlights the second finding
of this chapter: the ability of the graphical heuristic of Figure~\ref{fig:wpe_vs_mase_best} to flag
inappropriate models.  {\tt auto.arima} is an automated, mechanical procedure
for choosing modeling parameters for a given data set.  While
the tests and criteria employed by this algorithm
(Section~\ref{sec:arima}) are sophisticated, the results can still be
sub-optimal---if the initial space of models being searched is not
broad enough, for instance, or if one of the preliminary tests gives
an erroneous result.  Moreover, {\tt auto.arima} {\it always} returns a model, and it can
be very hard to detect when that model is bad.  The results discussed in this chapter suggest a
way to do so: if the 1-MASE score of an auto-fitted model like an {\tt auto.arima} result is
out of line with the WPE value of the data, that can be an indication
of inappropriateness in the order selection and parameter estimation
procedure.

The WPE of \svdfive ($0.677 \pm 0.006$) is higher than that of {\tt
  col\_major}.  This indicates that the rate of forward information
transfer of the underlying process is lower, but that 
observations from this system still contain a significant amount of
structure that can, in theory, be used to predict the future course of
the time series.
The 1-MASE scores of the \naive ~and \arima predictions for this system
are $20.870 \pm 0.192$ and $2.370 \pm 0.051$, respectively: that is,
20.87 and 2.37 times worse than a simple random-walk forecast\footnote{The \naive ~1-MASE
  score is large because of the bimodal nature of the distribution of
  the values of the signal, which makes guessing the mean a
  particularly bad strategy.  The same thing is true of the \svdthree
  signal.}of the
training set portions of the same signals. As before, the positions of these points on a WPE
vs. 1-MASE plot---significantly below and to the right of the shaded
region---should suggest to a practitioner that
a different method might do better.  Indeed, for \svdfive, \fnnLMA produces a 1-MASE score of $ 0.718\pm 0.048 $ and a cluster of
results that largely within the shaded region on the WPE-MASE plot.
This is consistent with the second finding of this chapter: the \fnnLMA method can capture
and reproduce the way in which the \svdfive system processes
information, but the \naive ~and \arima prediction methods cannot.

The WPE of \gcc is higher still: $0.943 \pm 0.001$.  This system
transmits very little information forward in time and provides almost
no structure for prediction methods to work with.  Here, the
random-walk predictor is the best of the methods used here.  This
makes sense; in a fully complex signal, where there is no predictive
structure to utilize, methods that depend on exploiting that
structure---like ARIMA and \fnnLMA---cannot get any traction.
%
%
Since fitting a hyperplane using least squares should filter out some
of the noise in the signal, the fact that \fnnLMA outperforms \arima
($1.530 \pm 0.021$ vs. $1.837 \pm 0.016$) may be somewhat
counterintuitive.  However, the small amount of predictive structure
that is present in this signal is nonlinear ({\it cf.,}~\cite{mytkowicz09}),
and \fnnLMA is designed to capture and exploit that kind of structure.
Note that all four \gcc clusters in Figure~\ref{fig:wpe_vs_mase_all}
are outside the shaded region; in the case of the random-walk
prediction, for instance, the 1-MASE value is $1.1381 \pm 0.011$.  This
is due to nonstationarity in the signal: in particular, differences
between the training and test signals.  The same effect is at work in
the \svdtwo results, for the same reasons---and visibly so, judging by the red
segment of Figure~\ref{fig:svd-ts-colored}, where the period and
amplitude of the oscillations are decreasing.

\svdone---the dark blue (first) segment of
Figure~\ref{fig:svd-ts-colored}---behaves very differently than the
other seven systems in this study.  Though its weighted permutation
entropy is very high ($0.957 \pm 0.016$), three of the four prediction
methods do quite well on this signal, yielding mean 1-MASE scores of
$0.714 \pm 0.075$ (\arima), $0.827 \pm 0.076$ (\fnnLMA), and $0.933 \pm 0.095$ (random walk).  This pushes the
corresponding clusters of points in Figure~\ref{fig:wpe_vs_mase_all}
well above the trend followed by the other seven signals.  The reasons
for this are discussed below.  The 1-MASE scores of
the predictions that are produced by the \naive ~method for this
system, however, are highly inconsistent.  The majority of the blue
diamond-shaped points on the top-right plot in
Figure~\ref{fig:wpe_vs_mase_all} are clustered near a 1-MASE score of
0.6, which is better than the other three methods.  In five of the 15
\svdone trials, however, there are step changes in the signal.  This
is a different nonstationarity than in the case of \col---large jump
discontinuities rather than small shifts in the baseline---and not one
that I am able to handle by simply forcing the ARIMA model to
include a particular term.  The \naive ~method has a very difficult
time with signals like this, particularly if there are multiple step
changes.  That raised the 1-MASE scores of these trials, pushing the
corresponding points\footnote{This includes the cluster
  of three points near 1-MASE $\approx 2.5$, as well as two points that
  are beyond the domain of the graph, at 1-MASE $\approx 11.2-14.8$.}to the right, and in turn raising both the mean and variance of this set of trials.

\begin{figure}[tb!]
\vspace{0.5cm}
  \centering
    \includegraphics[width=\textwidth]{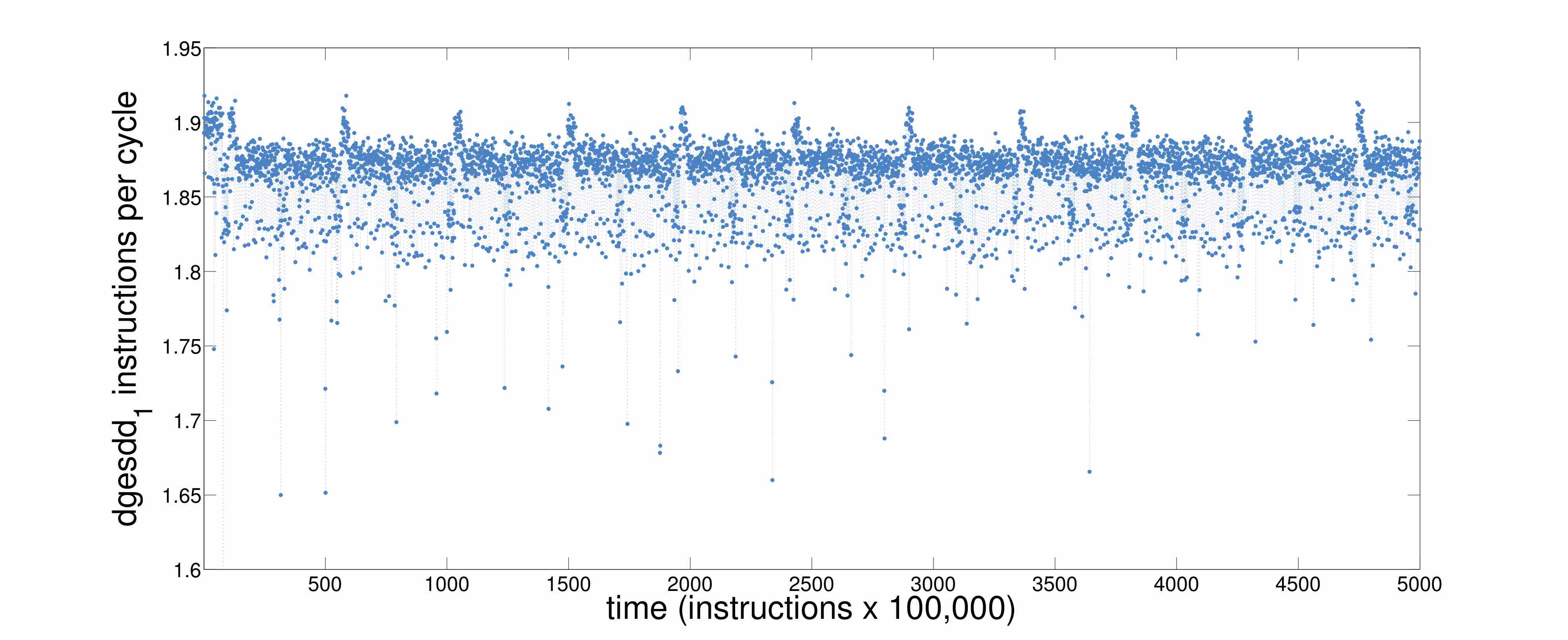}
\caption{ A small portion of the \svdone time series}\label{fig:svdone-ts}
\end{figure}

The effects described in the previous paragraph are also exacerbated
by the way 1-MASE is calculated.  Recall that 1-MASE scores are scaled
{\it relative to a random-walk forecast.} This creates several issues. 
Since random-walk prediction works very badly on signals with
frequent, large, rapid transitions (even simple periodic signals) this class of signals can exhibit low WPE, and high 1-MASE. This is because the
random-walk forecast of this signal will be 180 degrees out of phase
with the true continuation.  This effect can shift points leftwards
on a WPE vs. 1-MASE plot, and that is exactly why the \svdone clusters
in Figure~\ref{fig:wpe_vs_mase_all} are above the dashed curve.  This
time series, part of which is shown in closeup in
Figure~\ref{fig:svdone-ts},
is not quite the worst-case signal for a random-walk prediction,
but it still poses a serious challenge.  It
is dominated by a noisy regime (between $\approx$1.86 and
$\approx$1.88 on the vertical scale in Figure~\ref{fig:svdone-ts}),
punctuated by short excursions above 1.9.  In the former regime, which
makes up more than 80\% of the signal, there are frequent dips to 1.82
and occasional larger dips below 1.8.  These single-point dips are the
bane of random-walk forecasting.  In this particular case, roughly
40\% of the forecasted points are off by the width of the associated
dip, which skews the associated 1-MASE scores.  Signals like this are
also problematic for the \naive ~prediction strategy, since the
outliers have significant influence on the mean.  This compounds the
effect of the skew in the scaling factor and exacerbates the spread in
the \svdone 1-MASE values.

The second effect that can skew 1-MASE scores is nonstationarity.  Since
this metric is normalized by the error of a random-walk forecast
{\it on the training signal}, differences between the test signal and
training signal can create issues.  This is why the 1-MASE values in
Table~\ref{tab:wpeerror} are not identically one for every random-walk
forecast of every time series: the last 10\% of these signals is
significantly different from the first 90\%.  The deviation from 1.00
will depend on the process that generated the data---whether it has multiple regimes, what
those regimes look like, and how it switches between them---as well as
the experimental setup ({\it e.g.,} sensor precision and data length).  For
the processes studied here, these effects do not cause the 1-MASE values
to exceed 1.15, but pathological situations ({\it e.g.,} a huge switch in
scale right at the training/test signal boundary, or a signal that
simply grows exponentially) could produce higher values.  This
suggests another potentially useful heuristic: if the 1-MASE of a
random-walk prediction of a time series is significantly different
from 1, it could be an indication that the signal is nonstationary.

The curves in Figures~\ref{fig:wpe_vs_mase_best}
and~\ref{fig:wpe_vs_mase_all} are determined from finite sets of
methods and data.  I put a lot of thought and effort into making
these sets representative and comprehensive.  The forecast methods
involved ranged from the simple to the sophisticated; the time-series
data analyzed in this section are sampled from systems whose behavior
spans the dynamical behavior space.  While I am cautiously
optimistic about the generality of my conclusions, more exploration
will be required before I can make definitive or general conclusions.  However, my preliminary exploration  
 shows that data from the H\'enon map~\cite{henon}, the Lorenz 63 system \cite{lorenz}, the SFI dataset A~\cite{weigend-book}, and a random-walk process all fall within the
one-$\sigma$ volume of the fit in Figures~\ref{fig:wpe_vs_mase_best}
and~\ref{fig:wpe_vs_mase_all} region, as do various nonlinear
transformations of \svdtwo, \svdfive and \svdsix, so I am optimistic.

In this chapter I strictly used \fnnLMA as
the nonlinear exemplar to explore how prediction accuracy is related
to WPE. With this relationship established, applying this heuristic to assessing the appropriateness of \roLMA  for a given signal is quite trivial: one simply compares the accuracy of \roLMA, tabulated in Tables~\ref{tab:lorenz96error} and~\ref{tab:comp-perf-error}, with the WPE of that signal, given in Table~\ref{tab:wpeerror}, using the graphical heuristic in Figure~\ref{fig:wpe_vs_mase_best}.

Note that there has been prior work under a very similar title to
our paper on this topic~\cite{josh-pre}, but there are only superficial similarities
between the two research projects. Haven {\it et al.}~\cite{haven2005} utilize the
relative entropy to quantify the difference in predictability between
two distributions: one evolved from a small ensemble of past states
using the known dynamical system, and the other the observed
distribution. My work quantifies the predictability of a single
observed time series using weighted permutation entropy and makes no
assumptions about the generating process.

More closely related is the work of Boffetta {\it et
  al.}~\cite{boffetta02}, who investigated the scaling behavior of
finite-size Lyapunov exponents (FSLE) and $\epsilon$-entropy for a
wide variety of deterministic systems with known dynamics and additive
noise.  While the scaling of these measures acts as a general proxy
for predictability bounds, this approach differs from my work in a
number of fundamental ways.  First, \cite{boffetta02} is a theoretical
study that does not involve any actual predictions.  I focus on
real-world time-series data, where one does not necessarily have the
ability to perturb or otherwise interact with the system of interest,
nor can one obtain or manufacture the (possibly large) number of
points that might be needed to estimate the $\epsilon$-entropy for
small $\epsilon$.  Second, I do not require {\it a priori} knowledge
about the noise and its interaction with the system.  Third, I tie
information---in the form of the weighted permutation
entropy---directly to prediction error via calculated values of a
specific error metric.  Though FSLE and $\epsilon$-entropy allow for
the comparison of predictability between systems, they do not directly
provide an estimate of prediction error.  Finally, my approach also
holds for stochastic systems, where neither the FLSEs nor their
relationship to predictability are well defined.

\section{Summary}\label{sec:conc}

Forecast strategies that are designed to capture predictive structure
are ineffective when signal complexity outweighs information
redundancy.  This poses a number of serious challenges in practice.
Without knowing anything about the generating process, it is difficult
to determine how much predictive structure is present in a noisy,
real-world time series.  And even if predictive structure exists, a
given forecast method may not work, simply because it cannot exploit
the structure that is present ({\it e.g.,} a linear model of a nonlinear
process).  If a forecast model is not producing good results, a
practitioner needs to know why: is the reason that the data contain no
predictive structure---{\it i.e.,} that no model will work---or is the model
that s/he is using simply not good enough?

In this chapter, I have argued that redundancy is a useful proxy for
the inherent predictability of an empirical time series.  To
operationalize that relationship, I used an approximation of the
Kolmogorov-Sinai entropy, estimated using a weighted version of the
permutation entropy of~\cite{bandt2002per}.  This WPE technique---an
ordinal calculation of forward information transfer in a time
series---is ideal for my purposes because it works with real-valued
data and is known to converge to the true entropy value. Using a
variety of forecast models and more than 150 time-series data sets
from experiments and simulations, I have shown that prediction
accuracy is indeed correlated with weighted permutation entropy: the
higher the WPE, in general, the higher the prediction error.  The
relationship is roughly logarithmic, which makes theoretical sense,
given the nature of WPE, predictability, and 1-MASE.

An important practical corollary to this empirical correlation of
predictability and WPE is a practical strategy for assessing
appropriateness of forecast methods.  If the forecast produced by a
particular method is poor but the time series contains a significant
amount of predictive structure, one can reasonably conclude that that
method is inadequate to the task and that one should seek another
method.  \fnnLMA, for instance, performed better in
most cases because it is more general.  (This is particularly apparent
in the \col and \svdfive examples.)
The \naive ~method, which simply predicts the mean, can work very well
on noisy signals because it effects a filtering operation.  The simple
random-walk strategy outperforms \fnnLMA, \arima, and the \naive
~method on the \gcc signal, which is extremely complex---{\it i.e.,}
extremely low redundancy.

The curves and shaded regions in Figures~\ref{fig:wpe_vs_mase_best}
and~\ref{fig:wpe_vs_mase_all} operationalize the
discussion in the previous paragraph.  These geometric features are a
preliminary, but potentially useful, heuristic for knowing when a
model is not well-matched to the task at hand: a point that is below
and/or to the right of the shaded regions on a plot like
Figure~\ref{fig:wpe_vs_mase_all} indicates that the time series has
more predictive structure than the forecast model can capture and
exploit---and that one would be well advised to try another method. In the 
context of this thesis, this heuristic allows a practitioner to know if \roLMA is 
capturing all the available predictive structure in a time series or if another method such as \fnnLMA 
should be used instead.

These curves were determined empirically using a specific error metric
and a finite set of forecast methods and time-series traces.  If one
uses a different error metric, the geometry of the heuristic may be
different---and may not even make sense, if one uses a metric that
does not support comparison across different time series.  And while
the methods and traces used in this study were chosen to be
representative of the practice, they are of course not completely
comprehensive.  It is certainly possible, for instance, that the
nonlinear dynamics of computer performance is subtly different from
the nonlinear dynamics of other systems.  My preliminary results on
other systems, not shown here, {\it e.g.,} H\'enon, Lorenz 63, a random-walk process, SFI dataset A, and more, lead me to believe
that the results of this chapter will generalize beyond the examples presented.

%% file: concl.tex

\chapter{Conclusion and Future Directions}\label{ch:concandfuture}

Delay-coordinate embedding, the bedrock of nonlinear time-series analysis, has been the foundation of forecasting techniques for nonlinear dynamical systems.   
A significant hurdle in the application of these techniques in a real-time fashion or for nonstationary systems is proper estimation of the embedding dimension. 
The source of this difficulty is rooted in trying to obtain a diffeomorphic reconstruction of the observed system, {\it i.e.,} a topologically perfect representation. This thesis presented a paradigm shift away from that traditional approach, showing that for short-term delay-coordinate based forecasting of limited noisy data, {\it perfection is not  necessary, and can even be detrimental}. 

As a first step along this path, I introduced a novel forecasting schema, \roLMA, that sidesteps the difficult parameter estimation step by simply fixing $m=2$.  For a range of low- and high-dimensional synthetic and experimental systems, I showed that \roLMA produced short-term predictions on par with {\it or exceeding} the accuracy of traditional embeddings even though \roLMA employed {\it incomplete}
reconstructions of the dynamics---{\it i.e,} models that are not necessarily true embeddings.  This effected an experimental validation of the central premise of this thesis, {\it viz.,} the current paradigm in delay-coordinate embedding may be overly stringent for short-term forecasts of real-world data.

The utility of incomplete reconstructions is in stark contrast to traditional views of delay-coordinate embedding, so experimental validation of this bold claim is insufficient. In order to present a complete validation of this methodology, I provided two 
complementary theoretical frameworks, based in information theory and computational topology, to explain {\it why} and {\it how} this reduced-order modeling strategy works. 

The information-theoretic analysis focused on understanding how information about the future is stored in delay-coordinate vectors. Leveraging this knowledge, in collaboration with R.~G.~James, I constructed a novel metric, time-delayed active information storage (\mytau), for selecting forecast-specific reconstruction parameters. This approach to parameter selection is drastically different than standard approaches. Instead of focusing on the calculation of dynamical invariants as the end goal, \mytau maximizes the information about the future stored in each delay vector---explicitly optimizing the reconstruction for the purposes of forecasting. \mytau allows one to select parameter values that are tailored specifically to the quantity of data available, the signal-to-noise ratio of the time series and the required forecast horizon---and does so quickly, efficiently and directly from the data. \mytau independently corroborated the central claims of this thesis, showing that for noisy, limited datasets, often the state estimator used in \roLMA contained as much---or more---information about the near future as a full embedding. This result is counter-intuitive; one would think, up to some limit, each new dimension would add information to the model.  The fact that more information is not necessarily gained in each new dimension changes the way delay-coordinate based forecasting should be approached, and offers a fundamental explanation of why prediction in projection is effective in practice. 

The topological analysis in Section~\ref{sec:compTopo} questioned the fundamental need for a {\it diffeomorphic} reconstruction---especially when one is not interested in calculating dynamical invariants. In collaboration with J.~D.~Meiss, I conjectured that
it may be possible to ascertain information about the large-scale topology of the invariant set---specifically, the homology}---with a lower reconstruction dimension than that
needed to obtain an embedding.  
Using a simple canonical example, I showed that the witness complex correctly resolved the homology of the underlying 
invariant set, {\it viz.,} its Betti numbers, even if the reconstruction dimension was well below the thresholds for
which the embedding theorems assure smooth conjugacy between the true
and reconstructed dynamics. Since many
properties that one cares about for forecasting---the existence of periodic orbits, recurrence, entropy, etc.---depend only upon topology, the stabilization of large-scale topology at low dimensions effects an alternative validation of the central premise of this thesis. I further conjectured---but did not prove---that this unexpected resolution of large-scale topology at low dimensions may be due to the existence of a {\it homeomorphism} between the original and reconstructed dynamics.  Proving the broader claim that the delay-coordinate map is a homeomorphism at lower embedding dimensions than required for a diffeomorphism and that a homeomorphic reconstruction is sufficient for short-term forecasting will be a key future direction of research.

While I illustrated that \roLMA works for a broad spectrum of signals and provided sound theoretical evidence supporting why it works, it is important to remember that \roLMA---or any forecast model for that matter---will not be ideal for all tasks. Following this line of reasoning it was vital to understand {\it when} \roLMA  
would be effective. In this capacity, again in collaboration with R.~G.~James, I developed a model-free framework for quantifying when a time-series exhibits predictable features, {\it viz.,} bounded information production. This heuristic allows one to know {\it a priori} whether \roLMA---or any forecast method---is appropriate for forecasting a given time series.

There are a number of important avenues for future work associated with each topic proposed in this thesis; those avenues were all discussed individually in their respective chapters. However, the next frontier of this work, which draws upon all aspects of the research described in this dissertation, is  developing strategies for grappling with nonstationary time series---a serious challenge in any time-series modeling problem---in the context of delay coordinate based forecasting. 
We live in a nonlinear {\it and} nonstationary world. Real-world systems 
change over time: bearings 
break down, computer systems get updated, transistors wear out, the climate moves between glacial and interglacial periods, etc. 
The nonlinear and nonstationary nature of systems like these highlights the need for adaptive models, built on the fly, that require little to no human interaction.

\begin{figure}[tb!]
  \centering
  \includegraphics[width=0.65\textwidth]{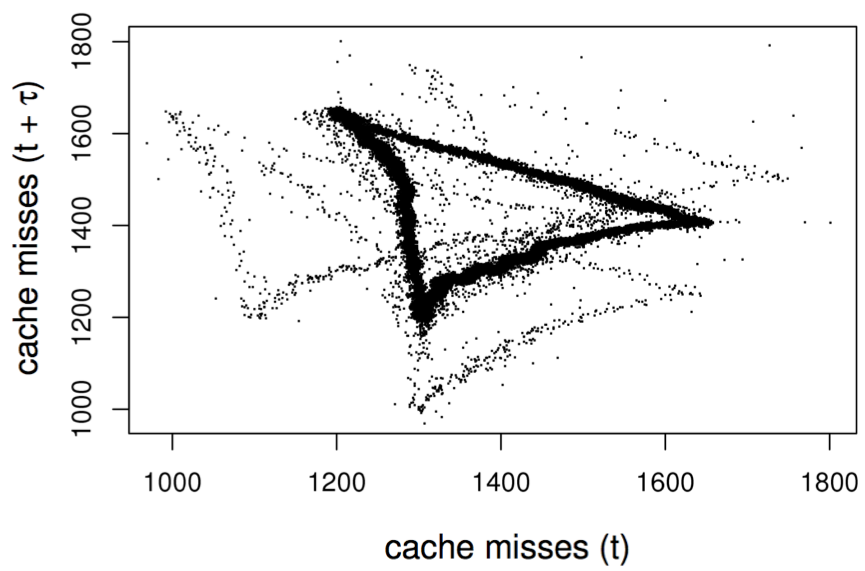}
  \caption{A two-dimensional reconstruction of L2 cache-misses of \row on an Intel Core Duo 
  with $\tau=100,000$ instructions.}
  \label{fig:ghost}
  \vspace{-0.75cm}
\end{figure}

\begin{figure}[bt!]
  \centering
  \vspace{1cm}
  \includegraphics[width=\textwidth]{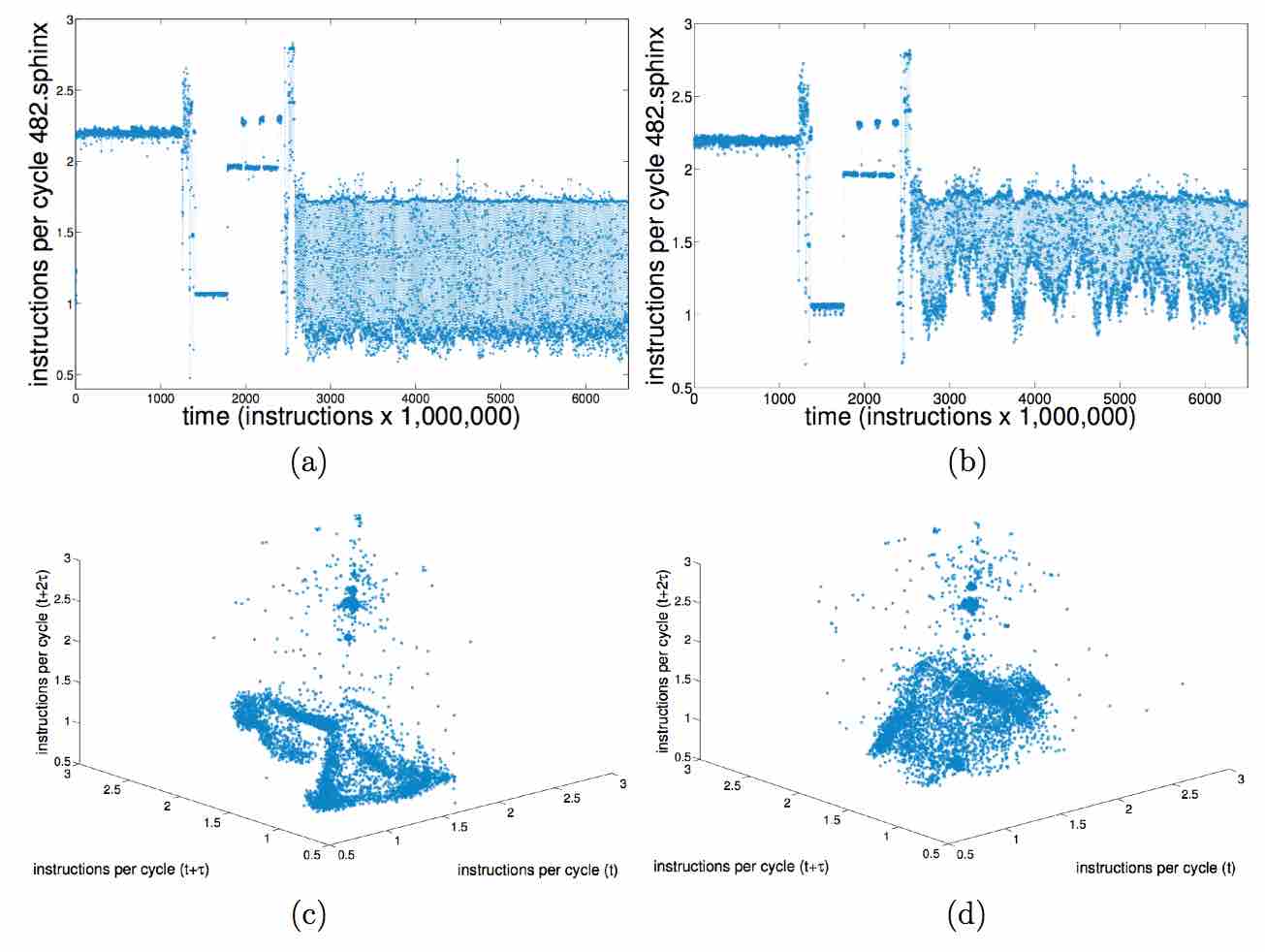}
  \caption{Processor load (IPC) traces [Top] and 3D projections of the respective reconstructed dynamics [Bottom] of \sphinx before [(a),(c)] and after [(b,d)] a kernel {\tt panic()} halt. }
  \label{fig:halt}
  \vspace{-0.75cm}
\end{figure}

My first experience with nonstationarity involved a performance trace of \row running on an Intel Core Duo that exhibited an interesting phenomenon termed ``ghost triangles" (as can be seen in Figure \ref{fig:ghost}) \cite{Alexander12}. After a routine (automatic) operating system update, not only were the ``ghost triangles" gone, but the entire triangular structure present in the two-dimensional reconstruction had been lost. 
I thought I had learned my lesson the hard way with unforeseen change, but the dynamics were even more sensitive than I originally thought. As a result of the auto-update fiasco, our group purchased a new computer and disabled its update scheduler. However, I soon learned this was not enough. When I went back to repeat some experiments, the computer crashed with a standard kernel {\tt panic()} halt. After this, the traces were never the same: the system halt caused a bifurcation in the performance dynamics. Figure~\ref{fig:halt} shows a before-and-after example involving another SPEC benchmark, 482.sphinx. According to computer design theory, this halt should {\it not} have changed anything. Nonetheless, it actually caused a fundamental shift in the performance  and a change in the dynamical structure.

These examples drive home the fact that real-world nonlinear 
systems routinely undergo bifurcations in their dynamics. This means that a forecast 
model should not be trained once and then used for all time. Instead, it should be constantly adapted to the current dynamics. For the \fnnLMA method, this is next to impossible due to the human-intensive parameter selection process. The agility of \roLMA (no need to estimate m) should make it possible to tackle delay-coordinate based forecasting of nonstationary time series. Detecting regime shifts---and adapting the time delay accordingly---is the first step in this important area of future work. 
 
 To accomplish this, I plan to study whether the methods described in this thesis can signal regime shifts in a time series. 
Fuzzy witness complexes may be a useful strategy in this vein of research by detecting and characterizing bifurcations \cite{BerwaldGV13}.
Suppose that, for example, the first part of the data set corresponds
trivially to an equilibrium---that is, to a set with $\beta_k = 0$ for
all $k>0$---but that this equilibrium undergoes a bifurcation to an
oscillatory regime midway through the data set.  In this case, a shift
in $\beta_1$ signals a regime change.  

Not all regime shifts are the result of a bifurcation, however. In cases like that, a change in information mechanics, {\it e.g.,} information storage (\mytau) or information production (WPE), could signal a regime shift---even if the topology of the new regime was too similar (or identical) to  the old regime. 
My WPE vs. 1-MASE results could be particularly powerful in this scenario, as their values could not only
help with regime-shift detection, but also suggest what kind of model
might work well in each new regime. Similarly, changes in information storage could also be quite useful. A change in \mytau  suggests a regime shift has occurred and simultaneously suggests a forecast-optimal parameter set for the new regime.  Even more importantly, it indicates whether new parameter values are even necessary in the new regime! 

Of particular interest in this new frontier of research will be
the class of so-called {\it hybrid systems}~\cite{hybrid}, which
exhibit discrete transitions between different continuous
regimes---{\it e.g.,} a lathe that has an intermittent instability or
traffic at an internet router, whose characteristic normal traffic
patterns shift radically during an attack.  Traded financial markets, too, are highly susceptible to jump processes.  Effective modeling and
prediction of these kinds of systems is quite difficult; doing so
adaptively and automatically is an important and interesting challenge.

The elements of this thesis could be combined to form a complete nonstationary forecasting framework that could address that challenge. The method  would work as follows: use \mytau to select forecast-optimal parameters and start forecasting using \roLMA. While forecasting on a small buffer of new data, monitor information production (WPE), information propagation and storage (\mytau) and the homology (witness complex). If any of these drastically change, a regime shift has occurred and the model should be rebuilt. Nicely, as \mytau and WPE are already being monitored on the new data buffer, the new regime is already modeled and forecasting can begin immediately.  

This thesis bridged the gap between rigorous nonlinear mathematical models---which are ineffective in real-time---and \naive methods that are agile enough for adaptive modeling of non-stationary processes. This in and of itself has real practical utility for a wide spectrum of forecasting tasks as a simple, agile, noise-resilient, forecasting strategy for nonlinear systems, but this thesis went far beyond practical optimizations at the sacrifice of theoretical rigor. Specifically, the theoretical analysis outlined here offered a deeper understanding of delay-coordinate embedding---an understanding that suggests (and justifies) the need for a new paradigm in delay-reconstruction theory that was the overall goal of this research project.

%% file: dissertation.bbl
\begin{thebibliography}{100}

\bibitem{Abdallah}
S.~A. Abdallah and M.~D. Plumbley.
\newblock A measure of statistical complexity based on predictive information
  with application to finite spin systems.
\newblock {\em Physics Letters A}, 376(4):275--281, 2012.

\bibitem{akaike1974}
H.~Akaike.
\newblock A new look at the statistical model identification.
\newblock {\em IEEE Transactions on Automatic Control}, 19(6):716--723, 1974.

\bibitem{alameldeen}
A.~Alameldeen and D.~Wood.
\newblock {IPC} considered harmful for multiprocessor workloads.
\newblock {\em IEEE Micro}, 26(4):8--17, 2006.

\bibitem{Alexander15}
Z.~Alexander, E.~Bradley, J.~D. Meiss, and N.~F. Sanderson.
\newblock Simplicial multivalued maps and the witness complex for dynamical
  analysis of time series.
\newblock {\em SIAM Journal on Applied Dynamical Systems}, 14(3):1278--1307,
  2015.

\bibitem{Alexander12}
Z.~Alexander, J.~D. Meiss, E.~Bradley, and J.~Garland.
\newblock Iterated function system models in data analysis: Detection and
  separation.
\newblock {\em Chaos: An Interdisciplinary Journal of Nonlinear Science},
  22(2):023103, 2012.

\bibitem{zach-IDA10}
Z.~Alexander, T.~Mytkowicz, A.~Diwan, and E.~Bradley.
\newblock Measurement and dynamical analysis of computer performance data.
\newblock In {\em Advances in Intelligent Data Analysis IX}, volume 6065.
  Springer Lecture Notes in Computer Science, 2010.

\bibitem{amigo2012permutation}
J.~Amig{\'o}.
\newblock {\em Permutation complexity in dynamical systems: Ordinal patterns,
  permutation entropy and all that}.
\newblock Springer, 2012.

\bibitem{lapack}
E.~Anderson, Z.~Bai, C.~Bischof, S.~Blackford, J.~Demmel, J.~Dongarra,
  J.~Du~Croz, A.~Greenbaum, S.~Hammarling, A.~McKenney, and D.~Sorensen.
\newblock {\em {LAPACK} Users' Guide}.
\newblock Society for Industrial and Applied Mathematics, Philadelphia, PA,
  third edition, 1999.

\bibitem{bandt2002per}
C.~Bandt and B.~Pompe.
\newblock Permutation entropy: A natural complexity measure for time series.
\newblock {\em Physical Review Letters}, 88(17):174102, 2002.

\bibitem{Bell03theco-information}
A.~J. Bell.
\newblock The co-information lattice.
\newblock In {\em Proceedings of 4th International Symposium on Independent
  Component Analysis and Blind Source Separation}, volume 2003, pages 921--926,
  2003.

\bibitem{BerwaldGV13}
J.~Berwald, M.~Gidea, and M.~Vejdemo-Johansson.
\newblock Automatic recognition and tagging of topologically different regimes
  in dynamical systems.
\newblock {\em Discontinuity, Nonlinearity, and Complexity}, 3(4):413--426,
  2015.

\bibitem{boffetta02}
G.~Boffetta, M.~Cencini, M.~Falcioni, and A.~Vulpiani.
\newblock Predictability: A way to characterize complexity.
\newblock {\em Physics Reports}, 356(6):367--474, 2002.

\bibitem{bollt2001}
{E.} Bollt, T.~Stanford, Y.~C. Lai, and K.~{\.Z}yczkowski.
\newblock {What symbolic dynamics do we get with a misplaced partition?: On the
  validity of threshold crossings analysis of chaotic time-series}.
\newblock {\em Physica D: Nonlinear Phenomena}, 154(3):259--286, 2001.

\bibitem{Holger-and-Liz}
E.~Bradley and H.~Kantz.
\newblock Nonlinear time-series analysis revisited.
\newblock {\em Chaos: An Interdisciplinary Journal of Nonlinear Science},
  25(9):097610, 2015.

\bibitem{davislinearts}
{P.} Brockwell and {R.} Davis.
\newblock {\em Introduction to Time Series and Forecasting}.
\newblock Springer-Verlag, New York, 2002.

\bibitem{papi}
S.~Browne, C.~Deane, G.~Ho, and P.~Mucci.
\newblock {PAPI}: A portable interface to hardware performance counters.
\newblock In {\em Proceedings of Department of Defense HPCMP Users Group
  Conference}, 1999.

\bibitem{Buzug92Comp}
Th. Buzug and G.~Pfister.
\newblock Comparison of algorithms calculating optimal embedding parameters for
  delay time coordinates.
\newblock {\em Physica D: Nonlinear Phenomena}, 58(1-4):127--137, 1992.

\bibitem{Buzugfilldeform}
Th. Buzug and G.~Pfister.
\newblock Optimal delay time and embedding dimension for delay-time coordinates
  by analysis of the global static and local dynamical behavior of strange
  attractors.
\newblock {\em Physical Review A}, 45(10):7073--7084, 1992.

\bibitem{Canova1995}
F.~Canova and {B.} Hansen.
\newblock Are seasonal patterns constant over time? a test for seasonal
  stability.
\newblock {\em Journal of Business \& Economic Statistics}, 13(3):237--252,
  1995.

\bibitem{Cao97Embed}
{L.} Cao.
\newblock Practical method for determining the minimum embedding dimension of a
  scalar time series.
\newblock {\em Physica D: Nonlinear Phenomena}, 110(1-2):43--50, 1997.

\bibitem{Carlsson14}
G.~Carlsson.
\newblock Topological pattern recognition for point cloud data.
\newblock {\em Acta Numerica}, 23:289--368, 2014.

\bibitem{Casdagli92dvsplots}
M.~Casdagli.
\newblock Chaos and deterministic versus stochastic non-linear modelling.
\newblock {\em Journal of the Royal Statistical Society, Series B},
  54:303--328, 1992.

\bibitem{casdagli-eubank92}
M.~Casdagli and S.~Eubank.
\newblock {\em Nonlinear Modeling and Forecasting}.
\newblock Addison Wesley, 1992.

\bibitem{Casdagli:1991a}
M.~Casdagli, S.~Eubank, {J. D.} Farmer, and {J. F.} Gibson.
\newblock State space reconstruction in the presence of noise.
\newblock {\em Physica D: Nonlinear Phenomena}, 51(1-3):52--98, 1991.

\bibitem{crutchfield2003}
J.~P. Crutchfield and D.~P. Feldman.
\newblock Regularities unseen, randomness observed: Levels of entropy
  convergence.
\newblock {\em Chaos: An Interdisciplinary Journal of Nonlinear Science},
  13(1):25--54, 2003.

\bibitem{PhysRevE.61.1353}
R.~L. Davidchack, Y.~C. Lai, E.~M. Bollt, and M.~Dhamala.
\newblock Estimating generating partitions of chaotic systems by unstable
  periodic orbits.
\newblock {\em Physical Review E}, 61(2):1353--1356, 2000.

\bibitem{deSilva04}
V.~de~Silva and E.~Carlsson.
\newblock Topological estimation using witness complexes.
\newblock In {\em Eurographics Symposium on Point-Based Graphics (2004)}, pages
  157--166. The Eurographics Association, Zurich, 2004.

\bibitem{simon-it}
S.~DeDeo.
\newblock Information theory for intelligent people.
\newblock Available at http://tuvalu.santafe.edu/~simon/it.pdf., 2015.

\bibitem{Dowker52}
C.~H. Dowker.
\newblock Homology groups of relations.
\newblock {\em Annals of Mathematics}, 56(1):84--95, 1952.

\bibitem{ELZ01}
H.~Edelsbrunner, D.~Letscher, and A.~Zomorodian.
\newblock Topological persistence and simplification.
\newblock In {\em {IEEE} Symposium on Foundations of Computer Science}, pages
  454--463, 2000.

\bibitem{eisele1999}
M.~Eisele.
\newblock Comparison of several generating partitions of the h\'{e}non map.
\newblock {\em Journal of Physics A}, 32(9):1533--1545, 1999.

\bibitem{fadlallah2013}
B.~Fadlallah, B.~Chen, A.~Keil, and J.~Pr{\'\i}ncipe.
\newblock Weighted-permutation entropy: A complexity measure for time series
  incorporating amplitude information.
\newblock {\em Physical Review E}, 87(2):022911, 2013.

\bibitem{fraser-swinney}
A.~Fraser and H.~Swinney.
\newblock Independent coordinates for strange attractors from mutual
  information.
\newblock {\em Physical Review A}, 33(2):1134--1140, 1986.

\bibitem{PhysRevLett.99.204101}
S.~Frenzel and B.~Pompe.
\newblock Partial mutual information for coupling analysis of multivariate time
  series.
\newblock {\em Physical Review Letters}, 99(20):204101, 2007.

\bibitem{my-masters}
J.~Garland.
\newblock {\em Prediction in projection: Computer performance forecasting, a
  dynamical systems approach}.
\newblock {M. S.} {T}hesis, Department of Applied Mathematics, University of
  Colorado at Boulder, 2011.

\bibitem{josh-IDA11}
J.~Garland and E.~Bradley.
\newblock Predicting computer performance dynamics.
\newblock In {\em Advances in Intelligent Data Analysis X}, volume 7014.
  Springer Lecture Notes in Computer Science, 2011.

\bibitem{josh-IDA13}
J.~Garland and E.~Bradley.
\newblock On the importance of nonlinear modeling in computer performance
  prediction.
\newblock In {\em Advances in Intelligent Data Analysis XII}, volume 8207,
  pages 210--222. Springer Lecture Notes in Computer Science, 2013.

\bibitem{joshua-pnp}
J.~Garland and E.~Bradley.
\newblock Prediction in projection.
\newblock {\em Chaos: An Interdisciplinary Journal of Nonlinear Science},
  25:123108, 2015.

\bibitem{josh-physicaD}
J.~Garland, E.~Bradley, and J.~D. Meiss.
\newblock Exploring the topology of dynamical reconstructions.
\newblock {\em Physica D: Nonlinear Phenomena}, 334:49--59, 2016.

\bibitem{ISIT13}
J.~Garland, R.~G. James, and E.~Bradley.
\newblock Determinism, complexity, and predictability of computer performance.
\newblock arXiv:1305.5408, 2013.

\bibitem{josh-pre}
J.~Garland, R.~G. James, and E.~Bradley.
\newblock Model-free quantification of time-series predictability.
\newblock {\em Physical Review E}, 90(5):052910, 2014.

\bibitem{josh-tdAIS}
J.~Garland, R.~G. James, and E.~Bradley.
\newblock Leveraging information storage to select forecast-optimal parameters
  for delay-coordinate reconstructions.
\newblock {\em Physical Review E}, 93(2):022221, 2016.

\bibitem{georges07}
A.~Georges, D.~Buytaert, and L.~Eeckhout.
\newblock Statistically rigorous {Java} performance evaluation.
\newblock In {\em Proceedings of the ACM SIGPLAN Conference on Object-Oriented
  Programming, Systems, Languages, and Applications (OOPSLA)}, pages 57--76,
  2007.

\bibitem{weigend93}
N.~Gershenfeld and A.~Weigend.
\newblock The future of time series.
\newblock In {\em Time Series Prediction: Forecasting the Future and
  Understanding the Past}. Santa Fe Institute Studies in the Sciences of
  Complexity, Santa Fe, NM, 1993.

\bibitem{Ghrist08}
R.~Ghrist.
\newblock Barcodes: {T}he persistent topology of data.
\newblock {\em Bulletin of the American Mathematical Society}, 45(1):61--75,
  2008.

\bibitem{Gibson92}
{J.} Gibson, {J.} Farmer, M.~Casdagli, and S.~Eubank.
\newblock An analytic approach to practical state space reconstruction.
\newblock {\em Physica D: Nonlinear Phenomena}, 57(1-2):1--30, 1992.

\bibitem{hybrid}
R.~Goebel, R.~G. Sanfelice, and A~Teel.
\newblock Hybrid dynamical systems.
\newblock {\em IEEE Control Systems Magazine}, 29(2):28--93, 2009.

\bibitem{ghkss}
P.~Grassberger, R.~Hegger, H.~Kantz, C.~Schaffrath, and T.~Schreiber.
\newblock On noise reduction methods for chaotic data.
\newblock {\em Chaos: An Interdisciplinary Journal of Nonlinear Science},
  3(2):127--141, 1993.

\bibitem{GrassbergerPhysicaD}
P.~{Grassberger} and I.~{Procaccia}.
\newblock {Measuring the strangeness of strange attractors}.
\newblock {\em Physica D: Nonlinear Phenomena}, 9(1-2):189--208, 1983.

\bibitem{Han1978133}
T.~S. Han.
\newblock Nonnegative entropy measures of multivariate symmetric correlations.
\newblock {\em Information and Control}, 36(2):133--156, 1978.

\bibitem{hasson2008influence}
C.~Hasson, R.~Van~Emmerik, G.~Caldwell, J.~Haddad, J.~Gagnon, and J.~Hamill.
\newblock Influence of embedding parameters and noise in center of pressure
  recurrence quantification analysis.
\newblock {\em Gait \& Posture}, 27(3):416--422, 2008.

\bibitem{haven2005}
K.~Haven, A.~Majda, and R.~Abramov.
\newblock Quantifying predictability through information theory: Small sample
  estimation in a non-{G}aussian framework.
\newblock {\em Journal of Computational Physics}, 206(1):334--362, 2005.

\bibitem{Hegger:1999yq}
R.~Hegger, H.~Kantz, and T.~Schreiber.
\newblock Practical implementation of nonlinear time series methods: The
  {TISEAN} package.
\newblock {\em Chaos: An Interdisciplinary Journal of Nonlinear Science},
  9(2):413--435, 1999.

\bibitem{spec2006}
J.~Henning.
\newblock {SPEC CPU2006} benchmark descriptions.
\newblock {\em SIGARCH Computer Architecture News}, 34(4):1--17, 2006.

\bibitem{henon}
M.~H\'{e}non.
\newblock A two-dimensional mapping with a strange attractor.
\newblock {\em Communications in Mathematical Physics}, 50(1):69--77, 1976.

\bibitem{autoARIMA}
R.~Hyndman and Y.~Khandakar.
\newblock Automatic time series forecasting: The forecast package for {R}.
\newblock {\em Journal of Statistical Software}, 27(3):1--22, 2008.

\bibitem{MASE}
{R.} Hyndman and {A.} Koehler.
\newblock Another look at measures of forecast accuracy.
\newblock {\em International Journal of Forecasting}, 22(4):679--688, 2006.

\bibitem{james2014many}
R.~G. James, J.~R. Mahoney, C.~J. Ellison, and J.~P. Crutchfield.
\newblock Many roads to synchrony: Natural time scales and their algorithms.
\newblock {\em Physical Review E}, 89(4):042135, 2014.

\bibitem{kantz97}
H.~Kantz and T.~Schreiber.
\newblock {\em Nonlinear Time Series Analysis}.
\newblock Cambridge University Press, Cambridge, 1997.

\bibitem{kydimension}
J.~L. Kaplan and J.~A. Yorke.
\newblock Chaotic behavior of multidimensional difference equations.
\newblock In {\em Functional Differential Equations and Approximation of Fixed
  Points}, volume 730 of {\em Lecture Notes in Mathematics}, pages 204--227.
  Springer Berlin Heidelberg, 1979.

\bibitem{KarimiL96}
A.~Karimi and {M.} Paul.
\newblock {Extensive chaos in the Lorenz-96 model}.
\newblock {\em Chaos: An Interdisciplinary Journal of Nonlinear Science},
  20(4):043105, 2010.

\bibitem{KBA92}
M.~Kennel, R.~Brown, and H.~Abarbanel.
\newblock Determining minimum embedding dimension using a geometrical
  construction.
\newblock {\em Physical Review A}, 45(6):3403--3411, 1992.

\bibitem{KSG}
A.~Kraskov, H.~St{\"o}gbauer, and P.~Grassberger.
\newblock Estimating mutual information.
\newblock {\em Physical Review E}, 69(6):066138, 2004.

\bibitem{Kugi96}
D.~Kugiumtzis.
\newblock State space reconstruction parameters in the analysis of chaotic time
  series---the role of the time window length.
\newblock {\em Physica D: Nonlinear Phenomena}, 95(1):13--28, 1996.

\bibitem{KPSSunit}
D.~Kwiatkowski, {P.} Phillips, P.~Schmidt, and Y.~Shin.
\newblock {Testing the null hypothesis of stationarity against the alternative
  of a unit root: How sure are we that economic time series have a unit root?}
\newblock {\em Journal of Econometrics}, 54(1-3):159--178, 1992.

\bibitem{lebeck02}
A.~Lebeck, J.~Koppanalil, T.~Li, J.~Patwardhan, and E.~Rotenburg.
\newblock A large, fast instruction window for tolerating cache misses.
\newblock In {\em Proceedings of the International Symposium on Computer
  Architecture (ISCA)}, pages 59--70, 2002.

\bibitem{liebert-wavering}
W.~Liebert, K.~Pawelzik, and H.~Schuster.
\newblock Optimal embeddings of chaotic attractors from topological
  considerations.
\newblock {\em Europhysics Letters}, 14(6):521--526, 1991.

\bibitem{Liebert89}
W.~Liebert and H.~Schuster.
\newblock Proper choice of the time delay for the analysis of chaotic time
  series.
\newblock {\em Physics Letters A}, 142(2-3):107--111, 1989.

\bibitem{lind95}
D.~Lind and B.~Marcus.
\newblock {\em An introduction to symbolic dynamics and coding}.
\newblock Cambridge University Press, 1995.

\bibitem{jidt}
J.~T. Lizier.
\newblock {JIDT}: An information-theoretic toolkit for studying the dynamics of
  complex systems.
\newblock {\em Frontiers in Robotics and Artificial Intelligence}, 1(11):1--20,
  2014.

\bibitem{lorenz}
E.~Lorenz.
\newblock Deterministic nonperiodic flow.
\newblock {\em Journal of the Atmospheric Sciences}, 20(2):130--141, 1963.

\bibitem{lorenz-analogues}
E.~Lorenz.
\newblock Atmospheric predictability as revealed by naturally occurring
  analogues.
\newblock {\em Journal of the Atmospheric Sciences}, 26(4):636--646, 1969.

\bibitem{lorenz96Model}
{E.} Lorenz.
\newblock Predictability: A problem partly solved.
\newblock In {\em Predictability of Weather and Climate}, pages 40--58.
  Cambridge University Press, 2006.

\bibitem{mantegna1994linguistic}
{R.} Mantegna, {S.} Buldyrev, {A.} Goldberger, S.~Havlin, {C.} Peng, M.~Simons,
  and {H.} Stanley.
\newblock Linguistic features of noncoding {DNA} sequences.
\newblock {\em Physical Review Letters}, 73(23):3169--3172, 1994.

\bibitem{martinerie92}
J.~Martinerie, A.~Albano, A.~Mees, and P.~Rapp.
\newblock Mutual information, strange attractors, and the optimal estimation of
  dimension.
\newblock {\em Physical Review A}, 45(10):7058--7064, 1992.

\bibitem{McGill-1954}
W.~J. McGill.
\newblock Multivariate information transmission.
\newblock {\em Psychometrika}, 19(2):97--116, 1954.

\bibitem{McNames98anearest}
J.~McNames.
\newblock A nearest trajectory strategy for time series prediction.
\newblock In {\em Proceedings International Workshop on Advanced Black-Box
  Techniques for Nonlinear Modeling}, pages 112--128, 1998.

\bibitem{rwMeese}
R.~Meese and K.~Rogoff.
\newblock Empirical exchange rate models of the seventies: Do they fit out of
  sample?
\newblock {\em Journal of International Economics}, 14(1):3--24, 1983.

\bibitem{meissDynamcis}
J.~D. Meiss.
\newblock {\em Differential dynamical systems}.
\newblock Monographs on Mathematical Modeling and Computation. Society for
  Industrial and Applied Mathematics, Philadelphia, PA, 2007.

\bibitem{mischaikow99}
K.~Mischaikow, M.~Mrozek, J.~Reiss, and A.~Szymczak.
\newblock Construction of symbolic dynamics from experimental time series.
\newblock {\em Physical Review Letters}, 82(6):1144--1147, 1999.

\bibitem{tippdirk}
T.~Moseley, {J.} Kihm, D.~Connors, and D.~Grunwald.
\newblock Methods for modeling resource contention on simultaneous
  multithreading processors.
\newblock In {\em Proceedings of the International Conference on Computer
  Design}, 2005.

\bibitem{rwCCE}
J.~Mu{\'c}k and P.~Skrzypczy{\'n}ski.
\newblock Can we beat the random walk in forecasting {CEE} exchange rates?
\newblock National Bank of Poland Working Papers 127, National Bank of Poland,
  Economic Institute, 2012.

\bibitem{mytkowicz09}
T.~Myktowicz, A.~Diwan, and E.~Bradley.
\newblock Computers are dynamical systems.
\newblock {\em Chaos: An Interdisciplinary Journal of Nonlinear Science},
  19(3):033124, 2009.

\bibitem{todd-phd}
T.~Mytkowicz.
\newblock {\em Supporting experiments in computer systems research}.
\newblock PhD thesis, University of Colorado, November 2010.

\bibitem{PhysRevE.87.022905}
C.~Nichkawde.
\newblock Optimal state-space reconstruction using derivatives on projected
  manifold.
\newblock {\em Physical Review E}, 87(2):022905, 2013.

\bibitem{nussbaum01}
S.~Nussbaum and J.~Smith.
\newblock Modeling superscalar processors via statistical simulation.
\newblock In {\em Proceedings of the 2001 International Conference on Parallel
  Architectures and Compilation Techniques (PACT)}, pages 15--24, 2001.

\bibitem{binding}
E.~Olbrich, N.~Bertschinger, N.~Ay, and J.~Jost.
\newblock How should complexity scale with system size?
\newblock {\em The European Physical Journal B}, 63(3):407--415, 2008.

\bibitem{Olbrich97}
E.~Olbrich and H.~Kantz.
\newblock Inferring chaotic dynamics from time-series: On which length scale
  determinism becomes visible.
\newblock {\em Physical Letters A}, 232(1-2):63--69, 1997.

\bibitem{packard80}
N.~Packard, J.~P. Crutchfield, J.~Farmer, and R.~Shaw.
\newblock Geometry from a time series.
\newblock {\em Physical Review Letters}, 45(9):712--716, 1980.

\bibitem{pecoraUnified}
{L.} Pecora, L.~Moniz, J.~Nichols, and {T.} Carroll.
\newblock A unified approach to attractor reconstruction.
\newblock {\em Chaos: An Interdisciplinary Journal of Nonlinear Science},
  17(1):013110, 2007.

\bibitem{pesin1977characteristic}
Y.~B. Pesin.
\newblock Characteristic {L}yapunov exponents and smooth ergodic theory.
\newblock {\em Russian Mathematical Surveys}, 32(4):55--114, 1977.

\bibitem{petersen1989}
{K.} Petersen.
\newblock {\em Ergodic theory}.
\newblock Cambridge Studies in Advanced Mathematics. Cambridge University
  Press, 1989.

\bibitem{pikovsky86-sov}
{A.} Pikovsky.
\newblock Noise filtering in the discrete time dynamical systems.
\newblock {\em Soviet Journal of Communications, Technology and Electronics},
  31(5):911--914, 1986.

\bibitem{Robins02}
V.~Robins.
\newblock Computational topology for point data: Betti numbers of
  $\alpha$-shapes.
\newblock In {\em Morphology of Condensed Matter}, volume 600 of {\em Lecture
  Notes in Physics}, pages 261--274. Springer Berlin Heidelberg, 2002.

\bibitem{rosenstein94}
{M.} Rosenstein, {J.} Collins, and {C.} De~Luca.
\newblock Reconstruction expansion as a geometry-based framework for choosing
  proper delay times.
\newblock {\em Physica D: Nonlinear Phenomena}, 73(1-2):82--98, 1994.

\bibitem{rossler76}
O.~R\"ossler.
\newblock An equation for continuous chaos.
\newblock {\em Physical Letters A}, 57(5):397--398, 1976.

\bibitem{rudolph-measurable-dynamics}
D.~J. Rudolph.
\newblock {\em {Fundamentals of measurable dynamics. Ergodic theory on Lebesgue
  spaces.}}
\newblock Oxford: Clarendon Press, 1990.

\bibitem{sauer-delay}
T.~Sauer.
\newblock Time-series prediction by using delay-coordinate embedding.
\newblock In {\em Time Series Prediction: Forecasting the Future and
  Understanding the Past}. Santa Fe Institute Studies in the Sciences of
  Complexity, Santa Fe, NM, 1993.

\bibitem{sauer91}
T.~Sauer, J.~Yorke, and M.~Casdagli.
\newblock Embedology.
\newblock {\em Journal of Statistical Physics}, 65(3--4):579--616, 1991.

\bibitem{PhysRevLett.85.461}
T.~Schreiber.
\newblock Measuring information transfer.
\newblock {\em Physical Review Letters}, 85(2):461--464, 2000.

\bibitem{Shalizi2008}
C.~R. Shalizi and J.~P. Crutchfield.
\newblock Computational mechanics: Pattern and prediction, structure and
  simplicity.
\newblock {\em Journal of Statistical Physics}, 104(314):817--879, 2001.

\bibitem{Shannon1951}
C.~Shannon.
\newblock {Prediction and entropy of printed English}.
\newblock {\em Bell Systems Technical Journal}, 30(1):50--64, 1951.

\bibitem{shannon64}
C.~Shannon.
\newblock {\em The Mathematical Theory of Communication}.
\newblock University of Illinois Press, 1964.

\bibitem{sherwood02}
T.~Sherwood, E.~Perelman, G.~Hamerly, and B.~Calder.
\newblock Automatically characterizing large scale program behavior.
\newblock In {\em Proceedings of the International Conference on Architectural
  Support for Programming Languages and Operating Systems (ASPLOS)}, pages
  45--57, 2002.

\bibitem{Small2004283}
M.~Small and C.~K. Tse.
\newblock Optimal embedding parameters: a modelling paradigm.
\newblock {\em Physica D: Nonlinear Phenomena}, 194(3--4):283--296, 2004.

\bibitem{smithdatabound}
{L.} {S}mith.
\newblock Intrinsic limits on dimension calculations.
\newblock {\em Physical Letters A}, 133(6):283--288, 1988.

\bibitem{Smith199250}
{L.} {S}mith.
\newblock Identification and prediction of low dimensional dynamics.
\newblock {\em Physica D: Nonlinear Phenomena}, 58(1--4):50--76, 1992.

\bibitem{sorenson}
H.~W. Sorenson.
\newblock {\em Kalman filtering: Theory and application}.
\newblock IEEE Press, 1985.

\bibitem{sprottBook}
J.~C. Sprott.
\newblock {\em Chaos and time-series analysis}.
\newblock Oxford University Press, 2003.

\bibitem{strelioff2014bayesian}
C.~C. Strelioff and J.~P. Crutchfield.
\newblock Bayesian structural inference for hidden processes.
\newblock {\em Physical Review E}, 89(4):042119, 2014.

\bibitem{Studen1998}
M.~Studen\'{y} and J.~Vejnarov\'{a}.
\newblock The multiinformation function as a tool for measuring stochastic
  dependence.
\newblock In {\em Proceedings of the NATO Advanced Study Institute on Learning
  in Graphical Models}, pages 261--297. Kluwer Academic Publishers, 1998.

\bibitem{sugihara90}
G.~Sugihara and R.~May.
\newblock Nonlinear forecasting as a way of distinguishing chaos from
  measurement error in time series.
\newblock {\em Nature}, 344(6268):734--741, 1990.

\bibitem{takens}
F.~Takens.
\newblock Detecting strange attractors in fluid turbulence.
\newblock In {\em Dynamical systems and turbulence}, pages 366--381. Springer,
  Berlin, 1981.

\bibitem{tausz2012javaplex}
A.~Tausz, M.~Vejdemo-Johansson, and H.~Adams.
\newblock Java{P}lex: {A} research software package for persistent
  (co)homology.
\newblock In {\em Proceedings of ICMS 2014}, volume 8592 of {\em Lecture Notes
  in Computer Science}, pages 129--136, 2014.

\bibitem{theiler-window}
J.~Theiler.
\newblock Spurious dimension from correlation algorithms applied to limited
  time series data.
\newblock {\em Physical Review E}, 34(3):2427--2432, 1986.

\bibitem{tsonisdatabound}
A.~A. Tsonis, J.~B. Elsner, and K.~P. Georgakakos.
\newblock Estimating the dimension of weather and climate attractors: Important
  issues about the procedure and interpretation.
\newblock {\em Journal of the Atmospheric Sciences}, 50(15):2549--2555, 1993.

\bibitem{halting-problem}
A.~Turing.
\newblock {On computable numbers with an application to the
  Entscheidungsproblem}.
\newblock {\em Proceeding of the London Mathematical Society}, 1936.

\bibitem{total-cor}
S.~Watanabe.
\newblock Information theoretical analysis of multivariate correlation.
\newblock {\em IBM Journal of Research and Development}, 4(1):66--82, 1960.

\bibitem{weigend-book}
A.~Weigend and N.~Gershenfeld.
\newblock {\em Time series prediction: Forecasting the future and understanding
  the past}.
\newblock Santa Fe Institute Studies in the Sciences of Complexity, Santa Fe,
  NM, 1993.

\bibitem{wolf}
A.~Wolf.
\newblock {\em Quantifying chaos with {Lyapunov} exponents}.
\newblock Princeton University Press, Princeton NJ, 1986.

\bibitem{yeung2012first}
R.~W. Yeung.
\newblock {\em A first course in information theory}.
\newblock Springer Science \& Business Media, 2012.

\bibitem{Yule27}
U.~Yule.
\newblock On a method of investigating periodicities in disturbed series, with
  special reference to {W}olfer's sunspot numbers.
\newblock {\em Philosophical Transactions of the Royal Society of London
  Series. Series A, Containing papers of a Mathematical or Physical Character},
  226(636-646):267--298, 1927.

\bibitem{Zomorodian05}
A.~Zomorodian and G.~Carlsson.
\newblock Computing persistent homology.
\newblock {\em Discrete Computational Geometry}, 33(2):249--274, 2005.

\end{thebibliography}
